%% file: making-v8.tex
\documentclass[12pt,leqno]{article}
\newcommand{\kmcomment}[1]{}
\newcommand{\ds}{\displaystyle }

\usepackage[a4paper,truedimen,scale={.85,.85},centering]{geometry}
\usepackage[all]{xy}
\usepackage[toc,page]{appendix}
\usepackage{newtxtext,newtxmath}
\usepackage{cases,array}
{\list{}{\leftmargin=#1 \topsep=0pt \parsep=0pt \itemsep=0pt}\item[]}%
{\endlist}
\RequirePackage{verbatim}
{\noindent\quote\endgraf\verbatim}%
{\endverbatim\endquote\endgraf\medskip}%
\makeatletter
\newenvironment{qverb*}%
{\noindent\quote\endgraf\@nameuse{verbatim*}}%
{\@nameuse{endverbatim*}\endquote\endgraf\medskip}%
\makeatother
\input{km_cmd.tex}

\usepackage{graphicx,color,verbatim,setspace,supertabular,ntheorem} 
\usepackage{epic,eepic}
\usepackage{amsmath,latexsym,amssymb,amscd,alltt}
\usepackage[vcentermath]{youngtab}
\newcommand{\fadj}{d^{\ast}}
\newcommand{\wt}{\operatorname{wt}}

{
\theorembodyfont{\rmfamily}%{\bfseries}{\rmfamily}{\sffamily}
\theoremstyle{plain}%{break}  
\newtheorem{thm}{Theorem}[section] 

\newtheorem{exam}{Example}[section]
\newtheorem{defn}{Definition}[section]
\newtheorem{kmProp}{Proposition}[section]
\newtheorem{Lemma}{Lemma}[section]
 
\newtheorem{kmCor}[kmProp]{Corollary}
\newtheorem{kmRemark}{\textbf{Remark}}[section]
} 

\renewcommand{\include}{\input}
\numberwithin{equation}{section}
\title{Cohomology groups of 
homogeneous Poisson structures  
}

\author{
Kentaro Mikami\thanks{ 
  Department of Computer Science and Engineering 
  Akita University, partially supported by
JSPS KAKENHI Grant Number  JP26400063, JP23540067 and JP20540059.}
 \and Tadayoshi Mizutani\thanks{
Professor Emeritus, Saitama University }
}
\date{May 28, 2017}

\begin{document}
\parindent=0pt
%\onehalfspacing
%\begin{spacing}{2} 
\allowdisplaybreaks 
        \setstretch{1.5}

\maketitle
\thispagestyle{empty}
\tableofcontents 

\include{Intro}

\include{Include-0v1}

\include{Contrib}

\include{Top-Betti}

\include{Euler-v4}

\include{Combi}
\include{Suppl-2}

%\include{MyAppend}

%$\kk{111}$

\bibliographystyle{plain}
\bibliography{km_refs}
\end{document}

%% file: km_cmd.tex
\newcommand{\kmqed}{\hfill\ensuremath{\blacksquare}}

\newcommand{\myB}{B}
 
\newcommand{\Oira}[2]{\overline{\chi}^{#1}_{\;\; #2}}

\newcommand{\T}[1]{T_{#1}}
\newcommand{\ovfrakg}{\overline{\frakg}}

\newcommand{\mydq}{d}

\newcommand{\mydov}{\overline{d}}
\newcommand{\phd}[1]{ \phi_{#1}} 
\newcommand{\md}[1]{ \text{md}_{#1}} 
\newcommand{\mmd}[1]{\tbinom{n-1+#1}{#1}} 
 
\newcommand{\futo}[1]{{\mathbf #1}}

\newcommand{\pdel}{\partial} 
\newcommand{\pdl}[1]{\:\pdel_{#1}}
\newcommand{\kd}[1]{k'{}_{#1}}
\newcommand{\last}{top }
\newcommand{\LS}[2]{ \Lambda^{#2} \frakSS{#1}}

\newcommand{\frakS}[1]{\mathfrak{S}_{#1}}
\newcommand{\frakSS}[1]{\overline{\mathfrak{S}}_{#1}}
\newcommand{\myS}[1]{S_{#1}}
\newcommand{\mySS}[1]{\overline{S}_{#1}}
\newcommand{\myC}[2]{\text{C}^{#1}_{#2}}
\newcommand{\myCC}[2]{\overline{\text{C}}^{#1}_{#2}}
\newcommand{\myH}[2]{\text{H}^{#1}_{#2}}
\newcommand{\myHH}[2]{\overline{\text{H}}^{#1}_{#2}}
\newcommand{\frakX}[1]{\mathfrak{X}_{#1}} 
\newcommand{\mN}{\ensuremath{\mathbb{N}}} %Use in Math-mode
\newcommand{\mR}{\ensuremath{\mathbb{R}}} %Use in Math-mode
\newcommand{\kk}[1]{\ensuremath{k_{#1}}}

\newcommand{\junban}{\text{od}} 
\newcommand{\zL}[2]{ z^{(#1)}_{\bullet < #2}}
\newcommand{\zR}[2]{ z^{(#1)}_{\bullet > #2}}
\newcommand{\zF}[1]{ z^{(#1)}_{\text{Full}}} 

\newcommand{\zM}[2]{ z^{(#1)}_{ < #2 >}}
\newcommand{\zBetw}[3]{ z^{(#1)}_{ #2< \bullet < #3}}

\newcommand{\frakg}{\mathfrak{g}}
\newcommand{\frakh}{\mathfrak{h}}

\newcommand{\fraksl}{\mathfrak{sl}}

\newcommand{\fraksp}{\mathfrak{sp}}

\newcommand{\Z}[2]{z^{#2}_{#1}}

\newcommand{\myP}[1]{\ensuremath{P_{#1}}}
\newcommand{\myPP}[1]{\ensuremath{R_{[#1]}}}
\newcommand{\tbdl}[1]{\mathrm{T}(#1)}

\newcommand{\Sbt}[2]{[#1,#2]}                
\newcommand{\SbtS}[2]{[#1,#2]_{\text{\footnotesize S}}} %Use in Math-ode 
\newcommand{\SbtR}[2]{\ensuremath{[#1,#2]_{\mR}}}                
\newcommand{\parity}[1]{(-1)^{#1}}

\newcommand{\Pkt}[2]{\{#1,#2\}}% 

\newcommand{\Bkt}[2]{\{#1,#2\}_{\pi}}

\newcommand{\w}[1]{\boldsymbol{x}^{#1}}
\newcommand{\myUnit}{\delta_{\futo{0}}} 
\newcommand{\polygon}[1]{\operatorname{Poly}(\mR^{#1})} 
\newcommand{\myCasim}[1]{\operatorname{Cent}_{\pi}#1}
\newcommand{\myCasimm}{\operatorname{Cent}_{\pi}}
\newcommand{\myCO}[3]{\Lambda^{#1}_{#2}(#3)}

\newcommand{\z}[1]{z_{#1}}
\newcommand{\eps}[1]{\epsilon_{#1}}
\newcommand{\mynabla}[2]{\nabla^{#1}_{#2}}
\newcommand{\wnabla}[2]{\widehat{\nabla}^{#1}_{#2}}

\newcommand{\remainder}{\mathbf{r}}
\newcommand{\myLM}[1]{\operatorname{LM}(#1)}

%% file: Intro.tex
\section{Introduction}% 
Researches on the  cohomology of  formal Hamiltonian vector fields on
symplectic $2n$-planes have made much progress after the work of
Gel'fand-Kalinin-Fuks (\cite{MR0312531}).    We see S.~Metoki's
work \cite{metoki:shinya} in which the author showed there exists a
non-trivial relative cocycle in degree 9 in the case  $n=1$.  

In 2009, D.~Kotschick and S.~Morita (\cite{KOT:MORITA}) gave some
important contribution to this research area.  One of the important
notion adopted there was the notion of weight (Metoki says it as ``type'' in 
\cite{metoki:shinya}).    
They decomposed the Gel'fand-Fuks cohomology groups
according to the weights 
and investigated the multiplication of cocycles.  

Also in a similar method,  the structure of
cochains were studied in the case of symplectic 2-plane in
(\cite{M:N:K}) and in (\cite{KM:affirm}), and further 
in the case of symplectic 4-plane in \cite{Mik:Nak} or   
in the case of symplectic 6-plane in \cite{KM:D6}.  

In this paper, we will study of Lie algebra cohomology groups of 
formal Hamiltonian vector fields of some kind of Poisson structures on $\ds
\mR^n$ 
including linear symplectic structures on $\ds\mR^{2n}$ 
by the unified way introducing the notion of weight defined
depending to homogeneity of Poisson structures. 
Using manipulation of Young diagrams, we get the recursive formulae for the
Euler characteristics depending on the homogeneity and the weight
(see Lemma \ref{lemma:euler:recur}).  Moreover, in Theorem
\ref{thm:euler:num},  
we state that for each 1-homogeneous Poisson structure, the   
Euler characteristics is always zero for any weight.  
And we will give some different view point for studying Poisson cohomology
groups.

\section{Preliminaries} \label{sect:preliminaries}
\subsection{Homology and cohomology of Lie algebra}
We recall the definition of homology and cohomology groups of Lie algebras.  
Given a Lie algebra $\ds\frakg$, we can consider cohomology groups or
homology groups 
%on a representation space $(V,\phi)$ 
of a Lie algebra $\frakg$.   
        
\subsubsection{Cohomology groups}
        Let $\myC{m}{} = \{ \sigma : \frakg \times \cdots \times \frakg
        \rightarrow \mR \mid \text{skew-symmetric }
        m\text{-multi-linear}\} =  \Lambda^m (\frakg^{*}) $.  
        \begin{align*}
                        (d\sigma)(u_0,\ldots, u_m) = & \sum_{0\leq
                        i<j\leq m} (-1)^{i+j}
                \sigma( [u_i,u_j],\ldots, \widehat{u_i},\ldots,
                \widehat{u_j},\ldots) \; . 
% \\& + \sum_{i=0}^{m} (-1)^i \phi(u_i) ( 
%                \sigma( \ldots, \widehat{u_i},\ldots) )\;.
        \end{align*}
\subsubsection{Homology groups}
%        Let $\myC{}{m} = V \otimes \Lambda^m (\frakg) $.  
        Let $\myC{}{m} = \Lambda^m (\frakg) $.  
        \begin{align*}
                \pdel ( \&^{\wedge}(u_1,\ldots, u_m)) = 
                &  \sum_{1\leq
                        i<j\leq m} (-1)^{i+j}
                \&^{\wedge} ( [u_i,u_j],\ldots, \widehat{u_i},\ldots,
                \widehat{u_j},\ldots)\;,  
% \\&
% + \sum_{i=1}^{m} (-1)^i \phi(u_i)(\futo{v}) \otimes \&^{\wedge}  
% ( \ldots, \widehat{u_i},\ldots) \;, \\
        \end{align*} 
        where the notation $\ds\&^{\wedge} (u_1,\ldots, u_m) $ means the
        wedge product of $\ds u_1,\ldots,u_m$, namely, 
        $\ds\&^{\wedge} (u_1,\ldots, u_m)  = u_1 \wedge \cdots \wedge u_m$.  
%\newpage
\subsubsection{Duality}
Hereafter, we only consider the trivial representation $\ds (V,\phi)$ of
Lie algebra, namely, $\ds V=\mR$ and $\phi$ is the zero map. 
        Let $\ds \futo{u}_i$ be vectors and 
        $\ds \futo{a}_i$ be covectors.  
        The natural coupling of 
        $\ds \futo{u}_1 \wedge \cdots  \wedge  \futo{u}_k $
        and 
        $\ds \futo{a}_1 \wedge \cdots  \wedge  \futo{a}_k $
        is given by 
        \begin{align*}
        \langle \futo{u}_1 \wedge \cdots  \wedge  \futo{u}_k ,
                \futo{a}_1 \wedge \cdots  \wedge  \futo{a}_k \rangle =&  
         \sum_{\sigma} \operatorname{sgn}(\sigma) 
        \langle \futo{u}_1 ,  \futo{a}_{\sigma(1)}\rangle 
        \langle \futo{u}_2 ,  \futo{a}_{\sigma(2)}\rangle 
        \cdots 
        \langle \futo{u}_k ,  \futo{a}_{\sigma(k)}\rangle \\
                =&  \begin{vmatrix}
                        \langle \futo{u}_1 ,  \futo{a}_{1}\rangle  &
                        \ldots &
                        \langle \futo{u}_1 ,  \futo{a}_{k}\rangle  \\
                        \vdots & \vdots & \vdots \\
                        \langle \futo{u}_k ,  \futo{a}_{1}\rangle  &
                        \ldots &
                        \langle \futo{u}_k ,  \futo{a}_{k}\rangle
                \end{vmatrix} 
        \end{align*}
        When we are interested in the behavior of 
        $\ds \futo{u}_1 \wedge \futo{u}_2$, we have the following
        formula:
        \begin{kmProp}\label{prop:two}
                \begin{align*}
                        &        \langle (\futo{u}_1 \wedge \futo{u}_2) \wedge 
                        \cdots  \wedge  \futo{u}_k ,
                \futo{a}_1 \wedge \cdots  \wedge  \futo{a}_k \rangle \\
                        =&  
                        \sum_{i<j} (-1)^{i+j+1} 
                        \langle \futo{u}_1 \wedge \futo{u}_2, 
                        \futo{a}_i \wedge \futo{a}_j \rangle  
                        \langle \futo{u}_3 \wedge  
                        \cdots  \wedge  \futo{u}_k ,
                \futo{a}_1 \wedge \cdots
                        %\wedge 
                        \widehat{ \futo{a}_i} 
                        %\wedge
                        \cdots
                        %\wedge 
                        \widehat{ \futo{a}_j} %\wedge 
                        \cdots
                        \wedge  \futo{a}_k \rangle   
                \end{align*}
        \end{kmProp}
\textbf{Proof:}
        \begin{align*}
                \text{LHS} =&  
         \sum_{\sigma} \operatorname{sgn}(\sigma) 
        \langle \futo{u}_1 ,  \futo{a}_{\sigma(1)}\rangle 
        \langle \futo{u}_2 ,  \futo{a}_{\sigma(2)}\rangle 
        \langle \futo{u}_3 ,  \futo{a}_{\sigma(3)}\rangle 
        \cdots 
        \langle \futo{u}_k ,  \futo{a}_{\sigma(k)}\rangle \\
                =& \sum_{i<j} \sum_{ \{\sigma(1),\sigma(2)\}=\{i,j\}} 
         \operatorname{sgn}(\sigma) 
        \langle \futo{u}_1 ,  \futo{a}_{\sigma(1)}\rangle 
        \langle \futo{u}_2 ,  \futo{a}_{\sigma(2)}\rangle 
        \langle \futo{u}_3 ,  \futo{a}_{\sigma(3)}\rangle 
        \cdots 
        \langle \futo{u}_k ,  \futo{a}_{\sigma(k)}\rangle \\
                =& \sum_{i<j} \sum_{
                        \{\tau(3),\ldots,\tau(k)\}=\{\hat{i},\hat{j}\}} 
                (
                \operatorname{sgn}\begin{pmatrix} 1 &2 & \\
                i & j & \tau \end{pmatrix} 
        \langle \futo{u}_1 ,  \futo{a}_{\sigma(1)}\rangle 
        \langle \futo{u}_2 ,  \futo{a}_{\sigma(2)}\rangle 
        \langle \futo{u}_3 ,  \futo{a}_{\tau(3)}\rangle 
        \cdots 
        \langle \futo{u}_k ,  \futo{a}_{\tau(k)}\rangle \\
                &
                        + 
                \operatorname{sgn}\begin{pmatrix} 1 &2 & \\
                j & i & \tau\end{pmatrix} 
        \langle \futo{u}_1 ,  \futo{a}_{\sigma(1)}\rangle 
        \langle \futo{u}_2 ,  \futo{a}_{\sigma(2)}\rangle 
        \langle \futo{u}_3 ,  \futo{a}_{\tau(3)}\rangle 
        \cdots 
                \langle \futo{u}_k ,  \futo{a}_{\tau(k)}\rangle )\\
                =& \sum_{i<j} \sum_{
                        \{\tau(3),\ldots,\tau(k)\}=\{\hat{i},\hat{j}\}} 
                (
                (-1)^{i+j+1} 
                \operatorname{sgn} (\tau) 
        \langle \futo{u}_1 ,  \futo{a}_{i}\rangle 
        \langle \futo{u}_2 ,  \futo{a}_{j}\rangle 
        \langle \futo{u}_3 ,  \futo{a}_{\tau(3)}\rangle 
        \cdots 
        \langle \futo{u}_k ,  \futo{a}_{\tau(k)}\rangle \\
                &
                        + 
                (-1)^{i+j+2} 
                \operatorname{sgn}( \tau )
        \langle \futo{u}_1 ,  \futo{a}_{j}\rangle 
        \langle \futo{u}_2 ,  \futo{a}_{i}\rangle 
        \langle \futo{u}_3 ,  \futo{a}_{\tau(3)}\rangle 
        \cdots 
                \langle \futo{u}_k ,  \futo{a}_{\tau(k)}\rangle )\\
                =& \sum_{i<j} 
                (-1)^{i+j+1} 
                (
        \langle \futo{u}_1 ,  \futo{a}_{i}\rangle 
        \langle \futo{u}_2 ,  \futo{a}_{j}\rangle 
                - 
        \langle \futo{u}_1 ,  \futo{a}_{2}\rangle 
                \langle \futo{u}_2 ,  \futo{a}_{1}\rangle ) 
                \\&\qquad 
        \langle \futo{u}_3 \wedge \cdots \wedge  \futo{u}_{k}, 
                \futo{a}_1 \wedge \cdots
                        %\wedge 
                        \widehat{ \futo{a}_i} 
                        %\wedge
                        \cdots
                        %\wedge 
                        \widehat{ \futo{a}_j} %\wedge 
                        \cdots
                        \wedge  \futo{a}_k \rangle   
        \end{align*}
\(\ds \langle d( \futo{u}), \futo{a}_1 \wedge \futo{a}_2 \rangle = - 
 \langle  \futo{u}, [\futo{a}_1 , \futo{a}_2] \rangle \quad 
 \text{and}\quad  
 \langle d( \futo{u}), \futo{a}_1 \wedge \futo{a}_2 \rangle =  
 \langle  \futo{u}, \fadj( \futo{a}_1 \wedge \futo{a}_2 )\rangle\)
 implies
\(\ds\fadj( \futo{a}_1 \wedge \futo{a}_2) = - [\futo{a}_1, \futo{a}_2 ]
 \)  
\begin{align*}
        & \langle d( \futo{u}_1 \wedge \futo{u}_2), \futo{a}_1 \wedge \futo{a}_2
  \wedge \futo{a}_3       \rangle =  
        \langle  \futo{u}_1 \wedge \futo{u}_2, \fadj( 
        \futo{a}_1 \wedge \futo{a}_2 \wedge \futo{a}_3 ) \rangle \\ 
        \noalign{and}
        & \langle d( \futo{u}_1 \wedge \futo{u}_2), \futo{a}_1 \wedge \futo{a}_2
  \wedge \futo{a}_3       \rangle =  
        \langle d( \futo{u}_1) \wedge \futo{u}_2 
       -  
         \futo{u}_1 \wedge d(\futo{u}_2 ), \futo{a}_1 \wedge \futo{a}_2
        \wedge \futo{a}_3       \rangle \\& 
        = \sum_{i<j} (-1)^{i+j+1} ( 
        \langle d( \futo{u}_1), \futo{a}_i \wedge \futo{a}_j \rangle
        \langle \futo{u}_2, \futo{a}_k \rangle  
        -
        \langle d( \futo{u}_2), \futo{a}_i \wedge \futo{a}_j \rangle
        \langle \futo{u}_1, \futo{a}_k \rangle  
        ) \\
        & 
        = \sum_{i<j} (-1)^{i+j+1} ( 
        \langle  \futo{u}_1, \fadj(\futo{a}_i \wedge \futo{a}_j) \rangle
        \langle \futo{u}_2, \futo{a}_k \rangle  
        -
        \langle  \futo{u}_2 , \fadj(\futo{a}_i \wedge \futo{a}_j) \rangle
        \langle \futo{u}_1, \futo{a}_k \rangle  
        ) \\ &
        = \sum_{i<j} (-1)^{i+j+1}  
        \langle  \futo{u}_1 \wedge \futo{u}_2 , \fadj(\futo{a}_i \wedge
        \futo{a}_j) \wedge \futo{a}_{k}  \rangle
         \\
         \noalign{thus}
         & 
         \fadj( 
        \futo{a}_1 \wedge \futo{a}_2 \wedge \futo{a}_3 ) = 
         \sum_{i<j} (-1)^{i+j+1}  
        \fadj(\futo{a}_i \wedge
        \futo{a}_j) \wedge \futo{a}_{k} =
         \sum_{i<j} (-1)^{i+j}  
        [ \futo{a}_i,  
        \futo{a}_j] \wedge \futo{a}_{k} 
\end{align*}
In general, 
\begin{align*} 
        & \langle 
        d( \futo{u}_1 \wedge \cdots \wedge \futo{u}_k ), 
        \futo{a}_1 \wedge \cdots \wedge \futo{a}_{1+k}\rangle =  
        \langle 
         \futo{u}_1 \wedge \cdots \wedge \futo{u}_k , \fadj( 
        \futo{a}_1 \wedge \cdots \wedge \futo{a}_{1+k})\rangle  \\
        =& \langle \sum_{i} (-1)^{1+i} 
        d( \futo{u}_i)  \wedge \cdots \wedge \widehat{\futo{u}_i} \wedge \cdots \wedge \futo{u}_k , 
        \futo{a}_1 \wedge \cdots \wedge \futo{a}_{1+k}\rangle \\=  & 
        \sum_{i} (-1)^{1+i} \sum_{p<q}(-1)^{1+p+q}\langle 
        d( \futo{u}_i), \futo{a}_p \wedge \futo{a}_q \rangle  
\langle \cdots \widehat{\futo{u}_i} \wedge \cdots \wedge \futo{u}_k , 
        \futo{a}_1 \wedge \cdots \wedge \widehat{\futo{a}_{p}} 
        \cdots\widehat{\futo{a}_{q}}    \wedge 
 \cdots \rangle \\
        =&  \sum_{p<q}(-1)^{1+p+q}
        \sum_{i} (-1)^{1+i}
        \langle 
         \futo{u}_i, \fadj(\futo{a}_p \wedge \futo{a}_q) \rangle  
\langle \cdots \widehat{\futo{u}_i} \wedge \cdots \wedge \futo{u}_k , 
        \futo{a}_1 \wedge \cdots \wedge \widehat{\futo{a}_{p}} 
        \cdots\widehat{\futo{a}_{q}}    \wedge 
 \cdots \rangle \\
        =&  \sum_{p<q}(-1)^{1+p+q} 
        \langle 
         \futo{u}_1, \wedge \cdots \wedge \futo{u}_k , 
\fadj(\futo{a}_p \wedge \futo{a}_q) \wedge   
        \futo{a}_1 \wedge \cdots \wedge \widehat{\futo{a}_{p}} 
        \cdots\widehat{\futo{a}_{q}}    \wedge 
\cdots \rangle \; . 
\end{align*}
Thus
\begin{align*}         
          \fadj( 
        \futo{a}_1 \wedge \cdots \wedge \futo{a}_{1+k}) & = 
        \sum_{p<q}(-1)^{1+p+q} 
\fadj(\futo{a}_p \wedge \futo{a}_q) \wedge   
        \futo{a}_1 \wedge \cdots \wedge \widehat{\futo{a}_{p}} 
        \cdots\widehat{\futo{a}_{q}}    \wedge 
 \cdots \\ 
         & = 
        \sum_{p<q}(-1)^{p+q} 
        [\futo{a}_p , \futo{a}_q] \wedge   
        \futo{a}_1 \wedge \cdots \wedge \widehat{\futo{a}_{p}} 
        \cdots\widehat{\futo{a}_{q}}    \wedge 
 \cdots  \; .
\end{align*}
We get a more general formula of Proposition \ref{prop:two} as follows:
\begin{kmProp}
\begin{align*} 
        &\langle (\futo{u}_1 \wedge \dots\wedge \futo{u}_k) \wedge
        \futo{u}_{k+1} \wedge \cdots  \wedge  \futo{u}_n ,
                \futo{a}_1 \wedge \cdots  \wedge  \futo{a}_n \rangle \\
                        =& 
                        \sum_{\#M=k} (-1) ^{ |M| + k(k+1)/2} 
\langle \futo{u}_1 \wedge \dots\wedge \futo{u}_k,  
                        \futo{a}_{M} \rangle
                        \langle 
                \futo{u}_{k+1} \wedge \cdots  \wedge  \futo{u}_n , 
                \futo{a}_{\widehat{M}} \rangle
                        \end{align*} 
where $M$ runs in the family of $k$-subsets of $1,2,\ldots,n$.   
        $\ds \futo{a}_{M} = \futo{a}_ {m_1} \wedge \cdots \wedge 
         \futo{a}_ {m_k}$ 
        with $\ds M = \{ m_1 < m_2< \cdots < m_k\}$ and $\ds |M| =
        \sum_{i=1}^k m_i$,  
         and  $\ds \widehat{M} =
        \{1,2,\ldots,n\}\smallsetminus M$. 
\end{kmProp}

%%%%%%%%%%%%%%%%%%%

%\subsubsection{Elementary or general fact} 
When the chain and the cochain spaces are finite dimensional, then  we
can express each Betti number by the dimensions of the kernel spaces
and the (co)chain spaces as follows: 

%\[\rightarrow C^{m -1}\rightarrow C^m\rightarrow C^{m+1}\rightarrow\]
Cochain complex case:
\begin{center}
\begin{tabular} {c|*{7}{c}}
& $\rightarrow$ &  
        $\ds \myC{m -1}{}$ & $\rightarrow$ &   
        $\ds \myC{m}{} $ &  $\rightarrow$ &   
        $\ds \myC{m+1}{}$ &  $\rightarrow$ \\\hline
$\dim$ & & $d_{m-1}$ && $d_{m}$ &&  $d_{m+1}$ &  \\
$\ker$ & & $k_{m-1}$ && $k_{m}$ &&  $k_{m+1}$ &  \\\hline
Betti  & &  && $k_{m}+k_{m-1}-d_{m-1}$ &&   &  \\\hline
 \end{tabular} 
\end{center}

\bigskip

Chain complex case: 
%\[\rightarrow C^{m -1}\rightarrow C^m\rightarrow C^{m+1}\rightarrow\]
\begin{center}
\begin{tabular} {c|*{7}{c}}
& $\leftarrow$ &  
        $\ds \myC{}{m -1}$ & $\leftarrow$ &   
        $\ds \myC{}{m} $ &  $\leftarrow$ &   
        $\ds \myC{}{m+1}$ &  $\leftarrow$ \\\hline
$\dim$ & & $d_{m-1}$ && $d_{m}$ &&  $d_{m+1}$ &  \\
$\ker$ & & $\ell_{m-1}$ && $\ell_{m}$ &&  $\ell_{m+1}$ &  \\\hline
Betti  & &  && $\ell_{m}+\ell_{m+1}-d_{m+1}$ &&   &  \\\hline
\end{tabular} 
\end{center}

\bigskip

Now assume that $\ds\dim \myC{m}{}$ and $\ds\dim \myC{}{m}$ are finite
dimensional and $\ds\dim \myC{m}{} = \dim \myC{}{m}$ for each $m$.  
%Relations between the two: 
Since $\ds \langle d_{p} \sigma , u \rangle = \langle  \sigma ,
\pdel_{p+1} u \rangle $, we have $\ds \left( \text{Ker} d_{p}\right)^{0}
= \text{Im} (\pdel_{p+1})$.  
Evaluating the dimensions of both sides, we have 
\[ \dim C^p - \dim \text{Ker} d_p = \dim C_{p+1} - \dim \text{Ker}
\pdel_{p+1}\;.\] 
Using the notations above, this is equivalent to
\[ \dim C^{p} - k_p = \dim C_{p+1} - \ell_{p+1}\;. \]
Thus, we claim that 
\begin{align*}
        \dim H_{m} =& \ell_{m} + \ell_{m+1} - \dim C_{m+1} \\
        =& k_{m-1} - \dim C^{m-1} + \dim C_{m} 
        \\& 
        + k_{m} - \dim C^{m} + \dim C_{m+1}  - \dim C_{m+1} \\ 
        =& k_{m-1} + k_{m} - \dim C^{m-1} \quad (\text{if} \; \dim
        C^j= \dim C_j \;)\\
        =& \dim H^{m} \;. 
\end{align*} 
Thus, we have the following statement and 
we may choose cochain complex or chain complex if we are
interested in only Betti numbers.   
\begin{kmProp}
        If the chain spaces and the cochain spaces are finite
        dimensional and $\ds\dim \myC{m}{} = \dim \myC{}{m}$ for each
        $m$, then  $\ds \dim H^{m} = \dim H_{m}$ for each $m$. 
\end{kmProp}

%%%%%%%%%%%%%%%%%%%%%%

\subsection{Schouten bracket and Poisson structure} \label{sub:schou} 
Let us recall the Schouten bracket on a $n$-dimensional smooth manifold
$M$.   Let $\ds\Lambda^j \tbdl{M}$ be the space of $j$-vector fields on
$M$. In particular, $\ds\Lambda^1 \tbdl{M}$ is ${\mathfrak X}(M)$, the
Lie algebra of smooth vector fields on $M$, and $\ds\Lambda^0 \tbdl{M}$
is $\ds C^{\infty}(M)$.  $\ds A = \sum_{j=0}^ n \Lambda^j \tbdl{M}$ is
the exterior algebra of multi-vector fields on $M$. 

%%%%%%%%%%%%%%%%%%%%%%%%%%%%%%%%%%%%%%%%%%%%

For 
$\ds P \in \Lambda^{p} \tbdl{M}$ and 
$\ds Q \in \Lambda^{q} \tbdl{M}$,  the Schouten bracket   
$\ds\SbtS{P}{Q}$ is defined to be an element in $\ds\Lambda^{p+q-1} \tbdl{M}$.  
This bracket satisfies the following formulas. 
\begin{align*}
\SbtS{Q}{P} & = - \parity{ (q+1)(p+1)} \SbtS{P}{Q}\quad
		\text{(symmetry)}\;, \\
 0 & = \mathop{\frakS{}}_{p,q,r} \parity{(p+1)(r+1)} \SbtS{P}{\SbtS{Q}{R}}   
 \quad \text{(the Jacobi identity)}\;,\\  
%
%\SbtS{P}{\SbtS{Q}{R}} &= 
% \SbtS{ \SbtS{P}{Q}}{R} 
%+ \parity{ (p+1)(q+1)} \SbtS{Q}{\SbtS{P}{R}}\\
\SbtS{P}{Q\wedge R} &= \SbtS{P}{Q}\wedge R + (-1)^{(p+1)q} Q \wedge
\SbtS{P}{R}\;,\\ 
\SbtS{P \wedge Q}{R} &= P \wedge \SbtS{Q}{R}  + (-1)^{q (r+1)} 
\SbtS{P}{R} \wedge Q\;,\\
\noalign{another expression of Jacobi identity is the next} 
\SbtS{P}{\SbtS{Q}{R}} &= 
\SbtS{ \SbtS{P}{Q}} {R} + (-1)^{(p+1)(q+1)} \SbtS{Q}{\SbtS{P}{R}} 
\;,\\ 
\SbtS{\SbtS{P}{Q}}{R} &= 
\SbtS{P}{\SbtS{Q}{R}} + (-1)^{(q+1)(r+1)} \SbtS{\SbtS{P}{R}}{Q}\;.   
\\\noalign{Note that these formulae above are valid without any change
		for the case of $\mR$-linear wedge product.}
\\\noalign{To define the usual Schouten bracket uniquely, we need 
also the following basic formulae concerning functions.}
\SbtS{X}{Y} &= \text{Jacobi-Lie bracket of}\ X\ \text{and}\ Y\;, \\
\SbtS{X}{f} &= \langle X, df \rangle\;.  
\end{align*} 
For vector fields $X_1,\ldots.X_p$ and $Y_1,\ldots,Y_q$ the Schouten bracket
of $X_1\wedge\cdots \wedge X_p$ and 
        $Y_1\wedge\cdots \wedge Y_q$  
is given by 
\begin{align*}
        \SbtS{X_1\wedge\cdots \wedge X_p}
        {Y_1\wedge\cdots \wedge Y_q} 
         = \sum_{i,j} (-1)^{i+j} \Sbt{X_i}{Y_j} \wedge 
        (X_1\wedge\cdots \widehat{X_i} \cdots X_{p}) 
        \wedge 
        (Y_1\wedge\cdots \widehat{Y_j} \cdots Y_{q}) \;.
\end{align*}

%%%%%%%%%%%%%%%%%%%%%%%%%%%%%%%%%%%%%%%%%%%%%%%%%%%%%%%%%%%%%%%%%%%%%%%

Let $\pi$ be a 2-vector fields on $M$. $\pi$ is a Poisson structure if
and only if $\ds\SbtS{\pi}{\pi} = 0$. Locally, let $\ds (x_1\ldots x_n)$
be a local coordinates and 
\[ \ds\pi = 
\frac{1}{2}\sum_{i,j} p_{ij} \pdel_i \wedge \pdel_j 
\quad \text{where}\quad \pdel_i = \frac{\pdel}{\pdel x_i} \ \text{and}\  
p_{ij} + p_{ji}=0\;. \]    
The Schouten bracket of $\pi$ and itself is calculated as follows 
\begin{align*}
2 \SbtS{\pi}{\pi} =&\sum_{i,j} \Sbt{\pi}{p_{ij}\pdel_i\wedge\pdel_j}\\
=& \sum_{i,j}(  \Sbt{\pi}{p_{ij}} \wedge \pdel_i \wedge \pdel_j  
+ p_{ij} \SbtS{\pi}{\pdel_i} \wedge \pdel_j  
- p_{ij}\pdel_i \wedge \SbtS{\pi}{\pdel_j}  ) \;, 
\\ 
4 \SbtS{\pi}{\pi} =& 
\sum_{k\ell} \sum_{i,j}(  
\Sbt{ p_{k\ell} \pdel_k \wedge \pdel_{\ell} }{p_{ij}} \wedge \pdel_i \wedge \pdel_j  
+ p_{ij} \SbtS{ p_{k\ell} \pdel_k \wedge \pdel_{\ell}
}{\pdel_i} \wedge \pdel_j  \\ & 
- p_{ij}\pdel_i \wedge \SbtS{ p_{k\ell} \pdel_k \wedge \pdel_{\ell} }{\pdel_j})
\\ =& 
\sum_{k\ell} \sum_{i,j}(  
p_{k\ell} (\pdel_{\ell} p_{ij})  \pdel_k  \wedge \pdel_i \wedge \pdel_j 
- p_{k\ell} (\pdel_{k} p_{ij})  \pdel_{\ell}\wedge \pdel_i \wedge \pdel_j 
\\ & 
- p_{ij} ( \pdel_i   p_{k\ell}) \pdel_k \wedge \pdel_{\ell}  
        \wedge \pdel_j  + p_{ij}(\pdel_j p_{k\ell})  \pdel_i \wedge 
        \pdel_k \wedge \pdel_{\ell} )
\\ =& 4 \sum_{ijk} \sum_{\ell} p_{i\ell} \frac{\pdel p_{jk}}{\pdel x_{\ell}}\;  
        \pdel_i \wedge \pdel_j \wedge \pdel_k
		\;. 
\end{align*} 
The Poisson bracket of $f$ and $g$  is given by $\Bkt{f}{g} = \pi(df,dg)$ so 
$\ds\Bkt{x_i}{x_j} = p_{ij}$ and 
\begin{align*}
\Bkt{\Bkt{x_i}{x_j}}{x_k} =& \Bkt{p_{ij}}{x_k} = 
\frac{1}{2} \sum_{\lambda,\mu} p_{\lambda,\mu} (
\pdel_{\lambda}p_{ij} \pdel_{\mu}x_k 
- \pdel_{\mu }p_{ij} \pdel_{\lambda}x_k )
= - \sum_{\lambda} p_{k\lambda} \pdel_{\lambda} p_{ij} \;, 
\end{align*}
$\pi$ is Poisson if and only if Jacobi identity for the bracket holds, and it
is equivalent to 
\begin{equation} 
\mathop{\mathfrak S}_{i,j,k} 
\sum_{\lambda} p_{k\lambda} \frac{ \pdel p_{ij}}{\pdel x_{\lambda} } 
                = 0 
		\quad 
		(\text{cyclic sum with respect to  } i,j,k) \;. 
				\label{bkt:jacobi}
\end{equation}

A given Poisson structure $\pi$ on a manifold $M$ yields two kinds of
Lie algebras: one is the space of smooth functions $\ds C^{\infty}(M)$
on $M$ with the Poisson bracket, the other is the space of Hamiltonian
vector fields defined by $\Bkt{f}{\cdot} = \pi(df,\cdot)$ with the
Jacobi-Lie bracket. The relation of these Lie algebras is  described by
the following short exact sequence 
\[ 0 \rightarrow \myCasim(M) \to C^{\infty}(M) \rightarrow \big\{
\Bkt{f}{\cdot}\in \mathfrak{X}(M) \big\} \rightarrow 0\] 
where $\ds \myCasim(M) = \big\{ f \in C^{\infty}(M) \mid \Bkt{f}{\cdot} =
0\big\}$, the center of $\ds C^{\infty}(M)$, whose element is  called a
Casimir function.  Obviously  $\ds \myCasim(M)$ is a ideal of Lie algebra.  Thus we
have an isomorphism 
\[  C^{\infty}(M) /\myCasim(M)  \cong \big\{
 \Bkt{f}{\cdot}\in \mathfrak{X}(M) \big\}\] 
 as Lie algebras.  Roughly speaking, the space $\ds \myCasim(M)$ shows how far
the Poisson structure of the manifold $M$ is from symplectic structure,
in the case of symplectic structure, $\ds \myCasim(M)$ is just the space of
constant functions.  

\begin{defn} 
We endow the space $\ds  \mR^n$ with the Cartesian coordinates $\ds
(x_1,\ldots, x_n)$. Then any 2-vector field $\pi$ is written as 
\[  \pi = \frac{1}{2}\sum_{i,j} p_{ij}(x) \pdel_i \wedge \pdel_j
\qquad \text{where\quad} \pdel_i = \frac{\pdel}{\pdel x_i}\quad
\text{and}\quad p_{ji}(x) = -p_{ij}(x)\;. \] 
Again, $\pi$ is a Poisson structure if and only if the
Schouten bracket $\ds  \SbtS{\pi}{\pi} $ vanishes or equivalently
satisfy (\ref{bkt:jacobi}).  
We say that a Poisson structure $\pi$ on $\ds\mR^n$ is
$h$-\textbf{homogeneous} when all the coefficients $p_{ij}(x) =
\Bkt{x_i}{x_j}$ are homogeneous polynomials of degree $h$ in 
$\ds x_1,\ldots, x_n$.  
\end{defn} 
The space of polynomials $\ds\mR[x_1,\ldots, x_n]$ is a Lie sub-algebra
with respect to the Poisson bracket defined by a ($h$-)homogeneous
Poisson structure and is a quotient Lie algebra of $\ds\mR[x_1,\ldots,
x_n]$ modulo Casimir polynomials as mentioned above.  
Thus, we may consider the two kinds of Lie algebra cohomology groups.

\begin{kmRemark}\label{rem:exactVShomog}
        There is a notion of
        subclass of Poisson structures, called {\em exact} Poisson
        structures or 
        {\em homogeneous} Poisson structures on a manifold $M$
        satisfying the condition that the Poisson tensor $\pi$ has a
        vector field $X$ with $\ds \SbtS{X}{\pi}= \pi$.  
As been discussed in \cite{Mik:Miz:Sat}, $h$-homogeneous Poisson
        structure on $\ds\mR^n$ in our sense with $h \ne 2$ is automatically 
        {\em exact} or {\em homogeneous}, but when $h=2$ there are some
        examples 2-homogeneous in our sense but not {\em exact} nor {\em
        homogeneous}.  It will be better to call our $h$-homogeneous
        Poisson structure by {\em naive} homogeneous Poisson structures,
        but sometimes call just homogeneous Poisson structures.  
\end{kmRemark}

\begin{kmRemark} \label{rem:LiePoisson}
It is well-known that a $1$-homogeneous Poisson structure is nothing but
a Lie Poisson structure and the space  $\mR^n$  is the dual space of a
Lie algebra and the Poisson bracket is the Lie algebra bracket defined 
on the space of linear functions on it, that is, the Lie algebra itself. 
In precise, let  $\ds\frakh$ be
        a finite dim Lie algebra.    
For $F:\frakh^{*}\rightarrow \mR$, 
define $\ds \frac{\delta F}{\delta \mu} \in \frakh^{**}$ by 
$\ds \frac{\delta F}{\delta \mu}(\nu) = \frac{d}{dt} F(\mu+t\nu)_{\mid t=0}$
and $\ds \Pkt{F}{H}(\mu) := \langle [ \frac{\delta F}{\delta \mu} ,
        \frac{\delta H}{\delta \mu}], \mu\rangle$.  In particular,
        $\Pkt{v}{w} = [v,w]$  for $v, w \in \frakh$.   
        \kmcomment{
Lie Poisson structures are 1-homog Poisson structures and conversely any
        1-homog Poisson structure is just some Lie Poisson structure.  
} 
\end{kmRemark}

\subsection{Examples} 
On $\ds\mR^n$ take a $h$-homogeneous 2-vector field $\pi$. Then $\pi$ is
written as $\ds\pi = \frac{1}{2}\sum_{i,j} p_{ij}(x) \pdel _i \wedge
\pdel_j$ where $\ds p_{ij}(x)+p_{ji}(x)=0$ and $\ds p_{ij}(x)$ are
$h$-homogeneous.  Poisson condition for $\pi$ is given globally by
(\ref{bkt:jacobi}).  

If $n=3$ then the condition is rather simple and is written as 
\begin{equation} 
p_{12} \left(\frac{\pdel p_{23}}{\pdel x_2} 
               -\frac{\pdel p_{31}}{\pdel x_1}\right)   
+p_{23} \left(\frac{\pdel p_{31}}{\pdel x_3} 
               -\frac{\pdel p_{12}}{\pdel x_2}\right)   
+p_{31} \left(\frac{\pdel p_{12}}{\pdel x_1} 
  -\frac{\pdel p_{23}}{\pdel x_3}\right)  = 0
  \label{eqn:Pcond:3} 
\end{equation} 
for $h$-homogeneous polynomials $\{ p_{12}(x), p_{23}(x), p_{31}(x) \}$. 

%%%%%%%%%%% 
We will try to find some of homogeneous Poisson structures on $\mR^3$.
2-vector fields on $\ds\mR^3$ are classified into 3 types by shape: 
\begin{align*}
\text{(1)}\quad & f \pdel_i \wedge \pdel _j\;,\\
\text{(2)}\quad & \pdel_i \wedge (f \pdel_j + g\pdel_k) \text{ with } fg\ne 0
\text{ and } \{i,j,k\}=\{1,2,3\}\;,\\
\text{(3)}\quad & \sum_{i=1}^3 f_{i} \pdel_{i+1}\wedge \pdel_{i+2} \text{ with }
f_1 f_2 f_3 \ne 0 \text{ and } \pdel_{i+3}= \pdel_{i}\;. 
\end{align*} 
It is obvious 0-homogeneous 2-vector field $\ds\pi_0 = \frac{1}{2}
\sum_{ij}^{n} c_{ij} \pdel_i \wedge \pdel_j $ satisfies the Poisson
condition automatically.  We examine the Poisson condition for $\ds f
\pi_0$.  Our classical computation is the following: now $p_{ij} = f
c_{ij}$ and we check the left-hand-side of (\ref{bkt:jacobi}): Then 
\begin{equation} 
\mathop{\mathfrak{S}}_{ijk} \sum_{\lambda=1}^{n} p_{k\lambda} 
\frac{\pdel p_{ij}}{\pdel x_{\lambda}} = 
\mathop{\mathfrak{S}}_{ijk} \sum_{\lambda=1}^{n} f c_{k\lambda} 
 c_{ij} \frac{\pdel f}{\pdel x_{\lambda}}  
 = f \sum_{\lambda=1}^{n}  \frac{\pdel f}{\pdel x_{\lambda}}  
 \mathop{\ds\mathfrak{S}}_{ijk} 
 c_{ij} c_{k\lambda}\;.  \label{eqn:classic:Pcond}
\end{equation}
When $n=3$, (\ref{eqn:classic:Pcond}) is 0 because 
$\ds 
\mathop{\ds\mathfrak{S}}_{123} c_{12}  c_{3\lambda} =0$ 
for each $\lambda = 1,2,3$, thus $\ds f\pi_0$ is Poisson for any
function $f$.  
\kmcomment{ 
  $\hrulefill$

It is obvious 0-homogeneous 2-vector field $\ds\pi_0 = 
c_{1} \pdel_2 \wedge \pdel_3 
+c_{2} \pdel_3 \wedge \pdel_1 
+c_{3} \pdel_1 \wedge \pdel_2 
$ satisfies the Poisson
condition automatically, and $\ds f \pi_0$ satisfies the Poisson condition for
any function $f$. Our proof is the following:
Since \begin{align*}
        \SbtS{f\pi_0}{f\pi_0} =& 2 f \SbtS{\pi_0}{f} \wedge \pi_0
        \\ \noalign{and so}
        \frac{1}{2}  \SbtS{\SbtS{\SbtS{\SbtS{f\pi_0}{f\pi_0} }{x_1}}{x_2}}{x_3} 
    =& \SbtS{\SbtS{\pi_0}{f}}{x_1} \SbtS{\SbtS{\pi_0}{x_2}}{x_3} 
      + \SbtS{\SbtS{\pi_0}{f}}{x_2} \SbtS{\SbtS{\pi_0}{x_3}}{x_1} 
      \\& 
      + \SbtS{\SbtS{\pi_0}{f}}{x_3} \SbtS{\SbtS{\pi_0}{x_1}}{x_2} 
      \\
      =& \SbtS{\SbtS{\pi_0}{f}}{x_1} ( - c_1) 
      + \SbtS{\SbtS{\pi_0}{f}}{x_2}  ( - c_2) 
      + \SbtS{\SbtS{\pi_0}{f}}{x_3}  ( - c_3) 
      \\ =& 
      \SbtS{ c_1 \SbtS{\pi_0}{x_1}
       + c_2 \SbtS{\pi_0}{x_2} + c_3 \SbtS{\pi_0}{x_3} 
      }{f}
      \\ = & \SbtS{0}{f} = 0 
\end{align*}
}%endOfkmcomment
Thus, 2-vector field of type (1) satisfies the Poisson condition
automatically, so we may choose $h$-homogeneous polynomial $f$.  

About type (2), we may assume $\ds\pdel_1 \wedge ( f \pdel_2 + g
\pdel_3) $.  Then the Poisson condition is 
\begin{equation} 
        \frac{\pdel f}{\pdel x_1} g - 
        f \frac{\pdel g}{\pdel x_1}  = 0\;.  \label{eqn:type:2}
\end{equation} 
If $f$ and $g$ satisfy (\ref{eqn:type:2}) then the commonly multiplied
$\phi f$ and $\phi g$ also satisfy (\ref{eqn:type:2}) because 
\[ (\phi f)' (\phi g) - (\phi f) (\phi g)' =
(\phi ' f + \phi f') (\phi g) - (\phi f) (\phi ' g + \phi g') = 0\]
where the dash of $f$, $f'$ means $\ds\frac{\pdel f}{\pdel x_1}$. 

If $f$ and $g$ are polynomials of the variables $x_2$ and $x_3$, then
(\ref{eqn:type:2}) holds and so $\ds\phi \pdel_1 \wedge ( f \pdel_2 + g
\pdel_3) $ give Poisson structures on $\mR^3$.  So, we 
  have found examples of Poisson structures: 
\[ 
\phi^{[i]}\pdel_1 \wedge ( f^{[j]}(x_2,x_3) \pdel_2 +  g^{[j]} (x_2,x_3)
\pdel_3 ) 
\] 
where $\ds f^{[i]}$, $\ds g^{[i]}$ and $\ds\phi^{[i]}$ mean $i$-th
homogeneous polynomials and $j>0$.  

As for the tensor of type (3), we have two kinds of examples of Poisson structures. 
\begin{quote}
(3-1) $\ds\sum_{i=1}^3 c_i x_i^p \pdel_i  \wedge x_{i+1}^p \pdel_{i+1}$ 
\quad 
(more general, 
$\ds\frac{1}{2} \sum_{i=1}^n c_{ij} x_i^p \pdel_i  \wedge x_{j}^p \pdel_{j}$
where $c_{ij}+c_{ji}=0$).  
\end{quote}
The reason is: 
\begin{align*}
& \SbtS{x_i^p \pdel_i \wedge x_{i+1}^p \pdel_{i+1}} 
        {x_j^p \pdel_j \wedge x_{j+1}^p \pdel_{j+1}} \\
        =& 
        \SbtS{x_i^p \pdel_i} {x_j^p \pdel_j}  
        \wedge x_{i+1}^p \pdel_{i+1} \wedge x_{j+1}^p \pdel_{j+1} 
        - \SbtS{x_i^p \pdel_i} {x_{j+1}^p \pdel_{j+1} }  
        \wedge x_{i+1}^p \pdel_{i+1} \wedge x_{j}^p \pdel_{j} 
        \\& 
        - \SbtS{x_{i+1}^p \pdel_{i+1}} {x_j^p \pdel_j}  
        \wedge x_{i}^p \pdel_{i} \wedge x_{j+1}^p \pdel_{j+1} 
        + \SbtS{x_{i+1}^p \pdel_{i+1}} {x_{j+1}^p \pdel_{j+1}}  
        \wedge x_{i}^p \pdel_{i} \wedge x_{j}^p \pdel_{j} 
        \\
        \noalign{and}
        & 
        \SbtS{x_i^p \pdel_i} {x_j^p \pdel_j} =
        x_i ^{p} p x_j ^{p-1} \delta_{ij} \pdel_j - 
        x_j ^{p} p x_i ^{p-1} \delta_{ij} \pdel_i = \delta_{ij} p (x_i x_j) ^{p-1} 
        ( x_i \pdel_j - x_j \pdel_i) = 0 \;. 
\end{align*}

\begin{quote}
(3-2) $\ds\sum_{i=1}^3 c_{i} x_{i}^{p} \pdel_{i+1}  \wedge \pdel_{i+2}$
where $p$ is a non-negative integer, $\ds c_i$ are constant and $\ds x_{i+3} = x_{i}$. 
\end{quote}
Reason is: 
\begin{align*}
& \SbtS{\sum_{i=1}^3 c_{i} x_{i}^{p} \pdel_{i+1}  \wedge \pdel_{i+2}} 
        {\sum_{j=1}^3 c_{j} x_{j}^{p} \pdel_{j+1}  \wedge \pdel_{j+2}} 
        \\
= & \sum_{i,j} c_{i} c_{j} 
         \SbtS { x_{i}^{p} \pdel_{i+1}  \wedge \pdel_{i+2}} 
        { x_{j}^{p} \pdel_{j+1}  \wedge \pdel_{j+2}} \\
= & \sum_{i,j} c_{i} c_{j} (
         \SbtS { x_{i}^{p} \pdel_{i+1}} {x_{j}^{p} \pdel_{j+1}  }  
        \wedge \pdel_{i+2} \wedge \pdel_{j+2}  
- \SbtS { x_{i}^{p} \pdel_{i+1}}{\pdel_{j+2}}  
        \wedge  \pdel_{i+2} \wedge  x_{j}^{p} \pdel_{j+1} 
        \\&
- \SbtS {  \pdel_{i+2}} { x_{j}^{p} \pdel_{j+1}  }  
        \wedge  x_{i}^{p} \pdel_{i+1} \wedge  \pdel_{j+2}  
+ \SbtS {  \pdel_{i+2}} {  \pdel_{j+2}}  
        \wedge  x_{i}^{p} \pdel_{i+1} 
        \wedge  x_{j}^{p} \pdel_{j+1} 
        )
        \\
= & \sum_{i,j} c_{i} c_{j} (
        ( x_{i}^{p} p  x_{j}^{p-1}  \delta_{i+1}^{j} \pdel_{j+1}    
        - p x_{i}^{p-1} x_{j}^{p}\delta_{j+1}^{i} \pdel_{i+1})     
        \wedge \pdel_{i+2} \wedge \pdel_{j+2}  
        \\& 
+ p x_{i}^{p-1} \delta_{j+2}^{i} \pdel_{i+1} 
        \wedge  \pdel_{i+2} \wedge  x_{j}^{p} \pdel_{j+1} 
-  p x_{j}^{p-1} \delta_{i+2}^{j} \pdel_{j+1}    
        \wedge  x_{i}^{p} \pdel_{i+1} \wedge  \pdel_{j+2}  ) 
        \\
        = & p \sum_{i,j} c_{i} c_{j}  (x_{i}x_{j})^{p-1} (
         x_{i} \delta_{i+1}^{j} \pdel_{j+1}    
        \wedge \pdel_{i+2} \wedge \pdel_{j+2}  
        -  x_{j}\delta_{j+1}^{i} \pdel_{i+1}     
        \wedge \pdel_{i+2} \wedge \pdel_{j+2}  
        \\& 
        +  x_{j} \delta_{j+2}^{i} \pdel_{i+1} 
        \wedge  \pdel_{i+2} \wedge  \pdel_{j+1} 
        -  x_{i} \delta_{i+2}^{j} \pdel_{j+1}    
        \wedge  \pdel_{i+1} \wedge  \pdel_{j+2}  
        ) 
        \\ 
        = & p \sum_{j=i+1} c_{i} c_{j}  (x_{i}x_{j})^{p-1} 
        x_{i}  \pdel_{i+2}    
        \wedge \pdel_{i+2} \wedge \pdel_{i+3}  
        -  p \sum_{j=i-1} c_{i} c_{j}  (x_{i}x_{j})^{p-1} 
        x_{i-1} \pdel_{i+1}     
        \wedge \pdel_{i+2} \wedge \pdel_{i+1}  
        \\& 
        +  p \sum_{j=i-2} c_{i} c_{j}  (x_{i}x_{j})^{p-1} 
        x_{i-2}  \pdel_{i+1} 
        \wedge  \pdel_{i+2} \wedge  \pdel_{i-1} 
        -  p \sum_{j=i+2} c_{i} c_{j}  (x_{i}x_{j})^{p-1} 
        x_{i} \pdel_{i+3}    
        \wedge  \pdel_{i+1} \wedge  \pdel_{i+4}  
        \\ 
        =& 0\;. 
\end{align*}
We show a direct computation to find Poisson structures for $h=1$.  
Take a general 1-homogeneous 2-vector field 
\[\pi = 
( c_1 x_1 + c_2 x_2 + c_3 x_3 ) \pdel_1 \wedge \pdel_2 
+( c_4 x_1 + c_5 x_2 + c_6 x_3 ) \pdel_2 \wedge \pdel_3 
+( c_7 x_1 + c_8 x_2 + c_9 x_3 ) \pdel_3 \wedge \pdel_1 \]  
where $\ds c_i$ are constant. 
Then the Poisson condition consists of 3 quadratic equations:  
\begin{align}
        & c_1 c_5 -c_2 c_4 + c_4 c_9 -c_6 c_7 =0\;,  \label{p:one}\\
        & c_1 c_8 -c_2 c_7 + c_5 c_9 -c_6 c_8 =0\;,  \label{p:two}\\
        & c_1 c_9 -c_2 c_6 + c_3 c_5 -c_3 c_7 =0\;.  \label{p:thr}
\end{align}
Solving the equations by the symbolic calculator Maple, we get 10 solutions:
\input mysoll.tex

If we transform some of variables $\{c_i\}$ by $\{ u_j\}$ as follows:
\[
c_{1}-c_{6} = u_{1}, c_{1}+c_{6} = u_{6}, c_{2}-c_{9} = u_{2}, c_{2}+c_{9} =
u_{9}, c_{5}-c_{7} = u_{5}, c_{5}+c_{7} = u_{7}\;. 
\]
The the Poisson condition (\ref{p:one}) -- (\ref{p:thr}) becomes 
\[ 
        -2 c_{4} u_{2}+u_{1} u_{7}+u_{6} u_{5} =0\;,\quad 
        2 c_{8} u_{1}+u_{5} u_{9}-u_{2} u_{7} = 0\;, \quad 
        2  c_{3} u_{5}-u_{2} u_{6}+u_{9} u_{1} =0\;. 
\]
Solving these equations, we get 4 solutions: 
\begin{align*}
        [1] \quad         u_{6} = & (2 c_{4} u_{2}-u_{1} u
		_{7})/u_{5}\;, 
u_{9} = (-2 c_{8} u_{1}+u_{2} u_{7})/u_{5}\;, 
c_{3} = (c_{4} u_{2}^2-u_{2} u_{1} u_{7}+c_{8} u_{1}^2)/u_{5}^2\, 
\\
\noalign{
where 
\( u_{1}, u_{2},  u_{5},  u_{7}, c_{4},  c_{8}\) are free. } 
[2] \quad 
u_{5} =& 0, 
        u_{7} = 2 c_{4} u_{2}/u_{1}, u_{9} = u_{2} u_{ 6}/u_{1}, 
c_{8} = u_{2}^2 c_{4}/u_{1}^2\quad 
\text{
where}\quad  
u_{1},  u_{2}, u_{6}, c_{3} , c_{4}\quad \text{are free.} \\ 
[3] \quad  
u_{1} = & 0, u_{2} = 0, u_{5} = 0\quad 
\text{where}\quad  
u_{6} ,  u_{7}, u_{9}, c_{3}, c_{4} , c_{8}\quad\text{are free.}\\ 
[4] \quad 
u_{1} =& 0, 
u_{5} = 0, u_{6} = 0, u_{7} = 0, 
c_{4} = 0\quad \text{where} \quad  
        u_{2} , u_{9} , c_{3} , c_{8} \quad \text{are free.} 
\end{align*}

% \young(123,455,6\null5)
% \young(\null\null\null,\null,\null)

%%%%%%%%%%%                      
If $n=4$, then the condition becomes more complicated as we get 4 equations 
from (\ref{bkt:jacobi}).  
\kmcomment{
If $n=4$, then the condition becomes more complicated as follows: 
\begin{align*} 
  p_{{12}}{\frac {\pdel}{\pdel x_{{1}}}}  p_{{13}} 
  & 
 +p_{{12}}{\frac {\pdel}{\pdel x_{{2}}}} p_{{23}} 
 -p_{{13}}{\frac {\pdel}{\pdel x_{{1}}}} p_{{12}}
 +p_{{13}}{\frac {\pdel}{\pdel x_{{3}}}} p_{{23}}
 +p_{{14}}{\frac {\pdel}{\pdel x_{{4}}}} p_{{23}}\\& 
 -p_{{23}}{\frac {\pdel}{\pdel x_{{2}}}} p_{{12}} 
 -p_{{23}}{\frac {\pdel}{\pdel x_{{3}}}} p_{{13}} 
 -p_{{24}}{\frac {\pdel}{\pdel x_{{4}}}} p_{{13}}
 +p_{{34}}{\frac {\pdel}{\pdel x_{{4}}}} p_{{12}} = 0 \\ 
 p_{{12}}{\frac {\pdel}{\pdel x_{{1}}}} p_{{14}} 
 & 
+p_{{12}}{\frac {\pdel}{\pdel x_{{2}}}} p_{{24}} 
+p_{{13}}{\frac {\pdel}{\pdel x_{{3}}}} p_{{24}}
-p_{{14}}{\frac {\pdel}{\pdel x_{{1}}}} p_{{12}}
+p_{{14}}{\frac {\pdel}{\pdel x_{{4}}}} p_{{24}}
\\ &
-p_{{23}}{\frac {\pdel}{\pdel x_{{3}}}} p_{{14}}  
-p_{{24}}{\frac {\pdel}{\pdel x_{{2}}}} p_{{12}}
-p_{{24}}{\frac {\pdel}{\pdel x_{{4}}}} p_{{14}} 
-p_{{34}}{\frac {\pdel}{\pdel x_{{3}}}} p_{{12}}= 0   \\ 
 p_{{12}}{\frac {\pdel}{\pdel x_{{2}}}} p_{{34}}
 & 
+p_{{13}}{\frac {\pdel}{\pdel x_{{1}}}} p_{{14}} 
+p_{{13}}{\frac {\pdel}{\pdel x_{{3}}}} p_{{34}} 
-p_{{14}}{\frac {\pdel}{\pdel x_{{1}}}} p_{{13}}
+p_{{14}}{\frac {\pdel}{\pdel x_{{4}}}} p_{{34}}\\&
+p_{{23}}{\frac {\pdel}{\pdel x_{{2}}}} p_{{14}} 
-p_{{24}}{\frac {\pdel}{\pdel x_{{2}}}} p_{{13}}
-p_{{34}}{\frac {\pdel}{\pdel x_{{3}}}} p_{{13}}
-p_{{34}}{\frac {\pdel}{\pdel x_{{4}}}} p_{{14}} = 0 \\ 
-p_{{12}}{\frac {\pdel}{\pdel x_{{1}}}} p_{{34}} 
& 
+p_{{13}}{\frac {\pdel}{\pdel x_{{1}}}} p_{{24}} 
-p_{{14}}{\frac {\pdel}{\pdel x_{{1}}}} p_{{23}}
+p_{{23}}{\frac {\pdel}{\pdel x_{{2}}}} p_{{24}} 
+p_{{23}}{\frac {\pdel}{\pdel x_{{3}}}} p_{{34}}\\ &  
-p_{{24}}{\frac {\pdel}{\pdel x_{{2}}}} p_{{23}}
+p_{{24}}{\frac {\pdel}{\pdel x_{{4}}}} p_{{34}} 
-p_{{34}}{\frac {\pdel}{\pdel x_{{3}}}} p_{{23}}
-p_{{34}}{\frac {\pdel}{\pdel x_{{4}}}} p_{{24}} = 0   \end{align*}
}%endOFkmcomment 

Any constant 2-vector field $\ds\pi_0 = \frac{1}{2} \sum_{ij}^ 4 c_{ij}
\pdel_i \wedge \pdel_j$ is Poisson, and about $f\pi_0$ if $\ds
\pi_0\wedge \pi_0 \ne 0$, i.e., symplectic then $\ds f\pi_0$ is Poisson
only for constant function $f$ and if $\ds\pi_0\wedge \pi_0 =0$, i.e.,
rank is 2 then $\ds f\pi_0$ is Poisson for any $f$ from
(\ref{eqn:classic:Pcond}).  Here, we used the relations $\ds
\mathop{\ds\mathfrak{S}}_{ijk} c_{ij}c_{k\ell} =0$ if $\ell\in
\{i,j,k\}$ and $\ds\left(\mathop{\ds\mathfrak{S}}_{ijk}
c_{ij}c_{k\ell}\right) \pdel_1\wedge \cdots \wedge\pdel_4  = \pm \pi_0
\wedge \pi_0$ if $\{i,j,k,\ell\}=\{1,2,3,4\}$.  
Thus, we see different situations even for 3 or 4 dimensional.

%% file: mysoll.tex
%% Created by Maple 16.02, Linux
%% Source Worksheet: Untitled (1)
%% Generated: Thu Oct 15 16:58:59 JST 2015
%\documentclass[14pt]{jsarticle}
%\usepackage{amsmath}
%\def\emptyline{\vspace{12pt}}
%\begin{document}
\begin{align*}           
        c_{{1}}=& {\frac {c_{{2}}c_{{7}}-c_{{5}}c_{{9}}+c_{{8}}c_{{6}}}{c_{{8}}}},
c_{{2}}= c_{{2}},
c_{{3}}=  {\frac {-c_{{9}}c_{{2}}c_{{7}}+c_{{5}}{c_{{9}}}^{2}-c_{{9}}c_{{8}}c_{{6}}+c_{{6}}c_{{2}}c_{{8}}}{c_{{8}} \left( c_{{5}}-c_{{7}} \right) }},
\\
& 
c_{{4}}= {\frac {c_{{5}}c_{{2}}c_{{7}}-{c_{{5}}}^{2}c_{{9}}+c_{{5}}c_{{8}}c_{{6}}
-c_{{7}}c_{{6}}c_{{8}}}{c_{{8}} \left( c_{{2}}-c_{{9}} \right) }},
c_{{5}}= c_{{5}},
c_{{6}}= c_{{6}},
c_{{7}}= c_{{7}},
c_{{8}}= c_{{8}},
c_{{9}}= c_{{9}}
\\
%%% 
        c_{{1}}=& {\frac {c_{{9}}c_{{7}}}{c_{{8}}}},
        c_{{2}}= c_{{9}},
        c_{{3}}={\frac {{c_{{9}}}^{2}}{c_{{8}}}},
        c_{{4}}=c_{{4}},
c_{{5}}= c_{{5}},
                c_{{6}}= {\frac {c_{{5}}c_{{9}}}{c_{{8}}}},
                c_{{7}}=c_{{7}},
                c_{{8}}= c_{{8}},
                c_{{9}}= c_{{9}} 
\\ 
        c_{{1}}=& c_{{6}},
c_{{2}}=c_{{9}},
c_{{3}}=c_{{3}},
c_{{4}}=c_{{4}},
c_{{5}}=c_{{5}},
c_{{6}}=c_{{6}},
c_{{7}}=c_{{5}},
c_{{8}}=c_{{8}},
c_{{9}}=c_{{9}} 
\\ 
        c_{{1}}= & {\frac {c_{{5}}c_{{2}}}{c_{{8}}}},
c_{{2}}=c_{{2}},
c_{{3}}=c_{{3}},
c_{{4}}={\frac {{c_{{5}}}^{2}}{c_{{8}}}},
c_{{5}}=c_{{5}},
c_{{6}}={\frac {c_{{5}}c_{{9}}}{c_{{8}}}},
c_{{7}}=c_{{5}},
c_{{8}}=c_{{8}},
c_{{9}}=c_{{9}} 
\\
        c_{{1}}=& c_{{1}},
c_{{2}}={\frac {c_{{5}}c_{{9}}}{c_{{7}}}},
c_{{3}}={\frac {c_{{9}} \left( c_{{1}}c_{{5}}c_{{7}}-c_{{4}}c_{{5}}c_{{9}}
                +c_{{1}}{c_{{7}}}^{2} \right)}{{c_{{7}}}^{3}}},
c_{{4}}=c_{{4}},
c_{{5}}=c_{{5}},
\\
& 
c_{{6}}={\frac {c_{{1}}c_{{5}}c_{{7}}
                -c_{{4}}c_{{5}}c_{{9}} +c_{{4}}c_{{9}}c_{{7}}}{{c_{{7}}}^{2}}},
c_{{7}}=c_{{7}},
c_{{8}}=0,
c_{{9}}=c_{{9}}
\\
        c_{{1}}=& c_{{1}},
c_{{2}}=c_{{9}},
c_{{3}}=c_{{3}},
c_{{4}}=c_{{4}},
c_{{5}}=c_{{7}},
c_{{6}}=c_{{1}},
c_{{7}}=c_{{7}},
c_{{8}}=0,
c_{{9}}=c_{{9}} 
\\
c_{{1}}=& {\frac {c_{{4}}c_{{2}}}{c_{{5}}}},
c_{{2}}=c_{{2}},
c_{{3}}={\frac {c_{{6}}c_{{2}}}{c_{{5}}}},
c_{{4}}=c_{{4}},
c_{{5}}=c_{{5}},
c_{{6}}=c_{{6}},
c_{{7}}=0,
c_{{8}}=0,
c_{{9}}=0
\\ 
        c_{{1}}=& {\frac {c_{{6}}c_{{2}}}{c_{{9}}}},
c_{{2}}=c_{{2}},
c_{{3}}=c_{{3}},
c_{{4}}=0,
c_{{5}}=0,
c_{{6}}=c_{{6}},
c_{{7}}=0,
c_{{8}}=0,
c_{{9}}=c_{{9}} 
\\
        c_{{1}}=& c_{{1}},
c_{{2}}=c_{{2}},
c_{{3}}=c_{{3}},
c_{{4}}=0,
c_{{5}}=0,
c_{{6}}=0,
c_{{7}}=0,
c_{{8}}=0,
c_{{9}}=0 \\
        c_{{1}}=& c_{{1}},
c_{{2}}=0,
c_{{3}}=c_{{3}},
c_{{4}}=c_{{4}},
c_{{5}}=0,
c_{{6}}=c_{{6}},
c_{{7}}=0,
c_{{8}}=0,
c_{{9}}=0 
\end{align*}

%% file: Include-0v1.tex
\section{Polynomial algebra including constants} 
When a Poisson structure $\ds \pi$ is given on a manifold $M$, we have a Lie
algebra $\ds (C^{\infty}(M), \Bkt{\cdot}{\cdot})$ and we may consider
Lie algebra cohomology groups in ``primitive'' sense as follows: the
coboundary operator is a derivation of degree $+1$ and for
each ``$k$-cochain'' $\sigma$, $\ds d \sigma$ is defined by  
\begin{align*}
        (d \sigma) (f_0,\ldots,f_k) &= \sum_{i<j} (-1)^{i+j}\sigma
        (\Bkt{f_i}{f_j},\ldots, \widehat{f_i}, \ldots,  \widehat{f_j},
        \ldots, f_k)\;,  
\end{align*}
where $\ds \widehat{f_i}$ means omitting $\ds f_i$.  
Or, the boundary operator is given by 
\begin{align*}
        \pdel \left( \&^{\wedge} (f_1,\ldots,f_k)\right) &
        = \sum_{i<j} (-1)^{i+j}
        \&^{\wedge}(\Bkt{f_i}{f_j},\ldots,\widehat{f_i},\ldots,\widehat{f_j},
        \ldots, f_k)\;.  
\end{align*}

But, it seems hard to handle ``primitive (co)chain spaces'' because
the (co)chain complexes are
huge spaces. Instead, in the following subsections we consider only 
        homogeneous Poisson structures and 
polynomial functions and reduce the spaces by the notion of ``weight''.
Also, we study about Lie algebra (co)homology groups of formal Hamiltonian
vector fields of the homogeneous Poisson structures.   
\kmcomment{ 
In Gel'fand-Fuks cohomology theory of the formal Hamiltonian vector
fields on the homogeneous Poisson manifold $\ds\mR^n$,
 we may ignore the constant
polynomials from the original algebra 
because the constant polynomials are central elements with respect the
Poisson bracket (called Casimir polynomials).  
        \kmcomment{
\(\ds\ovfrakg = \sum_{j=1}^{\infty} \mySS{j} \) where 
$\ds \frakSS{j}$ is the dual space of the space of $j$-th homogeneous
        polynomials. }

In this note,}%endOFkmcomment
        We first handle the algebra including constant polynomials, i.e.,  
\begin{equation} % \operatorname{Poly}(\mR^n) 
\polygon{n}
=  \sum^{\infty}_{j=\color{red}{0}}
\mySS{j} \; 
\text{, and the dual space is}\; \polygon{n}^{*}=  
\sum^{\infty}_{j=\color{red}{0}} \frakSS{j} \; 
\label{eqn:with0}
\end{equation} 
where $\ds\mySS{j}$ is the space of $j$-th homogeneous
        polynomials and  $\ds \frakSS{j}$ is the dual space of $\ds\mySS{j}$. 

The dimension of $\frakSS{j}$ or $\mySS{j}$ is $\ds \tbinom{n-1+j}{j}$,
and 
$\ds \dim \frakSS{0}=\dim \mySS{0} = 1$ in particular.  The $m$-th
cochain space or chain space 
of the Lie algebra $\ds\polygon{n}$ by the Poisson bracket
        with polynomials in \eqref{eqn:with0} are given by  
\begin{align*} \myCO{m}{}{\polygon{n}^{*}} &= 
\sum_{
{\color{red}0}\leq
i_1\leq \cdots\leq i_m} \frakSS{i_1} \wedge \cdots \wedge \frakSS{i_m} 
= \sum 
\Lambda^{k_0}\frakSS{0} \otimes 
\Lambda^{k_1}\frakSS{1} \otimes \cdots \otimes 
\Lambda^{k_{\ell}}\frakSS{\ell} 
\\
        \myCO{m}{}{\polygon{n}^{ }} &= 
\sum_{
{\color{red}0}\leq
i_1\leq \cdots\leq i_m} \mySS{i_1} \wedge \cdots \wedge \mySS{i_m} 
= \sum 
\Lambda^{k_0}\mySS{0} \otimes 
\Lambda^{k_1}\mySS{1} \otimes \cdots \otimes 
\Lambda^{k_{\ell}}\mySS{\ell} 
\end{align*}
where $\ds \sum_{j={\color{red}0}}^{\infty} k_j = m$.  
When $\ds k_j =0$ for some $j$, then $\ds\Lambda^{k_j}\frakSS{j} = \mR$
and  $\ds \Lambda^{k_j}\frakSS{j} \otimes
\Lambda^{k_{\ell}}\frakSS{\ell} = \Lambda^{0}\frakSS{j} \otimes
\Lambda^{k_{\ell}}\frakSS{\ell} = \Lambda^{k_{\ell}}\frakSS{\ell}$ and   
$\ds \Lambda^{k_j}\mySS{j} \otimes
\Lambda^{k_{\ell}}\mySS{\ell} = 
 \Lambda^{k_{\ell}}\mySS{\ell}$ in the similar way.    
When $m=0$, we have \(\ds \Lambda^0(\polygon{n}^{*}) = \mR\)
and \(\ds \Lambda^0(\polygon{n}^{ }) = \mR\), 
where $\ds\mR$ is the (coefficient)
field of algebra  $\ds\polygon{n}$.  As a basis of $\ds\mySS{j}$ we have
$\ds\{ \w{A} \mid A=(a_1,\ldots,a_n) \in\mN^n\;\text{ and }\; |A|=
a_1+\cdots+a_n=j\}$, and the dual basis $\ds\{\z{A}\}$ is given by 
\[ \z{A} = 
\myUnit \circ (\frac{1}{A!} \frac{\pdel}{\pdel x_1}^{a_1} \cdots 
        \frac{\pdel}{\pdel x_n}^{a_n})
       =  
\myUnit \circ (\frac{1}{A!} \pdel ^{A} )
        \quad\text{with}\quad |A|=j
\]
where $\ds A! = a_1 !\cdots a_n !$ and $\ds \myUnit$ is the
Dirac delta function, which evaluates the target at the origin. Thus, $\ds
\z{[0,\ldots,0]}= 
        \myUnit\in\frakSS{0}$ is the dual basis of $\ds
        1\in\mySS{0}=\mR=$ \{constant polynomials\}. 

\subsection{Weight decomposition}
\begin{defn}
We introduce the notion of weight associated with the $h$-homogeneous
Poisson structure on $\ds\mR^n$ as follows: 
\begin{align}
        \text{the weight of non-zero element of }\; & 
        \Lambda^{k_0}\frakSS{0} \otimes 
\Lambda^{k_1}\frakSS{1} \otimes \cdots \otimes 
\Lambda^{k_{\ell}}\frakSS{\ell} \; \text{ is }\; 
\sum_{j={\color{red}0}}^{\ell}  k_{j} ( j-2+h)\; , \\
        \noalign{and}
        \text{the weight of non-zero element of }\; & 
        \Lambda^{k_0}\mySS{0} \otimes 
\Lambda^{k_1}\mySS{1} \otimes \cdots \otimes 
\Lambda^{k_{\ell}}\mySS{\ell} \; \text{ is }\; 
\sum_{j={\color{red}0}}^{\ell}  k_{j} ( j-2+h)\; .  
        \end{align}
\end{defn}
Thus, the $m$-th cochain and chain space of the weight $w$ are given by  
\begin{align} 
        \myCO{m}{w} {\polygon{n}^{*}} &= \sum \Lambda^{k_0}\frakSS{0} \otimes 
\Lambda^{k_1}\frakSS{1} \otimes \cdots \otimes 
\Lambda^{k_{\ell}}\frakSS{\ell}  \label{eqn:one:ato}
\\ 
        \myCO{m}{w} {\polygon{n}^{ }} &= \sum \Lambda^{k_0}\mySS{0} \otimes 
\Lambda^{k_1}\mySS{1} \otimes \cdots \otimes 
\Lambda^{k_{\ell}}\mySS{\ell}  
\end{align}
with three conditions: 
\begin{align}
 & \sum_{j=0}^{\infty} k_j = m\;, \label{cond:deg}\\
 & \sum_{j=0}^{\infty} k_j (j-2+h) = w\;, \label{cond:wt}\\
 \noalign{and}
 & 0\leq k_j \leq \dim(\frakSS{j}) = \dim(\mySS{j}) = \tbinom{n-1+j}{j}
 \; \text{ for }\; j=0,1,\ldots\;. \label{cond:dim}
\end{align}
\kmcomment{
\eqref{cond:dim} yields $\ds k_0=0$ or $1$. Thus, we can decompose 
\eqref{eqn:one:ato} as follows: 
\begin{equation} \myCO{m}{w}  {\polygon{n}^{*}}
        = \sum_{k_0=0}  
\Lambda^{k_1}\frakSS{1} \otimes \cdots \otimes 
\Lambda^{k_{\ell}}\frakSS{\ell}  
\; \oplus\;
        \sum_{\kk{0}=1}  \frakSS{0} \otimes 
        \Lambda^{\kk{1}}\frakSS{1} \otimes \cdots \otimes 
\Lambda^{\kk{\ell'}}\frakSS{\ell'} 
\label{eqn:two}
\end{equation}
The first term of the direct sum is $\ds \myCC{m}{w}$ we already handle. 
About the second term, 
\eqref{cond:deg} implies 
$\ds  \sum_{j=1}^{\infty} \kk{j} = m -1 $, and  \eqref{cond:wt} implies  
$\ds \sum_{j=1}^{\infty} \kk{j} (j-2+h) = w + (2-h)$, thus,   
the second term coincides with $\ds \frakSS{0} \otimes
\myCC{m-1}{w+2-h}$, namely  
 \begin{equation}
         \Lambda ^m_{w}  (\polygon{n}) 
         =  \myCC{m}{w}\; 
\oplus\; \frakSS{0} \otimes 
 \myCC{m-1}{w+2-h} \;.
\label{eqn:three} 
\end{equation} 
}%endOFkmcomment

\subsection{Possibility of weight and range of degree}
\label{subsec:possible}
Using \eqref{cond:deg} and \eqref{cond:wt}, $m=0$ implies $w=0$, in
other words, $\ds \myCO{0}{w}{\polygon{n}^{*}} = 0$ if $w\ne 0$ and 
$\ds \myCO{0}{0}{\polygon{n}^{*}} = \mR$. 
\begin{align}
        \noalign{
We prepare  two equations, one is 
        $(2-h)$ times \eqref{cond:deg} :} 
& (2-h) \sum_{j=0}^{\infty} k_j  = (2-h)m  \label{jou:one}\\
        \noalign{the other is a reform of \eqref{cond:wt}:}
        & \sum_{j=0}^{\infty} jk_{j} = w + (2-h)m \label{jou:wt}\\
        \noalign{$\eqref{jou:wt}-\eqref{jou:one}$ implies} 
        & \sum_{j<2-h} (j-2+h)k_j + 
        \sum_{j>2-h} (j-2+h)k_j = w \;. \label{jou:three}
\end{align}
In \eqref{jou:three}, if $ h-2\geq 0$ then it turns out $w \geq 0$. 
When $h-2 < 0$, \eqref{jou:three} implies 
\[ 0 \leq 
        \sum_{j>2-h} (j-2+h)k_j = w 
        - \sum_{j<2-h} (j-2+h)k_j =\begin{cases} 
                w + k_{0} & (h=1) \\[-2mm]
        w + 2k_{0} + k_{1} &(h=0)\end{cases} 
        \]
\eqref{cond:dim} says $\ds 0\leq k_0 \leq 1$ and 
$\ds 0\leq k_1 \leq n$,  and so 
$w \geq -1 $ when $h=1$ and $w \geq -(2 + n)$ when $h=0$. 

Subtracting 
$(3-h)$ times \eqref{cond:deg} from  \eqref{jou:wt}, we have 
\begin{align*} 
        &    \sum_{j<3-h} (j-3+h)k_j + 
\sum_{j>3-h} (j-3+h)k_j = w -m  \\ 
        &    
0 \leq \sum_{j>3-h} (j-3+h)k_j = w -m - \sum_{j<3-h} (j-3+h)k_j  \\
        & 
        m \leq w  + \sum_{j<3-h} (3-h-j)k_j = \begin{cases}
                w & h \geq 3 \\[-2mm]
                w+ k_0 \leq w+1 & h=2 \\[-2mm]
                w+ 2k_0 +k_1\leq w+2+n & h=1 \\[-2mm]
        w+ 3k_0 +2 k_1 + k_2\leq w+ \frac{(n+2)(n+3)}{2} & h=0 \end{cases}\\ 
\end{align*} 

Subtracting 
$(1-h)$ times \eqref{cond:deg} from  \eqref{jou:wt}, we have 
\begin{align*} 
        &    \sum_{j<1-h} (j-1+h)k_j + 
\sum_{j>1-h} (j-1+h)k_j = w +m  \\ 
        &    
0 \leq \sum_{j>1-h} (j-1+h)k_j = w + m - \sum_{j<1-h} (j-1+h)k_j  \\
        & 
        -m \leq w + \sum_{j<1-h} (1-h-j) k_{j} =\begin{cases} 
                w & h\geq 1\\[-2mm]
        w+ k_{0} \leq w+1 & h =0\end{cases}
\end{align*}
We may summarize the discussion above as the following table.
\begin{center}
\begin{tabular}{|c|c|l|}
        \hline 
        homogeneity & weight & range of $m$ of $m$-th cochains \\\hline
$h=0$ & $ w \geq -2-n$ & $\max(0, -(w+1)) \leq m \leq w+(n+2)(n+3)/2$\\
$h=1$ & $ w \geq -1$ & $\max(0, -w)) \leq m \leq w+n+2$\\
$h=2$ & $ w \geq 0$ & $0 \leq m \leq w+1$\\
$h\geq 3$ & $ w \geq 0$ & $0 \leq m \leq w$\\\hline
\end{tabular}
\end{center}

\begin{exam}
        When $h=0$ and $w=-(2+n)$, we determine the cochain spaces $\ds
        \myCO{m}{(-2-n)} {\polygon{n}^{*}}$. 
        In this case, \eqref{cond:wt} is 
        \begin{align*}
                -2-n =w =& \sum_{j} k_j(j-2+h) = \sum_{j}k_j(j-2) =
                -2k_{0} -k_{1} +k_{3} + 2 k_{4} + \cdots \;. 
                \\
                2k_{0}-2 +k_{1}-n =& k_{3} + 2 k_{4} + \cdots \;.  
                \end{align*}

Since $\ds 0\leq k_0\leq 1$ and $\ds 0\leq k_1\leq n$,  the last
        equation implies $\ds k_{0}=1$,  $\ds k_{1}=n$  and $k_{j}=0$
        for $j >2$. Applying these facts to \eqref{cond:deg}, we have
        $\ds k_2 = m-(1+n)$. Thus 
        \[ \myCO{m}{(-2-n)} {\polygon{n}^{*}} = \frakSS{0}
        \otimes \Lambda^{n} \frakSS{1}
        \otimes \Lambda^{m-(1+n)} \frakSS{2}\;.\] 
And, $\ds 0\leq m-(1+n) \leq (n+1)n/2$, i.e., $\ds 1+n\leq m
        \leq w+(n+2)(n+3)/2$. 

When $w=-1-n$, we have 
        \begin{align}
                2k_{0}-1 +k_{1}-n = k_{3} + 2 k_{4} + \cdots \;.  
        \end{align}
        Then $\ds k_{0}=1$ and $\ds k_{1}=n$ or $\ds k_{1}=n-1$. 
        When $\ds k_{0}=1$ and $\ds k_{1}=n$, $\ds k_{3}=1,\; k_{j}=0\;
        (j>3)$ and using \eqref{cond:deg}, $\ds k_2 = m-2-n$. 

        When $\ds k_{0}=1$ and $\ds k_{1}=n-1$, $\ds k_j=0\; (j>2)$ and
        $\ds k_2 = m-n$ again using \eqref{cond:deg}.  
        \[ \myCO{m}{(-1-n)} {\polygon{n}^{*}} = \frakSS{0} \otimes \left( 
        \Lambda^{n} \frakSS{1}
        \otimes \Lambda^{m-(2+n)} \frakSS{2} \otimes \frakSS{3} \;
        \oplus\; + 
        \Lambda^{n-1} \frakSS{1}
        \otimes \Lambda^{m-n} \frakSS{2} \right)
        \;.\] 

When $w=-n$, we have 
        \begin{align}
                2k_{0} +k_{1}-n = k_{3} + 2 k_{4} + \cdots \;.  
        \end{align}
        If $\ds k_{0}=0$ then $\ds k_{1}=n, \; k_{j}=0\; (j>2)$ and so
        \eqref{cond:deg} says $\ds k_2= m-n$. 

        If $\ds k_{0}=1$, 
        then $\ds k_{1}=n,\; n-1, \text{ or }\; n-2$.  
        If $\ds k_{0}=1,\; k_{1}=n-2$, then 
        $k_{j}=0\;
        (j>2)$ 
        and using \eqref{cond:deg}, $\ds k_2 = m-n+1$.   
        If $\ds k_{0}=1,\; k_{1}=n-1$, then 
        $k_{3}=1\; k_{j}=0\; (j>3)$ 
        and using \eqref{cond:deg}, $\ds k_2 = m-n-1$.   
        If $\ds k_{0}=1,\; k_{1}=n$, then 
        $\; k_{3}=2,\; k_{j}=0\;
        (j>3)$ 
        \kmcomment{ and using \eqref{cond:deg}, $\ds k_2 = m-n-3$ }
        or 
        $\; k_{3}=0,\; k_{4}=1\; k_{j}=0\;
        (j>4)$.  \kmcomment{and using \eqref{cond:deg}, $\ds k_2 =
        m-n-2$.} 
        Using \eqref{cond:deg}, $\ds k_{2}= m-n-3$ in the former case
        and $\ds k_{2}= m-n-2$ in the last case. 
Thus, we have 
        \begin{align*} \myCO{m}{(-n)} {\polygon{n}^{*}} = & 
        \Lambda^{n} \frakSS{1}
                \otimes \Lambda^{m-n} \frakSS{2} \\ & 
\oplus\; \frakSS{0} \otimes ( 
 \Lambda^{n-2} \frakSS{1} \otimes \Lambda^{m-n+1} \frakSS{2}\;\oplus\;  
 \Lambda^{n-1} \frakSS{1} \otimes \Lambda^{m-n-1} \frakSS{2} \otimes
                \frakSS{3} \\& \qquad  \qquad 
                \;\oplus\;  
 \Lambda^{n} \frakSS{1} \otimes \Lambda^{m-n-3} \frakSS{2} \otimes
                \Lambda^{2}\frakSS{3} \;\oplus\;  
 \Lambda^{n} \frakSS{1} \otimes \Lambda^{m-n-2} \frakSS{2} \otimes
                \frakSS{4} )
        \;.\end{align*} 
\end{exam}

\begin{exam}
        When $h=1$ and $w=-1$, we determine the cochain spaces $\ds
        \myCO{m}{(-1)} {\polygon{n}^{*}}$. 
        In this case, \eqref{cond:wt} is 
        \begin{align*}
                -1 =w =& \sum_{j} k_j(j-2+h) = \sum_{j}k_j(j-1) =
                -k_{0}+k_{2} + 2 k_{3} + \cdots \;. 
                \end{align*}
   This implies $\ds k_{0}=1$ and $k_{j}=0$ for $j\geq 2$. Applying these
        facts to \eqref{cond:deg}, we have $\ds m = 1+k_1$, i.e., $\ds
        k_1 = m-1$. Thus 
        $\ds \myCO{m}{(-1)} {\polygon{n}^{*}} = \frakSS{0}\otimes
        \Lambda^{m-1} \frakSS{1}$. Thus, $\ds 0\leq m-1 \leq n$, i.e.,
        $\ds 1\leq m\leq w+n+2$. 

By the same way, we have 
        \begin{align*} 
                \myCO{m}{0} {\polygon{n}^{*}} =& 
        \Lambda^{m} \frakSS{1}\; \oplus \;
        \frakSS{0}\otimes
        \Lambda^{m-2} \frakSS{1} \otimes \frakSS{2}
\\
                \myCO{m}{1} {\polygon{n}^{*}} = &
                \Lambda^{m-1} \frakSS{1} \otimes \frakSS{2} \; \oplus \;
                \frakSS{0}\otimes \left(
        \Lambda^{m-2} \frakSS{3} \;\oplus 
        \Lambda^{m-3} \frakSS{1} \otimes 
                \Lambda^{2} \frakSS{2} \right) 
        \end{align*} 
\end{exam}

\subsection{About $\ds \mydov(\myUnit)$ } 
Since the coboundary operator $\mydov$ is a derivation of degree 1, it is known
that the behavior of $\mydov$ is under controlled by the action of
$\mydov$ for each 1-cochain $\sigma$. 

From the definition of coboundary operator of Lie algebra cohomology, we see that  
\begin{align} 
        \mydov \sigma =&  
        \frac{1}{2} \sum (\mydov \sigma ) (\w{A},\w{B}) \z{A} \wedge \z{B} 
= - \frac{1}{2} \sum \langle \sigma , \Bkt{\w{A}}{\w{B}}\rangle \z{A}
        \wedge \z{B}\;,  
        \label{eqn:one:cochain} \\
        \noalign{and, in particular} 
        (\mydov\myUnit)(\w{A},\w{B}) = & - \langle \myUnit , \Bkt{\w{A}}{\w{B}}\rangle
=0  \quad\text{if}\quad |A| |B| = 0 \; \text{or}\; |A|+|B|\ne 2-h\;,
\end{align}
where $h$ is the homogeneity of our Poisson structure.  
Thus, 
when $h>0$, if  
$\ds |A||B| \ne 0$ then $\ds |A|+|B| \geq 2 > 2-h$ and  
$\ds \mydov(\myUnit ) = 0$ holds as expected and we have  
$\ds \mydov(\frakSS{0}) = 0$.

Assume our 
Poisson structure $\pi$ is 0-homogeneous,  given by $\ds\Bkt{\w{\eps{i}}} {\w{\eps{j}}} =
p_{i,j}$ (constant). 
Let us denote $\ds \eps{i}$ whose
element is 1 if $i$-th position and 0 otherwise.  
Then 
\begin{align}
        & (\mydov \;\myUnit\;) ( 
\w{\eps{i}}, \w{\eps{i}})  =  - \langle \myUnit,  \Bkt{
        \w{\eps{i}}}{\w{\eps{i}}}  =  - p_{i,j}\;, \label{alig:one}\\
& \mydov \;\myUnit\; = - \frac{1}{2}\sum_{i,j} p_{i,j} \z{\eps{i}} \wedge
\z{\eps{j}} \cong  - \pi\;,
        \label{alig:two} \\
%\noalign{we see directry} 
        & \mydov \pi= 0 \;\text{ from
        \eqref{alig:two}}\; . \label{alig:three}
\end{align}  
\begin{kmRemark} When $h=0$, 
        the weight of $\ds\myUnit$ is $-2$ and the weight of
        $\ds \z{\eps{i}}$ is $-1$ and that of $\ds\pi$ is $-2$.
\end{kmRemark}
\begin{kmRemark}
Preparing formulae of $\ds \mydov  \z{\eps{i}}$, we reconfirm
\eqref{alig:three} by those properties of $\ds \mydov \z{\eps{i}}$: 
\begin{align*}
        \mydov \z{\eps{\ell}} 
        =& - \frac{1}{2} \sum_{A,B} \langle
        \z{\eps{\ell}}, \Bkt{ \w{A}} {\w{B}} \rangle  \z{A}\wedge\z{B} 
        = -  \sum_{|A|=1,|B|=2} \langle
        \z{\eps{\ell}}, \Bkt{ \w{A}} {\w{B}} \rangle \z{A}\wedge\z{B}
        \\ =& -  \sum_{i,j} \langle
        \z{\eps{\ell}}, \Bkt{ \w{\eps{i}}} {\w{2\eps{j}}} \rangle
        \z{\eps{i}}\wedge\z{2\eps{j}}
         - \frac{1}{2} \sum_{i}\sum_{j\ne k} \langle
        \z{\eps{\ell}}, \Bkt{ \w{\eps{i}}} {\w{\eps{j}+\eps{k}}} \rangle  
        \z{\eps{i}}\wedge\z{\eps{j}+\eps{k}} \\
        =& - \frac{1}{2} \sum_{i,j} \langle
        \z{\eps{\ell}}, \Bkt{ \w{\eps{i}}} {\w{2\eps{j}}} \rangle
        \z{\eps{i}}\wedge\z{2\eps{j}}
         - \frac{1}{2} \sum_{i,j,k} \langle
        \z{\eps{\ell}}, \Bkt{ \w{\eps{i}}} {\w{\eps{j}+\eps{k}}} \rangle  
        \z{\eps{i}}\wedge\z{\eps{j}+\eps{k}} 
        \\
        =& - \frac{1}{2} \sum_{i,j} \langle
        \z{\eps{\ell}}, 2 p_{i,j} \w{\eps{j}} \rangle
        \z{\eps{i}}\wedge\z{2\eps{j}}
         - \frac{1}{2} \sum_{i,j,k} \langle
        \z{\eps{\ell}}, 
        p_{i,j}  \w{\eps{k}}+  p_{i,k}  \w{\eps{j}} \rangle  
        \z{\eps{i}}\wedge\z{\eps{j}+\eps{k}} 
        \\
        =& -  \sum_{i,j} 
          p_{i,\ell} \z{\eps{i}}\wedge\z{2\eps{\ell}}
         - \frac{1}{2} \sum_{i,j} p_{i,j}  
        \z{\eps{i}}\wedge\z{\eps{j}+\eps{\ell}} 
         - \frac{1}{2} \sum_{i,k} 
 p_{i,k} 
        \z{\eps{i}}\wedge\z{\eps{\ell}+\eps{k}} \\
        =& Q_{\ell} \wedge \z{2\eps{\ell}} + \sum_{j} Q_{j}\wedge 
        \z{\eps{\ell}+\eps{j}} 
        \; \text{
                where}\;  \ds Q_{j} = \sum_{i} p_{j,i} \z{\eps{i}}\;. \\ 
%\end{align*}
%\begin{align*}
      2 \mydov\pi =& d (\sum_{i,j}p_{i,j} \z{\eps{i}} \wedge \z{\eps{j}}) 
        =  \sum_{i,j}p_{i,j} \left( (d\z{\eps{i}}) \wedge \z{\eps{j}} - 
        \z{\eps{i}} \wedge d\z{\eps{j}}\right) \\
        =& 2 \sum_{j} Q_{j}  \wedge d\z{\eps{j}}  
        = 2 \sum_{j} Q_{j}  \wedge (
        Q_{j} \wedge \z{2\eps{j}} + \sum_{k} Q_{k} \wedge
        \z{\eps{j}+\eps{k}}) 
        \\
        =&  2 \sum_{j,k} Q_{j}  \wedge  Q_{k} \wedge
        \z{\eps{j}+\eps{k}} = 0\;. 
\end{align*}
\end{kmRemark}

\subsection{Cohomology and homology groups with respect to weight}

\begin{kmProp} 
        The coboundary operator $\ds\mydov$ on $\ds\polygon{n}$ preserves the weight, namely, $\ds
        \mydov ( \myCO{m}{w}{\polygon{n}^{*}} ) \subset
        \myCO{m+1}{w}{\polygon{n}^{*}}$ holds.  Thus  we have the
well-defined cohomology group for each weight; 
        \[ \myHH{m}{w}(\polygon{n}) := 
        \text{Ker} ( \mydov : \myCO{m}{w}{\polygon{n}^{*}} 
        \rightarrow \myCO{m+1}{w} {\polygon{n}^{*}} 
        ) / 
        \mydov ( \myCO{m-1}{w}{\polygon{n}}^{*})\;. \] 
\end{kmProp}

\textbf{Proof:} From the linearity of $\mydov$, it is enough only to
check $\ds\mydov(\sigma)  \in \myCO{m+1}{w}{\polygon{n}^{*}}$ for any
generator $\ds \sigma = \sigma_1 \wedge \cdots \wedge \sigma_m \in
\myCO{m}{w}{\polygon{n}^{*}}$ where $\ds\sigma_i \in \frakSS{\phi(i)}$
for $i=1,\ldots,m$.  From the definition, $\ds w = \sum_{i=1}^m
\wt(\sigma_i)$ where we denote the weight of $\ds \sigma_i$
by  $\ds\wt( \sigma_i)$, namely $\ds\wt( \sigma_i) = \phi(i) +h-2 $.  
If
$\ds f\in\frakSS{a}$ and $\ds g\in\frakSS{b}$, then we have \(\ds \Bkt{f}{g} \in
\frakSS{a+b+h-2}\) because the Poisson structure $\pi$ is
$h$-homogeneous.  From (\ref{eqn:one:cochain}) we see that
\(\ds\mydov(\sigma_i) \in \mydov( \frakSS{\phi(i)}) \subset \sum_{a\leq
b, a+b= \phi(i)-h+2} \frakSS{a}\wedge \frakSS{b} \), and we have 
\begin{align*}
\wt( \mydov(\sigma_i)) =& (a+h-2) + (b+h-2) = \phi(i)+h-2 = \wt(\sigma_i)\ .\\  
        \noalign{Thus,} 
& \wt( \sigma_1 \wedge \cdots \wedge \sigma_{i-1} \wedge \mydov( \sigma_i) 
\wedge \sigma_{i+1} \wedge \cdots \wedge \sigma_{m}) \\ 
=  & 
\wt( \sigma_1 ) + \cdots + \wt(\sigma_{i-1}) +  \wt( \mydov( \sigma_i)) +  
\wt( \sigma_{i+1}) + \cdots + \wt(\sigma_{m}) 
\\
=  & 
\wt( \sigma_1 ) + \cdots + \wt(\sigma_{i-1}) +  \wt( \sigma_i) +  
\wt( \sigma_{i+1}) + \cdots + \wt(\sigma_{m}) 
 = \wt(\sigma) = w\end{align*} 
and we conclude \( \wt( \mydov( \sigma )) = \wt(\sigma)\).  \kmqed

Also, the weight is preserved by the boundary operator $\ds\pdel$ and 
we have the following Proposition.
\begin{kmProp} 
        The boundary operator $\ds\pdel $ on $\ds\polygon{n}$ preserves the weight, namely, $\ds
        \pdel ( \myCO{m}{w}{\polygon{n}^{}} ) \subset
        \myCO{m-1}{w}{\polygon{n}^{}}$ holds.  Thus  we have the
well-defined homology group for each weight; 
        \[ \myHH{ }{m,w}(\polygon{n}) := 
        \text{Ker} ( \pdel  : \myCO{m}{w}{\polygon{n}^{ }} 
        \rightarrow \myCO{m-1}{w} {\polygon{n}^{ }} 
        ) / 
        \pdel ( \myCO{m+1}{w}{\polygon{n}}^{ })\;. \] 
\end{kmProp} 

\section{Lie algebra (co)homology of Hamiltonian vector fields of
polynomial potentials} 
For a general Poisson structure $\pi$ on a manifold $M$, the Lie algebra
of Hamiltonian vector fields is identified with $\ds
C^{\infty}(M)/\myCasimm$ where $\ds\myCasimm$ is the space of Casimir
functions. The Lie bracket is given by 
\(\ds [ [f],[g] ] = [\Bkt{f}{g}]\) where $\ds f,g\in C^{\infty}(M)$ and $\ds 
[f] = f+ \myCasimm \in C^{\infty}(M)/\myCasimm$. 
\kmcomment{ 
In ``primitive sense'', the coboundary operator is 
 a derivation of degree $+1$ which is defined
 for ``$k$-cochain'' $\sigma$  by the following way:  
\begin{align*}  (\mydq \sigma)([f_0],\ldots,[f_k]) &= \sum_{i<j}(-1)^{i+j} \sigma( 
[ [f_i],[f_j]], \ldots, \widehat{ [f_i] },\ldots , \widehat{ [f_j]
        },\ldots, [f_k])\;, \\  \noalign{In particular,} 
        (\mydq\sigma)([f],[g]) &=  - \langle \sigma, [ [f],[g]]\rangle = 
        - \langle \sigma, [ \Bkt{f}{g} ]\rangle 
        = - \langle \sigma , \Bkt{f}{g} \rangle
        \\&
        = (\mydov \sigma) (f,g) 
        \; \text{if }\; \langle
\sigma , \myCasimm  \rangle = 0\;. \end{align*} 
More precisely, as mentioned in the previous section, 
}%endOFkmcomment

For a given  homogeneous Poisson
structure $\pi$ of $\ds\polygon{n}$, the
Lie algebra $\frakg$ of Hamiltonian vector fields with polynomial
potentials is identified as $\ds\frakg \cong \polygon{n}/\myCasimm$, 
where $\ds \myCasimm = \{ f\in \polygon{n}  \mid \Bkt{f}{\cdot} =
0\}$ is the space of Casimir polynomials.  
First we prepare a small lemma: 
\begin{Lemma}\label{lemma:one}
        Let $f$ be a Casimir polynomial and let $\ds f_j$ be the
              $j$-homogeneous part of $f$, i.e., $\ds f_j\in \mySS{j}$
              and so
              $\ds f = \sum_j f_j$. Then each $\ds f_j$ are Casimir
              polynomials. (The converse is obviously true.) 
\end{Lemma} 
The Lemma above implies $\ds\myCasim{j} = \{ f_j \mid f\in\myCasimm \}$. 
Since $\ds \mySS{0}\subset\myCasimm$, we have   
\begin{align*}
        \frakg \cong & \polygon{n}/\myCasimm = 
        \sum_{j=0}^{\infty} \mySS{j}/\myCasim{j} = 
        \sum_{j=1}^{\infty} \myS{j} 
        \kmcomment{
                \text{if}\; h=0\;,  \\
        \noalign{and} 
        \frakg  \cong  & \polygon{n}/\myCasimm = 
\sum_{j=1}^{\infty} \mySS{j}/\myCasimm 
= \ovfrakg /\myCasimm\;\text{if}\; h>0\; .}
\; \text{where}\; \myS{j} = \mySS{j}/\myCasim{j}\;.  
\end{align*} 

The dual space  
$\ds\frakg^{*}$ is decomposed as 
$\ds\frakg^{*} = \sum_{j=1}^{\infty} \frakS{j}$ where   
$\ds\frakS{j} = \{\; 
\sigma \in \frakSS{j} \mid \langle \sigma , \myCasimm  
\rangle =0\;\}$.      
\kmcomment{
be the space $\ds\{\; \sigma \in 
\polygon{n} ^{*} \mid  \langle \sigma, \myCasimm  
\rangle = 0\; \}$, where $\langle\, \cdot\,,\, \cdot\, \rangle$ is the 
natural pairing.  We have the decomposition 
$\ds\frakg^{*} = \sum_{k} \frakS{k}$ where   $\ds\frakS{k} = \{\; 
\sigma \in \frakSS{k} \mid \langle \sigma , \myCasimm  
\rangle =0\;\}$ and we can determine the dual space 
$\ds\myS{k}$ of $\ds\frakS{k}$  in $\ds\mySS{k}$.  
We see $\ds \frakS{0}=0$.  }%endOFkmcomment
Thus, the $m$-th cochain or chain spaces with the weight $w$ of $\frakg$
are given by
\begin{align}
        \myC{m}{w} (\frakg) =  \sum 
        \Lambda^{k_1}\frakS{1} \otimes 
        \Lambda^{k_1}\frakS{2} \otimes \cdots \otimes 
        \Lambda^{k_{\ell}}\frakS{\ell} 
        \\
        \myC{}{m,w} (\frakg) =  \sum 
        \Lambda^{k_1}\myS{1} \otimes 
        \Lambda^{k_1}\myS{2} \otimes \cdots \otimes 
        \Lambda^{k_{\ell}}\myS{\ell} 
\end{align} 
with \eqref{cond:deg} and \eqref{cond:wt}. 
About dimensional restriction, 
we only say that \[\ds 0\leq k_j \leqq \dim 
\frakS{j} = \dim\frakSS{j} - \#\{\text{linear independent } 
j\text{-homogeneous Casimir polynomials}\}\;. \]  

From the definition of coboundary operators, we see the next property: 
for each 1-cochain $\sigma$ of $\frakg$,  
\begin{align*} (\mydq \sigma) ([f],[g]) &= 
  - \langle \sigma, [\Bkt{f}{g}]\rangle 
  =  - \langle \sigma, \Bkt{f}{g} + \myCasimm  \rangle 
  \\ & 
  = - \langle \sigma, \Bkt{f}{g} \rangle 
  = (\mydov \sigma) (f,g)\;. 
\end{align*}
Using the derivation property of coboundary operator, we get 
\begin{equation}\label{eqn:the:same} (\mydq \sigma)
([f_0],[f_1],\ldots,[f_k]) = (\mydov \sigma) (f_0,f_1,\ldots,f_k)
\text{ for } \sigma \in \myC{k}{}\;.  \end{equation} In the sense above, 
we may say that $\ds \mydq$ is the restriction of 
$\ds\mydov$. Thus, we have the next Proposition.  

\begin{kmProp}  \label{kmProp::CD} 
$\ds \mydq$  is the natural restriction of $\ds\mydov$ and we have 
the following commutative diagram  
\begin{equation}
        \begin{CD}
\myC{m}{} = \Lambda^m \frakg^{*} @. \subset @.
                 \Lambda^m \left(\polygon{n}^{*}\right)\\
@Vd VV   @. @VV\overline{d}V \\ 
\myC{m+1}{} = \Lambda^{m+1} \frakg^{*} @. \quad \subset \quad @.  
                \Lambda^{m+1} \left(\polygon{n}^{*}\right) \;. 
\end{CD}
\label{Kont:CD}
\end{equation} 
The coboundary operator $\ds\mydq$ preserves the weight, namely, $\ds
\mydq ( \myC{m}{w} ) \subset \myC{m+1}{w}$ holds. 
        Thus  we have the
cohomology group 
        \[\ds \myH{m}{w}(\frakg) := \text{Ker} ( \mydq : \myC{m}{w}  \rightarrow
\myC{m+1}{w}) / \mydq ( \myC{m-1}{w})\;.\]   

\end{kmProp}
\kmcomment{
\textbf{Proof:} 
It is enough to show the proposition is valid only for 1-cochains,
because $\ds\mydov$ and $\ds\mydq$ are both derivations.  Let $\ds\{
\widehat{ \kappa}_i \}$ be a basis of $\ds \myCasimm $ and $\ds\{
\widehat{\kappa}_i , \widehat{\mu}_j \}$ be a basis of $\ds\ovfrakg$,
and let $\ds\{ \kappa_i , \mu_j \}$ be the dual basis.  Then $\ds\{
\mu_j \}$ is a basis of $\ds\frakg^{*}$.  For each 1-cochain $\sigma\in
\myC{\bullet}{}$ 
\begin{align*}
\mydov (\sigma) =& - \frac{1}{2} \langle \sigma, 
\Bkt{\widehat{\kappa}_i}{\widehat{\kappa}_j}\rangle \kappa_i \wedge \kappa_j 
 -  \langle \sigma, 
\Bkt{\widehat{\kappa}_i}{\widehat{\mu}_j}\rangle \kappa_i \wedge \mu_j 
 - \frac{1}{2} \langle \sigma, 
\Bkt{\widehat{\mu}_i}{\widehat{\mu}_j}\rangle \mu_i \wedge \mu_j \\
=& - \frac{1}{2} \langle \sigma, 
\Bkt{\widehat{\mu}_i}{\widehat{\mu}_j}\rangle \mu_i \wedge \mu_j 
= \mydq ( \sigma )
\end{align*} 
since $\ds\Bkt{\widehat{\kappa}_i}{\cdot} =0$ because $
\widehat{\kappa}_i$ are Casimir polynomials. \kmqed
}%endOFkmcomment
As a direct corollary of Proposition \ref{kmProp::CD}, we have
\begin{kmCor}
        $\ds\myH{\bullet}{w}(\frakg) = \left(\myC{\bullet}{w}\cap \ker ( \mydov:
        \myCO{\bullet}{w}{\polygon{n}^{*}} \rightarrow
        \myCO{\bullet+1}{w}{\polygon{n}^{*}} )\right) \big/ \mydq(
        \myC{\bullet-1}{w})$ . 
\end{kmCor} 

\bigskip

\begin{exam}
        In the case $\pi= \pdel_1 \wedge \pdel_2$ on $\ds\mR^3$, 
        i.e., when $h=0$, we
        see the decomposition of $\frakg$.  
        By direct calculation, we know \(\ds \myCasimm\cap \mySS{j} =
        \{\mR x_3{}^j\}\), and 
        \begin{align*}
                & \mySS{0}/\myCasimm = 0\;,\; 
                \mySS{1}/\myCasimm = \text{LSpan}(x_1,x_2)\;,\; 
                \mySS{2}/\myCasimm = \text{LSpan}(x_1{}^2,x_2{}^2, x_1 x_2 ,
                x_1 x_3, x_2 x_3)\;,\;  \\
                & 
                \mySS{j}/\myCasimm \cong \mySS{j} \setminus \mR x_3{}^j
                \; (j=1,2,\ldots)\;, \\ 
                \noalign{and so}
                & 
                \frakS{j} = \text{LSpan} ( \z{A} \mid |A|=j
                \;\text{and}\; a_3\ne j )\;. 
        \end{align*} 
\end{exam}

\subsection{Possibility of weight and range of degree in Hamiltonian
case}
\label{subsec:possib}
By the same argument in \S~\ref{subsec:possible}, we have 
\begin{center}
\begin{tabular}{|c|c|l|}
        \hline 
        homogeneity & weight & range of $m$ of $m$-th cochains \\\hline
        $h=0$ & $ w \geq -\dim\frakS{1}$ & $\max(0, -w) \leq m \leq w+
        2\dim \frakS{1} + \dim \frakS{2} $\\
        $h=1$ & $ w \geq 0$ & $0 \leq m \leq w+\dim\frakS{1}$\\
        $h\geq 2$ & $ w \geq (h-1)\dim\frakS{1}$ & $0 \leq m \leq w$\\\hline
\end{tabular} 
\end{center} 
where $\ds \dim\frakS{1} \leq n$ and $\dim\frakS{2} \leq n(n+1)/2$, and 
depend on the existence of Casimir polynomials.  

%%% HH

\subsection{Finding a basis of $\ds\mySS{j}/\myCasim{j}$ by 
Groebner basis theory} 
If we deal with the algebra of Hamiltonian vector fields, 
it seems hard to say which of homology or cohomology 
is easier, faster, lighter or better.  
We prepare the vector subspace $\ds\myCasim{j}$ of  
Casimir polynomials in $\ds\mySS{j}$. 
Let $\ds f_1,\ldots,f_k$ be a basis of $\ds \myCasim{j}$.  
Since $\ds \frakS{j} = \{ \sigma\in\frakSS{j} \mid \langle \sigma,
\myCasim{j}\rangle = 0\}$,  
take a generic element $\ds \sigma =  \sum_{A} c_A \z{A} \in
\frakSS{j}$ and  
solve the linear equations $\ds
\langle \sum_{A} c_A \z{A}, f_i \rangle =0$ for  $i=1,\ldots ,k$. 
Substituting the solution to $\sigma$, we get a basis of
$\ds\frakS{j}$. 

For $\ds\myS{j} = \mySS{j}/\myCasim{j}$, 
conceptually we know well about 
%$\ds\myS{m} = \mySS{m}/\myCasim{m}$, 
the structure of the quotient vector space.   
% Let $W\subset V$ be a subspace of $V$. Then 
% the quotient vector space $\ds V/W$ is given by \[ (v+W) + (v' + W):=
% v+v'+W\;,\;\text{and}\; \lambda(v+W):= \lambda v + W\;. \]
But it is not so easy to fix a basis of $\ds \mySS{j}/\myCasim{j}$ in general. 
It is well-known that 
if 
there is a linear map $\ds P:\mySS{j} \rightarrow \mySS{j}$ with $\ds P^2=P$
and $\ds P^{-1}(0)=\myCasim{j}$. Then the subspace $\ds
P(\mySS{j})\subset \mySS{j}$ is isomorphic with 
$\ds \mySS{j}/\myCasim{j}$ by $\ds
P(f) \mapsto f+\myCasim{j}$.

Now we make use of Groebner basis theory and observe that  
the normal form $\ds\phi$ with respect to 
some basis $\ds f_1,\ldots,f_k$ of $\ds \myCasim{j}$ satisfies the property
above and is the projection we want.  

\subsubsection{Remainder modulo Casimir polynomials} \label{koko}
Let $\ds\mySS{j}$ be the vector space of homogeneous $j$-polynomials of 
$x_1,\ldots,x_n$, i.e., 
$\ds\mySS{j}$ is generated by $\ds\{\w{A} \mid |A|=j\}$.  
Let $f_1,\ldots,f_k$ be a basis of $\ds \myCasim{j}$.  
To apply
Groebner basis theory, we fix a term order $\ds x_1> \cdots > x_n$, for
instance, the graded reverse lexicographic order.     
By the notation $\myLM{}$, we mean  
the leading monomial with respect to our term order.  
If $\ds g_1$ and $\ds g_2$ are linearly independent and  
$\ds \myLM{g_1} = \myLM{g_2}$, then 
$\ds h_1 = c_1 g_1 + c_2 g_2$ and  
$\ds h_2 = c_2 g_1 - c_1 g_2$ are again linearly independent, and  
$\ds \myLM{h_1} =  \myLM{g_1} = \myLM{g_2}$ and 
$\ds \myLM{h_2} <  \myLM{g_1} = \myLM{g_2}$, where 
$\ds g_1 = c_1 \myLM{g_1} + \cdots$ and   
$\ds g_2 = c_2 \myLM{g_2} + \cdots$.  
Thus, we may assume that $\myLM{f_1} > \cdots > \myLM{f_k}$ and each
$f_i$ contains only $\myLM{f_i}$ among  $\ds \{ \myLM{f_{\ell}} \mid
\ell=1..k\}$  after changing a basis by a suitable linear transformation
if necessary. This means we take the Groebner basis of $\ds\myCasim{j}$.  

For each $g\in\mySS{j}$, we recall the division algorithm with respect
to $f_1,\ldots,f_k$:  
\begin{equation}\label{eqn:algor}
        g = \sum_{i=1}^k a_i f_i + \remainder
\end{equation}
with the next two conditions:
\begin{itemize}
        \item[(a)]
if $\remainder \ne 0$, each term of $\remainder$ is not multiple of 
                $\myLM{f_i}$ for each $i$, and 
\item[(b)] 
if $a_i \ne 0$ then $\ds \myLM{a_i f_i} \leq \myLM{g}$. 
\end{itemize}
$\remainder$ is called the remainder or the normal form 
of $g$ with respect to $f_1,\ldots,f_k$.  

In our case, the situation is very simple and we have the next Lemma:

\begin{Lemma}\label{lemma:key}
The remainder $\ds\remainder$ of 
a given $j$-homogeneous polynomial $g$ with respect 
 to $f_1,\ldots,f_k$ is obtained as 
\begin{align}\label{eqn:lock}
\remainder =& \sum_{\w{A} \not\in \{\myLM{f_i}\mid i=1..k\} } c_A \w{A}
        \quad
\text{where}\; c_A \;\text{are constants}\;, \\ \noalign{and}
        \label{eqn:lock:two}
g =& \sum_{i=1}^k a_i f_i + \remainder\; \quad\text{where}\; a_i\;\text{are
        constants}\;. 
\end{align}
The constants $\ds a_i$ and $\ds c_A$ with 
$\ds \w{A} \not\in \{\myLM{f_i}\mid i=1..k\} $ are uniquely determined
by $g$ (depending on the term order).  If we denote the remainder
$\remainder$ by $\phi(g)$, then  $\phi^2 = \phi$ and $\phi^{-1}(0) =
        \text{LSpan}(f_1,\ldots,f_k)$.  
\end{Lemma}
\textbf{Proof:} 
In our case, since $g$ and $f_i$ are homogeneous with the same
degree $j$, the condition (b) says $\ds a_i$ must be constant. Thus
\eqref{eqn:algor} yields that $\remainder$ is also the same
$j$-homogeneous. The condition (a) yields that  
\eqref{eqn:lock}. Since in the space $\ds\mySS{j}$, $\ds\{ 
\myLM{f_i}\mid i=1..k\}$ 
and $\ds \{ 
\w{A} \not\in \{\myLM{f_i}\mid i=1..k\}\}$ give a basis, and we see that 
\[ g = \sum_{i} p_i \myLM{f_i} + \sum_{   
\w{A} \not\in \{\myLM{f_i}\mid i=1..k\}}
b_B \w{B} =
\sum_{i} a_i f_i + \sum_{   
\w{A} \not\in \{\myLM{f_i}\mid i=1..k\}}
c_B \w{B}\;. 
\]
\kmcomment{
About the second assertion, 
if $\ds \w{A}=\myLM{f_i}$ for some $i$, then $\ds \w{A} = a_i f_i + g$ with
$\ds \myLM{g}< \myLM{f_i}$. Now
        $\ds\phi(\w{A}) = \phi(g)$ and so  
\begin{align*}
        \myLM{\phi(\w{A})} =& \myLM{\phi(g)} \leq \myLM{g}  
         < \myLM{f_i} = \w{A}\;. %\hspace{20mm} \blacksquare
\end{align*} 

\vspace{-7mm} 
}%endOFkmcomment

\hfill $\blacksquare$

Lemma \ref{lemma:key} %in Section \ref{koko}.  
implies directly the next proposition. 
\begin{kmProp}\label{prop::fixing::nf}
        For each $\ds f\in \mySS{j}$, $\ds\phi(f)=0$ if and only if $\ds f\in
\myCasim{j}$.  Let $\ds\{ \w{A} \mid |A|= j\}$ be the natural basis of
$\ds\mySS{j}$. Then 
$\{ \w{A} \in \mySS{j} \mid \phi( \w{A} ) = \w{A} \}$ is a basis of
        $\ds\myS{j} = \phi(\mySS{j}) \subset \mySS{j}$.  
        Also $\ds \{\w{A} \mid |A|=j\} \smallsetminus 
        \{ \myLM{ f_i } \mid f_i \in \myCasim{j} \}$ gives the same basis of 
        $\ds\myS{j} = \phi(\mySS{j}) \subset \mySS{j}$, where we assumed
        that $\ds\{ f_i \}$ is a basis of $\ds\myCasim{j}$ with mutually
        distinct leading monimials $ \{ \myLM{ f_i } \}$.  
        % ({\color{red} So far, no clear proof for it!})
%We refer Lemma \ref{lemma:key} in Section \ref{koko}.  
\end{kmProp}

For $f\in \myS{i}$ and 
$g\in \myS{j}$,  the Lie bracket is given by $\ds [f,g]= \phi( \Bkt{f}{g} )$, 
where  $\phi$ is now 
the normal form with respect to 
$\ds \myCasim{i+j-2+h}$, here $h$ is the homogeneity of the Poisson
structure $\pi$.   Thus, the boundary operator $\pdel$ of $\frakg$ is
given by 
\begin{equation}
        \pdel \&^{\wedge}(f_1,\ldots,f_k) = 
        \sum_{i<j}  (-1)^{i+j} 
        \&^{\wedge}( \phi(\Bkt{f_i}{f_j}),\ldots,\widehat{f_i},
        \ldots,\widehat{f_j},\ldots)  
\end{equation}
where $\ds f_i \in \myS{\sigma(i)}\subset \mySS{\sigma(i)}$. 
\begin{kmRemark}
In the above, we do not miss $\phi$-projection. See 
         \eqref{align:hanrei} in the next example. 
\end{kmRemark}

\begin{exam} \label{exam:sl2} 
        Take a Lie Poisson structure $\fraksl(2)$ as a concrete
        example. 
        \[ \Bkt{p}{q}= r\;,\;   \Bkt{r}{p}= 2p\;,\;   \Bkt{r}{q}= -2q\;.  
        \] 
        Then we have \[\ds \myCasim{2j-k}= \begin{cases}\; \emptyset
                & \text{if}\; k=1\;, \\
       \; \mR \cdot (4 p q +r^2)^{j} & \text{if}\; k  = 0\; . \end{cases}\]  

Let $\ds \w{A} = p^{a_1}  q^{a_2}  r^{a_3} $ and $\ds \z{A}$ are the dual.  
Let $\phi$ be the normal form $\phi$ with respect to $\ds 4 p q + r^2$
        under the term order $\ds p>q>r$.  Then $\ds\phi( \w{1,1,0}) = -
        \frac{1}{4} \w{0,0,2}$ and $\ds\phi( \w{A}) = \w{A}$ for $A \ne
        [1,1,0]$. We get a basis $\w{A}$ ($A \ne [1,1,0]$) of
        $\ds\myS{2} = \phi( \mySS{2})$.   
        Since $\ds\myCasim{1}= 0$, we see $\ds\myS{1} = \mySS{1}$.  
Take $\ds f= x_1 \in \myS{1}$ and $\ds g= x_2 x_3 \in \myS{2}$, then  
        \begin{align} \label{align:hanrei}
                & \Bkt{f}{g} = x_3{}^2 - 2 x_1 x_2 \not\in
                \phi(\mySS{2})\;,\;\text{ and }\; 
        \phi(\Bkt{f}{g}) = \frac{3}{2} x_3{}^2 \in
        \phi(\mySS{2})\;.\end{align}  

About $\ds \sum_{|A|=2} c_A \z{A}\in \frakS{2} \subset \frakSS{2}$, 
\begin{align*}
0 =& \langle \sum_{|A|=2} c_A \z{A}, 4pq +r^2\rangle=4 c_{1,1,0}+c_{0,0,2} \\
\noalign{implies}
\sum_{|A|=2}  c_A \z{A} = & 
\sum_{a_1\ne a_2} c_A\z{A} -\frac{c_{0,0,2}}{4}\z{1,1,0}+c_{0,0,2}\z{0,0,2} \;
        \text{ or }\;  
\sum_{a_1\ne a_2} c_A\z{A} +c_{1,1,0}\z{1,1,0} - 4c_{1,1,0}\z{0,0,2} \; . 
\end{align*}
        Thus, (natural) bases of $\ds \frakS{2}$ are given by  
\[ \z{A}\; (a_1 \ne a_2)\; \;\text{ and }\; -\frac{1}{4}\z{1,1,0}+\z{0,0,2}\;
        \;\text{ or }\;\; 
\z{A}\; (a_1 \ne a_2)\; \;\text{ and }\; \z{1,1,0}-4 \z{0,0,2}\; . \]

\end{exam}

\subsubsection{Easy way to get the remainder}
Even though 
the most polite way is using the normal form of Groebner basis theory,  
in our case, 
we may use the following
formula due to Lemma \ref{lemma:key}: 
\begin{equation}\label{eqn:anchoku}
        \phi(g) = g - \sum_{i} \frac{ \pdel^{ A_i} g }{ \pdel^{ A_i}
f_i} f_i \end{equation} where 
$\ds A_i $ is given by $\ds \w{A_i} = \myLM{f_i}$ and $\ds \pdel^{A} =
(\frac{\pdel}{\pdel x_1})^{a_1}\cdots 
(\frac{\pdel}{\pdel x_n})^{a_n}$.  In practical calculation, 
\eqref{eqn:anchoku} seems faster than the normal form.

\begin{exam} 
        We deal with Example \ref{exam:sl2} with 3 variables  $\ds
        x_1,x_2, x_3$ with the graded reverse lexicographic order as the
        term order.     
        Consider 2-homogeneous space $\ds\mySS{2}$ and the subspace 
        $\ds \{ f_1 = 4 x_1 x_2 + x_3{}^2\}$. Now $\ds \myLM{f_1}= x_1
        x_2$ and so the normal form of $g$ is given by $\ds 
        \sum_{A \ne [1,1,0]} c_A \w{A} $ of 
        \begin{equation} \label{align:one:eg}
                g = a f_1 + \sum_{A \ne [1,1,0]} c_A \w{A}
                 =   
                4 a \w{[1,1,0]} 
                + ( a + c_{[0,0,2]}) \w{[0,0,2]} 
                + \sum_{A \ne [1,1,0], [0,0,2]} c_A \w{A} 
                \;. 
        \end{equation}
When $g= x_1 x_2$,  then
        \eqref{align:one:eg} implies  
        $4a=1$, $\ds c_A=0$ 
        ($\ds A \ne [1,1,0], [0,0,2]$) and $\ds 
                a + c_{[0,0,2]} = 0$, thus $\ds\phi(x_1 x_2) = -
                \frac{1}{4} x_3{}^2$.   

        When $g= x_3 {}^2$,  then
        \eqref{align:one:eg} implies  
        $4a=0$, $\ds c_A=0$ 
        ($\ds A \ne [1,1,0], [0,0,2]$) and $\ds 
                a + c_{[0,0,2]} = 1$, thus $\ds\phi(x_3 {}^2) = 
                 x_3{}^2$.   

        If we use the formula \eqref{eqn:anchoku}, $\ds\phi(g)$ is given
        by 
        \[ \phi(g) = g - \frac{ \pdel_1 \pdel _2 g }{4} 
        (4 x_1 x_2 + x_3{}^2) \]
        and $\ds \phi(x_1 x_2 ) = x_1 x_2 - \frac{1}{4}
        (4 x_1 x_2 + x_3{}^2)  = - \frac{ x_3{}^2}{4}$, 
        $\ds \phi(x_3{}^ 2 ) = x_3 {}^2 - 0 (4 x_1 x_2 + x_3{}^2)  = 
         x_3 {}^2$, and so on. 
\end{exam}

Another easy way to get the remainder in our case is the following:
\begin{enumerate}
        \item Take $\ds\myCasim{j}$.  
                Let $\ds \{ f_{\ell} \mid \ell=1..k\}$ be a basis
                with the properties described before.

        \item Get the leading monomials of $\ds\myCasim{j}$,  
                $\ds \text{LMCasim} = 
                \{ \myLM{f_{\ell}} \mid \ell=1..k\}$

                Now, $\phi(\w{A}) = \w{A}$ for $\ds \w{A}\not\in \text{LMCasim}$.  
        \item To get $\ds \phi(\myLM{f_i})$, 
                we use $\ds \phi(f_i)=0$. Since $\phi$ is
                $\mR$-linear and satisfies the property above, we get
                the expression of 
                $\phi( \myLM{f_i} )$  by $\ds\mR$-linear combination of 
                $\ds \w{A}\not\in \text{LMCasim}$.  
\end{enumerate}
\begin{exam} 
        Again, we deal with the Lie Poisson modeled by $\ds\fraksl(2)$ with 3 variables  $\ds
        x_1,x_2, x_3$ with the graded reverse lexicographic order as the
        term order.     
        Consider 2-homogeneous space $\ds\mySS{2}$ and the Casimir subspace 
        is given by 
        $\ds \{ f_1 = 4 x_1 x_2 + x_3{}^2\}$. Now $\ds \myLM{f_1}= x_1
        x_2$ and so $\ds \phi( \w{A}) = \w{A}$ unless $a_1=a_2=1$ where
        $\phi$ is the normal form with respect to the Casimir subspace.  
        Since \[\ds 0 = \phi(f_1) = 4 \phi(x_1 x_2) + \phi( x_3{}^2) = 
         4 \phi(x_1 x_2) +  x_3{}^2\;,\] 
we get  
        $\ds \phi(x_1 x_2) = - \frac{1}{4} x_3{}^2$.
\end{exam}

\subsection{Fixing the dual basis of $\ds\mySS{j}/\myCasim{j}$ in
$\ds\frakS{j}$}
Let $\phi$ be the normal form with respect to $\ds \myCasim{j}$ on
$\ds\mySS{j}$.   
Using the natural basis  $\ds\w{A}$ with $|A|=j$ in $\ds\mySS{j}$, we
consider $\ds \tau = \sum_{|A|=j} \phi(\w{A}) \otimes \z{A}$.  
We have known that $\ds \{ \w{B} \in \mySS{j} \mid \phi( \w{B}) = \w{B}
\}$ is a basis of  $\ds\myS{j} = \mySS{j}/\myCasim{j}$ in Proposition
\ref{prop::fixing::nf}, and   
the $\ds \w{B}$-coefficient of $\ds\tau$ 
with $\ds \{ \w{B} \in \mySS{j} \mid \phi( \w{B}) = \w{B} \}$ 
form a basis of $\ds\frakS{j}$.  

\begin{exam} Again, we deal with Example \ref{exam:sl2} with 3 variables  $\ds
        x_1,x_2, x_3$ with the graded reverse lexicographic order as the
        term order.     
        Consider 2-homogeneous space $\ds\mySS{2}$ and the subspace 
        $\ds \{ f_1 = 4 x_1 x_2 + x_3{}^2\}$. Now $\ds \myLM{f_1}= x_1
        x_2$ and the normal forms are given by 
                $\ds\phi(x_1 x_2) = - \frac{1}{4} x_3{}^2$ and 
                $\ds\phi(\w{B}) = \w{B}$ for $\ds |B|=2$ and $\ds B\ne [1,1,0]$. 
        Then we have $\ds \tau = \sum_{a_1 \ne a_2} \w{A} \otimes \z{A} 
        + (- \frac{1}{4} x_3{}^2 ) \otimes \z{1,1,0} + 
        ( x_3{}^2 ) \otimes \z{0,0,2}  
                 $ and we get a basis of $\ds\frakS{2}$ 
        \[ \z{A} \; (|A|=2\;,\; a_1 \ne a_2)\;, \;\; 
         - \frac{1}{4}  \z{1,1,0} + \z{0,0,2}\;. \]  
\end{exam}

\begin{exam}
Consider $\ds \fraksl(2)$ in $\ds \mR^4$, 
        i.e., the Poisson bracket is given by 
        \begin{align*}
                \Bkt{x_1}{x_2} & = x_3 = - \Bkt{x_2}{x_1} \;,\; 
         \Bkt{x_1}{x_3} = -2 x_1 = - \Bkt{x_3}{x_1} \;,\; 
        \Bkt{x_2}{x_3} = 2 x_2 = - \Bkt{x_3}{x_2} \;,\\
                \Bkt{x_i}{x_4} & =  
        \Bkt{x_4}{x_i} =  0 \quad (i=1,2,3)\; .\end{align*} 
Then $\ds \myCasimm$ is
        rather complicated and generated by $\ds f(4x_1x_2+x_3{}^2)
        g(x_4)$ where $f,g$ are polynomials of one variable.  
        Thus, $\ds \myCasim{j}$ are given as
        \begin{align*} 
                \myCasim{1} & = [ x_4]\;,\; 
        \myCasim{2} = [ 4x_1x_2+x_3{}^2 \;,\; x_4{}^2]\;,\;
        \myCasim{3} = [ (4x_1x_2+x_3{}^2)x_4 \;,\; x_4{}^3]\;,\;
                \\
                \myCasim{4} & = [ 
        (4x_1x_2+x_3{}^2)^2\;,\; 
        (4x_1x_2+x_3{}^2)x_4{}^2 \;,\; 
        x_4{}^4]\;,\;\ldots \end{align*} 
        The leading monomials of $\myCasim{j}$ are         
\begin{align*} 
        \myLM{\myCasim{1}} & = [ x_4]\;,\; 
        \myLM{\myCasim{2}} = [ x_1x_2 \;,\; x_4{}^2]\;,\;
        \myLM{\myCasim{3}} = [ x_1x_2x_4 \;,\; x_4{}^3]\;,\;
                \\
             \myLM{\myCasim{4}} & = [ 
        (x_1x_2)^2\;,\; 
        (x_1x_2)x_4{}^2 \;,\; 
        x_4{}^4]\;,\;\ldots \end{align*} 
and so we have bases of $\ds \myS{j}$ as follows: 
\begin{align*} 
\text{basis of }\myS{1} & = [ x_1, x_2, x_3]\;,\; 
\text{basis of }\myS{2} = [\w{A} \mid |A|=2] \smallsetminus 
        [ x_1x_2 \;,\; x_4{}^2]\;,\; \\
        \text{basis of }\myS{3} &=  [\w{A} \mid |A|=3] \smallsetminus 
        [ x_1x_2x_4 \;,\; x_4{}^3]\;,\;
                \\
  \text{basis of }\myS{4} & =
         [\w{A} \mid |A|=4] \smallsetminus 
        [ 
        (x_1x_2)^2\;,\; 
        (x_1x_2)x_4{}^2 \;,\; 
        x_4{}^4]\;,\;\ldots \end{align*} 
Concerning to find a basis of $\ds \frakS{j}$, we need information about
        the normal form $\phi$ with respect to $\ds\myCasim{j}$. 

        When $\ds \frakS{1}$, we see that 
        $\ds \phi(\w{\epsilon_i})= \w{\epsilon_i}$ for $i=1,2,3$ and 
        $\ds \phi(\w{\epsilon_4})= 0$. Thus, $\ds \tau = \sum_{i=1}^3
        \w{\epsilon_i} \otimes \z{ \epsilon_i}$ and we get a basis 
        $\ds \{ \z{\epsilon_i} \mid i=1,2,3\}$ on $\ds \frakS{1}$.   

When $\ds \frakS{2}$, we get $\ds \phi( x_1 x_2) = - \frac{1}{4}
        x_3{}^2$, $\phi( x_4{}^2)=0$ and the other monomials are
        invariant by $\phi$.   
        Now $\ds \tau = 
        \sum_{|A|=2} \phi(\w{A}) \otimes \z{A} = 
        \sum_{A \ne [1,1,0,0], [0,0,0,2]} \w{A}\otimes \z{A} 
        + (-\frac{1}{4} x_3{}^2 ) \otimes \z{[1,1,0,0]} $ and so 
        we get a basis of $\ds\frakS{2}$ 
        \[ \z{A} \; (A \ne[1,1,0,0],[0,0,2,0], [0,0,0,2] )\; \text{ and
        }\;
        -\frac{1}{4} \z{[1,1,0,0]} + 
        \z{[0,0,2,0]}\;.  \]

By the same way, about $\ds \frakS{3}$ we have a basis
        \[ \z{A} \; (A \ne [0,0,0,3], [1,1,0,1],[0,0,2,1])\; \text{ and
        }\;
        -\frac{1}{4} \z{[1,1,0,1]} + 
        \z{[0,0,2,1]}\;.  \]
As a basis of $\ds\frakS{4}$ we have
        \[ \z{A} \; (A \ne [0,0,0,4],[0,0,2,2], [1,1,0,2], [0,0,4,0],
        [1,1,2,0],[2,2,0,0])\] and
       \[ 
        -\frac{1}{2} \z{[2,2,0,0]} + \z{[1,1,2,0]}\;,\;  
        -\frac{1}{4} \z{[1,1,0,2]} + \z{[0,0,2,2]}\;,\;  
        -\frac{1}{16} \z{[2,2,0,0]} + \z{[0,0,4,0]}\;.  
        \]

\end{exam}

\begin{kmRemark} 
About some comparison of the effort of using our method or using the
        normal form in order to get the remainder, we refer to the
        subsection \ref{elapsed:time}. There, we got data of elapsed
        time for concrete Poisson structures.  
\end{kmRemark}

\kmcomment{
\subsubsection{Comparison of elapsed time}
\paragraph{Polynomials}
\begin{center}
\begin{tabular}{|c|*{2}{c}|*{2}{c}|}\hline 
& \multicolumn{2}{|c|}{sl2} & \multicolumn{2}{c|}{Heisen3} \\\hline \hline 
        weight & HomoL & Cohom & HomoL & Cohom \\\hline
        wt=0 & 0.021 & 0.021 & 0.024 & 0.019 \\
        wt=1 & 0.047 & 0.057 & 0.043 & 0.042 \\
        wt=2 & 0.221 & 0.269 & 0.136 & 0.163 \\
        wt=3 & 1.467 & 1.664 & 0.920 & 0.980 \\
wt=4 & 13.274 & 14.651 & 7.260 & 7.592\\\hline \end{tabular}
\end{center}

cf: \verb+../km/Odd_GKF/new-exec-(0,1,2)-v<>.mpl+ is for Cohomology of
Polynomials.   

cf: \verb+../km/Odd_GKF/Homol/homoL-exec-(0,1,2)-v<>.mpl+ is for Homology 
of Polynomials.

\paragraph{Hamiltonian}
\begin{center}
\begin{tabular}{|c|*{2}{c}|*{2}{c}|}\hline 
& \multicolumn{2}{|c|}{sl2} & \multicolumn{2}{c|}{Heisen3} \\\hline \hline 
        weight & HomoL & Cohom & HomoL & Cohom \\\hline
        wt=0 & 0.039(0.036) & 0.022 & 0.023({\color{red}0.032}) & 0.020 \\
        wt=1 & 0.197(0.186) & 0.062 & 0.086(0.080) & 0.044 \\
        wt=2 & 1.300(1.198) & 0.285 & 0.429(0.353) & 0.090 \\
        wt=3 & 8.236(7.745) & 1.512 & 2.335(1.865) & 0.312 \\
wt=4 & 44.539(40.497) & 11.453 & 11.645(8.873) & 1.552\\\hline \end{tabular}
\end{center}

The numbers in the parentheses are obtained by getting the remainder by 
\eqref{eqn:anchoku}. Except only one case printed in red color,  
using \eqref{eqn:anchoku} is faster than using the normal form. Even
though, computing cohomology groups is faster than homology computing in 
Hamiltonian vector fields cases.  

cf:\verb+../Odd_GKF/new-exec-2-v<2>.mpl+ is for Cohomology 
of Hamiltonian (just switch \texttt{myHamil}).

cf:\verb+../Odd_GKF/Homol/homoL-exec-3-v<>.mpl+ is for Homology 
of Hamiltonian.
}%%%%endOFkmcomment

\section{Subalgebra $\ds \ovfrakg$ excluding constant polynomials} 

\(\ds\sum_{j=\ell}^{\infty} \mySS{j} \) is a subalgebra of
$\ds\polygon{n}$ for each $\ell \geq  1$ 
when the homogeneity of Poisson structure $\pi$ is     
$\geq 1$, and   
when the homogeneity is 0,    
\(\ds \sum_{j=\ell}^{\infty} \mySS{j} \) is a subalgebra 
for ${\color{red} \ell \geq 2}$.    
The notion of weight is valid on those subalgebras and we have
the weight decomposition of cochain complex of \(\ds\ovfrakg 
= \sum_{j=1}^{\infty} \mySS{j} \) when $h>0$ and  \(\ds\ovfrakg 
= \sum_{j=2}^{\infty} \mySS{j} \) when $h=0$.  
We denote 
by $\ds \myCC{m}{w}$ the $m$-th cochain space with the weight $w$,
namely 
\begin{equation} 
        \myCC{m}{w} = \sum \Lambda^{k_0}\frakSS{0} \otimes 
\Lambda^{k_1}\frakSS{1} \otimes \cdots \otimes 
\Lambda^{k_{\ell}}\frakSS{\ell}  \label{eqn:sub:one}
\end{equation}
with the condition \eqref{cond:deg}, \eqref{cond:wt} and \eqref{cond:dim} with the
additional restriction $\ds k_0=k_1=0$ if $h=0$ and 
$\ds k_0=0$ if $h>0$.

\subsection{Possibility of weight and range of degree when excluding
constant polynomials}
We follow the discussion in the subsection 
\ref{subsec:possible} with 
the restriction 
$\ds k_0=0$ if $h>0$  
or 
$\ds k_0=k_1=0$ if $h=0$,  
we see that  
\begin{center}
\begin{tabular}{|c|c|l|}
        \hline 
        homogeneity & weight & range of $m$ of $m$-th cochains \\\hline
$h=0$ & $ w \geq 0$ & $0\leq m \leq w+n(n+1)/2$\\
$h=1$ & $ w \geq 0$ & $0\leq m \leq w+n$\\
$h>1$ & $ w \geq 0$ & $0\leq m \leq w$\\\hline
\end{tabular}
\end{center}

\subsection{$\ds\ovfrakg$ and $\ds\polygon{n}$ when homogeneity $h>0$}
The dimensional condition 
\eqref{cond:dim} yields $\ds k_0=0$ or $1$. Thus, we can decompose 
\eqref{eqn:one:ato} as follows: 
\begin{equation} \myCO{m}{w}{ \polygon{n}^{*}}
        = \sum_{k_0=0}  
\Lambda^{k_1}\frakSS{1} \otimes \cdots \otimes 
\Lambda^{k_{\ell}}\frakSS{\ell}  
\; \oplus\;
        \sum_{\kk{0}=1}  \frakSS{0} \otimes 
        \Lambda^{\kk{1}}\frakSS{1} \otimes \cdots \otimes 
\Lambda^{\kk{\ell'}}\frakSS{\ell'} 
\label{eqn:two}
\end{equation} 
The first term of the direct sum is the $m$-th cochain space
$\ds\myCC{m}{w}$ of $\ds \ovfrakg$.   About the second term, 
\eqref{cond:deg} implies 
$\ds  \sum_{j=1}^{\infty} \kk{j} = m -1 $, and  \eqref{cond:wt} implies  
$\ds \sum_{j=1}^{\infty} \kk{j} (j-2+h) = w + (2-h)$, thus,   
the second term coincides with $\ds \frakSS{0} \otimes
\myCC{m-1}{w+2-h}$, namely  
\begin{align}
\myCO{m}{w}{\polygon{n}^{*}} =&  \myCC{m}{w}\; 
\oplus\; \frakSS{0} \otimes \myCC{m-1}{w+2-h} \;.  \label{eqn:three}\\
\noalign{Similarly, we have} 
\myCO{m}{w}{\polygon{n}^{}} =&  \myCC{}{m,w}\; 
\oplus\; \mySS{0} \otimes \myCC{}{m-1,w+2-h} \;.  \label{eqn:three:sonoNi} 
\end{align} 
Since we know that $\ds d(\frakSS{0})=0$
and \( \pdel (\&^{\wedge} (1,u_2,\ldots,u_m))= - 1 \wedge 
\pdel (\&^{\wedge} (u_2,\ldots,u_m)) \), we have
\begin{kmProp} When $h>0$, the followings hold
\begin{align*} 
\myHH{m}{w} (\polygon{n}) \cong &  \myHH{m}{w}(\ovfrakg)\; \oplus\; 
        \myHH{m-1}{w+2-h}(\ovfrakg) \;, \\
\myHH{}{m,w} (\polygon{n}) \cong &  \myHH{}{m,w}(\ovfrakg)\; \oplus\; 
\myHH{}{m-1,w+2-h}(\ovfrakg) \;.  
\end{align*}
\end{kmProp}
The proposition above says that it is enough to handle $\ds
\myHH{m}{w}(\ovfrakg)$ to study $\ds 
         \text{H}^m_{w} (\polygon{n})$. Thus, from now on we only deal with
         $\ds\ovfrakg$ and the cochain complex $\ds
         \myCC{m}{w}(\ovfrakg)$ (sometimes simply denote $\ds  
         \myCC{m}{w}$) or the chain complex
         $\ds \myCC{}{m,w}(\ovfrakg)$ (or $\ds  
         \myCC{}{m,w}$).    
We rewrite Proposition 
 \ref{kmProp::CD} by using the words of $\ds\ovfrakg$ as follow:

\begin{kmProp} 
        When $h>0$, 
$\ds \mydq$  is the natural restriction of $\ds\mydov$ and we have 
the following commutative diagram  
\begin{equation}
        \begin{CD}
\myC{m}{} = \Lambda^m \frakg^{*} @. \subset @.
\myCC{m}{} = \Lambda^m \ovfrakg^{*}   \\
@Vd VV   @. @VV\overline{d}V \\ 
\myC{m+1}{} = \Lambda^{m+1} \frakg^{*} @. \quad \subset \quad @.  
\myCC{m+1}{} = \Lambda^{m+1} \ovfrakg^{*}   \;. 
\end{CD}
\label{Kont:CD:two}
\end{equation} 
The coboundary operator $\ds\mydq$ preserves the weight, namely, $\ds
\mydq ( \myC{m}{w} ) \subset \myC{m+1}{w}$ holds. 
        Thus  we have the
cohomology group 
        \[\ds \myH{m}{w}(\frakg) := \text{Ker} ( \mydq : \myC{m}{w}  \rightarrow
\myC{m+1}{w}) / \mydq ( \myC{m-1}{w})\;.\]   
\end{kmProp}
As a direct corollary of Proposition above, we have
\begin{kmCor} When $h>0$, 
        $\ds\myH{\bullet}{w}(\frakg) = \left(\myC{\bullet}{w}\cap \ker ( \mydov:
        \myCC{\bullet}{w}\rightarrow \myCC{\bullet+1}{w})\right) \big/ 
        \mydq ( \myC{\bullet-1}{w})$ .  
\end{kmCor}

\subsection{Young diagrams and cochain spaces}
\label{subsec:YD:cochain} 
The $m$-th cochain space of weight $w$ 
\[ 
\myCC{m}{w} = \sum
\Lambda^{k_1} \frakSS{1} \otimes 
\Lambda^{k_2} \frakSS{2} \otimes 
\Lambda^{k_3} \frakSS{3} \otimes \cdots  
\] 
satisfies the following three conditions \eqref{cond:deg},
\eqref{cond:wt} and \eqref{cond:dim} with $k_0=0$. 
\kmcomment{
\begin{align}
& k_1+ k_2 + \cdots+ k_j + \cdots 
= m\;, \label{cond:deg}\\
%\noalign{and}
& k_1(1+h-2) +  k_2(2+h-2)+\cdots 
+   k_j(j+h-2)+\cdots 
=w\;, \label{cond:wt} \\
%\noalign{and}
0\leq & k_j \leq \dim\frakSS{j} = (j+n-1)!/( j! (n-1)!) 
=\tbinom{j+n-1}{n-1}\;. \label{cond:dim} 
\end{align} 
}%endOFkmcomment
%%%%%%%%%%%%%%%%%%%%%%%%%%%%%%%%%%%%%%%%%%%%%%%%%%%%%%%%%%%%%%%%%%%%%%%%%%%%%%%%

Take a sequence \(\ds [k_1,k_2,\ldots]\) satisfying the three conditions
above.  Then there exists
an $\ell$ so that $\ds k_j=0$ ($j>\ell$) and we may write 
\(\ds [k_1,k_2,\ldots] = [k_1,k_2,\ldots,k_{\ell}] \) as a finite sequence.  
The conditions \eqref{cond:deg} and \eqref{cond:wt}
for \(\ds [k_1,k_2,\ldots,k_{\ell}] \) 
are equivalent to \eqref{cond:deg} and \eqref{jou:wt} with $\ds k_0=0$. 
\kmcomment{
\begin{align} 
& k_1+ k_2 + \cdots + k_j + \cdots + k_{\ell}
 = m 
 \label{eq::three:one}
 \\
 \noalign{and} 
 & k_1 + 2 k_2 + \cdots + j k_j + \cdots + \ell k_{\ell}
 = w + (2-h) m \;. 
 \label{eq::three:two}
\end{align} 
}%endOFkmcomment
As observed in \cite{M:N:K},
\( ( \underbrace{\ell,\ldots,\ell}_{k_{\ell}\text{-times}}, 
        \ldots, 
\underbrace{1,\ldots,1}_{k_{1}\text{-times}} ) \) is a Young diagram (in
traditional expression) whose  
length $m$ and consisting of $(w+(2-h)m)$ cells 
from \eqref{cond:deg} and \eqref{jou:wt}.   

\begin{center}
\setlength{\unitlength}{4mm}
\begin{picture}(12,11)(-2,0)
\path(0,0)(0,10)% (1,2)(0,2)(0,0)
%\put(0,0){$\mathbf O$}
\multiput(0,0)(0,1){3}{\line(1,0){1}}
\path(1,0)(1,2)
\put(1.5,1){\makebox(0,0)[l]{$k_{1}$}}
\multiput(0,2)(0,1){4}{\line(1,0){2}}
\path(2,2)(2,5)
\put(2.5, 3){$k_{2}$}
\multiput(0,8)(0,1){3}{\line(1,0){8}}
\path(8,8)(8,10)
\put(8.5, 9){\makebox(0,0)[l]{$k_{\ell}$}}
\put(4,10.5){$\ell$}
\put(3,10.5){\vector(-1,0){3}}
\put(5,10.5){\vector(1,0){3}}
%\put(5,0){\Large A} 
\put(-2,5){\makebox(0,0){$m$}}
\put(-2,6){\vector(0,1){4}}
\put(-2,4){\vector(0,-1){4}}
%\dashline{4}[0.8](6,6)(6,7)
\dottedline[$\cdot$]{0.2}(2,5)(3,5)
\dottedline[$\cdot$]{0.2}(3,5)(3,6)
\dottedline[$\cdot$]{0.2}(7,8)(7,6)
\dottedline[$\cdot$]{0.2}(3,6)(7,6)
%\put(5,0){\arc{1.2}{5.2}{6.2}}
\put(4,1.7){\makebox(0,0)[l]{Total area is}}
\put(4,0.5){\makebox(0,0)[l]{$w +(2-h)m$}}
\end{picture}
\end{center} 
Conversely, for a given Young diagram $\lambda$ (in traditional), define
\[ k_i := \#\{ j \mid \lambda_j =i\} = \text{ the number of rows width }i
        \text{ in } \lambda \; . 
\] 
Then \(\ds [k_1,k_2,\ldots] \) 
satisfies  
\begin{align*}
        & k_1+ k_2 + \cdots + k_j + \cdots + k_{\ell} =
        \ell(\lambda)=\text{length of }\lambda \;, \\
        \noalign{and}
        & k_1 + 2 k_2 + \cdots + j k_j + \cdots + \ell k_{\ell} = |\lambda| =
        \sum_{i=1}^{\ell(\lambda)}\lambda_i = \text{the total area of }\lambda
        \;. 
\end{align*}
Thus, our cochain space $\ds \myCC{m}{w}$ corresponds to the Young
diagrams of height $m$ and the weight $w+(2-h)m$ under the dimensional
condition. 
%(\ref{eq::three:one}) and (\ref{eq::three:two}). 
%
%
\subsection{Manipulation of Young diagrams} 
%\begin{kmRemark} \label{rmk:zero} 
Let $\ds\mynabla{A}{k}$ be the set of Young diagrams whose total area (the
number of cells) is $A$ and the height is $k$.  Especially, $\ds
\mynabla{A}{1} = \{ \yng(3)\cdots\yng(3) \}$ (the row of width $A$) and
$\ds\mynabla{A}{A} = \{ \yng(1,1,1,1,1) \}$ (the column of height $A$).
If $A\leq 0$ or $k\leq 0$ or $ A < k$ then $\ds\mynabla{A}{k} = \emptyset$.  
Denoting the element of 
$\ds\mynabla{A}{A} $  by $\ds\T{A}$, i.e.,      
$\ds\mynabla{A}{A} = \{ \T{A}\}$,  
we have the following recursive formula;
\begin{align} \mynabla{A}{k} =& 
\T{k} \cdot ( \mynabla{A-k}{0} \sqcup \mynabla{A-k}{1} \sqcup 
\cdots \sqcup \mynabla{A-k}{k} ) \label{rmk:recur}
\end{align} 
where ``$\cdot$'' above means distributive concatenating operation of
the tower $\ds\T{k}$ and other Young diagrams.  (In fact, the series of
$\sqcup$ stops at $\min(k,A-k)$, however, the above notation may not
cause any confusion.)  It is also convenient to regard $\ds\mynabla{0}{0}$
as the single set of the unital element and $\ds\mynabla{A}{0}$ ($A>0$) or
$\ds\mynabla{A}{k}$ ($A<k$ or $k<0$) as the single set of the null element
of the concatenation ``$\cdot$''. We see that  $\ds\mynabla{A}{1} = \{
\T{1}^A\}$.  Using this operation, we can list up the elements of the
set $\ds\mynabla{A}{k}$.   For example, we have the following:   
\begin{align*}
        \mynabla{5}{3} =& \T{3} \cdot ( \mynabla{2}{0} \sqcup 
         \mynabla{2}{1} \sqcup \mynabla{2}{2} \sqcup \mynabla{2}{3}) 
         = \T{3} \cdot ( \mynabla{2}{1} \sqcup \mynabla{2}{2} ) \\ 
 =& \T{3}\cdot \T{1} \cdot ( \mynabla{1}{0} \sqcup \mynabla{1}{1} ) \sqcup 
\T{3}\cdot \T{2} \cdot ( \mynabla{0}{0} \sqcup \mynabla{0}{1} \sqcup \mynabla{0}{2})
         \\
 =& \T{3}\cdot \T{1} \cdot \mynabla{1}{1}  \sqcup 
         \T{3}\cdot \T{2} \cdot \mynabla{0}{0} 
         = \{ \T{3}\cdot \T{1}^{2} , 
 \T{3}\cdot \T{2} \} = \{ \yng(3,1,1), \yng(2,2,1) \}.
\end{align*} 
If we decompose $A$ as $ A = ak+b$ where $a>0$ and $0\leq b <k$, then we
have 
\begin{align*} \mynabla{A}{k} =& 
\T{k} \coprod_{j_a\leq \cdots \leq j_1 \leq k}
\T{j_1}\cdots \T{j_{a-1}} \cdot 
 \mynabla{A-k- \sum_{s=1}^{a-1}j_s}{j_a} \;.
\end{align*} 
In (\ref{rmk:recur}), we replace $A$ by $A+1$ and $k$ by $k+1$, and we
have 
\[
\mynabla{A+1}{k+1} = 
\T{k+1} \cdot ( \mynabla{A-k}{0} \sqcup \mynabla{A-k}{1} \sqcup 
\cdots \sqcup \mynabla{A-k}{k} \sqcup \mynabla{A-k}{k+1} 
). \] 
If we rewrite (\ref{rmk:recur}) formally 
\[ 
 \mynabla{A-k}{0} \sqcup \mynabla{A-k}{1} \sqcup 
\cdots \sqcup \mynabla{A-k}{k}  = \T{k}^{-1}\cdot \mynabla{A}{k} \]
then we get another recursive formula; 
\begin{align}
        \mynabla{A+1}{k+1} &= \T{k+1}\cdot\T{k}^{-1} \cdot \mynabla{A}{k} \sqcup
        \T{k+1}\cdot \mynabla{A-k}{ k+1})\;. \label{eqn:like:binom}
        \\
\noalign{Denoting $\ds\T{k}^{-1}\cdot \mynabla{A}{k}$ by
        $\ds\wnabla{A}{k}$, we have another form}
        \wnabla{A+1}{k+1} &= 
        \wnabla{A}{k} \sqcup \T{k+1} \cdot  
        \wnabla{A-k}{k+1} \;.    
        \label{hatnabla:recur}
\end{align} 
$\wnabla{A}{k}$ satisfies  
\begin{equation} 
\wnabla{k}{k} =\{id\}\;, \; 
\wnabla{k}{1} =\{\T{1}^{k-1}\}\;,\; 
\wnabla{A}{k} =\{ 0 \} \ \text{if}\ A < k \;. 
\label{eqn:hatnabla:property}
\end{equation}

We see how the formula (\ref{hatnabla:recur}) works on the same example: 
\begin{align*}
        \wnabla{5}{3} =& 
        \wnabla{4}{2} \sqcup \T{3} \cdot 
        \wnabla{2}{3}  = 
        \wnabla{3}{1} \sqcup  
        \T{2}\cdot \wnabla{2}{2} 
        = \{ \T{1}^2 , \T{2} \} \\
        \mynabla{5}{3} =& \T{3}\cdot 
        \wnabla{5}{3} = \T{3}\cdot \{
\T{1}^2 ,\T{2} \} = \{ \T{3}\T{1}^2 , \T{3}\T{2}\}\;.  
\end{align*} 

In (\ref{eqn:like:binom}), $\ds \T{k+1}\cdot \T{k}^{-1}$ means adding one cell
under the bottom of the left-most column of a Young diagram, we may another
form of (\ref{eqn:like:binom}): 
\begin{equation} 
        \mynabla{A+1}{k+1} = \myB  \cdot \mynabla{A}{k} \sqcup
        \T{k+1}\cdot \mynabla{A-k}{ k+1} \label{eqn:add:one}
\end{equation}
where $\myB$ is the operation adding one cell
under the bottom of the left-most column of each Young diagram. 

Using 
(\ref{eqn:like:binom}), (\ref{hatnabla:recur}) or (\ref{eqn:add:one}),  
we may write each Young diagram as  
\(\ds  
\T{\ell_1}\cdot 
\T{\ell_2 } \cdots  
\T{\ell_{s-1}} \cdot  
\T{\ell_{s}}\) with 
$\ell _1 \geq \ell_2 \geq \cdots \geq \ell _{s-1} \geq \ell _{s} > 0$.  
We may call it the decomposition by towers (or icicles) $\{\ds\T{j}\}$.   

For a given Young diagram $\lambda$,    
to get its tower decomposition 
\begin{equation} \label{eqn:towers} 
\T{\ell_1}\cdot 
\T{\ell_2 } \cdots  
\T{\ell_{s-1}} \cdot  
\T{\ell_{s}}\end{equation}
with $\ell _1 \geq \ell_2 \geq \cdots \geq \ell _{s-1} \geq \ell _{s} > 0$, 
is visually obvious, just slicing $\lambda$ vertically.   
A mathematical formula is given by
\begin{align} 
        \label{trad:to:tower}
        \ell_{j} = \#\{ i \mid  \lambda_i-(j-1) > 0 \} = 
        \#\{ i \mid  \lambda_i \geq j \} \quad ( j=1,\ldots, \lambda_1)
        \;, 
        \end{align} 
        and this says that the tower decomposition $\ds (\ell_1,\ldots,
        \ell_s) $ is equal to the conjugate Young diagram of $\lambda$.      
\kmcomment{
Conversely, assume we are given a tower decomposition (\ref{eqn:towers}) with 
$\ell _1 \geq \ell_2 \geq \cdots \geq \ell _{s-1} \geq \ell _{s} > 0$. 

\begin{equation}\label{tower:to:trad}
\T{\ell_1}\cdot \T{\ell_2 } \cdots  \T{\ell_{s}}\mapsto 
         {}^t ( \ell _1, \ell_2,\ldots,\ell_{s}) 
        (= 
        \text{the conjugate of } ( \ell _1, \ell_2,\ldots,\ell_{s})) = \lambda 
        \quad\text{as Young diagrams.} 
\end{equation}
%is a Young diagram whose tower decomposition is (\ref{eqn:towers}). 
}%endOFkmcomment

\kmcomment{
Already we observed that if a Young diagram $\lambda $ is denoted by 
$\ds [k_1,\ldots, k_p]$, 
then $\ds T(m)\cdot \lambda$ is given by 
\begin{equation}\label{eqn:chuukan}
 k_1 = m - \sum_{j=1}^p k_j, \quad k_2 = k_1, \quad \ldots, k_{p+1} =
 k_p\end{equation}

We regard the right $(s-1)$ terms in (\ref{eqn:towers}) as $\lambda$,  and
apply (\ref{eqn:chuukan}), we get 
}%endOFkmcomment

\begin{Lemma}
The expression
$\ds [k_1,\ldots, k_p]$ 
        in our notation of the Young diagram \eqref{eqn:towers} 
is given by  
\begin{equation}
\label{tower:to:our} 
k_1 = \ell_{1} - \ell_{2} \;, \quad 
k_2 = \ell_{2} - \ell_{3} \;, \quad
\ldots\;, \quad  
k_{s-1} = \ell_{s-1} - \ell_{s} \;,\quad 
k_{s} = \ell_{s}\;.  
\end{equation}
\kmcomment{
\begin{align}
\notag k_1 &= \ell_{1} - \ell_{2} \\
\notag k_2 &= \ell_{2} - \ell_{3} \\
\label{tower:to:our} & \vdots \\
        \notag k_{s-1} &= \ell_{s-1} - \ell_{s} \\
\notag k_{s} &= \ell_{s} 
\end{align}
}%endOFkmcommet 
Proof: 
Let us prepare a general Young diagram and slice it vertically.  If we compare 
the height of the nearby two towers, then (\ref{tower:to:our}) is clear.  

If you prefer a proof by induction, then follow the next: 
When $s=1$, $\ds\T{\ell_1}$ is given by $k_1 = \ell_1$ and $k_j=0$ for
$j>1$ and (\ref{tower:to:our}) holds. 

When $s=2$, $\ds\T{\ell_1}\cdot\T{\ell_2}$ is given by $k_1 = \ell_1 -\ell_2$,
$k_2=\ell_2$ and $k_j=0$ for $j>2$, and (\ref{tower:to:our}) holds. 

Now assume  (\ref{tower:to:our}) holds for general width $s$ 
and consider width $s+1$ YD $\ds \T{\ell_1}\cdots \T{\ell_s} \T{\ell_{s+1}}$. 
Let us take 
$\ds \lambda = 
 \T{\ell_2}\cdots \T{\ell_s} \T{\ell_{s+1}}$. Then 
   (\ref{tower:to:our}) says that 
   \[ 
k_1 = \ell_{2} - \ell_{3}\;,   
k_2 = \ell_{3} - \ell_{4}\;, \cdots,\;   
        k_{s-1} = \ell_{s} - \ell_{s+1}\;,  
   k_{s} = \ell_{s+1}\] 
Then our notation
of $\ds
\T{\ell_{1}}\cdot\lambda$  is given by 
        $\ds \bar{k}_1 = \ell_{1} - \sum_{j=1}^s k_j  = \ell_{1} - \ell_{2}$ 
        and 
$\ds \bar{k}_j = k_{j-1} = \ell_{j} -\ell_{j+1}$ for ($j=2,\ldots, s$) and 
$\ds \bar{k}_{s+1} = k_{s} = \ell_{s+1} $.   
Thus (\ref{tower:to:our}) holds for $s+1$. \hfill \kmqed%Q.E.D.  
\end{Lemma}

Solving 
(\ref{tower:to:our}) reversely, we have 
\begin{equation}
        \label{our:to:tower}
        \ell_j = \sum_{i\geq j} k_i \quad (j=1,2,\ldots)
\end{equation}

\begin{kmCor}\label{cor:chuukan} 
If a Young diagram $\lambda $ is denoted by 
$\ds [k_1,\ldots, k_p]$, then $\ds T(m)\cdot \lambda$ is given by 
\begin{equation} 
        \label{eqn:chuukan}
        %\label{eqn:concat}
 \bar{k}_1 = m - \sum_{j=1}^p k_j, \quad \bar{k}_2 = k_1, \quad
 \ldots,\quad 
 \bar{k}_{p+1} = k_p\end{equation} 
\end{kmCor}
\textbf{Proof:}
Since $\ds \lambda=  [k_1,\ldots, k_p]= \T{\ell_1}\cdots \T{\ell_{p}}$ 
with (\ref{our:to:tower}), $\T{m}\cdot \lambda$ is given by  
$\ds \T{m}\cdot \T{\ell_1}\cdots \T{\ell_{p}}$. From (\ref{tower:to:our}), 
the notation is \begin{align*}
        \bar{k}_1 =& m - \ell_1 = m - \sum_{j\geq 1}k_j\;, 
        \quad 
        \bar{k}_2 = \ell_1  - \ell_2 - \sum_{j\geq 1}k_j 
        - \sum_{j\geq 2}k_j = k_1 \;, \\
        & \ldots , \quad 
        \bar{k}_{p+1} = \ell_p  - \ell_{p+1} - \sum_{j\geq p}k_j 
        - \sum_{j\geq {p+1} }k_j = k_p \; . %\hspace{50mm}\rule{1ex}{1.5ex}
        \hspace{30mm}\kmqed%Q.E.D.  
        \end{align*}
%\bigskip
\kmcomment{
The transformations among three notations: 
\begin{center}
\resizebox{0.6\textwidth}{!}{
\(
\xymatrix{
        \text{our } \ar[r]^{(\ref{our:to:trad})} & \text{trad}
        \ar[dl]^{(\ref{trad:to:tower})} \\
        \text{tower }\ar[u]^{(\ref{tower:to:our})} 
    &
    } 
    \hspace{10mm} 
\xymatrix{
\text{our } \ar[d]_{(\ref{our:to:tower})} & \text{trad } \ar[l]_{(\ref{trad:to:our})} \\
        \text{tower }\ar[ur]_{(\ref{tower:to:trad})} 
    &
    }
    \) 
    }
    \end{center}
%\bigskip 
}%endOFkmcomment
%\bigskip 
\kmcomment{
When a Young diagram $\lambda$ is given, in our notation by index $\ds
(k_{\ell},\dots,k_1)$, the decomposition of $\lambda$ into towers is the
following; The $j$-th tower from the left is $\ds\T{ \sum_{i=j}^{\ell}
k_i}$ ($j=1,\ldots,\ell$). Thus, the Young diagram $\ds\T{h}\cdot
\lambda$ is given by 
\begin{equation} 
        k'_1 = h - \sum_{i=1}^{\ell}k_i, \quad 
        k'_2 = k_1, \quad \cdots \quad ,k'_{\ell+1} = k_{\ell}\;.
        \label{eqn:concat}
\end{equation} 
}%endOFkmcomment
%\end{kmRemark}

\subsection{Examples of decomposition of cochain complex according to weights} 
We will now show some examples of decomposition of a cochain complex
according to weights.  
\subsubsection{case $h=1$}   
Since $h=1$, for a given weight $w$, the $m$-th cochain space $\ds
\myCC{m}{w}$ corresponds to the set $\ds\mynabla{w+m}{m}$ of Young diagrams
of area $w+m$ and the height $m$, and  we have  $m \leqq w + n$ as in
the table in \S~\ref{subsec:possible}.  

If $w=0$, $\ds\mynabla{m}{m} = \{\T{m}\}$, this means $\ds k_1=m$ and $\ds
k_j=0$ ($j>1$). Thus $\ds\myCC{m}{0} = \Lambda^m \frakS{1}$.  

If $w=1$, $\ds\mynabla{1+m}{m} = \{ \T{m}\cdot \T{1}\} $, this means $\ds
k_1=m-1$, $k_2=1$ and  $\ds k_j=0$ ($j>2$). Thus $\ds\myCC{m}{1} =
\Lambda^{m-1} \frakS{1}  \otimes \frakS{2}$.  

If $w=2$, $\ds\mynabla{2+m}{m} = \T{m}\cdot (\mynabla{2}{1} \sqcup
\mynabla{2}{2}
) = \T{m}\cdot \{\T{1}^2, \T{2}\} $, this means $\ds k_1=m-1$, $k_3=1$
and  $\ds k_j=0$ ($j \ne 1,3$) or $\ds k_1=m-2$, $k_2=2$ and  $\ds
k_j=0$ ($j > 2$).  
Thus $\ds\myCC{m}{2} = \Lambda^{m-1} \frakS{1}  \otimes \frakS{3}  +
\Lambda^{m-2} \frakS{1}  \otimes \Lambda^{2} \frakS{2}$.    

If $w=3$, 
\begin{align*} 
        \mynabla{3+m}{m} = & 
\T{m}\cdot (\mynabla{3}{1} \sqcup \mynabla{3}{2} \sqcup \mynabla{3}{3}) 
= \T{m}\cdot \{\T{1}^3 , \T{2}\cdot \T{1} , \T{3} \} 
= \{ \T{m} \cdot \T{1}^3, \T{m} \cdot \T{2}\cdot \T{1} ,\T{m} \cdot \T{3} 
\}
\end{align*}
this means
$\ds k_1=m-1$, $k_4=1$ and  $\ds k_j=0$ ($j \ne 1,4$),  
$\ds k_1=m-2$, $k_2=1$, $k_3=1$ and  $\ds k_j=0$ ($j > 3$) or 
$\ds k_1=m-3$, $k_2=3$ and  $\ds k_j=0$ ($j > 2$).  
Thus $\ds\myCC{m}{3} = 
\Lambda^{m-1} \frakS{1}  \otimes \frakS{4}  
+ \Lambda^{m-2} \frakS{1}  \otimes \frakS{2}  \otimes \frakS{3}
+ \Lambda^{m-3} \frakS{1}  \otimes \Lambda^{3} \frakS{2} $.    
Summarizing the expressions we got above, 
\begin{align} 
        \myCC{m}{0} =&  \Lambda^m \frakSS{1}\;. \label{h1:w0}\\  
        \myCC{m}{1} =&  
        \Lambda^{m-1} \frakSS{1}  \otimes \frakSS{2}\;, \label{h1:w1}\\   
        \myCC{m}{2} =&  
\Lambda^{m-1} \frakSS{1}  \otimes \frakSS{3}  
+ \Lambda^{m-2} \frakSS{1}  
\otimes \Lambda^{2} \frakSS{2}\;, \label{h1:w2}\\    
        \myCC{m}{3} =& 
\Lambda^{m-1} \frakSS{1}  \otimes \frakSS{4}  
+ \Lambda^{m-2} \frakSS{1}  \otimes \frakSS{2}  \otimes \frakSS{3}
+ \Lambda^{m-3} \frakSS{1}  \otimes \Lambda^{3} \frakSS{2}\;. \label{h1:w3}    
\end{align}

%OK \end{spacing} \end{document}

\subsubsection{case $h=1$ and $n=3$}
We restrict ourselves to $\ds  \mR^3$ and $h=1$. Then   
$m \leqq w + n = w+3$ as in \S~\ref{subsec:possib} and 
we obtain the following
\begin{exam}[n=3, h=1]
\begin{align*}
        \myCC{m}{0} =& \Lambda^{m}\frakSS{1} \quad
		(\tbinom{3}{m}\dim)\;, \\
\myCC{1}{1} =& \frakSS{2}\ (6\;\dim), \quad  \myCC{2}{1} = \frakSS{1} \otimes
\frakSS{2}\ (18\; \dim), \quad 
\myCC{3}{1} = \Lambda^2 \frakSS{1} \otimes \frakSS{2} \ (18\;\dim),\quad 
\\&\hspace{26mm}
\myCC{4}{1} = \Lambda^3 \frakSS{1} \otimes \frakSS{2}\ (6\;\dim),  
\\
\myCC{1}{2} =& \frakSS{3}\ (10\;\dim),  \quad  
\myCC{2}{2} = \frakSS{1} \otimes \frakSS{3} \oplus \Lambda^2 \frakSS{2}\
(45\;\dim), \quad  
\myCC{3}{2} = \Lambda^2 \frakSS{1} \otimes \frakSS{3}\oplus \frakSS{1}
\otimes \Lambda^2 \frakSS{2}\ (75\;\dim),  \\ 
\myCC{4}{2} =& \Lambda^3 \frakSS{1} \otimes \frakSS{3} \oplus \Lambda^2
\frakSS{1} \otimes \Lambda^2 \frakSS{2}\ (55\;\dim), \quad  
\myCC{5}{2} = \Lambda^3 \frakSS{1} \otimes \Lambda^2 \frakSS{2} \
(15\;\dim),  
\\ 
\myCC{1}{3} =& \frakSS{4} \ (15\;\dim),\quad  
\myCC{2}{3} = \frakSS{1} \otimes \frakSS{4} \oplus 
\frakSS{2} \otimes \frakSS{3}\  (105\;\dim),\\  
\myCC{3}{3} =& 
 \Lambda^2 \frakSS{1} \otimes \frakSS{4}
\oplus 
 \frakSS{1} \otimes \frakSS{2} \otimes \frakSS{3} 
 \oplus \Lambda^3 \frakSS{2}  
 \ (245\;\dim),  
 \\ 
\myCC{4}{3} =& 
 \Lambda^3 \frakSS{1}\otimes\frakSS{4}
\oplus 
\Lambda^2  \frakSS{1}\otimes\frakSS{2}\otimes\frakSS{3} 
\oplus
\frakSS{1}\otimes \Lambda^3 \frakSS{2} 
 \ (255\;\dim), \\  
\myCC{5}{3} =& 
\Lambda^3 \frakSS{1} \otimes  \frakSS{2} \otimes \frakSS{3}
\oplus  
\Lambda^2 \frakSS{1} \otimes \Lambda^3 \frakSS{2} \
(120\;\dim), \quad  
\myCC{6}{3} = \Lambda^3 \frakSS{1} \otimes \Lambda^3 \frakSS{2}  \ (20\;\dim) 
. 
\end{align*} 
\end{exam}

%OK \end{spacing} \end{document}

\subsection{Concrete Poisson structure when $\ds n=3$ and $h=1$: 
 \label{Conc::Eg}}  
We choose a specified Poisson structure, namely, we consider the Lie
algebra $\fraksl(2)$ and consider the Lie Poisson structure.  The
Poisson bracket is defined by 
$ \Bkt{p}{q}=r, 
\Bkt{r}{p}=2p, 
\Bkt{r}{q}=-2q$, where $p,q,r$ are the coordinates in $\mR^3$, namely,     
\[  \Bkt{F}{H} = \begin{vmatrix} 2q & 2p & r \\
        F_{p} & F_{q} & F_{r} \\
        H_{p} & H_{q} & H_{r} \end{vmatrix}
= \frac{ \pdel ( 2pq+r^2/2,F,H)}{ \pdel( p, q, r)}
\;. 
\]  

We use the notations $\ds \w{A} = p^{a_1} q^{a_2} r^{a_3}$ for each triple
$\ds A=(a_1,a_2,a_3)$ of non-negative integers. Then $\ds\{ \w{A} \mid
|A|(:= a_1+a_2+a_3) = k\}$ is a basis of $\mySS{k}$ and the dual
basis is given by $\ds\{ \z{A} \mid |A|=k\}$.  The
coboundary operator $\mydov$ for 1-cochains is defined as 
\[  \mydov (\z{C} )  = - \frac{1}{2} \sum _{A,B} \langle \z{C},
\Bkt{\w{A}}{\w{B} }\rangle
\z{A}\wedge \z{B}\;.\]
Since 
\begin{align*}
        \Bkt{\w{A}}{\w{B}} &=  
 \begin{vmatrix} 2q & 2p & r \\
         \ds\frac{a_1}{p} \w{A}  & \ds\frac{a_2}{q}\w{A}
         &\ds\frac{a_3}{r}\w{A} \\
\ds      \frac{b_1}{p} \w{B} & \ds\frac{b_2}{q}\w{B}
         &\ds\frac{b_3}{r}\w{B} 
 \end{vmatrix} = \frac{ \w{A+B} }{pqr } 
 \begin{vmatrix} 2pq & 2pq & r^2 \\
         a_1  & a_2  & a_3 \\
         b_1  & b_2  & b_3 \end{vmatrix} 
 \\ & 
=  2 \w{A+B-\epsilon_3} ( 
 \begin{vmatrix} a_2 & a_3 \\ b_2 & b_3\end{vmatrix}
 + \begin{vmatrix} a_3 & a_1 \\ b_3 & b_1\end{vmatrix}) 
 +  \w{A+B-\epsilon_1 -\epsilon_2 +\epsilon_3}   
 \begin{vmatrix} a_1 & a_2 \\ b_1 & b_2\end{vmatrix} 
\end{align*}
where $\ds\epsilon_1=(1,0,0),\epsilon_2=(0,1,0),\epsilon_3=(0,0,1)$, 
we have 
\begin{align}
        \mydov (\z{C} ) = - \sum_{A+B=C+\epsilon_3} ( 
 \begin{vmatrix} a_2 & a_3 \\ b_2 & b_3\end{vmatrix} + 
         \begin{vmatrix} a_3 & a_1 \\ b_3 & b_1\end{vmatrix}) \z{A}
                 \wedge \z{B} 
 - \frac{1}{2} \sum_{A+B = C 
 +\epsilon_1 +\epsilon_2 -\epsilon_3}    
        \begin{vmatrix} a_1 & a_2 \\ b_1 & b_2\end{vmatrix} \z{A} \wedge
                \z{B}\;.
\label{one:two} 
\end{align} 
For example, let $C=(1,0,0)$. In the first term of (\ref{one:two}), $C +
\epsilon_3 = (1,0,1)$ and so we find two summands corresponding to $A=(1,0,0),
B=(0,0,1)$ or  $A=(0,0,1), B=(1,0,0)$.   In the second term, $C+
\epsilon_1+  \epsilon_2-  \epsilon_3 = (2,1,-1)$ and no summand
corresponding to 
$B$ and $C$ with $B+C = (2,1,-1)$. Thus, we get 
\[ \ds\mydov (z_{1,0,0}) = -  
( \begin{vmatrix} 0 & 1 \\ 0 & 0\end{vmatrix} + 
 \begin{vmatrix} 1 & 0 \\ 0 & 1\end{vmatrix}
) z_{0,0,1} \wedge z_{1,0,0}
- 
( \begin{vmatrix} 0 & 0 \\ 0 & 1\end{vmatrix} + 
 \begin{vmatrix} 0 & 1 \\ 1 & 0\end{vmatrix}
 ) z_{1,0,0} \wedge z_{0,0,1} = -2 z_{0,0,1} \wedge z_{1,0,0}\ . 
\] 
Similarly we get $\ds\mydov (z_{0,1,0}) = 2 z_{0,0,1} \wedge
z_{0,1,0}$. For $\ds C = (0,0,1)$, in the first term $C+\epsilon _3 =
(0,0,2)$ and so $A=B=(0,0,1)$ and $\ds z_A\wedge z_B = 0$. In the second
term, since $A+B=(1,1,0)$ we have $A=(1,0,0)$ and  $B=(0,1,0)$, or
$A=(0,1,0)$ and  have $B=(1,0,0)$, and so 
\[\ds\mydov ( z_{0,0,1} ) = 
- \frac{1}{2} 
\begin{vmatrix} 0 & 1 \\ 1 & 0\end{vmatrix}  z_{0,1,0} \wedge z_{1,0,0}  
- \frac{1}{2} 
\begin{vmatrix} 1 & 0 \\ 0 & 1\end{vmatrix}  z_{1,0,0} \wedge z_{0,1,0}  
= z_{0,1,0} \wedge z_{1,0,0} \ .\]     
As a summary, we obtain
\begin{align}
\mydov (z_{1,0,0}) &= -2 z_{0,0,1} \wedge z_{1,0,0}\;, \label{d:one}
\\
\mydov (z_{0,1,0}) &= 2 z_{0,0,1} \wedge z_{0,1,0}\;, \label{d:two}
\\ 
\mydov ( z_{0,0,1} ) &= z_{0,1,0} \wedge z_{1,0,0}\;.  \label{d:three}
\end{align}

By a similar argument, for $\ds\mydov ( z_A )$ for $|A|=2$  we obtain
the following result: 
\begin{align} 
\mydov \left( z_{{2,0,0}} \right) \, = &  \,-4\,
  z_{{0,0,1}} \wedge z_{{2,0,0}}  +2\,
  z_{{1,0,0}} \wedge z_{{1,0,1}} \;, 
\\
\mydov \left( z_{{0,2,0}} \right) \, = & \,4\,
          z_{{0,0,1}}\wedge z_{{0,2,0}}  -2\,
  z_{{0,1,0}}\wedge z_{{0,1,1}} \;,   
 \\
\mydov \left( z_{{0,0,2}} \right) \, = &  \,
  z_{{0,1,0}}\wedge z_{{1,0,1}}  +
  z_{{0,1,1}} \wedge z_{{1,0,0}} \;,  
\\
\mydov \left( z_{{0,1,1}} \right) \, = &  \,2\,
  z_{{0,0,1}}\wedge z_{{0,1,1}}  +4\,
 z_{{0,0,2}} \wedge z_{{0,1,0}}  + 
 z_{{0,1,0}} \wedge z_{{1,1,0}}  +2\, 
z_{{0,2,0}} \wedge z_{{1,0,0}} \;,  
\\ 
\mydov \left( z_{{1,0,1}} \right) \, = & \,-2\, 
                z_{{0,0,1}} \wedge z_{{1,0,1}}  -4\,
                 z_{{0,0,2}}\wedge z_{{1,0,0}}  +2\,
  z_{{0,1,0}}\wedge z_{{2,0,0}}  -
  z_{{1,0,0}} \wedge z_{{1,1,0}} \;,  
 \\ 
\mydov \left( z_{{1,1,0}} \right) \, = & \,-2\, 
z_{{0,1,0}}\wedge z_{{1,0,1}}  -2\,  z_{{0,1,1}}\wedge z_{{1,0,0}}  
\;. 
\end{align}
Now we can compute the cohomology groups $\ds 
\myHH{\bullet}{1} $ for the weight 1 by examining the kernel  
of $\ds \mydov: \myCC{m}{1} \rightarrow \myCC{m+1}{1}$.  
Take a 1-cochain $\sigma \in \myCC{1}{1} = \frakSS{2}$ and examine
$\mydov( \sigma)=0$.  
\begin{align*}
0 & = 
 c_{2,0,0} (-4
  z_{0,0,1} \wedge z_{2,0,0}  +2
  z_{1,0,0} \wedge z_{1,0,1})  
+ c_{0,2,0} (4
          z_{0,0,1}\wedge z_{0,2,0}  -2
  z_{0,1,0}\wedge z_{0,1,1}  ) 
  \\ & 
+ c_{0,0,2} (
  z_{0,1,0}\wedge z_{1,0,1}  +
  z_{0,1,1} \wedge z_{1,0,0} ) 
\\& 
+ c_{0,1,1} (  2
  z_{0,0,1}\wedge z_{0,1,1}  +4
 z_{0,0,2} \wedge z_{0,1,0}  + 
 z_{0,1,0} \wedge z_{1,1,0}  +2 
z_{0,2,0} \wedge z_{1,0,0} ) 
\\& 
+  c_{1,0,1} (-2 
                z_{0,0,1} \wedge z_{1,0,1}  -4
                 z_{0,0,2}\wedge z_{1,0,0}  +2
  z_{0,1,0}\wedge z_{2,0,0}  -
  z_{1,0,0} \wedge z_{1,1,0} ) 
  \\& 
 + c_{1,1,0} ( -2 
        z_{0,1,0}\wedge z_{1,0,1}  -2  z_{0,1,1}\wedge z_{1,0,0})  
		\;. 
\end{align*}
Taking the interior product by $\ds \w{\epsilon_i}$ ($i=1,2,3$) , we have 
\begin{align*}
& 
4 c_{1, 0, 1} z_{0, 0, 2}+(2 c_{1, 1, 0}-c_{0, 0, 2})z_{0, 1, 1}-2c_{0, 1, 1} z_{0,
2, 0}+2 c_{2, 0, 0} z_{1, 0, 1}-c_{1, 0, 1} z_{1, 1, 0} = 0\;, \\
& 
-4 c_{0, 1, 1} z_{0, 0, 2}-2 c_{0, 2, 0} z_{0, 1, 1}+(-2 c_{1, 1, 0}+c_{0, 0,
2}) z_{1, 0, 1}+c_{0, 1, 1} z_{1, 1, 0}+2 c_{1, 0, 1} z_{2, 0, 0}
=0\;, \\
&
-2 c_{1, 0, 1} z_{1, 0, 1}+2 c_{0, 1, 1} z_{0, 1, 1}+4 c_{0, 2, 0} z_{0, 2,
0}-4 c_{2, 0, 0} z_{2, 0, 0}=0\;. 
\end{align*} 
Thus, we get 14 linear equations of $\ds c_A$  ($|A|=2$) and solving
them we see that $\ds c_{0,0,2} = 2 c_{1,1,0}$ and $\ds c_A =0$ ($A\ne
(0,0,2) , (1,1,0)$ and $|A|=2$), i.e., $\ds\sigma = c_{1,1,0}(
z_{1,1,0}+ 2 z_{0,0,2})$, namely,  the kernel of $\ds\mydov :
\myCC{1}{1} \rightarrow \myCC{2}{1}$ is spanned by $\ds z_{1,1,0}+ 2
z_{0,0,2}$, thereby the first Betti number is 1, and  the rank of $\ds
\mydov : \myCC{1}{1} \rightarrow \myCC{2}{1}$ is 5 because of $\ds\dim
\myCC{1}{1} = \dim \frakSS{2} = 6$. 

Next, we consider the second cochain space.  Take a 2-cochain $\ds
\sigma = \sum_{|\alpha|=1, |A|=2} c_{\alpha,A} z_{\alpha} \wedge z_{A}$.  
Applying the interior products $\ds i_{ \w{\alpha}} \circ i_{ \w{\beta}}
\circ i_{ \w{A}}$  with $ | \alpha | = |\beta|=1$ and $\ds |A|=2$ to 
$\ds\sum_{|\alpha|=1, |A|=2} c_{\alpha,A}\; \mydov( z_{\alpha} \wedge
z_{A}) = 0$, we get 17 linear equations of 18 variables.  Solving this
equations, the kernel of $\ds\mydov : \myCC{2}{1} \rightarrow
\myCC{3}{1}$ is linearly spanned by the following 5 terms.  
\begin{align*} 
& -2 z_{0,0,1} \wedge  z_{2,0,0}+ z_{1,0,0} \wedge z_{1,0,1},\quad 
z_{0,1,0} \wedge  z_{0,1,1}-2 z_{0,0,1} \wedge  z_{0,2,0},\quad 
-z_{1,0,0}\wedge z_{0,1,1}+z_{0,1,0} \wedge z_{1,0,1}, 
\\& 
z_{1,0,0} \wedge ( -4 z_{0,0,2}+   
z_{1,1,0})-2 z_{0,1,0}\wedge z_{2,0,0}+2 z_{0,0,1}\wedge z_{1,0,1}, 
\\& 
-2 z_{1,0,0}\wedge z_{0,2,0}+ z_{0,1,0}\wedge (-4 z_{0,0,2}+
\wedge z_{1,1,0}) +2 z_{0,0,1} \wedge  z_{0,1,1}
\;. 
\end{align*} 
Thus the kernel is 5-dimensional and the rank of $\ds\mydov :
\myCC{2}{1} \rightarrow \myCC{3}{1}$ is $13 (= 18-5)$. Thereby, the second
Betti number is 0.  
By the same method, the kernel of $\ds\mydov : \myCC{3}{1} \rightarrow
\myCC{4}{1}$ is 13 dimensional and so the third Betti number is 0 and
the rank is 5.  

Indeed, take an arbitrary 3-cochain 
$\ds\sigma = \sum_{\alpha,\beta,A} 
c_{\alpha,\beta,A} z_{\alpha} \wedge z_{\beta} \wedge z_{A}$  ( $\ds
|\alpha|=|\beta|=1$, $\ds |A|=2$). Then we have 
\begin{align*}
\mydov(\sigma) = &  ( c_{1,3,[0, 1, 1]}+ c_{2,3,[1, 0, 1]})  
 z_{1,0,0}\wedge z_{0,0,1}\wedge z_{0,1, 0}\wedge ( -4 z_{0,0,2}+ 
z_{1,1,0}) 
\\& 
+(-2 c_{1,2,[0,1,1]}-2 c_{1,3,[0,2,0]}-2 c_{2,3,[1,1,0]}+c_{2,3, [0,0,2]}) 
(z_{1,0,0}\wedge z_{0,0,1}\wedge z_{0,1,0}\wedge z_{0,1,1}) 
\\&
+(-4 c_{1,2,[0, 2, 0]}+2 c_{2,3,[0, 1, 1]}) 
(z_{1,0,0}\wedge z_{0,0,1}\wedge z_{0,1,0}\wedge z_{0,2,0}) 
\\&
+(2 c_{1,2,[1,0,1]} -2 c_{1,3,[1,1,0]}+c_{1,3,[0,0,2]}-2 c_{2,3,[2,0,0]})  
(z_{1,0,0}\wedge z_{0,0,1 }\wedge z_{0,1,0}\wedge z_{1,0,1}) 
\\&
+(4 c_{1,2,[2, 0, 0]}+2 c_{1,3,[1, 0, 1]}) 
(z_{1,0,0}\wedge z_{0 ,0,1}\wedge z_{0,1,0}\wedge z_{2,0,0})\;. 
\end{align*}
Since $\ds\myCC{5}{1}= (0)$, the kernel of 
$\ds\mydov : \myCC{4}{1} \rightarrow \myCC{5}{1}$ has the dimension 
$\ds\dim \myCC{4}{1} = \dim (\Lambda^3 \frakSS{1} \otimes \frakSS{2}) = 6$, so 
\[\ds\myHH{4}{1} \cong \text{LSpan}(z_{1,1,0},z_{0,0,2})/ \text{LSpan}(
z_{1,1,0}-4 z_{0,0,2})\] and  
the fourth  Betti number is 1. 
We summarize the discussion above into the table below.  
\begin{center}
\begin{tabular}{|c|*{9}{c}|} \hline 
wt=1 & $\ds\myCC{1}{1}$ & $\ds\mathop{\rightarrow}^{\mydov}$  &  
 $\ds\myCC{2}{1}$ & $\ds\mathop{\rightarrow}^{\mydov}$  &
 $\ds\myCC{3}{1}$ & $\ds\mathop{\rightarrow}^{\mydov}$  & 
 $\ds\myCC{4}{1}$ & $\ds\mathop{\rightarrow}^{\mydov}$  & 0\\\hline 
$\dim$ & 6 && 18 && 18 && 6  && \\
Ker $\dim$ & 1 && 5 && 13 && 6 && \\
rank && 5 && 13 && 5 && 0 & \\
Betti & 1 && 0 && 0 && 1 && \\\hline
\end{tabular} 
\end{center} 
It is well-known that $\ds\sum_{m>0} (-1)^m \dim \myCC{m}{w} =
\sum_{m>0} (-1)^m \dim \myHH{m}{w} $  holds and this number is called
the Euler characteristic . In our case above, the number is 0. Later, we will
show that this is true for 1-homogeneous Poisson structures in general.
\kmcomment{
In this paper, as the definition of Euler characteristic we sum up from
degree 1 and ignore the 0-th cochain complex $\ds\myCC{0}{w} = \mR$.
}%endOFkmcomment
If we want to continue  studying $\ds\myHH{\bullet}{w}$ for this
Poisson structure, we have to prepare $\ds\mydq ( z_A )$ further for $|A|\leq
w+1$.    

%OK \end{spacing} \end{document}
%%%%%%%%%%%%%%%%%%%%%%%%%%%%%%%%%%%%%%%%%%%%%%%%%%%%%%%%%%%
\medskip

Here we compute the Casimir polynomials: For $\ds F = \sum _{A} c_A
p^{a_1} q^{a_2} r^{a_3}$ where $A= (a_1,a_2,a_3)$ is a triple of
non-negative integers, and $\ds c_A$ are constant, we see that   
\[\ds\Bkt{r}{F} = \sum_{A} 2 c_A ( a_1 - a_2 ) p^{a_1} q^{a_2}
r^{a_3}\] 
so any Casimir polynomial should be $\ds F = \sum c_{i,i,k} (pq)^i r^k$.
By using the other condition $\ds\Bkt{p}{F}=0$, we see that 
\[ F = \sum c_k (4 p q + r^2)^k \] where $c_k$ are constant.  
Thus, $\ds\mySS{2k}$ contains a Casimir polynomial $\ds (4 pq +r^2)^k$
and $\ds\mySS{2k-1}$ does not contain any Casimir polynomials. 

For the weight =2, all the Betti numbers are trivial, and for the
weight=3 or 4, we see non-trivial Betti numbers in the tables below: 

\begin{center}
\tabcolsep=3pt
\begin{tabular}{|c|*{12}{c}c|}
\hline
$wt=3$
& $\ds  \myCC{1}{3}$ & $\rightarrow$ &
$\ds  \myCC{2}{3} $& $ \rightarrow $ &
$\ds  \myCC{3}{3} $& $ \rightarrow $ &
$\ds  \myCC{4}{3} $& $ \rightarrow $ &
$\ds  \myCC{5}{3} $& $ \rightarrow $ &
$\ds  \myCC{6}{3} $& $ \rightarrow $ & 
$\mathbf{0}$ \\
\hline
$\dim$  & 15  && 105 &&245 && 255 && 120 && 20 &&   \\
rank & & 14 && 91 && 154 && 100 && 20 &&0 & \\
$\dim$(ker) &  1 && 14 && 91 && 155 && 100 && 20 && \\
%$\dim$(image) &   && 0 && 14 && 91 && 154 && 100&& 20 &&\\
Betti num & 1 && 0 && 0 && 1 && 0 && 0 &&   \\
\hline
\end{tabular}

\bigskip

\tabcolsep=3pt
\begin{tabular}{|c|*{14}{c}c|}
\hline
$wt=4$
& $\ds  \myCC{1}{4}$ & $\rightarrow$ &
$\ds  \myCC{2}{4} $& $ \rightarrow $ &
$\ds  \myCC{3}{4} $& $ \rightarrow $ &
$\ds  \myCC{4}{4} $& $ \rightarrow $ &
$\ds  \myCC{5}{4} $& $ \rightarrow $ &
$\ds  \myCC{6}{4} $& $ \rightarrow $ &
$\ds  \myCC{7}{4} $& $ \rightarrow $ &
$\mathbf{0}$ \\
\hline
$\dim$  & 21  && 198 && 618 && 891 && 630 && 195&& 15 && \\
rank & & 21 && 176 && 441 && 450 && 179 &&15 && 0& \\
        $\dim$(ker) &   0 && 22 && 177 && 441 && 451 && 180 && 15 && \\
%$\dim$(image) &   && 0 && 14 && 91 && 154 && 100&& 20 &&\\
Betti num & 0 && 1 && 1 && 0 && 1 && 1 && 0 && \\
\hline
\end{tabular}
\end{center}

\medskip

%%%%%%%%%%%%%%%%%%%%%%%%%%%%%%%%%%%%%%%%%%%%%%%%%%%%%%%%%%%%%%%%%%%%%%%%%%%%%%%%
\begin{kmRemark} \label{rem:LieAlgItself} 
Given a finite dimensional Lie algebra $\frakh$, we have the so called
        Lie Poisson structure on $\ds\frakh^{*}$ as mentioned in Remark
        \ref{rem:LiePoisson}.  As Lie algebras,  $(\frakh, \Sbt{}{})$ is
        a subalgebra of $\ds (\sum_{k} \text{Sym}^k \frakh,
        \Pkt{\cdot}{\cdot})$.  
The cochain complex of Lie algebra $\frakh$ coincides with the weight
        zero cochain complex of $(\myCC{\bullet}{0}, \mydov)$ of Lie
        Poisson structure $\ds\frakh^{*}$.  The Euler characteristic of
        Lie algebra cohomology groups of $\frakh$ is zero and the Euler
        characteristic of $\ds \myHH{\bullet}{0}$ of  1-homogeneous
        Poisson structure with the weight 0 is also zero.  
\end{kmRemark}

%% file: Contrib.tex
\section{Contributions of Poisson structures}

We would like to know the concrete behavior of $\mydov$ for each weight.  
Since  $\mydov$ preserve weights, in the case of 
  $ h$-homogeneous structure, we have 
\[  \mydov( \frakSS{g} )
\subset \sum \frakSS{a} \wedge \frakSS{b}\]  
 where $g-2+h = (a-2+h) + (b-2+h)$, thus $a+b = g-h+2$.  

\subsection{$h=1$} 
When $h=1$, we see that 
\begin{align*}
  \mydov(\frakSS{1}) \subset &  \frakSS{1} \wedge \frakSS{1}, \qquad  
  \mydov(\frakSS{2}) \subset   \frakSS{1} \wedge \frakSS{2}, \qquad 
  \mydov(\frakSS{3}) \subset  \frakSS{1} \wedge \frakSS{3} \oplus  
                            \frakSS{2} \wedge \frakSS{2}
							\;. 
\end{align*}                           
These come from $\ds  \Bkt{\mySS{i}}{\mySS{j}} \subset
\mySS{i+j-1}$ in general. For instance,  
\begin{align*}
\Bkt {\mySS{1}}{ \mySS{1}} &\subset \mySS{1}\;, \quad  
\Bkt{\mySS{1}}{\mySS{2}} \subset \mySS{2}\;, \quad 
\Bkt{\mySS{1}}{\mySS{3}} \subset \mySS{3}\;, \quad 
\Bkt{\mySS{2}}{\mySS{2}} \subset \mySS{3}\;. 
\end{align*}

\subsubsection{$\pi = x_1 \pdel_2 \wedge \pdel_3+ x_2 \pdel_3 \wedge
\pdel_1+ x_3 \pdel_1 \wedge \pdel_2$, i.e., Lie-Poisson of
$\mathfrak{so}(3)$} As Lie algebras, 
$\mathfrak{so}(3)$ is isomorphic to $\ds\fraksl(2,\mR)$ and we have some data of
cohomology groups of lower weights as stated before.   

\subsubsection{$\pi = x_3 \pdel_1 \wedge \pdel_2$, i.e., Lie-Poisson of
the Heisenberg Lie algebra} 
The Casimir polynomials are $\ds\{ x_3^k\}$ ($k=1,2,\ldots$) and the
cohomology groups of two kinds for lower weights are follows:

\begin{center}
\tabcolsep=3pt 
\setlength{\extrarowheight}{-3pt}
\begin{tabular}{|c|*{6}{c}c|}
\hline
        $wt=0$ & $\ds\myCC{0}{0}$ & $\rightarrow $ &
         $\ds\myCC{1}{0}$ & $\rightarrow $ &
         $\ds\myCC{2}{0}$ & $\rightarrow $ &
         $\ds\myCC{3}{3}$ \\ \hline
         $\dim$       &    1 && 3 &&  3 && 1     \\
         $\dim(\ker)$&     1 &&  2  && 3 &&1     \\ 
         Betti num    &     1 &&  2  && 2  && 1    \\ \hline
\end{tabular}
\hfil 
\begin{tabular}{|c|*{6}{c}c|}
\hline
        $wt=0$ & $\ds\myC{0}{0}$ & $\rightarrow $ &
         $\ds\myC{1}{0}$ & $\rightarrow $ &
         $\ds\myC{2}{0}$ & $\rightarrow $ &
         $\ds\myC{3}{3}$ \\ \hline
         $\dim$       &    1 && 2 &&  1 && 0     \\
         $\dim(\ker)$&     1 &&  2  && 1 &&0     \\ 
         Betti num    &     1 &&  2  && 1  && 0    \\ \hline
\end{tabular} 
\end{center}

 \newcommand{\kmTemp}[2]{\bullet=#1}

\begin{center}
\tabcolsep=5pt
\setlength{\extrarowheight}{-3pt}
\begin{tabular}{|c|*{4}{c}|}
\hline
$wt=1$ &  
 $\ds\kmTemp{1}{1} $&
 $\ds\kmTemp{2}{1} $&
 $\ds\kmTemp{3}{1} $&
 $\ds\kmTemp{4}{1} $ \\ \hline
$\dim$ of $\myCC{\bullet}{1}$    & 6 & 18 & 18 & 6   \\
        $\dim(\ker)$     & 3 & 10 & 13 & 6   \\
Betti num       & 3 & 7  &  5 & 1   \\\hline
$\dim$ of $\myC{\bullet}{1}$ & 5 & 10 & 5 & 0   \\
        $\dim(\ker)$    & 3 & 7  & 5 & 0   \\
Betti num       & 3 & 5  &  2 & 0   \\\hline
\end{tabular}
\hfil 
%\medskip
% \end{center} \begin{center}
\tabcolsep=4pt
\begin{tabular}{|c|*{5}{c}|}
\hline
$ wt=2$ & 
$\ds  \kmTemp{1}{2} $&
$\ds  \kmTemp{2}{2} $&
$\ds  \kmTemp{3}{2} $&
$\ds  \kmTemp{4}{2} $&
$\ds  \kmTemp{5}{2} $ \\ \hline
$\dim$ of $\myCC{\bullet}{2}$ & 10 & 45 & 75 & 55 & 15 \\
$\dim(\ker)$  & 4 & 18 & 41 & 42 & 15 \\
Betti num   & 4 & 12 & 14 &  8 &  2\\ \hline
$\dim$ of $\myC{\bullet}{2}$     & 9 & 28 & 29 & 10 & 0 \\
$\dim(\ker)$  & 4 & 14 & 24 & 10 & 0 \\
Betti num &   4 & 9 & 10 &  5 &  0\\ \hline 
\end{tabular}
%\end{center} \begin{center} YY

\medskip

\tabcolsep=1pt
\begin{tabular}{|c|*{6}{c}|}
\hline
$wt=3$ &  
$\ds  \kmTemp{1}{3} $&
$\ds  \kmTemp{2}{3} $&
$\ds  \kmTemp{3}{3} $&
$\ds  \kmTemp{4}{3} $&
$\ds  \kmTemp{5}{3} $&
$\ds  \kmTemp{6}{3} $ \\ \hline
$\dim$ of $\myCC{\bullet}{3}$ & 15 & 105 & 245 & 255 & 120& 20  \\
$\dim(\ker)$ & 5 & 33 & 109 & 165 & 104&  20  \\
Betti num & 5 & 23 & 37 &   29 &  14 &  4  \\ \hline
$\dim$ of $\myC{\bullet}{3}$ & 14 & 73 & 114 & 65 & 10& 0 \\
$\dim(\ker)$ & 5 & 28 & 72 & 55 & 10 &   0  \\
Betti num & 5 & 19 & 27 &   13 &  0 &  0  \\ \hline 
\end{tabular}
%\end{center} \begin{center}
\hfil
%\medskip 
\tabcolsep=2pt
\begin{tabular}{|c|*{7}{c}|}
\hline
$wt=4$ & 
$\ds  \kmTemp{1}{4} $&
$\ds  \kmTemp{2}{4} $&
$\ds  \kmTemp{3}{4} $&
$\ds  \kmTemp{4}{4} $&
$\ds  \kmTemp{5}{4} $&
$\ds  \kmTemp{6}{4} $&
$\ds  \kmTemp{7}{4} $ \\ \hline
$\dim$ of $\myCC{\bullet}{4}$ & 21 & 198 & 618 & 891 & 630& 195 & 15  \\
        $\dim(\ker)$ & 6 & 49 & 221 & 477 & 464&  180 & 15  \\
Betti num & 6 & 34 & 72 & 80 &  50 &  14 & 0  \\ \hline
$\dim$ of $\myC{\bullet}{4}$ & 20 & 146 & 322 & 291 & 100& 5 & 0 \\
        $\dim(\ker)$ & 6 & 43 & 160 & 204 & 95&  5 & 0  \\
Betti num & 6 & 29 & 57 &   42 &  8 &  0 & 0  \\ \hline 
\end{tabular}
\end{center} 

\subsubsection{Solvable but not nilpotent cases}
The Lie algebra of upper triangle matrices of type (2,2) is an example
of not nilpotent but solvable Lie algebra of dimension 3. 
Let 
$\ds A_1 = \begin{bmatrix} 1 & 0 \\ 0 & 0\end{bmatrix}$, 
$\ds A_2 = \begin{bmatrix} 0 & 0 \\ 0 & 1\end{bmatrix}$, 
$\ds A_3 = \begin{bmatrix} 0 & 1 \\ 0 & 0\end{bmatrix}$. We get 
        \( [A_1,A_3]=A_3\;, [A_2,A_3]=-A_3 \) and the other brackets are
        trivial. The adjoint representation is given by 
\[ \operatorname{ad}_{A_1} = 
\begin{bmatrix} 0 & 0 & 0 \\ 0 & 0 & 0 \\ 0 & 0 & 1 \end{bmatrix}\;, \; 
 \operatorname{ad}_{A_2} = 
\begin{bmatrix} 0 & 0 & 0 \\ 0 & 0 & 0 \\ 0 & 0 & -1 \end{bmatrix}\;,\;  
 \operatorname{ad}_{A_3} = 
\begin{bmatrix} 0 & 0 & 0 \\ 0 & 0 & 0 \\ -1 & -1 & 0 \end{bmatrix}\;,  
\] 
and those show the algebra is not nilpotent, but the derived algebra is
1-dimensional and we see that the algebra is solvable.  We consider Lie
Poisson manifold whose Casimir polynomials are $\ds (x_1+x_2)^k$ where
we identity $\ds A_i$ by $\ds x_i$.   

\begin{center}
\tabcolsep=3pt 
\setlength{\extrarowheight}{-3pt}
\begin{tabular}{|c|*{6}{c}c|}
\hline
$wt=0$ & $\ds\myCC{0}{0}$ & $\rightarrow $ &
         $\ds\myCC{1}{0}$ & $\rightarrow $ &
         $\ds\myCC{2}{0}$ & $\rightarrow $ &
         $\ds\myCC{3}{3}$ \\ \hline
         $\dim$      &    1 && 3 && 3 && 1     \\
         $\dim(\ker)$&    1 && 2 && 2 && 1     \\ 
         Betti num   &    1 && 2 && 1 && 0    \\ \hline
\end{tabular}
\hfil 
\begin{tabular}{|c|*{6}{c}c|}
\hline
        $wt=0$ & $\ds\myC{0}{0}$ & $\rightarrow $ &
         $\ds\myC{1}{0}$ & $\rightarrow $ &
         $\ds\myC{2}{0}$ & $\rightarrow $ &
         $\ds\myC{3}{3}$ \\ \hline
         $\dim$       &   1 && 2 && 1 && 0 \\
         $\dim(\ker)$&    1 && 1 && 1 && 0 \\ 
         Betti num    &   1 && 1 && 0 && 0 \\ \hline
\end{tabular} 
\end{center}

\begin{center}
\tabcolsep=7pt
\setlength{\extrarowheight}{-3pt}
\begin{tabular}{|c|*{4}{c}|}
\hline
$wt=1$ &    
 $\ds\myCC{1}{1} $&  $\ds\myCC{2}{1} $&  $\ds\myCC{3}{1} $&  $\ds\myCC{4}{1} $   \\ \hline
$\dim$ &    6 & 18 & 18 & 6   \\
        $\dim(\ker)$     & 3 & 9 & 12 & 6   \\
Betti num       & 3 & 6  &  3 & 0   \\\hline
\end{tabular}
\hfil
\tabcolsep=7pt
\begin{tabular}{|c|*{4}{c}|}
\hline
$wt=1$ &    
 $\ds\myC{1}{1} $&  
 $\ds\myC{2}{1} $&  
 $\ds\myC{3}{1} $&  
 $\ds\myC{4}{1} $   \\ \hline
$\dim$     & 5 & 10 & 5 & 0   \\
        $\dim(\ker)$     & 2 & 5  & 5 & 0   \\
Betti num    &    2 & 2  &  0 & 0   \\\hline
\end{tabular}
\end{center}
\medskip

\begin{center} 
\tabcolsep=2pt
\begin{tabular}{|c|*{5}{c}|}
\hline
$ wt=2$ & 
$\ds  \kmTemp{1}{2} $&
$\ds  \kmTemp{2}{2} $&
$\ds  \kmTemp{3}{2} $&
$\ds  \kmTemp{4}{2} $&
$\ds  \kmTemp{5}{2} $ \\ \hline
$\dim$ of $\myCC{\bullet}{2}$ &     10 & 45 & 75 & 55 & 15 \\
$\dim(\ker)$ & 4 & 17 & 38 & 40 & 15 \\
Betti num &    4 & 11 & 10 &  3 &  0\\ \hline
$\dim$ of $\myC{\bullet}{2}$ &     9 & 28 & 29 & 10 & 0 \\
$\dim(\ker)$  & 3 & 10 & 19 & 10 & 0 \\
Betti num &    3 & 4 & 1 &  0 &  0\\ \hline 
\end{tabular}
%\end{center} \begin{center}
\hfil
%\medskip 
\tabcolsep=2pt
\begin{tabular}{|c|*{6}{c}|}
\hline
$wt=3$&    
$\ds  \kmTemp{1}{3} $&
$\ds  \kmTemp{2}{3} $&
$\ds  \kmTemp{3}{3} $&
$\ds  \kmTemp{4}{3} $&
$\ds  \kmTemp{5}{3} $&
$\ds  \kmTemp{6}{3} $ \\ \hline %XXX
$\dim$ of $\myCC{\bullet}{3}$ &    15 & 105 & 245 & 255 & 120& 20  \\
$\dim(\ker)$    & 5 & 32 & 103 & 156 & 100&  20  \\
Betti num      & 5 & 22 & 30 &   14 &  1 &  0  \\
\hline
$\dim$ of $\myC{\bullet}{3}$    & 14 & 73 & 114 & 65 & 10& 0  \\
$\dim(\ker)$    & 4 & 20 & 59 & 55 & 10 &   0  \\
Betti num      & 4 & 10 & 6 &   0 &  0 &  0  \\ \hline 
\end{tabular}
%\end{center} \begin{center}

\medskip

\tabcolsep=3pt
\begin{tabular}{|c|*{8}{c}|}
\hline
$wt=4$ &   
$\ds  \kmTemp{1}{4} $&
$\ds  \kmTemp{2}{4} $&
$\ds  \kmTemp{3}{4} $&
$\ds  \kmTemp{4}{4} $&
$\ds  \kmTemp{5}{4} $&
$\ds  \kmTemp{6}{4} $&
$\ds  \kmTemp{7}{4} $& \\ \hline
$\dim$ of $\myCC{\bullet}{4}$ &  21 & 198 & 618 & 891 & 630& 195 &
15 & \\
        $\dim(\ker)$ & 6 & 48 & 210 & 453 & 450&  180 & 15 & \\
Betti num &  6 & 33 & 60 & 45 &  12 &  0 & 0 & \\
\hline
$\dim$ of $\myC{\bullet}{4}$ & 20 & 146 & 322 & 291 & 100& 5 & 0 & \\
        $\dim(\ker)$ & 5 & 31 & 129 & 196 & 95&  5 & 0 & \\
Betti num & 5 & 16 & 14 &   3 &  0 &  0 & 0  &\\ \hline 
\end{tabular}
\end{center}

\subsection{$h=2$} 
Even though the normal form of analytic Poisson structures are studied
by J.F.~Conn (\cite{MR744864}), it is not clear what is the typical
2-homogeneous Poisson structure in our context.   

Here, we show some concrete examples: When $h=2$, we see $\ds
\Bkt{\mySS{i}}{\mySS{j}} \subset \mySS{i+j}$ and so we have 
\begin{align*}
\mydov(\frakSS{1}) =& (0)\;, \quad 
\mydov(\frakSS{2}) \subset \Lambda^2 \frakSS{1}\;, \quad 
\mydov(\frakSS{3}) \subset \frakSS{1} \otimes \frakSS{2}\;, \quad 
\mydov(\frakSS{4}) \subset \frakSS{1}\otimes \frakSS{3} \oplus
\Lambda^2 \frakSS{2}\;, 
\\
\mydov(\frakSS{5}) \subset& \frakSS{1}\otimes \frakSS{4} \oplus
\frakSS{2}\otimes \frakSS{3}\;. 
\end{align*}

We deal with the following 3 cases of 2-homogeneous Poisson structures
on $\ds\mR^3$.  
\begin{align*}
\text{case1:}\quad  & 
\pi =\frac{1}{2} ( x_1{}^2  \pdel_2 \wedge \pdel_3
+ x_2{}^2  \pdel_3 \wedge \pdel_1  + x_3{}^2  \pdel_1 \wedge \pdel_2), 
&& \text{Casimirs are } (x_1{}^3 + x_2{}^3 + x_3{}^3)^{k} 
\;. 
\\
\text{case2:}\quad & 
\pi = x_1 x_2  \pdel_1 \wedge \pdel_2 + x_2 x_3   \pdel_2 \wedge \pdel_3
+ x_3 x_1  \pdel_3 \wedge \pdel_1 , &&
 \text{Casimirs are } (x_1 x_2 x_3)^{k} 
 \;. 
\\
\text{case3:}\quad & 
\pi = x_1{}^2   \pdel_2 \wedge \pdel_3
+ x_3 x_1   \pdel_3 \wedge \pdel_1  + x_1 x_2  \pdel_1 \wedge \pdel_2 , 
& &\text{Casimirs are } (x_1{ }^2 + 2 x_2 x_3)^{k}\;. \\
\end{align*}

The cochain complex of 
weight 0 is trivial and only $\ds\frakSS{1}$ for weight 1. 
\subsubsection{weight=2} 
The three 2-homogeneous Poisson structures have the same table when
weight =2 as below.  
\begin{center}
\tabcolsep=3pt 
\setlength{\extrarowheight}{-3pt}
\begin{tabular}{|c|*{6}{c}c|}
\hline
$wt=2$
 & $\ds  
 \mathbf{0}$ & $\rightarrow$ &
 $\ds  \myCC{1}{2} $& $ \rightarrow $ &
 $\ds  \myCC{2}{2} $& $ \rightarrow $ &
$\mathbf{0} $ \\ \hline
$\dim$       &   && 6 && 3 &&    \\
%rank         &&  0 && 3 && 0 & \\
$\dim(\ker)$  &   && 3 && 3 &&   \\
%$\dim$(image)&   && 0 && 3  &&   \\
Betti num    &   && 3 && 0  &&   \\
\hline
\end{tabular}
\end{center}
However, $\ds\myH{\bullet}{2}$ differ as follows:

\tabcolsep=3pt 
\setlength{\extrarowheight}{-3pt}
\begin{tabular}{|c|*{2}{c}c|}
\hline
case1  & $\ds  \myC{1}{2}$ & $\rightarrow$ & $\ds\myC{2}{2}$   \\ \hline
$\dim$       &    6 && 3     \\
$\dim$(ker)  &    3 && 3    \\
Betti num    &    3 && 0     \\
\hline
\end{tabular}
\hfil
\begin{tabular}{|c|*{2}{c}c|}
\hline
case2  & $\ds  \myC{1}{2}$ & $\rightarrow$ & $\ds\myC{2}{2}$ \\\hline
$\dim$       &    6 && 3     \\
$\dim$(ker)  &    3 && 3    \\
Betti num    &    3 && 0     \\
\hline
\end{tabular}
\hfil 
\begin{tabular}{|c|*{2}{c}c|}
\hline
case3  & $\ds  \myC{1}{2}$ & $\rightarrow$ & $\ds\myC{2}{2}$ \\\hline
$\dim$       &    5 && 3     \\
$\dim$(ker)  &    2 && 3    \\
Betti num    &    2 && 0     \\
\hline
\end{tabular}

\medskip

\subsubsection{weight=3 }

\begin{center}
\tabcolsep=3pt 
\setlength{\extrarowheight}{-3pt}
\begin{tabular}{|cc|*{4}{c}c|}
\hline
& $wt=3$ & $\ds  \myCC{1}{3}$ & $\rightarrow $ & $\ds  \myCC{2}{3}$ & $\rightarrow $ &
$\ds  \myCC{3}{3}$   \\ 
& $\dim$       &    10 && 18 &&  1      \\\hline
case1 & Betti num    &     2 &&  9  && 0     \\ \hline
case2 & Betti num    &     4 &&  11  && 0     \\ \hline
case3 & Betti num    &     5 &&  12  && 0     \\ \hline
\end{tabular}
\end{center}

\setlength{\extrarowheight}{-3pt}
\tabcolsep=4pt
\begin{tabular}{|c|*{4}{c}c|}
\hline
 case1  & $\ds\myC{1}{3}$ & $\rightarrow $ & $\ds\myC{2}{3}$ & $\rightarrow $ &
$\ds\myC{3}{3}$\\\hline 
 $\dim$       &    9 && 18 &&  1\\
 Betti &     1 &&  9  && 0     \\ \hline
\end{tabular}
\hfil 
\begin{tabular}{|c|*{4}{c}c|}
\hline
case2  & $\ds\myC{1}{3}$ & $\rightarrow $ & $\ds\myC{2}{3}$ & $\rightarrow $ &
$\ds\myC{3}{3}$   \\\hline 
 $\dim$       &    9 && 18 &&  1 \\
 Betti &     3 &&  11  && 0     \\ \hline
\end{tabular}
\hfil
\begin{tabular}{|c|*{4}{c}c|}
\hline
case3  & $\ds\myC{1}{3}$ & $\rightarrow $ & $\ds\myC{2}{3}$ & $\rightarrow $ &
$\ds  \myC{3}{3}$   \\\hline 
 $\dim$       &    10 && 15 &&  1 \\
 Betti &     5 &&  9  && 0     \\ \hline
\end{tabular}

\medskip

\subsubsection{weight=4}

\begin{center}
\tabcolsep=3pt 
\setlength{\extrarowheight}{-3pt}
\begin{tabular}{|cc|*{4}{c}c|}
\hline
& $wt=4$
 &  
 $\ds\myCC{1}{4} $& $ \rightarrow $ &
 $\ds\myCC{2}{4} $& $ \rightarrow $ &
 $\ds\myCC{3}{4} $  \\ 
& $\dim$       &    15 && 45 &&  18   \\\hline
%case-1 &rank  &&     12 && 18 &&   \\
case1 &Betti num    &     3 &&  15 && 0     \\ \hline
%case-2 &rank  &&    12 && 18 &&   \\
case2 &Betti num    &     3 &&  15 && 0     \\ \hline
%case-3 &rank  &&   10 && 18 &&   \\
case3 &Betti num    &     5 && 17 && 0    \\ \hline
\end{tabular}
\end{center}
About $\ds\myH{\bullet}{4}$ we see that

\medskip

\tabcolsep=4pt 
\setlength{\extrarowheight}{-3pt}
\begin{tabular}{|c|*{4}{c}c|}
\hline
case1 &  
 $\ds\myC{1}{4} $& $ \rightarrow $ &
 $\ds\myC{2}{4} $& $ \rightarrow $ &
 $\ds\myC{3}{4} $  \\ \hline
$\dim$  &    15 && 42 &&  18   \\
Betti &     3 &&  12 && 0     \\ \hline
\end{tabular}
\hfil
\begin{tabular}{|c|*{4}{c}c|}
\hline
case2 &  
 $\ds\myC{1}{4} $& $ \rightarrow $ &
 $\ds\myC{2}{4} $& $ \rightarrow $ &
 $\ds\myC{3}{4} $  \\ \hline
$\dim$  &    15 && 42 &&  18   \\
Betti &     3 &&  12 && 0     \\ \hline
\end{tabular}
\hfil
\begin{tabular}{|c|*{4}{c}c|}
\hline
case3 &  
 $\ds\myC{1}{4} $& $ \rightarrow $ &
 $\ds\myC{2}{4} $& $ \rightarrow $ &
 $\ds\myC{3}{4} $  \\ \hline
$\dim$  &    14 && 40 &&  15   \\
Betti &     4 &&  15 && 0     \\ \hline
\end{tabular}

\medskip

\subsection{Comparison of elapsed time}
\label{elapsed:time}
Here, we compare the elapsed time for concrete Lie Poisson structures 
by $\fraksl(2,\mR)$ and 3-dim Heisenberg Lie algebra by Symbolic
Calculus software, Maple. Two classes mean 
the polynomial algebra and the algebra of Hamiltonian vector fields of
polynomials. 

\paragraph{Polynomials}
\begin{center}
\begin{tabular}{|c|*{2}{c}|*{2}{c}|}\hline 
& \multicolumn{2}{|c|}{sl2} & \multicolumn{2}{c|}{Heisen3} \\\hline \hline 
        weight & HomoL & Cohom & HomoL & Cohom \\\hline
        0 & 0.021 & 0.021 & 0.024 & 0.019 \\
        1 & 0.047 & 0.057 & 0.043 & 0.042 \\
        2 & 0.221 & 0.269 & 0.136 & 0.163 \\
        3 & 1.467 & 1.664 & 0.920 & 0.980 \\
4 & 13.274 & 14.651 & 7.260 & 7.592\\\hline \end{tabular}
\end{center}

\kmcomment{
cf: \verb+../km/Odd_GKF/new-exec-(0,1,2)-v<>.mpl+ is for Cohomology of
Polynomials.   

cf: \verb+../km/Odd_GKF/Homol/homoL-exec-(0,1,2)-v<>.mpl+ is for Homology 
of Polynomials.  
}

\paragraph{Hamiltonian}
\begin{center}
\begin{tabular}{|c|*{2}{c}|*{2}{c}|}\hline 
& \multicolumn{2}{|c|}{sl2} & \multicolumn{2}{c|}{Heisen3} \\\hline \hline 
        weight & HomoL & Cohom & HomoL & Cohom \\\hline
        0 & 0.039(0.036) & 0.022 & 0.023({\color{red}0.032}) & 0.020 \\
        1 & 0.197(0.186) & 0.062 & 0.086(0.080) & 0.044 \\
        2 & 1.300(1.198) & 0.285 & 0.429(0.353) & 0.090 \\
        3 & 8.236(7.745) & 1.512 & 2.335(1.865) & 0.312 \\
4 & 44.539(40.497) & 11.453 & 11.645(8.873) & 1.552\\\hline \end{tabular}
\end{center}

The numbers in the parentheses are obtained by calculating the remainder
by our method \eqref{eqn:anchoku}. Except only one case printed in red color,  
using \eqref{eqn:anchoku} is faster than using the normal form. Even
though, computing cohomology groups is faster than homology computing in 
Hamiltonian vector fields case.    

\kmcomment{
cf:\verb+../Odd_GKF/new-exec-2-v<2>.mpl+ is for Cohomology 
of Hamiltonian (just switch \texttt{myHamil}).

cf:\verb+../Odd_GKF/Homol/homoL-exec-3-v<>.mpl+ is for Homology 
of Hamiltonian.
}

%% file: Top-Betti.tex
\section{Top Betti number} 
As before, 
let $\ovfrakg $ be the Lie algebra defined by a non-trivial $h$-homogeneous
Poisson structure on $\ds\mR^n$.  
For a given weight $w$, we have a sequence of
cochain groups $\{\myCC{m}{w}\}$. 
As remarked in Subsection \ref{subsec:possible}, the sequence is finitely bounded, so   
let $\ds m_{0}(w) = \max \{ m \mid
        \myCC{m}{w} \ne 0 \}$.  
        We call $\ds\myCC{m_0(w)}{w} $,   
        $\ds\myHH{m_0(w)}{w} $ or    
        $\ds\dim\myHH{m_0(w)}{w} $ as the \last cochain group, the \last
        cohomology group or the \last Betti number.    

\subsection{Multi-index manipulation} 
As already seen in concrete examples, in order to deal with homogeneous
polynomials of $\ds x_1,\ldots,x_n$, we use monomials as a basis, and 
we use multi-index notation.  For each
positive integer $k$, let $\ds\mathfrak{M}[k] = \{ A \in \mN^{n} \mid
|A| = k\}$.  We denote the dual basis of $\ds \w{A}$ by $\ds z_{A}$ or
$\ds\Z{|A|}{A}$ emphasizing the degree of $A$.  
 
We put a total order in $\ds\mathfrak{M}[k]$ and assign natural number
$j$ to $A$ if $A$ is the $j$-th element in $\ds\mathfrak{M}[k]$ and
denote this assignment by $\junban(A) = j$.  
For $\ds\mathfrak{M}[1]$, it is natural to define as follows: for $\ds
E_j = (0,\ldots, \mathop{1}_{j},0\ldots)$,     $\junban( E_j ) = j$.    

For each positive $j$, we use the notation 
\begin{align*} \zF{j} & = \Z{j}{ \junban^{-1}(1)} \wedge \cdots \wedge
\Z{j}{ \junban( \# \mathfrak{M}[j])} && (\text{all}\ A\in
\mathfrak{M}[j])\;, \\
\zL{j}{A} & = \Z{j}{ \junban^{-1}(1)} \wedge \cdots \wedge
		\Z{j}{ \junban{}^{-1} ( \junban(A)-1)}&& (\text{all before}\
		A)\;,\\
\zR{j}{A} &= \Z{j}{ \junban^{-1}(\junban(A)+1)} \wedge \cdots \wedge
		\Z{j}{ \junban( \# \mathfrak{M}[j])}&& (\text{all after}\ A) \;,\\ 
\zM{j}{A} &= \zL{j}{A} \wedge \zR{j}{A}
		&& (\text{all except}\ A) \;,
\\ \noalign{in particular, 
if $A<B$} 
\zBetw{j}{A}{B} &= \Z{j}{\junban ^{-1}(\junban(A)+1)} \wedge \cdots\wedge 
\Z{j}{ \junban ^{-1}(\junban(B)-1)} 
&&(\text{between}\ A\ \text{and}\ B) \;. 
\end{align*}

%%%%%%%%%%%%%%%%%%%%%% 
\subsection{Case of polynomial algebra}
In the concrete examples of 2-homogeneous Poisson structures in previous
subsections, all the \last Betti number is zero, but we have the next
example of non-zero \last Betti number.    
\begin{exam} 
Let us consider a small example, $n=3$ and the 2-homogeneous Poisson bracket
given by \[ \Bkt{x_1}{x_2} = x_3{}^2, \quad \Bkt{x_1}{x_3} = 0, \quad
\Bkt{x_2}{x_3} = 0\] 
and study of weight 2 cohomology groups.  
1-cochain complex is $\ds\myCC{1}{2} = \frakSS{2}$ (6 dim) and  
2-cochain complex is $\ds\myCC{2}{2} = \Lambda^2 \frakSS{1}$ (3 dim), and     
$\ds\myCC{m}{2} = 0$ for $m >2$. $ \mydov ( \Z{1}{j} ) - 0$ ($j=1,2,3$) and 
\begin{align*}
        \mydov( \Z{2}{A} ) =& 
        - \sum_{i<j} \langle \Z{2}{A},
        \Bkt{x_i}{x_j}\rangle  \Z{1}{i} \wedge \Z{1}{j} 
        = 
        - \langle \Z{2}{A},
        x_3 {}^2 \rangle  \Z{1}{1} \wedge \Z{1}{2} = \begin{cases} 
                - \Z{1}{1} \wedge \Z{1}{2} & \text{if} \quad
				A=(0,0,2)\;, \\
                0 & \text{otherwise}\;.\end{cases} 
\end{align*} 
Thus, we have the table below left which tells that the \last Betti number is
not zero: 
\begin{center}        
\tabcolsep=3pt 
\setlength{\extrarowheight}{-3pt}
\begin{tabular}{|c|*{3}{c}|}
\hline
&  
 $\ds\myCC{1}{2} $& $ \rightarrow $ &
 $\ds\myCC{2}{2} $  
 \\\hline 
$\dim$  &    6 && 3   \\\hline
$\ker \dim $  &    5 && 3   \\\hline
Betti &     5 &&  2    \\ \hline
\end{tabular}
\hfil 
\begin{tabular}{|c|*{3}{c}|}
\hline
&  
 $\ds\myC{1}{2} $& $ \rightarrow $ &
 $\ds\myC{2}{2} $  
 \\\hline 
$\dim$  &    5 && 1   \\\hline
$\ker \dim $  &    4 && 1   \\\hline
Betti &     4 &&  0    \\ \hline 
\end{tabular} 
\end{center} 
The table above right of $\ds\frakg = \ovfrakg/\myCasim $ tells the 
\last Betti number is zero in this case. 
\end{exam} 
\kmcomment{
Our cochain complexes correspond to the Young diagrams with the
dimensional condition, and we always have the special Young diagram
satisfying the dimensional condition in the extreme sense, namely
$\text{YD}_0 = [ \md{j} \mid j=1\ldots \ell]$, where $\ds\md{j}=\dim
\frakSS{j} = \tbinom{n-1+j}{j}$.  
}%endOFkmcomment
%%%%%%%%%%%%%%%%%%%%%%%%%%%%%%%%%%%%%%%%%%%%% 
Each  of our cochain group corresponds to a Young diagram 
with our dimensional condition and we have special Young diagram 
which satisfies the extremal dimensional condition in the sense that 
each index of the Young diagram is maximal.  
Namely, for each fixed $\ell$, 
it is the one given by $\text{YD}_0 = [ \md{j} 
\mid j=1\ldots \ell]$, where $\ds \md{j}=\dim \frakSS{j} = \tbinom{n-1+j} 
{j}$.   
%%%%%%%%%%%%%%%%%%%%%%%%%%%%%%%%%%%%%%%%%%%%%
The corresponding factor of the
cochain complex is 
\[ \Lambda^{\md{1}} \frakSS{1} \otimes 
\Lambda^{\md{2}} \frakSS{2} \otimes \cdots 
\otimes \Lambda^{\md{\ell}} \frakSS{\ell}\] 
whose degree is \begin{equation}\ds 
         \label{eqn:three:one}
        m_{0} = \sum_{j=1}^{\ell}
        \md{j}\end{equation} and the weight is 
\begin{equation} \label{eqn:three:two}
        w_0 = \sum_{j=1}^{\ell} (j-2+h) \md{j}
\end{equation} 
from the definition, where $h$ the homogeneity of our Poisson bracket.
Concerning to these weight $\ds w_0$ and degree $\ds m_0$, we have the
following propositions.   

\begin{kmProp}\label{prop:high} For the given $\ds w_0$ and $\ds m_0$, 
        \begin{align*}
                \ds\myCC{m}{w_0} &= ( 0 ) \quad \text{if}\quad m>m_0 
                \\
                \noalign{and} 
         \myCC{m_0}{w_0} &= 
 \Lambda^{\md{1}} \frakSS{1} \otimes 
\Lambda^{\md{2}} \frakSS{2} \otimes \cdots 
\otimes \Lambda^{\md{\ell}} \frakSS{\ell} \quad \text{(which is 1 dimensional)}
\;. 
\end{align*} 
\end{kmProp}
\textbf{Proof:} 
To show the first claim, suppose $\ds\myCC{m}{w_0} \ne 0$ for some $m
>m_0$.  Then  there is a Young diagram $\ds [k_i \mid j=1\ldots s ]$
with 
\begin{equation} \label{eqn:three:three}
        \text{the weight is}\quad 
        w_0 = \sum_{i=1}^s (i-2+h) k_i \quad\text{and the height is}\quad  
        m =  \sum_{i=1}^s  k_i\end{equation} and also
satisfying the dimensional condition $\ds  0\leq k_i \leq \md{i}$, where $\ds
\md{i} = \mmd{i}$.  
If $s\leq \ell$ then $m \leq m_0$ contradicting to the assumption, so we may
assume $s > \ell$. Then, comparing the first equation of  
(\ref{eqn:three:three}) and 
(\ref{eqn:three:two}), we have 
\begin{equation}\label{eqn:four}
        \sum_{i>\ell} (i-2+h) k_i = \sum_{j=1}^{\ell} (j-2+h)
		(\md{j}-k_j)\;. 
\end{equation} On the other hand, 
\begin{align}
        0 < \Delta = m-m_{0} 
        =& \sum_{j=1}^{\ell} (k_j-\md{j}) + \sum_{i>\ell} k_i\;, \notag \\
\sum_{i>\ell} k_i =& 
\Delta + \sum_{j=1}^{\ell} (\md{j} - k_j )\;.  \label{eqn:five} 
\end{align}
Extracting 
$(\ell -2 +h)$ times (\ref{eqn:five}) from (\ref{eqn:four}), we have 
\[ \sum_{i>\ell} (i- \ell)k_i = - (\ell-2+h)\Delta + \sum_{j=1}^{\ell} (
j-\ell) ( \md{j}-k_j)\;. \] 
The left hand side is positive and the right hand side is non-positive, this
is impossible. Therefore, we conclude that $\ds\myCC{m}{w_0} = 0$ for $m>
m_0$. 

Starting from $m = m_0$, 
we claim that there is no other Young
diagram with the same weight and the
same height, namely, .  
we conclude that $s=\ell$ and $k_j = \md{j}$ and so 
\[ \myCC{m_0}{w_0} = 
 \Lambda^{\md{1}} \frakSS{1} \otimes 
\Lambda^{\md{2}} \frakSS{2} \otimes \cdots 
\otimes \Lambda^{\md{\ell}} \frakSS{\ell}\ . \]  
\kmqed 
\begin{kmProp}\label{prop:low} We use the same values $\ds w_0$ and $\ds m_0$
        for $h$.  
When $h=1$ the Young diagram defined by 
$\ds k_1= \md{1}-1$ , $\ds k_j=\md{j}$ ($j=2\ldots \ell$) is a factor, and 
the Young diagram defined by 
$\ds k_j=\md{j}$ ($j=1\ldots \ell-1$), 
$\ds k_{\ell}= \md{\ell}-2$, $\ds k_{2\ell-1}=1$ is 
another factor  of $\ds\myCC{m_0 -1}{w_0}$.  
When $h>1$, 
if $\ell>1$ then 
$\ds k_1= \md{1}-1$, $\ds k_j=\md{j}$
($j=2\ldots \ell-1$),  $\ds k_{\ell} = \md{\ell}-1$ and $\ds k_{\ell+h-1} =1$  
else if $\ell=1 $ then $k_1 = \md{1}-2$, $k_h =1$ 
is a factor of the direct sum of $\ds\myCC{m_0 -1}{w_0}$.   
\end{kmProp}
\textbf{Proof:} It is just calculation that the Young diagrams in the Proposition satisfy the dimensional
condition and its height is $m_0 -1$ and its weight is equal to $w_0$.  
But a visual understanding is the following:  

\begin{center}
        \setlength{\unitlength}{4mm} 

\begin{picture}(11,12)(-2,0)
        \put(-3,11){$h=1$ case 1}
\path(0,0)(0,10)% (1,2)(0,2)(0,0)
%\put(0,0){$\mathbf O$}
%\multiput(0,0)(0,1){3}{\line(1,0){1}}
\path(0,0)(1,0)(1,2)(2,2)
%\put(1.5,1){\makebox(0,0)[bl]{$k_{1}$}}
%\multiput(0,2)(0,1){4}{\line(1,0){2}}
\path(2,2)(2,5)
\put(2.5, 3){$k_{2}$}
%\multiput(0,8)(0,2){3}{\line(1,0){8}}
\path(7,7)(8,7)(8,10)(0,10)
%\put(8.5, 9){\makebox(0,0)[l]{$k_{\ell}$}}
\put(4,10.5){$\ell$}
\put(3,10.5){\vector(-1,0){3}}
\put(5,10.5){\vector(1,0){3}}
%\put(5,0){\Large A} 
\put(-2,5){\makebox(0,0){$m$}}
\put(-2,6){\vector(0,1){4}}
\put(-2,4){\vector(0,-1){4}}
%\dashline{4}[0.8](6,6)(6,7)
\dottedline[$\cdot$]{0.2}(2,5)(3,5)
\dottedline[$\cdot$]{0.2}(3,5)(3,6)
\dottedline[$\cdot$]{0.2}(7,7)(7,6)
\dottedline[$\cdot$]{0.2}(3,6)(7,6)
%\put(5,0){\arc{1.2}{5.2}{6.2}}
%\put(4,1.7){\makebox(0,0)[l]{Total area is}}
%\put(4,0.5){\makebox(0,0)[l]{$w +(2-h)m$}}
\put(4,7){\makebox(0,0){$w$}}
\shade{\path(0,0)(1,0)(1,10)(0,10)(0,0)}
\put(-1,1){\vector(1,0){4}}
\put(4,1){\makebox(0,0){cut}}
%\put(8,10){\makebox(0,0){X}}
%\shade{\path(1,2)(2,2)(2,3)(1,3)(1,2)}%
% \path(1,2)(2,2)(2,3)(1,3)(1,2)
%\put(1.5, 2.5){\makebox(0,0){cut}} 
%\shade{\path(7,7)(8,7)(8,8)(7,8)(7,7)}%
%% \path(7,7)(8,7)(8,8)(7,8)(7,7)
% \put(7.5, 7.5){\makebox(0,0){\scriptsize cut}}
%\path(8,10)(11,10)(11,9)(8,9)
%\put(9.5, 9.5){\makebox(0,0){put}}
\end{picture}
\hfil 
\begin{picture}(14,12)(-2,0)
        \put(-3,11){$h=1$ case 2}
\path(0,0)(0,10)% (1,2)(0,2)(0,0)
%\put(0,0){$\mathbf O$}
%\multiput(0,0)(0,1){3}{\line(1,0){1}}
\path(0,0)(1,0)(1,2)(2,2)
%\put(1.5,1){\makebox(0,0)[bl]{$k_{1}$}}
%\multiput(0,2)(0,1){4}{\line(1,0){2}}
\path(2,2)(2,5)
\put(2.5, 3){$k_{2}$}
%\multiput(0,8)(0,2){3}{\line(1,0){8}}
\path(7,8)(8,8)(8,10)(0,10)
%\put(8.5, 9){\makebox(0,0)[l]{$k_{\ell}$}}
\put(4,10.5){$\ell$}
\put(3,10.5){\vector(-1,0){3}}
\put(5,10.5){\vector(1,0){3}}
%\put(5,0){\Large A} 
\put(-2,5){\makebox(0,0){$m$}}
\put(-2,6){\vector(0,1){4}}
\put(-2,4){\vector(0,-1){4}}
%\dashline{4}[0.8](6,6)(6,7)
\dottedline[$\cdot$]{0.2}(2,5)(3,5)
\dottedline[$\cdot$]{0.2}(3,5)(3,6)
\dottedline[$\cdot$]{0.2}(7,7)(7,6)
\dottedline[$\cdot$]{0.2}(3,6)(7,6)
%\put(5,0){\arc{1.2}{5.2}{6.2}}
%\put(4,1.7){\makebox(0,0)[l]{Total area is}}
%\put(4,0.5){\makebox(0,0)[l]{$w +(2-h)m$}}
\put(4,7){\makebox(0,0){$w$}}
\shade{\path(0,0)(1,0)(1,10)(0,10)(0,0)}
\put(-1,1){\vector(1,0){4}}
\put(4,1){\makebox(0,0){cut}}
%\put(8,10){\makebox(0,0){X}}
%\shade{\path(1,2)(2,2)(2,3)(1,3)(1,2)}%
\path(1,2)(2,2)(2,3)(1,3)(1,2)
\put(1.5, 2.5){\makebox(0,0){cut}} 
%\shade{\path(7,7)(8,7)(8,8)(7,8)(7,7)}%
\path(7,7)(8,7)(8,8)(7,8)(7,7)
\put(7.5, 7.5){\makebox(0,0){\scriptsize cut}}
\path(8,10)(11,10)(11,9)(8,9)
\put(9.5, 9.5){\makebox(0,0){put}}
\end{picture} 
\\ 
\begin{picture}(18,13)(-5,0)
        \put(-5,11){$h>1$ case}
\path(0,0)(0,10)% (1,2)(0,2)(0,0)
%\put(0,0){$\mathbf O$}
%\multiput(0,0)(0,1){3}{\line(1,0){1}}
\path(0,0)(1,0)(1,2)(2,2)
%\put(1.5,1){\makebox(0,0)[bl]{$k_{1}$}}
%\multiput(0,2)(0,1){4}{\line(1,0){2}}
\path(2,2)(2,5)
\put(2.5, 3){$k_{2}$}
%\multiput(0,8)(0,2){3}{\line(1,0){8}}
\path(7,8)(8,8)(8,10)(0,10)
\put(8.5, 8.4){\makebox(0,0)[lt]{$k_{\ell}$}}
\put(4,10.5){$\ell$}
\put(3,10.5){\vector(-1,0){3}}
\put(5,10.5){\vector(1,0){3}}
%\put(5,0){\Large A} 
\put(-1.5,5){\makebox(0,0){$m$}}
\put(-1.5,6){\vector(0,1){4}}
\put(-1.5,4){\vector(0,-1){4}}
%\dashline{4}[0.8](6,6)(6,7)
\dottedline[$\cdot$]{0.2}(2,5)(3,5)
\dottedline[$\cdot$]{0.2}(3,5)(3,6)
\dottedline[$\cdot$]{0.2}(7,8)(7,6)
\dottedline[$\cdot$]{0.2}(3,6)(7,6)
%\put(5,0){\arc{1.2}{5.2}{6.2}}
%\put(4,1.7){\makebox(0,0)[l]{Total area is}}
%\put(4,0.5){\makebox(0,0)[l]{$w +(2-h)m$}}
\put(4,5){\makebox(0,0)[l]{$w$ is the sum of two areas}}
%\shade{\path(0,0)(1,0)(1,10)(0,10)(0,0)}
\put(-5.5,1){\vector(1,0){8}}
\put(4,1){\makebox(0,0){cut}}
%%%%%%%%%%%%%%%%%%%%%%%%%%%%%%%%%%%%%%%%%%%%%%%%%%%%%%%%
\path(-5,0)(-3,0)(-3,10)(-5,10)(-5,0)
\put(-4,-1){\makebox(0,0){$h-2$}}
\dottedline[$\cdot$]{0.2}(9,10)(15,10)(15,9)(9,9)(9,10)
\put(12,9.5){\makebox(0,0)[c]{gather cells}}
\end{picture}
\end{center} 

\kmqed

Our observation on the \last Betti number is the next theorem.

\begin{thm}\label{thm:last:Betti} 
Let $\ovfrakg $ be the Lie algebra 
$\ds\mR[x_1,\ldots,x_n]/\mR$ 
defined by a non-trivial $h$-homogeneous
Poisson structure on $\ds\mR^n$.  For weight $w$ given by $\ds w =
\sum_{j=1}^{\ell} (j-2+h) \tbinom{ n-1+j }{j}$ 
consider the \last cochain group of dimension 
$\ds m_{0}(w) = \sum_{j=1}^{\ell}  \tbinom{ n-1+j}{j}$.
Then the Betti
number $\ds\dim \myHH{m_0 (w)}{w} = 0$ for each $h$.    
\end{thm} 

\textbf{Proof:} 
We use the basis $\ds\{ \Z{j}{ A }\}$ of $\frakSS{j}$ where $\ds A\subset
\mN^n$ with $|A|=j$, i.e., $\ds A\in \mathfrak{M}[j]$. 

First suppose $h=1$. As we saw in Proposition \ref{prop:low}, we have at least
two factors of $\ds\myCC{m_0 -1}{w_0}$ and one of them is $\ds k_j= \md{j}$
for $j=1\ldots \ell-1$, $\ds k_{\ell} = \md{\ell}-2$ and 
$\ds k_{2\ell -1} = 1$. Take a cochain \[ \ds\sigma   
= \zF{1} \wedge \cdots \wedge \zF{\ell-1} \wedge (\zL{\ell}{A}\wedge
 \zBetw{\ell}{A}{B} \wedge \zR{\ell}{B}) \wedge \Z{2\ell -1}{C}.\] 
Homogeneity 1 implies 
 $\ds\mydov( \frakSS{1}) \subset \Lambda^2 \frakSS{1}$, 
 and in general  
 $ \mydov( \frakSS{k}) \subset \sum_{i\leq j,\; i+j=k+1} \frakSS{i} \wedge
 \frakSS{j}$, and   
the following holds 
 \begin{align*}
& \mydov( \zF{1} ) = 0\;, \\
& \zF{1} \wedge \mydov( \zF{2}) = 0\;, \\
& \vdots \\
& 
 \zF{1}  \wedge \cdots \wedge \zF{\ell-1} \wedge \mydov(\zL{\ell}{A}\wedge
 \zBetw{\ell}{A}{B} \wedge \zR{\ell}{B}) = 0 \;. 
 \end{align*}
Thus, 
 \begin{align*}
\mydov ( \sigma ) = &  \pm  
 \zF{1}  \wedge \cdots \wedge \zF{\ell-1} \wedge (\zL{\ell}{A}\wedge
 \zBetw{\ell}{A}{B} \wedge \zR{\ell}{B})\wedge \mydov( \Z{2\ell-1}{C}) 
 \\
 = & \pm 
 \zF{1}  \wedge \cdots \wedge \zF{\ell-1} \wedge (\zL{\ell}{A}\wedge
 \zBetw{\ell}{A}{B} \wedge \zR{\ell}{B}) 
 \\ & 
 \wedge 
 \sum_{ |P|+|Q|=1+|C|,\; \junban{P}\leq \junban{Q}} 
 \langle  \Z{2\ell-1}{C}, 
         \Bkt{\w{P}}{\w{Q}} \rangle \Z{|P|}{P} \wedge  \Z{|Q|}{Q}   
\\
 = & \pm 
 \zF{1}  \wedge \cdots \wedge \zF{\ell} 
 \langle  \Z{2\ell-1}{C}, 
         \Bkt{\w{A}}{\w{B}} \rangle \;. 
 \end{align*} 
Since the Poisson bracket is not trivial, there are some $i,j$ such that 
$\ds\Bkt{x_i}{x_j} \ne 0$.  Let us take $\ds A, B\in \mathfrak{M}[\ell]$ such
that 
$\ds \w{A} = x_i{ }^{\ell}$ and  
$\ds \w{B} = x_j{ }^{\ell}$.  Then  
$\ds\Bkt{\w{A}}{\w{B}} = \ell^2  (x_i x_j)^{\ell-1}  \Bkt{x_i}{x_j} \ne 0$ and we can find 
some $\ds C\in \mathfrak{M}[2\ell-1]$ satisfying 
$\ds\langle \Z{3}{C}, \Bkt{\w{A}}{\w{B}}\rangle  \ne 0$. This means $\ds
\mydov(\sigma ) \ne 0$ 
and also that   
$ \mydov :\myCC{m_{0}(w)-1}{w}\to\myCC{m_{0}(w)}{w} $ is surjective 
in this case,   
thus $\ds\dim\myHH{m_0 (w_0)}{w_0} = 0$.

Now assume that $h>1$.  Then Proposition \ref{prop:high} says 
$\ds\myCC{m_{0}(w)}{w} 
= \prod_{j=1}^{\ell}\Lambda
^{\md{j}} \frakSS{j}$ and 
$\ds\dim \myCC{m_{0}(w)}{w} =1$ and  Proposition \ref{prop:low} says  
one factor of 
$\ds\myCC{m_{0}(w)-1}{w}$ is given by  
\begin{equation} \ds k_1=\md{1}-1,\quad k_2= \md{2},\quad\ldots,\quad  
k_{\ell-1} = \md{\ell-1},\quad  
k_{\ell} = \md{\ell}-1,\quad  
k_{\ell+h-1} = 1\;. \end{equation}    

Since our Poisson bracket is non-trivial, for some $i_0$ and $j_0$ we may
assume $\ds\Bkt{x_{i_0}}{x_{j_0}} \ne 0$.  

If  $\ell >1 $ then 
we take a $\ds (m_{0}(w)-1)$-cochain $\sigma$ defined by 
\[\sigma = \zM{1}{i_0} \wedge
\zF{2} \wedge \cdots \wedge \zF{\ell-1} \wedge 
\zM{\ell}{A} \wedge \Z{\ell+h-1}{B} 
\]  
else if $\ell=1$ then we take 
$\ds\sigma = \zL{1}{i_0} \wedge \zBetw{1}{i_0}{j_0} \wedge 
\zR{1}{j_0} \wedge \Z{h}{B}$.

By the homogeneity $h>1$, 
the coboundary operator $\mydov$ has the following properties 
\begin{equation} \mydov (\frakSS{i}) = 0 \quad (i<h)\quad\text{and}\quad \mydov
        (\frakSS{k}) \subset \sum_{i+j=k+2-h,\ i\leqq j}\frakSS{i} \wedge
        \frakSS{j}\;. 
\end{equation} 

When $\ell=1$, \begin{align*}\mydov( \sigma ) =& 
\pm 
\zL{1}{i_0} \wedge \zBetw{1}{i_0}{j_0} \wedge \zR{1}{j_0} 
\wedge \mydov( \Z{h}{B}) \\
= & 
\pm 
\zL{1}{i_0} \wedge \zBetw{1}{i_0}{j_0} \wedge \zR{1}{j_0} 
\wedge \langle \Z{h}{B},\Bkt{x_i}{x_j}\rangle
\Z{1}{i}\wedge \Z{1}{j}  \\
= & 
\pm 
\zF{1}
\langle \Z{h}{B},\Bkt{x_{i_0}}{x_{j_0}}\rangle
\ne 0 \quad\text{for some}\  B.  
  \end{align*} 

When $1 < \ell < h$ then it holds 
\begin{align*} \mydov (\sigma) =& 
\pm \zM{1}{i_0}  \wedge  \zF{2}  \wedge \cdots
\wedge \zF{\ell-1} \wedge \zM{\ell}{A} \wedge
\mydov(\Z{\ell+h-1}{B}) 
\\ = & 
\pm \zM{1}{i_0} \wedge 
\zF{2} \wedge  
\zF{3}
\wedge \cdots \wedge \zF{\ell-1} \wedge
\zM{\ell}{A} \wedge \sum_{|P|+|Q|=\ell+1}
        \langle \Z{\ell+h-1}{B},\Bkt{\w{P}}{\w{Q}}\rangle z_{P}\wedge z_{Q} \\
= & 
\pm \zM{1}{i_0} \wedge 
\zF{2} \wedge  
\zF{3}
\wedge \cdots \wedge \zF{\ell-1} \wedge
\zM{\ell}{A} \wedge \sum_{|P|=1,\ |Q|=\ell}
        \langle \Z{\ell+h-1}{B},\Bkt{\w{P}}{\w{Q}}\rangle \Z{1}{P}\wedge \Z{\ell}{Q} 
\\
= & 
\pm 
        \langle \Z{\ell+h-1}{B},\Bkt{\w{E_{i_0}}} {\w{A}}\rangle  
\zF{1} \wedge \zF{2} \wedge  \zF{3}
\wedge \cdots \wedge \zF{\ell-1} \wedge
\zF{\ell}\;.  
\end{align*}
Since $\ds\Bkt{x_{i_0}}{ x_{j_0}{ }^{\ell}} = 
\ell { x_{j_0}}^{\ell-1} 
\Bkt{x_{i_0}}{ x_{j_0}} \ne 0$, there is some $B$ with $|B|=\ell+h-1$ such
that  
$\ds 
\langle \Z{\ell+h-1}{B},\Bkt{x_{i_0}}{ x_{j_0}{ }^{\ell}}  
\rangle  \ne 0 $. Take $A$ with $|A|= A[j_0] = \ell$. Then $\ds\mydov (\sigma)
\ne 0$.  

When $\ell=h$, 
\begin{align*} \mydov (\sigma) =& 
\pm \zM{1}{i_0}  \wedge  \zF{2}  \wedge \cdots
\wedge \zF{\ell-1} \wedge \mydov( \zM{\ell}{A} \wedge
\Z{\ell+h-1}{B}) 
\\ = & 
\pm \zM{1}{i_0}  \wedge  \zF{2}  \wedge \cdots
\wedge \zF{\ell-1} \wedge 
\left( \mydov( \zM{\ell}{A}) \wedge \Z{\ell+h-1}{B} \pm 
\zM{\ell}{A} \wedge \mydov( \Z{\ell+h-1}{B}) \right) 
\\
\noalign{then, because of $\ds 
\mydov( \zM{\ell}{A}) \subset \frakSS{1} \wedge \frakSS{1}\wedge
\Lambda^{\md{\ell}-1} \frakSS{\ell} $, we see}
=& 
\pm \zM{1}{i_0}  \wedge  \zF{2} \wedge \cdots
\wedge \zF{\ell-1} \wedge 
\zM{\ell}{A} \wedge \mydov( \Z{\ell+h-1}{B}) 
\\ 
\ne & 0 \quad \text{by the same argument when $\ell < h$.} 
\end{align*}
When $h < \ell$, we see that 
\[ \mydov (\sigma) = 
\pm \zM{1}{i_0}  \wedge  \zF{2}  \wedge \cdots
\wedge \zF{h-1} \wedge \mydov( 
\zF{h} \wedge \cdots 
\wedge \zF{\ell-1} \wedge 
\zM{\ell}{A} \wedge \Z{\ell+h-1}{B}) 
\] 
and since 
\begin{align*} 
& 
\zM{1}{i_0}  \wedge  \zF{2}  \wedge \cdots
\wedge \zF{h-1} \wedge \mydov( 
\zF{h}) = 
\zM{1}{i_0}  \wedge  \zF{2}  \wedge \cdots
\wedge \zF{h-1} \wedge \Lambda^2 \frakSS{1} 
\wedge \Lambda^{\md{h}-1} \frakSS{h} 
= 0 \;,\\
& 
\zM{1}{i_0}  \wedge  \zF{2}  \wedge \cdots
\wedge \zF{h-1} \wedge  
\zF{h} \wedge 
\mydov( \zF{h+1}) \\&\quad =   
\zM{1}{i_0}  \wedge  \zF{2}  \wedge \cdots
\wedge \zF{h-1} \wedge  
\zF{h} \wedge  \frakSS{1}\wedge \frakSS{2} 
\wedge \Lambda^{\md{h+1}-1} \frakSS{h+1} 
= 0\;,\\
&\quad \vdots\\
& 
\zM{1}{i_0}  \wedge  \zF{2}  \wedge \cdots
\wedge \zF{h-1} \wedge  
\zF{h} \wedge \cdots \wedge 
\mydov( \zF{\ell-1}) \\&\quad =   
\zM{1}{i_0}  \wedge  \zF{2}  \wedge \cdots 
\wedge \zF{\ell-2} \wedge 
\left(
\frakSS{1}\wedge \frakSS{\ell-h} 
+ \frakSS{2}\wedge \frakSS{\ell-h-1} + \cdots \right) 
\wedge \Lambda^{\md{\ell-1}-1} \frakSS{\ell-1} 
= 0\;, \\
& 
\zM{1}{i_0}  \wedge  \zF{2}  \wedge \cdots
\wedge \zF{h-1} \wedge  
\zF{h} \wedge \cdots 
\wedge 
\zF{\ell-1} \wedge   
\mydov( \zM{\ell}{A}) \\&\quad =   
\zM{1}{i_0}  \wedge  \zF{2}  \wedge \cdots 
\wedge \zF{\ell-1} \wedge 
\left(
\frakSS{1}\wedge \frakSS{\ell-h+1} 
+ \frakSS{2}\wedge \frakSS{\ell-h} + \cdots \right) 
\wedge \Lambda^{\md{\ell}-2} \frakSS{\ell} 
= 0\;, 
\end{align*}
we have 
\[ \mydov ( \sigma ) = 
\pm \zM{1}{i_0}  \wedge  \zF{2}  \wedge \cdots
\wedge \zF{h-1} \wedge 
\zF{h} \wedge \cdots 
\wedge \zF{\ell-1} \wedge 
\zM{\ell}{A} \wedge \mydov ( \Z{\ell+h-1}{B}) \]  
and $\ds\mydov ( \sigma )\ne 0$ again by the same argument when $\ell < h$.  
\[ 
\mydov ( \sigma ) = \beta 
\zF{1} \wedge \cdots \wedge \zF{\ell} \ne 0\] 
this means $\ds\mydov : \myCC{m_{0}(w)-1}{w}\rightarrow \myCC{m_{0}(w)}{w}$
is surjective and therefore $\ds\myHH{m_{0}(w)}{w} = 0$.  
\kmqed

\bigskip

\subsection{Case of algebra of Hamiltonian vector fields} 
For a given non-trivial
$h$-homogeneous Poisson structure on $\ds\mR^n$, let   $\ds\ovfrakg$ be
the polynomial algebra of $\mR^n$ by the Poisson structure.  And we can
consider the subalgebra $\ds\frakg = \ovfrakg/\myCasim $, which is
modulo Casimir polynomials $\ds \myCasim $.  Now we have to handle
$\ds\dim\frakS{j}$ and so cochain complexes $\ds\myC{m}{w}$ carefully.
Let us denote $\ds\dim \frakS{j}$ by $\phi_j$  which is dependent on the
Poisson structure.  In general, $\ds\phi_j = \dim\frakS{j} \leq
\dim\frakSS{j} = \md{j} = \tbinom{n-1+j}{j}$. Since the key discussion
in the proofs of Proposition \ref{prop:high}, \ref{prop:low} or Theorem
\ref{thm:last:Betti} is to check the dimensional condition, only by
replacing $\ds\md{j}$ by $\phd{j}$, we have analogous results  of
Proposition \ref{prop:high}, \ref{prop:low} and  Theorem
\ref{thm:last:Betti} as follows: 

Let 
\begin{equation}\ds 
         \label{phi:three:one}
        m_{1} = \sum_{j=1}^{\ell}
        \phd{j}\end{equation} and 
\begin{equation} \label{phi:three:two}
        w_1 = \sum_{j=1}^{\ell} (j-2+h) \phd{j}
\end{equation}
where $h$ the homogeneity of our Poisson structure(tensor). 
Then we have 

\begin{kmProp}\label{phi:high} Let $\ds w_1$ and $\ds m_1$ be as above, then
        we have the following
        \begin{align*}
                \ds\myC{m}{w_1} &= ( 0 ) \quad \text{if}\quad m>m_1 
                \\
                \noalign{and} 
         \myC{m_1}{w_1} &= 
 \Lambda^{\phd{1}} \frakS{1} \otimes 
\Lambda^{\phd{2}} \frakS{2} \otimes \cdots 
\otimes \Lambda^{\phd{\ell}} \frakS{\ell} \quad \text{(which is 1
dimensional)}\;. 
\end{align*} 
\end{kmProp}
Also we have 
\begin{kmProp}\label{phi:low} For the same $h$, $\ds w_1$ and $\ds m_1$, 
        when $h=1$ 
the Young diagram defined by $\ds k_1= \phd{1}-1$ and $\ds k_j=\phd{j}$
($j=2\ldots \ell$) is a direct summand of 
 $\ds\myC{m_1 -1}{w_1}$. The Young diagram defined by  
$\ds k_j=\phd{j}$ ($j=1\ldots \ell-1$),  
$\ds k_{\ell}= \phd{\ell}-2$ and 
$\ds k_{2\ell-1}= 1$ is another direct summand of $\ds\myC{m_1 -1}{w_1}$. 

When $h>1$, if $\ell>1$ then $\ds k_1= \phd{1}-1$, $\ds k_j=\phd{j}$
($j=2\ldots \ell-1$),  $\ds k_{\ell} = \phd{\ell}-1$ and $\ds
k_{\ell+h-1} =1$.  If $h>1$ and $\ell=1 $ then $k_1 = \phd{1}-2$, $k_h =1$ is
a summand of the direct sum of $\ds\myC{m_1 -1}{w_1}$.   
\end{kmProp}

Combining the two Propositions above, we get the following theorem. 
\begin{thm}\label{thm:phi:last:Betti} 
Let $\frakg$ be the Lie algebra of polynomials defined by a non-trivial
$h$-homogeneous Poisson structure on $\ds\mR^n$.  For the weight $w$
given by $\ds w = \sum_{j=1}^{\ell} (j-2+h) \phd{j}$, the degree of the
last cochain complex is $\ds m_{1}(w) = \sum_{j=1}^{\ell}  \phd{j}$ and
the Betti number $\ds\dim\myH{m_1 (w)}{w} = 0$ for each $h$.    
\end{thm}

%% file: Euler-v4.tex
\section{Euler characteristic of homogeneous Poisson structures} 
\label{sec:Euler} 
We deal with Lie algebra cohomology groups of Lie algebra $\ds\ovfrakg =
\mR[x_1,\ldots,x_n]/\mR $ with a $h$-homogeneous Poisson
structure on   $\ds\mR^n$.  

For given non-negative integers $w$ and $m$,
$m$-th cochain space with the weight $w$ is given by 
\[ \myCC{m}{w} := \sum \Lambda^{k_1}\frakSS{1} \otimes
\Lambda^{k_2}\frakSS{2} \otimes  \cdots \otimes
\Lambda^{k_\ell}\frakSS{\ell}\] 
with the conditions
\begin{equation} k_1+k_2+\cdots = m \quad \text{and} \quad
\sum_{j=1}^{\infty} (j+h-2) k_j = w\;. \label{eqns} 
\end{equation}
In the subsection \ref{subsec:YD:cochain}, we know the sequences \(\ds [k_1,k_2,\ldots] \)
satisfying the above two conditions are equal to $\ds\mynabla{w+(2-h)m}{m}$ 
(= the set of Young diagrams with area $w+(2-h)m$ and of length $m$).   
Since 
the Euler characteristic of \(\ds \{\myHH{\bullet}{w}\}\)
is equal to 
\(\ds \sum_{m}(-1)^{m} \dim \myCC{m}{w} \), or  
that of 
\(\ds\{ \myH{m}{\bullet}\}\) is equal to  
\(\ds \sum_{m} (-1)^{m}\dim \myC{m}{w} \),  
we have to manipulate $\ds \dim\myCC{m}{w} 
= \sum_{\lambda\in\mynabla{w+(2-h)m}{m}} \dim \lambda $. 
Let $\lambda\in\mynabla{w}{m}$ and $\ds \lambda=[k_1,\ldots,k_p] $ with $k_p>0$
in our notation. The dimension of $\lambda$ as a direct summand of cochain
space $\ds\myCC{m}{w}$ is 
\[\ds \dim\lambda = \tbinom{\md{1}}{k_1} \cdots \tbinom{\md{p}}{k_p}\] 
where $\ds \md{j} := \dim \frakSS{j} = \mmd{j}$.

\kmcomment{
Then $\ds B\cdot \lambda = [1+k_1,k_2,\ldots, k_p]$ and 
\begin{equation} 
        \label{eqn:Bcdot}
        \dim(B\cdot\lambda) = \tbinom{\md{1}}{1+k_1}\tbinom{\md{2}}{k_2} 
        \cdots \tbinom{\md{p}}{k_p}
        = (-1+ \frac{1+n}{1+k_1}) \dim\lambda \end{equation} 

\begin{equation}
\dim \lambda = \binom{\md{1}}{k_1} \cdots \binom{\md{\ell}}{k_{\ell}}
        \quad \text{for}\quad  
        \lambda = [k_1,\ldots,k_{\ell}]\quad \text{in our notation.}
        \end{equation}
}%endOFkmcomment
%\subsection{Treatment for general $h$}
Before applying recursive formula (\ref{eqn:add:one}), we study about 
        $B\cdot\lambda$ and 
        $\T{m}\cdot\lambda$ for $\lambda\in\mynabla{w+(1-h)m}{m}$. 
        Since $\ds B\cdot \lambda = [1+k_1,k_2,\ldots,k_{p}]$, 
\begin{align} 
        \notag 
\dim (B\cdot\lambda) &= \tbinom{\md{1}}{1+k_1} \cdot\tbinom{\md{2}}{k_2} 
        \cdots \tbinom{\md{p}}{k_{p}}\;, \\
        \notag 
\frac{\dim(B\cdot\lambda)}{\dim\lambda} &= 
\frac{\tbinom{n}{1+k_1}}{\tbinom{n}{k_1}}  
= \frac{n-k_1}{1+k_1} = \frac{1+n}{1+k_1} -1 \;,  
\\
\noalign{thus, we have }
\label{B:dimen}
\dim(B\cdot \lambda) & 
= ( \frac{1+n}{1+k_1} -1 ) \dim\lambda\;, \text{ where }
\lambda=[k_1,k_2,\ldots] \text{ in our notation.} 
\end{align}
From Corollary \ref{cor:chuukan}, we see 
\begin{equation}
\label{Tm:dimen}
\dim(\T{m}\cdot \lambda) =
\tbinom{\md{1}}{0}
\prod_{i=2}^{p+1} \tbinom{\md{i}}{k_i}\;,   
 \text{ where }
 k_i = k'_{i-1} \;,\; 
\lambda=[k'_1,k'_2,\ldots,k'_p] \text{ in our notation.  } 
\end{equation}

Applying recursive formula (\ref{eqn:add:one})
        for 
$\ds\myCC{m}{w} = \mynabla{w+(2-h)m}{m}$, we have 
\begin{align}
        \myCC{m}{w} % &= \mynabla{w+(2-h)m}{m} 
        & = B\cdot \mynabla{w+(2-h)m-1}{m-1} \sqcup 
        \T{m}\cdot\mynabla{w+(1-h)m}{m} 
        = B\cdot \myCC{m-1}{w+1-h} \sqcup  
        \T{m}\cdot\mynabla{w+(1-h)m}{m} \;,   \label{eqn:rec:again}
\end{align}
        here $\ds\mynabla{w+(1-h)m}{m} = 0$ if $w+(1-h)m < m$.  
So, using (\ref{B:dimen}) we see that 
\begin{align*}
        \sum_{m} (-1)^m \dim \myCC{m}{w} =& 
        \sum_{m} (-1)^m \sum_{\lambda \in \myCC{m-1}{w+1-h}} \dim
        (B\cdot\lambda) + 
        \sum_{m} (-1)^m \dim 
        (\T{m}\cdot\mynabla{w+(1-h)m}{m})\\ 
        =& 
        \sum_{m} (-1)^m \sum_{\lambda=(k_1,k_2,\ldots)\in \myCC{m-1}{w+1-h}} 
        ( -1 + \frac{1+n}{1+k_1}) \dim \lambda  
        + \sum_{m} (-1)^m \dim 
        (\T{m}\cdot\mynabla{w+(1-h)m}{m}) 
        \\
        =& 
        \sum_{m} (-1)^{m-1} \dim \myCC{m-1}{w+1-h} +  
        \sum_{m} (-1)^m \sum_{\lambda=(k_1,k_2,\ldots)\in \myCC{m-1}{w+1-h}} 
         \frac{1+n}{1+k_1} \dim \lambda \\& 
        + \sum_{m} (-1)^m \dim 
        (\T{m}\cdot\mynabla{w+(1-h)m}{m}) \;. 
\end{align*}
Thus, we have the next lemma. 
\begin{Lemma} 
        \label{lemma:euler:recur} 
        Let $\ds \Oira{h}{w}$ be the Euler characteristic of 
        $\ds \myHH{\bullet}{w}$ for $h$-homogeneous Poisson structure. Then 
        we have the following recursive formula: 
\begin{align}
        \label{rec:general}
        \Oira{h}{w}  =& 
        \Oira{h}{w+1-h} +  
        \sum_{m} (-1)^m \sum_{\lambda=(k_1,k_2,\ldots)\in 
        %\mynabla{w-1+(2-h)m}{m-1}
        \myCC{m-1}{w+1-h} 
        } 
         \frac{1+n}{1+k_1} \dim \lambda \\& 
        \notag
        + 
        \sum_{m} (-1)^m \sum_{ \lambda=(k_1,k_2,\ldots)\in \mynabla{w+(1-h)m}{m}} 
        \prod_{i} \tbinom{\md{i+1}}{k_i} 
\;.  
\end{align}
\end{Lemma}

\begin{kmRemark} 
The above (\ref{rec:general}) shows $h=1$ is special and
(\ref{rec:general}) tells us nothing about $\Oira{h=1}{w}$.  
\kmcomment{
When $h=0$,  (\ref{rec:general}) looks
like ``backward'' recursive formula. 
}%endOFkmcomment
\end{kmRemark}

\kmcomment{
\begin{kmRemark} 
In (\ref{rec:general}) let us assume $\ds k_1=0$ for each $\lambda$, 
        namely we only deal
        with the subalgebra \(\sum_{k>1} \mySS{k}\).  
        Then 
        (\ref{rec:general}) implies 
        \begin{align}
                \Oira{h}{w}  =& 
                - n\; \Oira{h}{w+1-h}   
        + \sum_{m} (-1)^m \dim 
        (\T{m}\cdot\mynabla{w+(1-h)m}{m}) 
\end{align}
In particular, when $h=1$ then above yields 
\begin{equation} 
        \label{h:one:kOne:zero} 
        (1+n) 
        \Oira{h=1}{w}  = 
         \sum_{m} (-1)^m \dim 
         (\T{m}\cdot\mynabla{w}{ m})\;. \end{equation} 
 and we already know the right-hand side is 0
 in the proof of Euler characteristic is 0 for $h-1$
in the section  
\ref{sec:another:h1}.   
\end{kmRemark}
}%endOFkmcomment

When $h=0$, we have a very clear property. 
\begin{thm}
For each weight $w$, \(\ds 
        \Oira{h=0}{w}  = 
        \Oira{h=0}{w+1}  \) holds. In particular \(\ds  
        \Oira{h=0}{w}  = 0\).  
\end{thm}
\textbf{Proof:}
We simply show that 
the two terms in 
        \eqref{rec:general} of 
Lemma \ref{lemma:euler:recur} are zero.   
\(\ds
        \lambda=(k_1,k_2,\ldots)\in 
        \myCC{m-1}{w+1} \) says that 
\( \sum k_{j} = m-1 \) and 
\( \sum (j-2) k_{j} = w+1 \), so   
the first term is zero as below.  
\begin{align*} 
& 
\sum_{\lambda=(k_1,k_2,\ldots)\in \myCC{m-1}{w+1-h} } 
  (-1)^{1+ \sum k_{j} } \frac{1+n}{1+k_1} \prod_{j} \tbinom{\md{j}}{k_{j}} 
  \\
= &   \sum_{ \sum (j-2) k_{j} = w+1} (-1)^{1+ \sum k_{j} } 
         \frac{1+n}{1+k_1} \prod_{j} \tbinom{\md{j}}{k_{j}}
         \\
= &  {\color{red}(} \sum_{k_{2}} (-1)^{k_{2}}  \tbinom{\md{2}}{k_{2}}
{\color{red})}
\sum_{ \sum (j-2) k_{j} = w+1} (-1)^{1+ \sum_{j\ne 2} k_{j} } 
         \frac{1+n}{1+k_1} \prod_{j\ne 2} \tbinom{\md{j}}{k_{j}}
         = 0 \;. 
         \end{align*}
Concerning to the second term, 
        \(\ds \lambda=(k_1,k_2,\ldots)\in \mynabla{w+m}{m}\) means 
\(\ds \sum k_{j} = m \) and 
\(\ds \sum j k_{j} = w + m  = w + \sum k_{j} \), i.e., 
\(\ds \sum (j-1) k_{j} = w \). Thus, 
\begin{align*}
\sum_{m} (-1)^m \sum_{\lambda=(k_1,k_2,\ldots)\in \mynabla{w+m}{m} } 
        \prod_{i} \tbinom{\md{i+1}}{k_i}  
        = & 
\sum_{\sum (i-1) k_{i} = w} (-1)^{ \sum k_{i} } 
        \prod_{i} \tbinom{\md{i+1}}{k_i}  \\
= & {\color{red}(}
\sum_{k_{1}} (-1)^{k_{1}}  \tbinom{\md{2}}{k_1}  {\color{red})}
\sum_{\sum (i-1) k_{i} = w} (-1)^{ \sum_{i\ne 1} k_{i} } 
        \prod_{i\ne 1} \tbinom{\md{i+1}}{k_i}   = 0\;. 
\end{align*}
Here \(\ds k_{1}\) is an element of only Young diagram so not required
the dimensional restriction.  

Since we have known that \(\ds \Oira{h=0}{w}  = \Oira{h=0}{w+1}  \) for
each weight $w$, we show \(\ds  \Oira{h=0}{w}  = 0\) for some special
$w$.  Now we will show this at the minimum weight $-n$ where the space
dimension is $n$.  

From the definition of \(\ds \myCC{m}{-n}\), \([k_{1},k_{2},\ldots]\)
satisfies \(\ds \sum_{j} k_{j} = m\) and 
\(\ds \sum_{j} (j-2) k_{j} = -n\), and   
\(\ds k_{j} \leq \md{j} = \tbinom{n-1+j}{n-1} \) for each $j$.    
We see directly that \(\ds k_{1}=n\) and 
 \(\ds k_{j}=0\) for $j>2$, and \(\ds k_{2} = m-n\). Thus
 \( \ds \myCC{m}{-n} = \Lambda^{n} \frakSS{1} \otimes 
  \Lambda^{m-n} \frakSS{2} \).  
Therefore, 
\(\ds \sum_{m} (-1)^m \dim \myCC{m}{-n} = \sum_{m} (-1)^{m} 
\dim (\Lambda^{n} \frakSS{1})\dim  (\Lambda^{m-n} \frakSS{2}) 
= 
 \sum_{m} (-1)^{m} 
\dim (\Lambda^{m-n} \frakSS{2}) = 0\).  
\kmcomment{
Applying recursive formula \eqref{eqn:rec:again} for 
$\ds\mynabla{a}{m}$ two times, we have 
\begin{align*}
\mynabla{a}{m} 
=& B\cdot \mynabla{a-1}{m-1} \sqcup 
        \T{m}\cdot\mynabla{a-m}{m}  
= B\cdot ( B\cdot  
        \mynabla{a-2}{m-2} \sqcup \T{m-1}
        \mynabla{a-m}{m-1}) 
        \sqcup 
        \T{m}\cdot(
        B\cdot \mynabla{a-m-1}{m-1}  
        \sqcup \T{m}\cdot \mynabla{a-2m}{m} ) \\
=&  B^2 \cdot   
        \mynabla{a-2}{m-2} \sqcup 
      B  \T{m-1} \cdot \mynabla{a-m}{m-1} 
        \sqcup 
         \T{m} B\cdot \mynabla{a-m-1}{m-1}  \sqcup  \T{m}^2 
        \cdot \mynabla{a-2m}{m}  \\
=&  B^2 \cdot   
        \mynabla{a-2}{m-2} \sqcup 
       \T{m} \cdot \mynabla{a-m}{m-1} 
        \sqcup 
         \T{m} B\cdot \mynabla{a-m-1}{m-1}  \sqcup  \T{m}^2 
        \cdot \mynabla{a-2m}{m}  
\end{align*}
where we used a property \(\ds B\cdot\T{m-1}= \T{m}\).  
\kmcomment{three times : 
\begin{align*}
         \mynabla{a}{m} 
         =& B\cdot \mynabla{a-1}{m-1} \sqcup 
        \T{m}\cdot\mynabla{a-m}{m}  \\
         =& B\cdot ( B\cdot  
        \mynabla{a-2}{m-2} \sqcup \T{m-1}
        \mynabla{a-m}{m-1}) 
        \sqcup 
        \T{m}\cdot(
        B\cdot \mynabla{a-m-1}{m-1}  
        \sqcup \T{m}\cdot \mynabla{a-2m}{m} ) \\
         =& B^2 \cdot ( 
        B\cdot \mynabla{a-3}{m-3}\sqcup
        \T{m-2}\cdot \mynabla{a-m}{m-2}) 
        \sqcup B\T{m-1}( 
        B\cdot \mynabla{a-m-1}{m-2} 
        + \T{m-1} \cdot \mynabla{a-2m+1}{m-1}) 
        \\& 
        \sqcup 
        \T{m} B \cdot(
         B \cdot \mynabla{a-m-2}{m-2} 
        \sqcup 
         \T{m-1} \mynabla{a-2m}{m-1} 
         ) 
        \sqcup \T{m}^2 \cdot  ( 
        B\cdot\mynabla{a-2m-1}{m-1} 
        \sqcup 
        \T{m}\cdot \mynabla{a-3m}{m} 
        ) 
        \\
         =& B^3 \cdot  
         \mynabla{a-3}{m-3}  \sqcup
        \T{m}\cdot \mynabla{a-m}{m-2} 
        \sqcup  
        \T{m}B\cdot \mynabla{a-m-1}{m-2} 
        \sqcup  \T{m} \T{m-1} \cdot \mynabla{a-2m+1}{m-1} 
        \\& 
        \sqcup 
        \T{m} B^2 \cdot 
          \mynabla{a-m-2}{m-2} 
        \sqcup  \T{m}^2  \cdot
         \mynabla{a-2m}{m-1} 
        \sqcup 
         \T{m}^2 B\cdot\mynabla{a-2m-1}{m-1} 
        \sqcup \T{m}^3 \cdot  
         \mynabla{a-3m}{m} 
\end{align*}
where we used a property \(\ds B\cdot\T{m-1}= \T{m}\). 
}%endOFkmcomment 
Apply the above
equation for \(\ds a= -n + 2m\). Then we have 
\begin{align*}
\mynabla{-n+2m}{m} 
=& B^2 \cdot  \mynabla{-n+2m-2}{m-2}  \sqcup
        \T{m}\cdot \mynabla{-n + m}{m-1} \sqcup  
        \T{m}B\cdot \mynabla{-n+m-1}{m-1} \sqcup  
        \T{m}^2  \cdot \mynabla{-n}{m} 
        \\ 
=& B^2 \cdot  \mynabla{-n+2m-2}{m-2}  \sqcup
        \T{m}\cdot \mynabla{-n + m}{m-1}
\end{align*}
here we checked the relation of length and area of Young diagrams.  
\kmcomment{three times: 
\begin{align*}
         \mynabla{-n+2m}{m} 
         =& B^3 \cdot  
         \mynabla{-n+2m-3}{m-3}  \sqcup
        \T{m}\cdot \mynabla{-n +m}{m-2} 
        \sqcup  
        \T{m}B\cdot \mynabla{-n+m-1}{m-2} 
        \sqcup  \T{m} \T{m-1} \cdot \mynabla{-n+1}{m-1} 
        \\& 
        \sqcup 
        \T{m} B^2 \cdot 
          \mynabla{-n+m-2}{m-2} 
        \sqcup  \T{m}^2  \cdot
         \mynabla{-n}{m-1} 
        \sqcup 
         \T{m}^2 B\cdot\mynabla{-n-1}{m-1} 
        \sqcup \T{m}^3 \cdot  
         \mynabla{-n-m}{m} 
\\ =& B^3 \cdot  \mynabla{-n+2m-3}{m-3}  
\sqcup \T{m}\cdot \mynabla{-n +m}{m-2} 
\sqcup  \T{m}B\cdot \mynabla{-n+m-1}{m-2} 
\sqcup \T{m} B^2 \cdot \mynabla{-n+m-2}{m-2}\;.  
\end{align*}
}%endOFkmcomment
Thus, 
\begin{align} 
\Oira{h=0}{w=-n} =& 
\sum_{m} (-1)^{m} \sum_{ \lambda \in
\mynabla{-n+2m-2}{m-2}} \dim ( B^2\cdot \lambda ) + 
\sum_{m} (-1)^{m} \sum_{ \lambda \in
\mynabla{-n+m}{m-1}} \dim( \T{m}\cdot \lambda ) 
\label{h0:w:-n}
\end{align} 
Again, we calculate the two terms in \eqref{h0:w:-n}
separately.  
Take an element \(\ds \lambda = [k'_{1}, k'_{2},\ldots] \in 
\mynabla{-n+2m-2}{m-2}\). They satisfy 
\(\ds \sum (j-2) k'_{j} = -n +2\) and \(\ds k'_{2}\) is free.  
\(\ds B^2\cdot \lambda = [k_{1},k_{2},
\ldots]\) with \(\ds k_{1} = 2+ k'_{1}\) and 
\(\ds k_{j} =  k'_{j}\)  for \(j > 1\).  Thus, \(\ds k_{2}\) is free. 
\begin{align*}
\text{The first term of \eqref{h0:w:-n} } =& \sum_{  
\sum (j-2) k'_{j} = -n +2 } (-1)^{\sum k_j} \prod
\tbinom{\md{j}}{k_j} \\
=& \sum_{k_2} (-1)^{k_2} \tbinom{\md{2}}{k_2} 
\sum_{  
\sum (j-2) k'_{j} = -n +2 } (-1)^{\sum_{j\ne 2} k_j} \prod_{j \ne 2} 
\tbinom{\md{j}}{k_j} = 0\; .
\end{align*}
About the second term in \eqref{h0:w:-n}, 
it makes sense only when $n=1$ and it is
impossible to deal with 2-vector field.  
Even though, we continue formal discussion. 
Take an element \(\ds \lambda = [k'_{1}, k'_{2},\ldots] \in 
\mynabla{-n+m}{m-1}\). They satisfy 
\(\ds \sum (j-1) k'_{j} = -n +1\) and \(\ds k'_{1}\) is free.  
We put 
\(\ds \T{m}\cdot \lambda = [k_{1},k_{2}, \ldots]\). Then 
\(\ds k_{1} = m - \sum k'_{j} = 1\) and 
\(\ds k_{j} =  k'_{j-1}\)  for \(j > 1\).  Thus, \(\ds k_{2}\) is free.  
\begin{align*}
\text{The second term of \eqref{h0:w:-n} } =& \sum_{  
\sum (j-1) k'_{j} = -n +1 } (-1)^{\sum k_j} \prod
\tbinom{\md{j}}{k_j} \\
=& \sum_{k_2} (-1)^{k_2} \tbinom{\md{2}}{k_2} 
\sum_{  
\sum (j-2) k'_{j} = -n +2 } (-1)^{\sum_{j\ne 2} k_j} \prod_{j \ne 2} 
\tbinom{\md{j}}{k_j} = 0\; .
\end{align*} 
}%endOFkmcomment
\kmqed 

\kmcomment{
%%%%%%%%%%%%
        \sum_{m} (-1)^m \sum_{\lambda=(k_1,k_2,\ldots)\in 
        \myCC{m-1}{w+1-h} 
        } 
         \frac{1+n}{1+k_1} \dim \lambda \\& 
        \notag
        + 
        \sum_{m} (-1)^m \sum_{ 
        \lambda=(k_1,k_2,\ldots)\in \mynabla{w+(1-h)m}{m}
        } 
        \prod_{i} \tbinom{\md{i+1}}{k_i} 
%%%%%%%%%%%%%%%%%%%%%%%%%%%%%%%%%%%%%%%%%%%%%%%%%%%
}%endOFkmcomment

\subsection{Our theorem and its proof}
In the previous section, we have dealt with some concrete examples and we 
expected the Euler characteristic of Lie algebra cohomology of Lie Poisson 
structure is zero.  We have the following theorem for Lie Poisson structures. 
\begin{thm} \label{thm:euler:num} 
On $\ds\mR^n$, consider a Poisson structure of homogeneity 1. Then for
each given weight $w$, the alternating sum of $\ds  \dim\; \myCC{m}{w}$
is 0, namely, the Euler characteristic of $\ds\{\myHH{\bullet}{w}\}$ is 0.
Also, the alternating sum of $\ds  \dim\; \myC{m}{w}$ is 0, namely,
the Euler characteristic of $\ds\{\myH{\bullet}{w}\}$ is 0.  
\end{thm}

\textbf{Proof:}
Since $h=1$, $m$-th cochain space with weight $w$ is 
$\ds  \myCC{m}{w} = \mynabla{w+m}{m}$. 
When $w=0$, then $\ds\myCC{m}{0} = \mynabla{m}{m}=\{\T{m}\} = \{(k_1=m)\}$
       and $\ds \dim\myCC{m}{0} = \dim\Lambda^{m}\frakSS{1} = \tbinom{n}{m}$
       and so $\ds \sum_{m}(-1)^m \dim\myCC{m}{0} = 
       \sum_{m}(-1)^m \tbinom{n}{m} = 0$. 

When $w>0$, applying  
(\ref{rmk:recur}) to $\ds\myCC{m}{w}$, we have 
\[
        \myCC{m}{w} = \mynabla{w+m}{m} = \T{m} \cdot ( \mynabla{w}{1}+\cdots +
\mynabla{w}{w})\;.  \]
\kmcomment{
We remark that 
the right-hand side term above appeared in (\ref{h:one:kOne:zero}).
}%endOFkmcomment
In Corollary \ref{cor:chuukan}, we see that 
if a Young diagram $\lambda = [k_1,k_2\ldots]$ then 
$\ds\T{m}\cdot\lambda = [m-\sum_{i}k_i, k_1,k_2,\ldots]$. 
Thus if $\lambda\in\mynabla{w}{j}$ then 
$\ds\T{m}\cdot\lambda = [m-j, k_1,k_2,\ldots]$ and 
\begin{align*}
        \dim(\T{m}\cdot\lambda)  &= 
\tbinom{n}{m-j} 
\tbinom{\md{2}}{k_1} 
\tbinom{\md{3}}{k_2} \cdots\;. \\ 
\noalign{Therefore,}
%\begin{align*}
%        \myCC{m}{w} &= \mynabla{w+m}{m} = \T{m} \cdot ( \mynabla{w}{1}+\cdots +
%        \mynabla{w}{w}) \\
        \sum_{m}(-1)^m \dim \myCC{m}{w} &= \sum_{m} \sum_{j=1}^ w 
(-1)^{m}
\dim (
        \T{m}\cdot \mynabla{w}{j}) 
        = \sum_{j=1}^w \sum_{\lambda\in\mynabla{w}{j}} \sum_{m} 
(-1)^{m}
        \dim (
        \T{m}\cdot \lambda) \\ 
        &= \sum_{j=1}^w \sum_{\lambda=(k_1,\ldots,k_p)\in\mynabla{w}{j}} 
        (\sum_{m}
(-1)^{m} 
        \tbinom{n}{m-j}) \tbinom{\md{2}}{k_1} \cdots 
       \tbinom{\md{p+1}}{k_p} 
        = \sum_{j=1}^w \sum_{\lambda\in\mynabla{w}{j}} 0  = 0.  
       \end{align*} 
\kmqed

\kmcomment{ 
\textbf{Proof:}
When $w=0$, from (\ref{h1:w0}), $\ds\myCC{m}{0} = \Lambda^m \frakS{1}$ and 
\( \sum_{m} (-1)^m \dim\; \myCC{m}{0} = \sum_{m} (-1)^m \dim\;( \Lambda^{m}
\frakS{1}) = 0\).   
When $w=1$, from (\ref{h1:w1}) 
we have \(\ds  \myCC{m}{1} = \Lambda^{m-1}\frakSS{1} \otimes \frakSS{2}\).  
Thus, 
\[ \sum_{m} (-1)^m \dim\; \myCC{m}{1} = \sum_{m} (-1)^m \dim\;( \Lambda^{m-1}
\frakSS{1}) \cdot \dim\; \frakSS{2} = (\sum_{m} (-1)^m \dim\;( \Lambda^{m-1}
\frakSS{1})) \cdot \dim\; \frakSS{2} = 0\;. \]  

When $w=2$, from  (\ref{h1:w2}) we have
\( \myCC{m}{2} = \Lambda^{m-1} \frakSS{1}\otimes \frakSS{3} \oplus  
\Lambda^{m-2} \frakSS{1}\otimes \frakSS{2}\)  
where $\ds  \Lambda^{k}\frakSS{1} = (0)$ for $k<0$. 
Thus, we see 
\begin{align*} \sum_{m}(-1)^m \dim\; \myCC{m}{2} =&
\sum_{m} (-1)^m 
\dim\; \Lambda^{m-1} \frakSS{1} \cdot \dim\;  \frakSS{3} 
+ 
\sum_{m} (-1)^m \dim \; 
\Lambda^{m-2} \frakSS{1} \cdot \dim\; \frakSS{2}
=  0\; .
\end{align*} 
When $w=3$, in a similar way, from (\ref{h1:w3}) we have
\[ \myCC{m}{3} = \Lambda^{m-1} \frakSS{1}\otimes \frakSS{4} \oplus  
\Lambda^{m-2} \frakSS{1}\otimes \frakSS{2} \otimes \frakSS{3} \oplus  
\Lambda^{m-3} \frakSS{1}\otimes \Lambda^{3} \frakSS{2} 
\;. 
\] 
Thus, we see 
\begin{align*}
       &  \sum_{m} (-1)^m \dim\;   \myCC{m}{3} \\
=&
 \sum_{m} (-1)^m \dim\;   (\Lambda^{m-1} \frakSS{1}\otimes \frakSS{4}) 
 + \sum_{m} (-1)^m \dim\;   
( \Lambda^{m-2} \frakSS{1}\otimes \frakSS{2} \otimes \frakSS{3}) \\& 
+ \sum_{m} (-1)^m \dim\;   
( \Lambda^{m-3} \frakSS{1}\otimes \Lambda^{3} \frakSS{2} )
\\=& 
 \sum_{m} (-1)^m \dim\;   
 (\Lambda^{m-1} \frakSS{1}) \dim (\frakSS{4}) 
 + \sum_{m} (-1)^m \dim\;   (
\Lambda^{m-2} \frakSS{1}) \dim( \frakSS{2} \otimes \frakSS{3}) \\&
+ \sum_{m} (-1)^m \dim\;   (
\Lambda^{m-3} \frakSS{1}) \dim (\Lambda^{3} \frakSS{2} )
\\
=& 0\;. 
\end{align*} 
In general, for a given weight $w$, from (\ref{rmk:recur}) we have
\begin{align*}
\mynabla{w+1}{1} &= \T{1}\cdot \mynabla{w}{1}\;, \\
\mynabla{w+2}{2} &= \T{2}\cdot( \mynabla{w}{1} \sqcup \mynabla{w}{2})\;, \\
\vdots & \\
\mynabla{w+w-1}{w-1} &= \T{w-1}\cdot( \mynabla{w}{1} \sqcup \cdots \sqcup
\mynabla{w}{w-1})\;, \\
\noalign{and if $w \leq m$}
\mynabla{w+m}{m} &= \T{m}\cdot( \mynabla{w}{1} \sqcup \cdots \sqcup \mynabla{w}{w})
\;. 
\end{align*} 
We point out that the second factor of the right-hand-side of the last
equation is independent of $m$. Now, we expand the equations above and
compute the dimension vertically. 

About the alternating sum of the first terms (sum along the most left
vertical line), we get the alternating sum $\ds\sum_{m} (-1)^m
\sum_{\lambda\in \mynabla{w}{1}} \dim (\T{m}\cdot \lambda)$, here for a
Young diagram \(\lambda\) we overuse the notation \(\ds\dim (\lambda) \) which
means \(\ds\dim (\Lambda^{k_1}\frakS{1} \otimes \Lambda^{k_2}\frakS{2}
\otimes \cdots)\) when \(\lambda \) corresponds to \( (k_1,k_2,\ldots) \).  

As we have seen in (\ref{eqn:concat}), for $\ds\lambda\in\mynabla{w}{1}$
$\T{m}\cdot \lambda$ corresponds to $(k_1=m-1, k_{w+1}=1)$ with $k_j =
0$ ($j\ne 1, w+1$).  Thus,  
\begin{align*}
\sum_{m} (-1)^m \sum_{\lambda\in \mynabla{w}{1}} \dim (\T{m}\cdot \lambda) 
        =& \sum_{m} (-1)^m \dim ( \Lambda^{m-1}\frakS{1} \otimes \frakS{w+1} )
        \\
        =& \left(\sum_{m} (-1)^m 
\dim ( \Lambda^{m-1}\frakS{1})\right)  \dim \frakS{w+1} = 0\;.\end{align*}  

About the $j$-th vertical line from the left ($j\leq w$), we consider
the alternating sum 
\[\ds\sum_{m\geq j} (-1)^m \sum_{\lambda\in \mynabla{w}{j}} \dim (\T{m}\cdot
\lambda)\;. \]

By (\ref{eqn:concat}), for $\ds\lambda\in\mynabla{w}{j}$ with $(k_1,\ldots)$ 
$\T{m}\cdot \lambda$ corresponds to $(\kd{1}=m-j, \kd{2}= k_{1}, 
 \kd{3}= k_{2}, \ldots )$.  Thus,  
\begin{align*}
        & \sum_{m\geq j } (-1)^m \sum_{\lambda\in \mynabla{w}{j}} \dim (\T{m}\cdot
        \lambda) \\
        =& \sum_{m\geq j} (-1)^m 
 \sum_{\lambda\in \mynabla{w}{j}}
\dim 
( \Lambda^{m-j}\frakS{1} \otimes \Lambda^{k_1} \frakS{2} \otimes 
 \Lambda^{k_2} \frakS{3} \otimes \cdots )
        \\
        =& \left(\sum_{m\geq j } (-1)^m 
\dim ( \Lambda^{m-j}\frakS{1})\right)  
\sum_{\lambda\in \mynabla{w}{j}} 
\dim( \Lambda^{k_1} \frakS{2} \otimes \Lambda^{k_2} \frakS{3} \otimes 
\cdots )  = 0\; .
\end{align*}  
Those say that the alternating sum of the dimension of the young
diagrams on each  vertical line is zero and since the alternating sum of
$\ds\dim \myCC{m}{w}$ is the sum of them, so the alternating sum of
$\ds\dim \myCC{m}{w}$ is zero. 

For $\ds\myC{m}{w}$, 
$\ds\dim \frakS{k}$ may differ from  that of 
$\ds\frakSS{k}$,  but in the discussion above we only used the fact that 
$\ds\sum_{k}(-1)^k  \dim ( \Lambda^k \frakSS{1}) =0$ and still 
$\ds\sum_{k}(-1)^k  \dim ( \Lambda^k \frakS{1}) =0$ holds.  
\kmqed
}%endOFkmcomment

%\subsection{Examples with $h=2$ where our theorem \ref{thm:euler:num} fails} 
\subsection{Examples with $h=2$} 
Since the Euler characteristic is given by \(\ds \sum_{m}
(-1)^{m}\dim\myCC{m}{w}\), and each \(\ds \myCC{m}{w}\) depends on $n$,
$w$, $m$ and only the homogeneity $h$ of Poisson structures, if we pick
up some $h$-homogeneous Poisson structure and calculate all the Betti
numbers and get the Euler characteristic, then that Euler characteristic
is common for the same $n$ and $h$.  

Here, we handle homogeneity 2 cases and show that the Euler
characteristic is not necessarily zero.  Since the normal form of
analytic Poisson structures are studied by J.F.~Conn  (\cite{MR744864})
and these are locally Lie Poisson structures, we have several cases of
2-homogeneous Poisson structures on $\ds\mR^3$.  
\begin{align*}
\text{case-1:}\quad  & 
\pi =\frac{1}{2} ( x_1{}^2  \pdel_2 \wedge \pdel_3
+ x_2{}^2  \pdel_3 \wedge \pdel_1  + x_3{}^2  \pdel_1 \wedge \pdel_2)
\;, \\
\text{case-2:}\quad & 
\pi = x_1 x_2  \pdel_1 \wedge \pdel_2 + x_2 x_3   \pdel_2 \wedge \pdel_3
+ x_3 x_1  \pdel_3 \wedge \pdel_1 \;, \\
\text{case-3:}\quad & 
\pi = x_1{}^2   \pdel_2 \wedge \pdel_3
+ x_3 x_1   \pdel_3 \wedge \pdel_1  + x_1 x_2  \pdel_1 \wedge \pdel_2
\;. 
\end{align*} 
Since $h=2$, we see $\ds  \Bkt{\mySS{i}}{\mySS{j}} \subset
\mySS{i+j}$ and so we have
\begin{align*}
\mydov(\frakSS{1}) =& (0), \quad 
\mydov(\frakSS{2}) \subset \Lambda^2 \frakSS{1}, \quad 
\mydov(\frakSS{3}) \subset \frakSS{1} \otimes \frakSS{2}, \quad 
\mydov(\frakSS{4}) \subset \frakSS{1}\otimes \frakSS{3} \oplus
\Lambda^2 \frakSS{2}\;, 
\\
\mydov(\frakSS{5}) \subset& \frakSS{1}\otimes \frakSS{4} \oplus
\frakSS{2}\otimes \frakSS{3}\;. 
\end{align*} 
The list below are the cochain complexes of n=3 (3 variables),
homogeneity is 2.  The subindex of $\ds  \myCC{\bullet}{w}$ is the
weight.     
\begin{align*} 
 \myCC{1}{1}=&  \frakSS{1}\ (3\;\dim)\;, \\ 
 \myCC{1}{2}=&  \frakSS{2}\ (6\;\dim)\; , \quad 
 \myCC{2}{2}=   \LS{1}{2}\ (3\;\dim)\;,  \\
 \myCC{1}{3}=&  \frakSS{3}\ (10\;\dim)\;, \quad 
 \myCC{2}{3}=   \frakSS{1} \otimes  \frakSS{2}\ (18\;\dim)\;, \quad 
 \myCC{3}{3}=   \LS{1}{3}\ (1\;\dim)\;,  \\
 \myCC{1}{4}=&  \frakSS{4}\ (15\;\dim)\;, \quad 
 \myCC{2}{4}= \frakSS{1} \otimes \frakSS{3} + \LS{2}{2} \ (45\;\dim)\;, \quad 
 \myCC{3}{4}=   \LS{1}{2} \otimes  \frakSS{2}\ (18\;\dim)\;,  \\
 \myCC{1}{5}=&  \frakSS{5}\ (21\;\dim)\; , \quad 
 \myCC{2}{5}=  (\frakSS{1} \otimes \frakSS{4})+ (\frakSS{2} \otimes \frakSS{3}) \ (105\;\dim) \\ 
 \myCC{3}{5}=&   (\LS{1}{2} \otimes \frakSS{3}) + (\frakSS{1} \otimes
 \LS{2}{2}) \ (75\;\dim)\; , \quad 
 \myCC{4}{5}=  \LS{1}{3} \otimes \frakSS{2}\ (6\;\dim)\;, \\
%%%%%%%%%%%%%%%%%%%%%%%%%%%%%%%%%%%%%%%%%%%%%%%%%%%%%%%%%%%%%%%%%%%%%%
\myCC{1}{6}=&\frakSS{6} \; (28\; \dim)\;, \quad 
\myCC{2}{6}=
 \frakSS{1}\otimes \frakSS{5} 
+ \frakSS{2}\otimes \frakSS{4} 
+ \LS{3}{2} 
\;(198\;\dim)\;,\quad 
\\
\myCC{3}{6}=& 
\LS{1}{2}\otimes \frakSS{4}
+ \frakSS{1}\otimes \frakSS{2}\otimes \frakSS{3} 
+ \LS{2}{3} 
\;(245\;\dim)\;, \\  
\myCC{4}{6}=&  \LS{1}{3}\otimes \frakSS{3} + \LS{1}{2}\otimes \LS{2}{2} 
\; (55\;\dim)\;, \quad 
% 6,2,"Eulernumberis",-20
\\
\myCC{1}{7}=&  \frakSS{7} \; (36\;\dim)\;, \quad 
\myCC{2}{7}= \frakSS{1}\otimes \frakSS{6} 
+ \frakSS{2}\otimes \frakSS{5} 
+ \frakSS{3}\otimes \frakSS{4} \; 
(360\;\dim)\;, \quad 
\\
\myCC{3}{7}=&
  \LS{1}{2}\otimes \frakSS{5} 
 + \frakSS{1}\otimes \frakSS{2}\otimes \frakSS{4} 
 + \frakSS{1}\otimes \LS{3}{2} 
 + 
 \LS{2}{2}\otimes \frakSS{3} 
 \; (618\;\dim)\;, \quad 
\\
\myCC{4}{7}=&
 \LS{1}{3}\otimes \frakSS{4}
 + \LS{1}{2}\otimes \frakSS{2}\otimes \frakSS{3}
 + \frakSS{1}\otimes \LS{2}{3} 
 \; (255\;\dim)\;, \quad 
\myCC{5}{7}= \LS{1}{3}\otimes \LS{2}{2} \; (15\;\dim)\;. 
%7,2,Eulernumberis-54 
 \end{align*}
From the list above, we see that the Euler characteristic for each weight
varies as follows when $n=3$: 
\begin{center} 
\tabcolsep=6pt 
\setlength{\extrarowheight}{-1pt}
\begin{tabular}{c|*{7}{c}}
weight & 1 & 2 & 3 & 4 & 5 & 6 & 7\\
\hline
Euler number   & $-3$ & $-3$ & 7 & 12 & 15 & $-20$ & $-54$
\end{tabular}
\end{center}
%\begin{kmRemark}
On the symplectic space $\ds\mR^2$, namely for homogeneity 0 non-degenerate Poisson
structure, we know some subalgebras whose Euler characteristics are not necessarily
zero (cf.\   
\cite{KOT:MORITA}, \cite{M:N:K}).  
%\end{kmRemark}
\subsection{Abstract way to get Euler characteristic for $h=2$} 
Since the 
cochain space $\ds\myCC{m}{w}$ with $h$ corresponds to   
$\ds\mynabla{ w+(2-h)m}{m}$, when $h=2$  
the cochain space \[\ds\myCC{m}{w} = \mynabla{ w}{m} \quad \text{for}\quad  m
=0,\ldots, w.\] 

Concerning to $m=0$, 
$\ds\myCC{0}{w} = 0$ 
when $w>0$,   
and  
$\ds\myCC{0}{0} = \mR$ as an unique exception.

\subsubsection{$h=2$ and $w=0$}
$w=0$ is exceptional and $\ds \myCC{0}{0} = \mR$, thus \(\Oira{h=2}{w=0} = 1\). 

\subsubsection{$h=2$ and $w=1$}
$\ds\myCC{0}{1}= 0$ and 
$\ds\myCC{1}{1}= \mynabla{1}{1} = \{\T{1}\} = \{(k_1=1)\}$. Thus 
\(\Oira{h=2}{w=1} = - \dim \frakSS{1} = -n\). 

\kmcomment{
\subsubsection{$h=2$ and $w=2$}
\begin{align*}
\myCC{1}{2} &= \mynabla{2}{1} = \T{1}^2 = ( k_2=1) = \frakSS{2} \\
\myCC{2}{2} &= \mynabla{2}{2} = \T{2} = ( k_1=2) = \Lambda^2\frakSS{1} 
\\ 
\Oira{}{w=2} & = - \dim \frakSS{2} + \dim \Lambda^2\frakSS{1} 
= - \tbinom{n-1+2}{2} + \tbinom{n}{2} = -n\;. 
\end{align*}

\subsubsection{$h=2$ and $w=3$}
\begin{align*}
\myCC{1}{3} &= \mynabla{3}{1} = \T{1}^3 = ( k_3=1) = \frakSS{3} \\
\myCC{2}{3} &= \mynabla{3}{2} = \T{2}\T{1} = ( k_1=1, k_2=1) = 
\frakSS{1} \otimes \frakSS{2} \\ 
\myCC{3}{3} &= \mynabla{3}{3} = \T{3} = ( k_1=3) = 
\Lambda^{3} \frakSS{1} \\ 
\Oira{}{w=3} & = - \dim \frakSS{3} + 
\dim \frakSS{1} \dim \frakSS{2} 
- \dim \Lambda^3\frakSS{1} 
\\&
= - \tbinom{n-1+3}{3} + \tbinom{n-1+1}{1} \tbinom{n-1+2}{2} 
- \tbinom{n}{3} 
= \frac{n}{6} (n-1)(n+4) \;. 
\end{align*}

\subsubsection{$h=2$ and $w=4$}
\begin{align*}
\myCC{1}{4} &= \mynabla{4}{1} = \T{1}^4 = ( k_4=1) = \frakSS{4} \\
\myCC{2}{4} &= \mynabla{4}{2} = \T{2}\T{1}^2+ \T{2}^2 = ( k_1=1, k_3=1),(
k_2=2) \\&
= \frakSS{1} \otimes \frakSS{3} + \Lambda^2 \frakSS{2} \\ 
\myCC{3}{4} &= \mynabla{4}{3} = \T{3}\T{1} = ( k_1=2, k_2=1 ) = 
\Lambda^{2} \frakSS{1} \otimes \frakSS{2}\\ 
\myCC{4}{4} &= \mynabla{4}{4} = \T{4} = ( k_1=4 ) = 
\Lambda^{4} \frakSS{1} \\ 
\Oira{}{w=4} & = - \dim \frakSS{4} +
(
\dim \frakSS{1} \dim \frakSS{3} + \dim  \Lambda^2 \frakSS{2})  
- \dim (\Lambda^2\frakSS{1} \otimes \frakSS{2}) 
+ \dim \Lambda^{4} \frakSS{1}  
\\&
= - \tbinom{n-1+4}{4} + \tbinom{n-1+1}{1} \tbinom{n-1+3}{3} 
+ \binom{ \tbinom{n-1+2}{2}}{2} - \tbinom{n}{2} \tbinom{n-1+2}{2} 
+ \tbinom{n}{4} 
\\& 
= \frac{1}{24} n(n-1)(n^2+7n+18) \;. 
\end{align*}
}%endOFkmcomment
\subsubsection{$h=2$ and general $w$}
From Lemma \ref{lemma:euler:recur}, if we put $h=2$ then we have a 
relation between 
$\ds \Oira{}{[w]}$ and 
$\ds \Oira{}{[w-1]}$ as follows:  
        \begin{equation}
                \label{eqn:recursive}
                \Oira{h=2}{w} = 
                \Oira{h=2}{w-1} 
        + 
        \sum_{m=1}^{w} (-1)^m \sum_{\lambda\in \mynabla{w-1}{m-1}}
        \frac{1+n}{1+k_1} \dim\lambda  
+ 
\sum_{m=1}^{w} (-1)^m \sum_{\lambda\in\mynabla{w-m}{m}}  
        \tbinom{\md{2}}{k_1} \cdots \tbinom{\md{p+1}}{k_p}
        \end{equation}
        \kmcomment{
We already have 
$\ds \Oira{}{[0]}$, $\ds \Oira{}{[1]}$, $\ds \Oira{}{[2]}$, 
$\ds \Oira{}{[3]}$, and $\ds \Oira{}{[4]}$. 
So, we experimental for lower $w$ if (\ref{eqn:recursive}) is true. 
}%endOFkmcomment
\paragraph{$w=2$:}
\begin{align*}
        \text{2nd term of } (\ref{eqn:recursive}) =&
        \sum_{m=1}^{2} (-1)^m \sum_{\lambda\in \mynabla{1}{m-1}}
        \frac{1+n}{1+k_1} \dim\lambda  
= (-1)^2 \sum_{\lambda\in \mynabla{1}{1}}
\frac{1+n}{1+k_1} \dim\lambda  = \frac{1+n}{2} n \\ 
\text{3rd term of } (\ref{eqn:recursive}) =& 
\sum_{m=1}^{2} (-1)^m \sum_{\lambda\in\mynabla{2-m}{m}}  
        \tbinom{\md{2}}{k_1} \cdots \tbinom{\md{p+1}}{k_p}
= (-1)^1 \sum_{\lambda\in\mynabla{1}{1}}  
        \tbinom{\md{2}}{k_1} \cdots \tbinom{\md{p+1}}{k_p}
        \\
        =&  -  \tbinom{\md{2}}{1} = - \md{2} = - \tbinom{n-1+2}{2}
        \end{align*}
        Thus, $\ds \Oira{h=2}{2} = \Oira{h=2}{1} + 
\text{2nd term} + 
\text{3rd term} = \Oira{h=2}{1} + \frac{1+n}{2}n - \tbinom{n+1}{2} = 
\Oira{h=2}{1} = -n$.

\paragraph{$w=3$:}
\begin{align*}
        \text{2nd term of } (\ref{eqn:recursive}) =&
        \sum_{m=1}^{3} (-1)^m \sum_{\lambda\in \mynabla{2}{m-1}}
        \frac{1+n}{1+k_1} \dim\lambda  
        \\
        = &
(-1)^2 \sum_{\lambda\in \mynabla{2}{1}} \frac{1+n}{1+k_1} \dim\lambda  
+(-1)^3 \sum_{\lambda\in \mynabla{2}{2}} \frac{1+n}{1+k_1} \dim\lambda  
\\
= &
\frac{1+n}{1+0}( \dim\T{1}^2)  
- \frac{1+n}{1+2} \dim\T{2} 
= (1+n) \mmd{2} - \frac{1+n}{3} \tbinom{n}{2}
\\ 
\text{3rd term of } (\ref{eqn:recursive}) =& 
\sum_{m=1}^{3} (-1)^m \sum_{\lambda\in\mynabla{3-m}{m}}  
        \tbinom{\md{2}}{k_1} \cdots \tbinom{\md{p+1}}{k_p}
= (-1)^1 \sum_{\lambda\in\mynabla{2}{1}}  
        \tbinom{\md{2}}{k_1} \cdots \tbinom{\md{p+1}}{k_p}
        \\
        =&  -  \tbinom{\md{2}}{0} \tbinom{\md{3}}{1} 
        = - \md{3} = - \mmd{3}
        \end{align*}
        Thus, 
        \begin{align*} \Oira{h=2}{3} =& \Oira{h=2}{2} + 
\text{2nd term} + 
\text{3rd term} = -n  
+(1+n) \mmd{2} - \frac{1+n}{3} \tbinom{n}{2}
- \mmd{3}\\ =& \frac{1}{6} n(n-1)(n+4) \;.
\end{align*}

\paragraph{$w=4$:}
\begin{align*}
        \text{2nd term of } (\ref{eqn:recursive}) =&
        \sum_{m=1}^{4} (-1)^m \sum_{\lambda\in \mynabla{3}{m-1}}
        \frac{1+n}{1+k_1} \dim\lambda  
        \\
        = &
 (-1)^2 \sum_{\lambda\in \mynabla{3}{1}} \frac{1+n}{1+k_1} \dim\lambda  
+(-1)^3 \sum_{\lambda\in \mynabla{3}{2}} \frac{1+n}{1+k_1} \dim\lambda  
+(-1)^4 \sum_{\lambda\in \mynabla{3}{3}} \frac{1+n}{1+k_1} \dim\lambda  
\\
= &
\frac{1+n}{1+0}( \dim\frakSS{3} )  
- \frac{1+n}{1+2} \dim\frakSS{1}  \dim\frakSS{2} 
+ \frac{1+n}{1+3} \binom{\dim\frakSS{1}}{3}  
\\
=& (1+n) \mmd{3} - \frac{1+n}{3}\mmd{1}\mmd{2} + \frac{1+n}{4}\tbinom{n}{3}
\\=& - \frac{1}{24}n^4 + \frac{1}{12} n^3 + \frac{13}{24}n^2 + \frac{5}{12} n
\\ 
\text{3rd term of } (\ref{eqn:recursive}) =& 
\sum_{m=1}^{4} (-1)^m \sum_{\lambda\in\mynabla{4-m}{m}}  
        \tbinom{\md{2}}{k_1} \cdots \tbinom{\md{p+1}}{k_p}
        \\
        =& (-1)^1 \sum_{\lambda\in\mynabla{3}{1}}  
        \tbinom{\md{2}}{k_1} \cdots \tbinom{\md{p+1}}{k_p}
        + 
(-1)^2 \sum_{\lambda\in\mynabla{2}{2}}  
        \tbinom{\md{2}}{k_1} \cdots \tbinom{\md{p+1}}{k_p}
        \\
        =&  -  \tbinom{\md{4}}{1} + \tbinom{\md{2}}{2} 
        = \frac{1}{12}n^4 - \frac{7}{12} n^2 - \frac{1}{2} n
        \end{align*}
        Thus, $\ds \Oira{h=2}{4} =  \Oira{h=2}{3} + 
\text{2nd term} + 
\text{3rd term} = 
\frac{1}{24}n^4 + \frac{1}{4} n^3 + \frac{11}{24}n^2 - \frac{3}{4} n
= \frac{1}{24} n(n-1)(n^2 + 7 n + 18) $. 
\begin{kmRemark}
        It seems to be interesting to find the generating function of
        $\ds \{ \Oira{h=2}{w}\}_{w}$.
\end{kmRemark}
\subsection{Euler characteristic for $h=3$} 
We put $h=3$ in ( \ref{rec:general}) 
of Lemma \ref{lemma:euler:recur}, we have a recursive formula for $h=3$:
\begin{align*}\Oira{h=3}{w} = \Oira{h=3}{w-{\color{red}2}}  & + 
        \sum_{m} (-1)^m \sum_{\lambda=(k_1,k_2,\ldots)\in 
        \mynabla{w-1-m}{m-1}} 
         \frac{1+n}{1+k_1} \dim \lambda 
         \\&
        + 
        \sum_{m} (-1)^m \sum_{ \lambda=(k_1,k_2,\ldots)\in \mynabla{w-2m}{m}} 
        \prod_{i} \tbinom{\md{i+1}}{k_i} 
\;.  
\end{align*} 
We use the above recursive formula or original definition 
\(\ds \dim\myCC{m}{w} = \sum_{\lambda \in \mynabla{w-m}{m}} \dim \lambda\) 
and 
we see the Euler characteristic for some $w$ as follows: 
\begin{align*}
        \Oira{h=3}{0} =& 1\;, \qquad  
        \Oira{h=3}{1} =  0 \;, \qquad  
        \Oira{h=3}{2} =  -n\;,   \qquad 
        \Oira{h=3}{3} =  - \md{2}\;, \\ 
        \Oira{h=3}{4} =&  - \md{3} + \tbinom{n}{2}\;,  \qquad 
        \Oira{h=3}{5} =  - \md{4} + \md{1} \md{2}\;,  \\
        \Oira{h=3}{6} =&  - \md{5} + \md{1} \md{3} + \tbinom{\md{2}}{2} -
        \tbinom{n}{3}\;.  
\end{align*} 

Still we are 
interested in the generating function of \(\ds\{ \Oira{h=3}{w}\}_{w}\).

%% file: Combi.tex
\section{Combinatorial approach to Poisson cohomology}% 
In Poisson geometry, the Poisson cohomology group is well-known as follows.  
\begin{defn}
For each natural number $m$, let us 
consider the vector space $\ds\myC{m}{} = \Lambda ^{m} \tbdl{M}$ of all
$m$-vector fields on $M$. Then we 
define a linear map $\ds\mydq: \myC{m}{} \rightarrow \myC{m+1}{}$ by the
Schouten bracket as
$\mydq( U ) = \SbtS{\pi}{U}$.  Then $\mydq\circ\mydq =0$ follows due to the
property 
$\ds\SbtS{\pi}{\pi} =0$ and the Jacobi identity of Schouten bracket.  
Thus, we have the Poisson cohomology group
$\ds\{ U \in \Lambda ^{m} \tbdl{M} \mid \SbtS{\pi}{U} = 0\}/
\{ \SbtS{\pi}{W} \mid W \in \Lambda ^{m-1} \tbdl{M} \}$.  
\end{defn} 
Although the definition of Poisson cohomology is clear,  calculation is
not easy in general because cochain complexes are huge.  
\kmcomment{
\subsection{Schouten bracket and Poisson cohomology} \label{sub:schou}
Let us recall the Schouten bracket and the Poisson cohomology group.
Let $M$ be a $n$-dimensional smooth manifold and  $\ds\Lambda^j \tbdl{M}$ 
the space of smooth $j$-vector fields on $M$. In particular, $\ds
\Lambda^1 \tbdl{M}$ is ${\mathfrak X}(M)$, the Lie algebra of smooth
vector fields on $M$, and $\ds\Lambda^0 \tbdl{M}$ is $\ds C^{\infty}(M)$.  
$\ds A = \sum_{j=0}^ n \Lambda^j \tbdl{M}$ is the exterior algebra of
multi-vector fields on $M$. 

%%%%%%%%%%%%%%%%%%%%%%%%%%%%%%%%%%%%%%%%%%%%

For 
$\ds P \in \Lambda^{p} \tbdl{M}$ and 
$\ds Q \in \Lambda^{q} \tbdl{M}$,   
$\ds\SbtS{P}{Q} \in \Lambda^{p+q-1} \tbdl{M}$ holds, 
\begin{align*}
        \SbtS{Q}{P} & = - \parity{ (q+1)(p+1)} \SbtS{P}{Q}\quad \text{(symmetry)} \\
 0 & = \mathop{\frakS{}}_{p,q,r} \parity{(p+1)(r+1)} \SbtS{P}{\SbtS{Q}{R}}   
 \quad \text{(the Jacobi identity)}\\  
%
%\SbtS{P}{\SbtS{Q}{R}} &= 
% \SbtS{ \SbtS{P}{Q}}{R} 
%+ \parity{ (p+1)(q+1)} \SbtS{Q}{\SbtS{P}{R}}\\
\SbtS{P}{Q\wedge R} &= \SbtS{P}{Q}\wedge R + (-1)^{(p+1)q} Q \wedge
\SbtS{P}{R}\\ 
\SbtS{P \wedge Q}{R} &= P \wedge \SbtS{Q}{R}  + (-1)^{q (r+1)} 
\SbtS{P}{R} \wedge Q \\
\noalign{another expression of Jacobi identity is the next} 
\SbtS{P}{\SbtS{Q}{R}} &= 
\SbtS{ \SbtS{P}{Q}} {R} + (-1)^{(p+1)(q+1)} \SbtS{Q}{\SbtS{P}{R}} 
\\ 
\SbtS{\SbtS{P}{Q}}{R} &= 
\SbtS{P}{\SbtS{Q}{R}} + (-1)^{(q+1)(r+1)} \SbtS{\SbtS{P}{R}}{Q} 
\\\noalign{so far, formulae are valid by using $\mR$-module properties of 
the Schouten bracket and the wedge product}
\SbtS{X}{Y} &= \text{Jacobi-Lie bracket of}\ X\ \text{and}\ Y \\
\SbtS{X}{f} &= \langle X, df \rangle 
\end{align*}

%%%%%%%%%%%%%%%%%%%%%%%%%%%%%%%%%%%%%%%%%%%%%%%%%%%%%%%%%%%%%%%%%%%%%%%

Let $\pi$ be a 2-vector fields on $M$. $\pi$ is a Poisson structure if and only if $\ds
\SbtS{\pi}{\pi} = 0$.  
Locally, let $\ds (x_1\ldots x_n)$ be a local coordinates and 
\[ \ds\pi = 
\frac{1}{2}\sum_{i,j} p_{ij}
\pdel_i \wedge 
\pdel_j 
\quad \text{where}\quad \pdel_i = \frac{\pdel}{\pdel x_i}, \ \text{and}\  
p_{ij} + p_{ji}=0\ .\]    
The Schouten bracket of $\pi$ and itself is calculated as follows 
\begin{align*}
2 \SbtS{\pi}{\pi} =& \sum_{i,j} \Sbt{\pi}{p_{ij}\pdel_i \wedge \pdel_j}  \\
=& 
\sum_{i,j}(  
\Sbt{\pi}{p_{ij}} \wedge \pdel_i \wedge \pdel_j  
+ p_{ij} \SbtS{\pi}{\pdel_i} \wedge \pdel_j  
- p_{ij}\pdel_i \wedge \SbtS{\pi}{\pdel_j}  
)
\\ 
4 \SbtS{\pi}{\pi} =& 
\sum_{k\ell} \sum_{i,j}(  
\Sbt{ p_{k\ell} \pdel_k \wedge \pdel_{\ell} }{p_{ij}} \wedge \pdel_i \wedge \pdel_j  
+ p_{ij} \SbtS{ 
        p_{k\ell} \pdel_k \wedge \pdel_{\ell}
}{\pdel_i} \wedge \pdel_j  \\ & 
- p_{ij}\pdel_i \wedge \SbtS{ 
        p_{k\ell} \pdel_k \wedge \pdel_{\ell}
}{\pdel_j}  
)
\\
=& 
\sum_{k\ell} \sum_{i,j}(  
p_{k\ell} (\pdel_{\ell} p_{ij})  \pdel_k  \wedge \pdel_i \wedge \pdel_j 
- p_{k\ell} (\pdel_{k} p_{ij})  \pdel_{\ell}\wedge \pdel_i \wedge \pdel_j 
\\ & 
- 
p_{ij} 
    ( \pdel_i   p_{k\ell}) \pdel_k \wedge \pdel_{\ell}  
        \wedge \pdel_j  + p_{ij}(\pdel_j p_{k\ell})  
        \pdel_i \wedge 
        \pdel_k \wedge \pdel_{\ell} 
)
\\
=& 4 \sum_{ijk} \sum_{\ell} p_{i\ell} \frac{\pdel p_{jk}}{\pdel x_{\ell}}\;  \pdel_i \wedge \pdel_j
\wedge \pdel_k
\end{align*} 
The Poisson bracket of $f$ and $g$ is given by $\Bkt{f}{g} = \pi(df,dg)$ so 
$\ds\Bkt{x_i}{x_j} = p_{ij}$ and 
\begin{align*}
        \Bkt{\Bkt{x_i}{x_j}}{x_k} =& \Bkt{p_{ij}}{x_k} = 
        \frac{1}{2} \sum_{\lambda,\mu} p_{\lambda,\mu} (
        \pdel_{\lambda}p_{ij} \pdel_{\mu}x_k 
        - \pdel_{\mu }p_{ij} \pdel_{\lambda}x_k )
        = - \sum_{\lambda} p_{k\lambda} \pdel_{\lambda} p_{ij}
\end{align*}
$\pi$ is Poisson if and only if Jacobi identity for the bracket holds if and only if
\begin{equation} \mathop{\mathfrak S}_{i,j,k} 
        \sum_{\lambda} p_{k\lambda} \pdel_{\lambda} p_{ij} = 0 \label{bkt:jacobi}
\end{equation}
}%endOFkmcomment

% \subsection{$\mR$-module Poisson-like cohomology}  
\subsection{Poisson-like cohomology}  

Here, we discuss ``Poisson-like'' cohomology for a given homogeneous
``Poisson-like'' structure restricting cochain spaces to vector fields
with polynomial coefficients, and also the notion of ``weight'' to
reduce our discussion in finite dimensional vector spaces. 

\begin{defn}
Let 
$\ds\frakX{pol}$ be $\ds   \{ X\in {\mathfrak X}(\mR^n) \mid \langle
dx_j, X\rangle \ 
\text{are polynomials for each }\ j \}$, and 
let 
$\ds\frakX{\ell}$ be $\{ X\in \frakX{pol} \mid \langle d
x_j,X\rangle \  
\text{are $\ell$-homogeneous polynomials for each}\ j \}$ 
for each non-negative integer $\ell$.   
\end{defn}

We see that $\ds\frakX{pol} = \mathop{\oplus}_{\ell} \frakX{\ell}$ as
$\mR$-vector space. Thus,  the exterior 2-power of $\ds\frakX{pol}$ is
\[ \ds\Lambda^{2} \frakX{pol} =  \frakX{pol} \wedge   \frakX{pol} =
        \sum_{i\leq j} \frakX{i}\wedge\frakX{j} \quad \text{(direct
sum)} \] as $\ds\mR$-modules, for instance.   

\begin{kmRemark} \label{rmk:linear} We have    
$\ds x_1 \pdel_1 \in \frakX{1}$ and $\ds x_1{ }^2 \pdel_1 \in
\frakX{2}$.  $\ds (x_1 \pdel_1)\wedge (x_1{ }^2 \pdel_1) \ne 0$ as
$\ds\mR$-modules but as $\ds C^{\infty}(\mR^n)$-modules we see that $\ds
(x_1 \pdel_1)\wedge (x_1{ }^2 \pdel_1) = x_1{ }^3\;  \pdel_1 \wedge
\pdel_1  = 0$.  
\end{kmRemark} 
For each natural number $m$, we consider $\ds\Lambda^m \frakX{pol}$ and
have natural $\mR$-module decomposition 
\[ \ds\Lambda^m \frakX{pol} = 
\sum_{m = \kk{0}+ \kk{1} + \cdots } 
\Lambda ^{k_0} \frakX{0} \otimes 
\Lambda ^{k_1} \frakX{1} \otimes \cdots 
\Lambda ^{k_{\ell}} \frakX{\ell}\;. \]  
Since $\ds\dim\frakX{j} = \tbinom{n-1+j}{j} n$, we have restrictions  
\begin{equation}
0\leq k_j \leq \dim \frakX{j} = \tbinom{n-1+j}{j} n\;.
\label{dim:cond}
\end{equation}

\begin{defn} 
Let us fix a non-negative integer $h$ (which plays a role of the
homogeneity of homogeneous Poisson-like 2-vector later).  We define the
weight $w$ of a non-zero element of 
$\ds\Lambda ^{k_0} \frakX{0} \otimes \Lambda ^{k_1} \frakX{1} \otimes
\cdots \Lambda ^{k_{\ell}} \frakX{\ell} $  to be 
\begin{equation}  
k_0 \; (0+1-h)  + k_1\; (1+1-h) 
+ \cdots + k_{\ell}\; (\ell+1-h)\;. \label{yd:area}\end{equation}
\end{defn}

\begin{defn}
For each $m$ and $w$, define a vector subspace 
\[ \ds\myC{m}{w} := 
\sum_{\text{"our cond"}} 
\Lambda ^{k_0} \frakX{0} \otimes 
\Lambda ^{k_1} \frakX{1} \otimes \cdots 
\Lambda ^{k_{\ell}} \frakX{\ell}\;. \] 
The "our cond" are (\ref{dim:cond}), (\ref{yd:area}) and 
\begin{equation}
         k_0+ k_1 + \cdots + k_{\ell} = m\;. \label{yd:deg} 
\end{equation}
\end{defn}

Now, we restrict the Schouten bracket of 
$\ds\oplus \Lambda^{\bullet} \tbdl{\mR^n}$ to 
$\ds\oplus \Lambda^{\bullet} \frakX{pol}$ and have a new bracket
$\SbtR{\cdot}{\cdot}$ and we call the $\mR$-Schouten bracket.    
\begin{defn} 
The $\mR$-Schouten bracket is characterized as follows (almost the same
in the subsection \ref{sub:schou}):   
For 
$\ds P \in \Lambda^{p} \frakX{pol}$ and 
$\ds Q \in \Lambda^{q} \frakX{pol}$,   
$\ds\SbtR{P}{Q} \in \Lambda^{p+q-1} \frakX{pol}$ holds, and  
\begin{align*}
        \SbtR{Q}{P} & = - \parity{ (q+1)(p+1)} \SbtR{P}{Q}\quad
		\text{(symmetry)}\;, \\
 0 & = \mathop{\frakS{}}_{p,q,r} \parity{(p+1)(r+1)} \SbtR{P}{\SbtR{Q}{R}}   
 \quad \text{(the Jacobi identity)}\;,\\  
%
%\SbtR{P}{\SbtR{Q}{R}} &= 
% \SbtR{ \SbtR{P}{Q}}{R} 
%+ \parity{ (p+1)(q+1)} \SbtR{Q}{\SbtR{P}{R}}\\
\SbtR{P}{Q\wedge R} &= \SbtR{P}{Q}\wedge R + (-1)^{(p+1)q} Q \wedge
\SbtR{P}{R}\;,\\ 
\SbtR{P \wedge Q}{R} &= P \wedge \SbtR{Q}{R}  + (-1)^{q (r+1)} 
\SbtR{P}{R} \wedge Q\;, \\
\noalign{another expression of Jacobi identity is the next} 
\SbtR{P}{\SbtR{Q}{R}} &= 
\SbtR{ \SbtR{P}{Q}} {R} + (-1)^{(p+1)(q+1)} \SbtR{Q}{\SbtR{P}{R}} 
\;, 
\\ 
\SbtR{\SbtR{P}{Q}}{R} &= 
\SbtR{P}{\SbtR{Q}{R}} + (-1)^{(q+1)(r+1)} \SbtR{\SbtR{P}{R}}{Q} 
\;,
\\
% \noalign{so far, formulae are valid by using 
% $\mR$-module properties of the Schouten bracket and the wedge product}
\SbtR{X}{Y} &= \text{Jacobi-Lie bracket of}\ X\ \text{and}\ Y 
%\SbtR{X}{f} &= \langle X, df \rangle 
\;. 
\end{align*} 
\end{defn}

The $\mR$-Schouten bracket also has an explicit expression given by 
        \begin{equation}
                \SbtR{u_1\wedge \cdots \wedge u_p}
                {v_1\wedge \cdots \wedge v_q}
                = \sum_{i,j}(-1)^{i+j} \SbtR{u_i}{v_j}\wedge 
                (u_1\wedge \cdots  \widehat{u_i} \cdots \wedge u_p) 
                \wedge 
                ( v_1\wedge \cdots \widehat{v_j}\cdots \wedge v_q)
                \label{eqn:SbtR}
        \end{equation} where
        $\ds u_i, v_j \in \frakX{pol}$ and  
        $\ds\widehat{u_i}$ means omitting $\ds u_i$. 

We have the same property that the ordinary Schouten bracket has.  
\begin{kmProp} 
Let $\ds\pi\in\Lambda^{2}\frakX{pol}$ and $\ds P\in
\Lambda^{p}\frakX{pol}$. Then 
\begin{equation}
2\SbtR{\pi}{\SbtR{\pi}{P}} + \SbtR{P}{\SbtR{\pi}{\pi}} =0 
\end{equation} 
holds and so if $\ds\SbtR{\pi}{\pi} = 0$ then $\ds
\SbtR{\pi}{\SbtR{\pi}{\cdot}} = 0$ holds on $\ds
\Lambda^{\bullet}\frakX{pol}$.   
\end{kmProp}
\begin{defn} \label{defn:PoissonLike} 
We call a 2-vector $\ds\pi\in\Lambda^{2}\frakX{pol}$ is Poisson-like if
$\pi$ satisfies $\ds\SbtR{\pi}{\pi} = 0$.   

A Poisson-like 2-vector $\pi$ is $h$-homogeneous if $\ds\pi \in
\mathop{\oplus}_{h = i+j,\ i\leq j} \frakX{i} \wedge \frakX{j}$.  
\end{defn}

\begin{kmProp} 
Let $\pi$ be a $h$-homogeneous Poisson-like 2-vector as in Definition
\ref{defn:PoissonLike}.  We see that $\ds\SbtR{\pi}{\myC{m}{w}} \subset
\myC{m+1}{w} $ and we get a sequence of cochain complexes: $\ds\mydq :
\myC{m}{w} \rightarrow \myC{m+1}{w}$, and the cohomology groups.  We
call them the Poisson-like cohomology groups of homogeneous Poisson-like
2-vector on $\mR^n$.  
\end{kmProp}

\textbf{Proof:} 
Since $\ds\pi \in \mathop{\oplus}_{h = i+j,\ i\leq j} \frakX{i} \wedge
\frakX{j}$,   
we have  
\[\SbtR{\pi}{\frakX{\ell}} \subset \mathop{\oplus}_{h=i+j,\ i\leq j} 
( \frakX{i} \wedge \frakX{j+\ell-1} \oplus  \frakX{i+\ell-1} \wedge
\frakX{j} ) \]  
and the weight of $\ds\frakX{i} \wedge \frakX{j+\ell -1}$ is  $(i+1-h)
+ (j+\ell-1+1-h) = \ell+1-h$ and the weight of
$\ds\frakX{i+\ell-1}\wedge \frakX{j} $ is $(i+\ell-1 +1-h) + (j+1-h) =
\ell + 1 -h$.  
Thus, after applying $\mydq$, the degree changes to $m+1$ but
the weight is invariant.  \kmqed

\medskip

\begin{kmRemark}\label{rmk:Phi}
We have two kinds of Schouten brackets, $\ds\SbtS{\cdot}{\cdot}$ and
$\ds\SbtR{\cdot}{\cdot}$.  
Let $\Phi: \oplus \Lambda^{\bullet}_{\mR}\frakX{pol} \rightarrow 
\oplus \Lambda^{\bullet}\frakX{pol}$ be the natural map relaxing the
$\mR$-linearity of the wedge product to the ordinary tensorial product. 
As commented in Remark 
 \ref{rmk:linear},  
$\ds (x_1 \pdel_1)\wedge (x_1{ }^2 \pdel_1) \ne 0$ 
as $\ds\mR$-modules and  
$\ds\Phi\left((x_1 \pdel_1)\wedge (x_1{ }^2 \pdel_1)\right) 
= x_1{ }^3\;  \pdel_1 \wedge \pdel_1  = 0$. From the definition of
$\ds\SbtR{\cdot}{\cdot}$, 
we see quickly that $\ds\SbtS{\Phi(P)}{\Phi(Q)} = \Phi( \SbtR{P}{Q})$ for each  
$P, Q \in \oplus \Lambda^{\bullet}_{\mR}\frakX{pol}$.  
Thus, for any Poisson-like 2-vector $\pi$, $\ds\Phi(\pi)$  is a Poisson
2-tensor.  But, the converse is not true. Namely, 
let $\overline{\pi}$ be an ordinary Poisson structure. There is no guarantee
that 
each inverse element $\mathbf{u}\in \Phi^{-1}(\overline{\pi})$ is 
Poisson-like.  
Let $\ds\overline{\pi} = x_{3} \pdl{1} \wedge \pdl{2}  
-2 x_{1} \pdl{1} \wedge \pdl{3} + 
2 x_{2} \pdl{2} \wedge \pdl{3} 
$, which is a Poisson structure due to $\fraksl(2)$   
we have already used.   
Let $\ds\phi = \pdl{1} \wedge ( x_{3} \pdl{2} ) 
-2 \pdl{1} \wedge (x_{1} \pdl{3}) + 
2 \pdl{2} \wedge (x_{2} \pdl{3}) 
$  so that $\ds\Phi(\phi) = \overline{\pi}$. 
But, $\ds\SbtR{\phi}{\phi}/4$ is calculated to be non-zero as follows. 
\[  
\pdl{1} \wedge 
\pdl{3} \wedge
(x_{3} \pdl{2})
-2 
\pdl{2} \wedge
\pdl{3} \wedge 
( x_{2} \pdl{3})
-2
\pdl{1} \wedge 
\pdl{3} \wedge 
( x_{1} \pdl{3}) 
+
\pdl{1} \wedge 
\pdl{2} \wedge
(x_{3} \pdl{3} ) 
-4
\pdl{1}\wedge 
\pdl{2} \wedge 
(x_{2} \pdl{2} ) \;. 
\] 
\end{kmRemark}

\bigskip

We rewrite (\ref{yd:area}) 
and  (\ref{yd:deg})
as follows:
\begin{align}
& 0\; k_0  + 1\; k_1
+ \cdots + \ell\; k_{\ell} = w + (h-1) \; m \;,\label{yd:area:second} \\
& k_1 + \cdots + k_{\ell} = m - k_0 \;. \label{yd:height} 
\end{align}
The first equation means the total area and  
the second equation means the height of the Young diagram 
$\ds (k_j \mid j=1\ldots \ell)$.   

Hereafter, we will show some concrete examples and difference between this
cohomology and the cohomologies in previous sections.

\subsubsection{case h=0} 
(\ref{yd:area}) and (\ref{yd:deg}) say that if the degree  \(m=0\) then
the weight \(w=0\), in other words, if the weight \(w\ne 0\) then
\(\ds\myC{0}{w} = \emptyset\).  We denote the left-hand-sides of
(\ref{yd:area:second}) and (\ref{yd:height}) by $A$ and $H$
respectively, then using \(h=0\) we have \( A=w-m, H=m-k_0\), and so
\(\ds\myC{m}{w} = \sum_{H} \Lambda^{m-H} \frakX{0} \otimes
\nabla(w-m,H)\), where we suppose \( \nabla(0,0) \) is the singleton of
the trivial Young diagram but \( \nabla(A,0) \) with \(A>0\) is the
empty set and do not sum up the terms containing these direct summands.
We may regard $H$ as a parameter with the restrictions \( H \leq w/2,\
m-n \leq H \leq m \) because by  adding \( m-H \geq 0, w-m\geq H \) we
see \( H\leq w/2\).  
Thus, when $w=1$ we have 
\[\ds\myC{m}{1} = \sum_{H \leq 1/2 } 
 \Lambda^{m-H} \frakX{0} \otimes \nabla(1-m,H) 
 = \Lambda^{m} \frakX{0} \otimes \nabla(1-m,0) 
 = \begin{cases} 
         \frakX{0} & \text{if}\ m = 1\;, \\[-2mm]
 \emptyset  & \text{otherwise.}\end{cases} 
\] 
The Euler characteristic for \(w=1\), 
we may denote it by \(\ds \Oira{h=0}{w=1}= -n\).

\begin{align*}\myC{m}{2} &= \sum_{H \leq 2/2 } 
 \Lambda^{m-H} \frakX{0} \otimes \nabla(2-m,H) 
 = \Lambda^{m-1} \frakX{0} \otimes \nabla(2-m,1) 
 + \Lambda^{m} \frakX{0} \otimes \nabla(2-m,0) 
 \\&
 = \begin{cases} 
         \frakX{1} & \text{if}\ m = 1\;, \\[-2mm]
         \Lambda^{2}\frakX{0} & \text{if}\ m = 2\;, \\[-2mm]
 \emptyset  & \text{otherwise. } \end{cases}
         \end{align*} 
         \(\ds\Oira{h=0}{w=2}= -n(n+1)/2 \). 

\[\ds\myC{m}{3} = 
  \Lambda^{m-1} \frakX{0} \otimes \nabla(3-m,1) 
 + \Lambda^{m} \frakX{0} \otimes \nabla(3-m,0) 
 =  \begin{cases} 
         \frakX{2} & \text{if}\ m = 1\; ,\\[-2mm]
         \frakX{0}\otimes\frakX{1} & \text{if}\ m = 2\; ,\\[-2mm]
         \Lambda^{3}\frakX{0} & \text{if}\ m = 3\;, \\[-2mm]
 \emptyset  & \text{otherwise.} \end{cases}
\] 
\(\ds\Oira{h=0}{w=3}= (n-1)n(n+1)/3 \).

\begin{align*} \myC{m}{4}  &=  
  \Lambda^{m} \frakX{0} \otimes \nabla(4-m,0) 
 + \Lambda^{m-1} \frakX{0} \otimes \nabla(4-m,1) 
 + \Lambda^{m-2} \frakX{0} \otimes \nabla(4-m,2) 
 \\
 & 
 = \begin{cases} 
         \frakX{3} & \text{if}\ m = 1\;, \\[-2mm]
         \frakX{0}\otimes\frakX{2} + \Lambda^{2} \frakX{1} &
         \text{if}\ m = 2\;, \\[-2mm]
         \Lambda^{2} \frakX{0} \otimes \frakX{1}  & \text{if}\ m =
         3\;, \\[-2mm]
         \Lambda^{4}\frakX{0} & \text{if}\ m = 4\;, \\[-2mm]
 \emptyset  & \text{otherwise.} \end{cases}
 \end{align*} 
 \(\ds\Oira{h=0}{w=4}= - (n-3) (n-1)n(n+2)/8 \).

Since the 0-homogeneous Poisson-like 2-vectors are of the form  
$\ds\sum_{i=1}^{r} \pdel_{2i-1} \wedge \pdel_{2i} $
after a suitable change of coordinates,  we take $\pi = \pdel_{1} \wedge
\pdel_{2}$ when $n=3$.  Then \[\ds\SbtR{\pi}{ w^A \pdel_j} = 
A_{1} ( w^{A-E_1}\pdel_j) \wedge \pdel_2 - 
A_{2} ( w^{A-E_2}\pdel_j) \wedge \pdel_1 \]
where $A\in\mathfrak{M}[k]$ and $j,k=1,\ldots,n$. We get Betti numbers as
follows.  
\begin{center}
\setlength{\extrarowheight}{-3pt}
\begin{tabular}{|c|*{3}{c}|} \hline 
&  $\ds\myC{1}{2} $& $ \rightarrow $ &
    $\ds\myC{2}{2} $  \\\hline 
$\dim$  &    9 && 3   \\\hline
$\ker \dim $  &    6 && 3   \\\hline
Betti &     6 &&  0    \\ \hline 
\end{tabular} 
\hfil 
\begin{tabular}{|c|*{5}{c}|} \hline 
&  $\ds\myC{1}{3} $
& $ \rightarrow $ & $\ds\myC{2}{3} $ 
& $ \rightarrow $ & $\ds\myC{3}{3} $ 
    \\\hline 
$\dim$  &    18 && 27 && 1   \\\hline
$\ker \dim $  &  3 && 26  && 1  \\\hline
Betti &     3 &&  11  && 0   \\ \hline 
\end{tabular} 
\hfil 
\begin{tabular}{|c|*{5}{c}|} \hline 
&  $\ds\myC{1}{4} $
& $ \rightarrow $ & $\ds\myC{2}{4} $ 
& $ \rightarrow $ & $\ds\myC{3}{4} $ 
    \\\hline 
$\dim$  & 30 && 90 && 27   \\\hline
$\ker \dim $  &  3 && 63  && 27  \\\hline
Betti &     3 &&  36  && 0   \\ \hline 
\end{tabular} 
\end{center} 

\subsubsection{case h=1} 
In this case, (\ref{yd:area:second}) says the weight $w$ is just the
total area of Young diagram and (\ref{yd:height}) says its height is $m
-k_0$. Thus, $\ds m-k_0 \leq w$ and $\ds m \leq w + n$ from
(\ref{dim:cond}). 

If the weight $w=0$, then $\ds k_j=0$ ($j>0$) and $\ds k_0=m$ and so 
\[ \ds\myC{m}{0} = \Lambda^{m} \frakX{0}\quad \text{for}\quad
m=0,\ldots,n.\]  

If the weight $w=1$, we know $\ds\nabla(1, m-k_0) = \begin{cases}
        \emptyset & \text{if} \quad m-k_0 \le 0\;, \\[-2mm]
        \{\T{1}\} & \text{if}\quad
m-k_0=1\;, \end{cases}$ 
\\
we see $\ds k_0 = m-1$, $\ds k_1=1$ and 
\(\ds\myC{m}{1} = \Lambda^{m-1} \frakX{0} \otimes \frakX{1}
\quad\text{for}\quad m=1,\ldots,n+1\).    

In the same way, 
if when the weight $w=2$, we know $\ds\nabla(2, m-k_0) = \begin{cases} 
        \{\T{2}\} & \text{if}\quad   m-k_0=2\;, \\[-2mm]
        \{\T{1}^ 2\} & \text{if}\quad   m-k_0=1\;, \\[-2mm] 
        \emptyset & \text{otherwise,}
\end{cases}$ 
we see ($\ds k_0 = m-1$, $\ds k_2=1$ ) or  
($\ds k_0 = m-2$, $\ds k_1=2$), and so    
\[ \ds\myC{m}{2} = \Lambda^{m-1} \frakX{0} \otimes \frakX{2}
\oplus 
\Lambda^{m-2} \frakX{0} \otimes \Lambda^{2} \frakX{1}
\ . \]   
If the weight $w=3$, we know $\ds\nabla(3, m-k_0) = \begin{cases} 
        \{\T{3}\} & \text{if}\quad   m-k_0=3\;, \\[-2mm]
        \{\T{2}\cdot \T{1}\} & \text{if}\quad   m-k_0=2\;, \\[-2mm]
        \{\T{1}^3 \} & \text{if}\quad   m-k_0=1\;, \\[-2mm] 
        \emptyset & \text{otherwise,}
\end{cases}$ 
\\
we see ($\ds k_0 = m-3$, $\ds k_1=3$ ),   
($\ds k_0 = m-2$, $\ds k_1=1$, $\ds k_2=1$) or 
($\ds k_0 = m-1$, $\ds k_3=1$), thus we have  
\[ \ds\myC{m}{3} = \Lambda^{m-3} \frakX{0} \otimes \Lambda^{3}\frakX{1}
\oplus \Lambda^{m-2} \frakX{0} \otimes  \frakX{1} \otimes \frakX{2}
\oplus \Lambda^{m-1} \frakX{0} \otimes  \Lambda^{3}\frakX{1} 
\ . \]

Assume $n=3$ now. Then we get     
\begin{align*}
        \myC{0}{0} &= \mR , \quad  
\myC{1}{0} = \frakX{0}, \quad  
\myC{2}{0} = \Lambda^{2} \frakX{0}, \quad    
\myC{3}{0} = \Lambda^{3} \frakX{0}\;,\\  
\myC{1}{1} & = \frakX{1}.\quad   
\myC{2}{1} = \frakX{0} \otimes \frakX{1}. \quad   
\myC{3}{1} = \Lambda^2 \frakX{0} \otimes
\frakX{1}.\quad   \myC{4}{1} = \Lambda^3 \frakX{0}
\otimes \frakX{1}\;, \\ 
\myC{1}{2} &= \frakX{2},\quad   
\myC{2}{2} = \frakX{0}\otimes \frakX{2}+ \Lambda^2 \frakX{1},\quad   
\myC{3}{2} = \Lambda^2 \frakX{0}\otimes \frakX{2}+ 
\frakX{0} \otimes \Lambda^2 \frakX{1}, \quad  
\\ 
\myC{4}{2} &= \Lambda^3 \frakX{0}\otimes \frakX{2}+ 
\Lambda^2 \frakX{0} \otimes  
\Lambda^2 \frakX{1}, \quad   
\myC{5}{2} = \Lambda^3 \frakX{0}\otimes
\Lambda^2 \frakX{1}\;, \\ 
\myC{1}{3} & = \frakX{3},  \quad 
\myC{2}{3}  = \frakX{0} \otimes \frakX{3} + \frakX{1}\otimes \frakX{2}
\;, 
\\
\myC{3}{3} & = 
\Lambda^2 \frakX{0} \otimes \frakX{3} + 
\frakX{0} \otimes \frakX{1}\otimes \frakX{2} 
+ \Lambda^3 \frakX{1} \;,
\\
\myC{4}{3} & = 
\Lambda^3 \frakX{0} \otimes \frakX{3} + 
\Lambda^2 \frakX{0} \otimes \frakX{1}\otimes \frakX{2} 
+ \frakX{0} \otimes \Lambda^3 \frakX{1} \;,
\\
\myC{5}{3} & = 
\Lambda^3 \frakX{0} \otimes \frakX{1}\otimes \frakX{2} 
+ \Lambda^2 \frakX{0} \otimes \Lambda^3 \frakX{1},  
\quad 
\myC{6}{3} = 
\Lambda^3 \frakX{0} \otimes \Lambda^3 \frakX{1} \;. 
\end{align*}

\begin{kmRemark} 
From the above several examples, we expect the Euler characteristic is 0
when $h=1$ likewise as the section \ref{sec:Euler}.  
\end{kmRemark}

In order to find concrete 1-homogeneous Poisson-like 2-vectors on
$\ds\mR^3$, we prepare a candidate in general form, say  \[\ds u :=
\sum_{i,j,k}^{3} c_{i, E_j, k} \pdl{i} \wedge ( w^{E_j} \pdl{k}) \] with
the condition $\ds\SbtR{u}{u} = 0$.  Then we have a system of
2-homogeneous polynomials of $\ds c_{i, E_j, k}$.  

It seems hard to know the whole solutions, but  one of many solutions is
$\pi = (\pdl{1} - \pdl{3}) \wedge ( x_{1} \pdl{3} + x_{3} \pdl{3})$.  We
show the tables of Betti numbers of Poisson-like cohomologies for weight
from 0 to 3.  
\begin{center}
\tabcolsep=3pt 
\setlength{\extrarowheight}{-3pt}
\begin{tabular}{|c|*{5}{c}|} \hline 
        wt=0 
&  $\ds\myC{1}{0} $
& $ \rightarrow $ & $\ds\myC{2}{0} $ 
& $ \rightarrow $ & $\ds\myC{3}{0} $ 
    \\\hline 
$\dim$    & 3 && 3 && 1  \\\hline
$\ker\dim$& 2 && 2 && 1  \\\hline
Betti &     2 && 1 && 0  \\ \hline 
\end{tabular} 
\hfil 
\begin{tabular}{|c|*{7}{c}|} \hline 
        wt = 1
&  $\ds\myC{1}{1} $
& $ \rightarrow $ & $\ds\myC{2}{1} $ 
& $ \rightarrow $ & $\ds\myC{3}{1} $ 
& $ \rightarrow $ & $\ds\myC{4}{1} $ 
    \\\hline 
$\dim$    & 9 && 27 && 27 && 9 \\\hline
$\ker\dim$& 3 && 15 && 20 && 9 \\\hline
Betti &     3 && 9  && 8  && 2  \\ \hline 
\end{tabular} 
\hfil
\begin{tabular}{|c|*{9}{c}|} \hline 
        wt=2
&  $\ds\myC{1}{2} $
& $ \rightarrow $ & $\ds\myC{2}{2} $ 
& $ \rightarrow $ & $\ds\myC{3}{2} $ 
& $ \rightarrow $ & $\ds\myC{4}{2} $ 
& $ \rightarrow $ & $\ds\myC{5}{2} $ 
    \\\hline 
$\dim$    & 18 && 90 && 162 && 126 && 36 \\\hline
$\ker\dim$& 4 && 35 && 93 && 98 && 36 \\\hline
Betti &     4 && 21 && 38  && 29 && 8 \\ \hline 
\end{tabular} 

\begin{tabular}{|c|*{11}{c}|} \hline 
        wt=3
&  $\ds\myC{1}{3} $
& $ \rightarrow $ & $\ds\myC{2}{3} $ 
& $ \rightarrow $ & $\ds\myC{3}{3} $ 
& $ \rightarrow $ & $\ds\myC{4}{3} $ 
& $ \rightarrow $ & $\ds\myC{5}{3} $ 
& $ \rightarrow $ & $\ds\myC{6}{3} $ 
    \\\hline 
$\dim$    & 30 && 252 && 660 && 768 && 414 && 84  \\\hline
$\ker\dim$& 3 && 74 && 315 && 504 && 344 && 84 \\\hline
Betti &     3 && 47 && 137 && 159 && 80 && 14 \\ \hline 
\end{tabular} 
\end{center}

\subsubsection{case h=2} 

In this subsection we deal with homogeneous Poisson-like 2-vectors with
$h=2$.      Comparing (\ref{yd:area:second}) and (\ref{yd:height}), we
see that $\ds -w \leq k_0 $ and so $\ds w \geq -n$.  If $\ell >1$ in
(\ref{yd:area:second}) or (\ref{yd:height}),  then
$(\ref{yd:area:second}) - 2( \ref{yd:height})$ implies $\ds m \leq w +
2k_0+ k_1 \leq w + 2 n + n^2$.  If $\ell=1$, then $w = -k_0$ and $m =
k_0+ k_1$.  If $\ell=0$, $m=k_0$ and $ w = - k_0$.  Finding
$\ds\myC{m}{w}$ is equivalent to finding the Young diagrams of height
$\ds m-k_0$ and the area $w+m$ for each $\ds k_0$. 

Since $\ds w\geq -n$, we put $\ds w = -n +j$ with non-negative integer
$j$.  Then $\ds k_0 \geq n-j$.  When $j=0$, i.e., the weight $\ds w =
-n$, we have $\ds k_0=n$ and $\ds\nabla(-n+m,m-n) = \{ \T{m-n} \}$, this
says $\ds k_1 = m-n$ and $\ds k_{\ell} = 0$ for $\ell >1$. Thus, we get
\begin{equation} \ds\myC{m}{-n} = \Lambda^{n}\frakX{0} \otimes
\Lambda^{m-n}\frakX{1}\ . \end{equation} 

When $j=1$, i.e., the weight $\ds w = -n + 1$, we see that $\ds k_0=n-1$
or $\ds k_0=n$.   If $\ds k_0=n-1$, $\ds\nabla(-n+1+m,m-n+1) = \{
\T{m-n+1} \}$, this says $\ds k_1 = m-n+1 $ and $\ds k_{\ell} = 0$ for
$\ell >1$.  If $\ds k_0=n$,  then $\ds\nabla(-n+1+m,m-n) =
\T{m-n}\cdot\{ \nabla(1,1) \}$, this says $\ds k_1 = m-n-1$,  $\ds k_2 =
1$,  and $\ds k_{\ell} = 0$ for $\ell >2$.  
Thus, 
\begin{equation}\ds\myC{m}{1-n} = 
\Lambda^{n-1}\frakX{0} \otimes \Lambda^{m-n+1}\frakX{1} \oplus 
\Lambda^{n}\frakX{0} \otimes \Lambda^{m-n-1}\frakX{1} \otimes \frakX{2}\ .
\end{equation}
When $j=2$, i.e., $w=-n+2$, possibilities of $\ds k_0$ are 
$\ds k_0=n-2$, $\ds k_0=n-1$ or $\ds k_0=n$.   
If $\ds k_0=n-2$, $\ds\nabla(-n+m+2,m-n+2) = \{\T{m-n+2}\}$, we see 
$\ds k_1= m-n+2, \; k_{\ell} =0\; (\ell>1)$. 
If $\ds k_0=n-1$, $\ds\nabla(-n+m+2,m-n+1) = \T{m-n+1}\cdot \nabla(1,1)$, i.e., 
$\ds k_1= m-n, \; k_2 = 1,\; k_{\ell} = 0\; (\ell>2)$. 
If $\ds k_0=n$, $\ds\nabla(-n+m+2,m-n) = \T{m-n}\cdot ( \nabla(2,1) +
\nabla(2,2))
= \T{m-n}\cdot \T{1}^ 2 + \T{m-n}\cdot \T{2}
$, i.e., 
$\ds k_1= m-n-1, \; k_{3} = 1,\; k_{\ell} = 0\; (\ell \ne 1,3)$ or  
$\ds k_1= m-n-2, \; k_{2} = 2,\; k_{\ell} = 0\; (\ell > 2)$.     
Combining those, we have 
\begin{align}\ds\myC{m}{2-n} = & 
\Lambda^{n-2}\frakX{0} \otimes \Lambda^{m-n+2}\frakX{1} \oplus 
\Lambda^{n-1}\frakX{0} \otimes \Lambda^{m-n}\frakX{1} \otimes \frakX{2}
\notag \\
& \oplus  
\Lambda^{n}\frakX{0} \otimes \Lambda^{m-n-1}\frakX{1} \otimes \frakX{3}
\oplus 
\Lambda^{n}\frakX{0} \otimes \Lambda^{m-n-2}\frakX{1} \otimes \Lambda^{2}\frakX{2} 
\ .
\end{align}
By the same discussion for $j=3,4$, we get 
\begin{align*}\ds\myC{m}{3-n} = & 
\Lambda^{n-3}\frakX{0} \otimes \Lambda^{m-n+3}\frakX{1} 
\oplus \Lambda^{n-2}\frakX{0} \otimes \Lambda^{m-n+1}\frakX{1}\otimes\frakX{2}
\\&
\oplus 
\Lambda^{n-1}\frakX{0} \otimes (
\Lambda^{m-n}\frakX{1} \otimes \frakX{3} \oplus 
\Lambda^{m-n-1}\frakX{1} \otimes \Lambda^{2}\frakX{2}) 
\\&
\oplus 
\Lambda^{n}\frakX{0} \otimes (
\Lambda^{m-n-1}\frakX{1} \otimes \frakX{4} \oplus 
\Lambda^{m-n-2}\frakX{1} \otimes \frakX{2} \otimes \frakX{3} \oplus 
\Lambda^{m-n-3}\frakX{1} \otimes \Lambda^{3} \frakX{2} 
)\;, \\ 
\myC{m}{4-n} = & 
\Lambda^{n-4}\frakX{0} \otimes \Lambda^{m-n+4}\frakX{1} 
\oplus \Lambda^{n-3}\frakX{0} \otimes \Lambda^{m-n+2}\frakX{1} \otimes \frakX{2}
\\&
\oplus 
\Lambda^{n-2}\frakX{0} \otimes (
\Lambda^{m-n+1}\frakX{1} \otimes \frakX{3} \oplus 
\Lambda^{m-n}\frakX{1} \otimes \Lambda^{2}\frakX{2}) 
\\&
\oplus 
\Lambda^{n-1}\frakX{0} \otimes (
\Lambda^{m-n}\frakX{1} \otimes \frakX{4} \oplus 
\Lambda^{m-n-1}\frakX{1} \otimes \frakX{2} \otimes \frakX{3} \oplus 
\Lambda^{m-n-2}\frakX{1} \otimes \Lambda^{3} \frakX{2} 
) \\& 
\oplus 
\Lambda^{n}\frakX{0} \otimes (
\Lambda^{m-n-1}\frakX{1} \otimes \frakX{5} \oplus 
\Lambda^{m-n-2}\frakX{1} \otimes \frakX{2} \otimes \frakX{4} 
\\ &\qquad\qquad\quad   \oplus  
\Lambda^{m-n-2}\frakX{1} \otimes \Lambda^{2} \frakX{3} 
\oplus \Lambda^{m-n-3}\frakX{1} \otimes \Lambda^{2} \frakX{2} \otimes
\frakX{3} 
\oplus \Lambda^{m-n-4}\frakX{1} \otimes \Lambda^{4} \frakX{2} 
) \;. 
\end{align*} 
Now assume $n=3$. Then we have  
\begin{align*}
\myC{m}{-3} = & \Lambda^3 \frakX{0} \otimes \Lambda^{m-3} \frakX{1}  
		\;, \\
\myC{m}{-2} = & \Lambda^2 \frakX{0} \otimes \Lambda^{m-2} \frakX{1} + 
        \Lambda^3 \frakX{0} \otimes \Lambda^{m-4} \frakX{1} \otimes \frakX{2}
		\;, \\
\myC{m}{-1} = &  
\frakX{0} \otimes \Lambda^{m-1} \frakX{1} + 
\Lambda^2 \frakX{0} \otimes \Lambda^{m-3} \frakX{1} \otimes \frakX{2} + 
\Lambda^3 \frakX{0} \otimes \Lambda^{m-4} \frakX{1} \otimes \frakX{3} + 
\Lambda^3 \frakX{0} \otimes \Lambda^{m-5} \frakX{1} \otimes \Lambda^2
\frakX{2} \;, \\ 
\myC{m}{0} =& \Lambda^m \frakX{1} + 
        \frakX{0} \otimes \Lambda^{m-2} \frakX{1} \otimes \frakX{2} +
        \Lambda^2 \frakX{0} \otimes \Lambda^{m-3} \frakX{1} \otimes \frakX{3}
        + \Lambda^2 \frakX{0} \otimes \Lambda^{m-4} \frakX{1} \otimes
        \Lambda^2 \frakX{2} 
        \\&
        + \Lambda^3 \frakX{0} \otimes \Lambda^{m-4} \frakX{1} \otimes
         \frakX{4} 
        + \Lambda^3 \frakX{0} \otimes \Lambda^{m-5} \frakX{1} \otimes
         \frakX{2} \otimes \frakX{3} 
        + \Lambda^3 \frakX{0} \otimes \Lambda^{m-6} \frakX{1} \otimes
        \Lambda^3 \frakX{2} 
		\;,
\\ 
        \myC{m}{1} =& 
        \Lambda^{m-1} \frakX{1} \otimes \frakX{2}  
        +  \frakX{0} \otimes ( \Lambda^{m-2} \frakX{1} \otimes \frakX{3}  
        +  \Lambda^{m-3} \frakX{1} \otimes \Lambda^{2}\frakX{2} ) 
        \\ &
        +  \Lambda^{2} \frakX{0} \otimes ( \Lambda^{m-3} \frakX{1} \otimes \frakX{4}  
        +  \Lambda^{m-4} \frakX{1} \otimes \frakX{2} \otimes \frakX{3} 
        +  \Lambda^{m-5} \frakX{1} \otimes \Lambda^{3}\frakX{2} 
        ) \\
        & 
        +  \Lambda^{3} \frakX{0} \otimes ( \Lambda^{m-4} \frakX{1} \otimes \frakX{5}  
        +  \Lambda^{m-5} \frakX{1} \otimes \frakX{2} \otimes \frakX{4} 
        +  \Lambda^{m-5} \frakX{1} \otimes \Lambda^{2}\frakX{3} 
        \\&\qquad \qquad 
        +  \Lambda^{m-6} \frakX{1} \otimes \Lambda^{2}\frakX{2}\otimes
        \frakX{3}  
        +  \Lambda^{m-7} \frakX{1} \otimes \Lambda^{4}\frakX{2}
        ) \;. 
\end{align*}

\begin{kmRemark} 
In homogeneous Poisson case, we have Theorem \ref{thm:euler:num} which
says the Euler characteristic of $1$-homogeneous Poisson structure is
always zero. On the other hand, we have a concrete example of
$2$-homogeneous Poisson structure which Euler characteristic is not
zero. Contrarily, in the case of homogeneous Poisson-like cohomology
groups we expect all the Euler characteristic may be zero by looking at
several concrete cochain complexes.  
\end{kmRemark}

We take the following one as a 2-homogeneous Poisson-like 2-vector 
\begin{align}\pi = & 
- \pdl{3} \wedge ( x_2 x_3 \pdl{1})+ 
(x_2 \pdl{1} ) \wedge ( x_2 \pdl{2} )
+( x_3 \pdl{1} ) \wedge ( x_3 \pdl{3} )+
\pdl{2} \wedge ( x_{3}{ }^2 \pdl{1} ) \notag \\& -
( x_2 \pdl{2}) \wedge ( x_3 \pdl{1} )+
( x_2 \pdl{1} ) \wedge (x_3 \pdl{3} ) -
\pdl{3} \wedge (x_3{ }^2 \pdl{1}) +
\pdl{2} \wedge (x_2 x_3 \pdl{1} ) 
\ . \label{eqn:P:T:like}
\end{align}  
We show some examples of Betti numbers of the cohomologies defined by the
above $\pi$. 

\begin{center}
        \setlength{\tabcolsep}{0.7pt}
        \setlength{\extrarowheight}{-3pt}
\begin{tabular}{|c|*{19}{c}|} \hline 
        wt=$-3$
& $\ds\myC{3}{-3}$
& $\rightarrow $& $ \myC{4}{-3} $ 
& $\rightarrow $& $ \myC{5}{-3} $ 
& $\rightarrow $& $ \myC{6}{-3} $ 
& $\rightarrow $& $ \myC{7}{-3} $ 
& $\rightarrow $& $ \myC{8}{-3} $ 
& $\rightarrow $& $ \myC{9}{-3} $ 
& $\rightarrow $& $ \myC{10}{-3} $ 
& $\rightarrow $& $ \myC{11}{-3} $ 
& $\rightarrow $& $ \myC{12}{-3} $ 
    \\\hline 
$\dim$    & 1 && 9 && 36 && 84 && 126 && 126 && 84 && 36 && 9 &&1  \\\hline
$\ker\dim$& 0 && 1 && 8 && 28 && 56 && 70 && 56 && 28 && 8 && 1 \\\hline
Betti &     0 && 0 && 0 && 0 && 0 && 0 && 0 && 0&& 0 && 0 \\ \hline 
\end{tabular} 

\begin{tabular}{|c|*{23}{c}|} \hline 
        wt=$-2$
& $\ds\myC{2}{-2}$
& $\rightarrow $& $ \myC{3}{-2} $ 
& $\rightarrow $& $ \myC{4}{-2} $ 
& $\rightarrow $& $ \myC{5}{-2} $ 
& $\rightarrow $& $ \myC{6}{-2} $ 
& $\rightarrow $& $ \myC{7}{-2} $ 
& $\rightarrow $& $ \myC{8}{-2} $ 
& $\rightarrow $& $ \myC{9}{-2} $ 
& $\rightarrow $& $ \myC{10}{-2} $ 
& $\rightarrow $& $ \myC{11}{-2} $ 
& $\rightarrow $& $ \myC{12}{-2} $ 
& $\rightarrow $& $ \myC{13}{-2} $ 
    \\\hline 
$\dim$    & 3 && 27 && 126 && 414 && 1026 && 1890 && 2520 && 2376 &&
1539 && 651 && 162 && 18  \\\hline
$\ker\dim$& 0 && 4 && 26 && 103 && 315 && 722 && 1183 && 1346 && 1032 &&
507 && 144 && 18\\\hline
Betti &     0 && 1 && 3 && 3 && 4 && 11 && 15 && 9&& 2 && 0 &&0 && 0 \\ \hline 
\end{tabular} 

\setlength{\tabcolsep}{0.1truept} 
\begin{tabular}{|c|*{27}{c}|} \hline 
        $w=-1$
& $\ds\myC{1}{w}$
& $\rightarrow $& $ \myC{2}{w} $ 
& $\rightarrow $& $ \myC{3}{w} $ 
& $\rightarrow $& $ \myC{4}{w} $ 
& $\rightarrow $& $ \myC{5}{w} $ 
& $\rightarrow $& $ \myC{6}{w} $ 
& $\rightarrow $& $ \myC{7}{w} $ 
& $\rightarrow $& $ \myC{8}{w} $ 
& $\rightarrow $& $ \myC{9}{w} $ 
& $\rightarrow $& $ \myC{10}{w} $ 
& $\rightarrow $& $ \myC{11}{w} $ 
& $\rightarrow $& $ \myC{12}{w} $ 
& $\rightarrow $& $ \myC{13}{w} $ 
& $\rightarrow $& $ \myC{14}{w} $ 
    \\\hline 
$\dim$    & 3 && 27 && 162 && 768 && 2745 && 7371 && 15084 && 23544 &&
    27621 && 23745 && 14418 && 5832 && 1407 && 153 \\\hline
$\ker\dim$& 1 && 5 && 25 && 142 && 663 && 2228 && 5481 && 10124 && 13974 &&
    14039 && 9872 && 4579 && 1254 && 153\\\hline
Betti &     1 && 3 && 3 && 5 && 37 && 146 && 338 && 521 && 554 && 392 && 166
    && 33 && 1 && 0 \\ \hline 
\end{tabular} 
\end{center} 

\subsection{Poisson cohomology of polynomial modules}
Since $\ds\Phi(\Lambda^m \frakX{pol}) = \mR[x_1,\ldots,x_n] 
\otimes \Lambda^m \frakX{0}$, we have a decomposition 
$\ds\Phi(\Lambda^m \frakX{pol}) = \oplus_{p} \Delta^m_p$
where 
the subspace $\ds\Delta^m_p$ is given by 
$\ds\Delta^m_p = p\text{-polynomials}
\otimes \Lambda^m \frakX{0}$.   
\begin{defn}

For a given non-negative integer $h$, the weight of each non-zero
element of $\ds\Delta^m_p$ is defined as $p-(h-1)m$.  We define the
space of the elements of degree $m$ and of weight $w$, $\ds\myCC{m}{w}$
by 
\begin{equation} \ds\myCC{m}{w} = 
\left(w+(h-1)m\right)\text{-polynomials} \otimes \Lambda^m
\frakX{0}\ . \label{defn:cochain:poly} \end{equation}  
\end{defn}
We see easily the next Proposition. 
\begin{kmProp} 
If $\ds\pi\in \Delta^2_{h}$, then $\ds\SbtS{\pi}{\myCC{m}{w}} \subset
\myCC{m+1}{w}$. Furthermore, if $\ds\SbtS{\pi}{\pi}=0$ then for each fixed
weight $w$, $\ds\{ \myCC{m}{w} \}_{m} $ with $ u \rightarrow \SbtS{\pi}{u}$
forms  a cochain complexes.  
\end{kmProp} 
We may call the cohomology groups of the cochain complexes above as
homogeneous Poisson polynomial cohomology groups. 

Using 
$\Phi: \oplus \Lambda^{\bullet}_{\mR}\frakX{pol} \rightarrow 
\oplus \Lambda^{\bullet}\frakX{pol}$ in Remark \ref{rmk:Phi},  
we have a commutative diagram: 
\begin{equation}
        \begin{CD}
                \myC{m}{w}  @> \SbtR{\pi}{\cdot} >>  @. \myC{m+1}{w}    \\
@V\Phi VV   @. @VV\Phi V \\ 
\myCC{m}{w}  @>\SbtS{\Phi(\pi)}{\cdot}>> @.  \myCC{m+1}{w}    
\end{CD}
\label{Konta:CD}
\end{equation} 
We remark that if $m>n$ then $\ds\Phi(\myC{m}{w}) = 0$ even though $\ds 
\myC{m}{w} \ne 0$. 

If $h=1$ in (\ref{defn:cochain:poly}), we have directly the next
proposition.   
\begin{kmProp} \label{prop:euler:hOne} 
On $\mR^n$, for each weight $w$ and for each 1-homogeneous Poisson
structure, the Euler characteristic of Poisson polynomial cohomology
groups is always zero.  
\end{kmProp}
\textbf{Proof:} 
$\ds\dim \myCC{m}{w} = 
\dim (w\text{-polynomials})\: \dim( \Lambda^m \frakX{0}) = \tbinom{n-1+w}{n-1}
\tbinom{n}{m}$, and \[\ds\sum_{m=0}^n (-1)^m \ \dim 
\myCC{m}{w} = 
\tbinom{n-1+w}{n-1} \sum_{m=0}^n (-1)^m  \tbinom{n}{m} = 0 \ .\] \kmqed

The $\Phi$-image of the 2-homogeneous Poisson-like 2-vector
(\ref{eqn:P:T:like}) in the previous subsection, is just 
\[\overline{\pi} =  
( x_2{}^2 
- x_{3}{}^2 ) \pdl{1} \wedge  \pdl{2} 
% + x_2  x_3 \pdl{1} \wedge \pdl{2} 
% - x_2 x_3 \pdl{1} \wedge  \pdl{2}  
+ 2 ( x_2 x_3 
+  x_3{}^2) \pdl{1}  \wedge  \pdl{3} 
\] 
and satisfies $\ds\SbtS{\overline{\pi}}{\overline{\pi}} =0$, namely
$\ds\overline{\pi}$ is a usual Poisson 2-vector field.  
In the following, we show several examples of the Poisson polynomial
cohomology groups of $\ds\overline{\pi}$ on $\mR^3$.  
\begin{center} 
        \setlength{\tabcolsep}{2truept}
        \setlength{\extrarowheight}{-3pt}
\begin{tabular}{|c|*{1}{c}|} \hline 
wt=$-3$ &  $\ds\myCC{3}{-3} $
    \\\hline 
$\dim$    & 1 \\\hline
$\ker\dim$& 1 \\\hline
Betti &     1 \\ \hline 
\end{tabular} 
\hfil 
\begin{tabular}{|c|*{3}{c}|} \hline 
wt=$-2$ &  $\ds\myCC{2}{-2} $
& $ \rightarrow $ & $\ds\myCC{3}{-2} $ 
    \\\hline 
$\dim$    & 3 && 3    \\\hline
$\ker\dim$& 2 && 3   \\\hline
Betti &     2 && 2  \\ \hline 
\end{tabular} 
\hfil 
\begin{tabular}{|c|*{5}{c}|}
        \hline 
        wt=$-1$
&  $\ds\myCC{1}{-1} $
& $ \rightarrow $ & $\ds\myCC{2}{-1} $ 
& $ \rightarrow $ & $\ds\myCC{3}{-1} $ 
    \\\hline 
$\dim$    & 3 && 9 && 6   \\\hline
$\ker\dim$& 1 && 6 && 6 \\\hline
Betti &     1 && 4 && 3 \\ \hline 
\end{tabular} 

\begin{tabular}{|c|*{7}{c}|}
        \hline 
        wt=0
&  $\ds\myCC{0}{0} $
& $ \rightarrow $ & $\ds\myCC{1}{0} $ 
& $ \rightarrow $ & $\ds\myCC{2}{0} $ 
& $ \rightarrow $ & $\ds\myCC{3}{0} $ 
    \\\hline 
$\dim$    & 1 && 9 && 18 && 10  \\\hline
$\ker\dim$& 1 && 4 && 12 && 10 \\\hline
Betti &     1 && 4 && 7 && 4 \\ \hline 
\end{tabular} 
\hfil 
\begin{tabular}{|c|*{7}{c}|} \hline 
wt=1 &  $\ds\myCC{0}{1} $
& $ \rightarrow $ & $\ds\myCC{1}{1} $ 
& $ \rightarrow $ & $\ds\myCC{2}{1} $ 
& $ \rightarrow $ & $\ds\myCC{3}{1} $ 
    \\\hline 
$\dim$    & 3 && 18 && 30 && 15  \\\hline
$\ker\dim$& 0 && 7 && 20 && 15 \\\hline
Betti &     0 && 4 && 9 && 5 \\ \hline 
\end{tabular} 
\hfil 
\begin{tabular}{|c|*{7}{c}|} \hline 
wt=2 &  $\ds\myCC{0}{2} $
& $ \rightarrow $ & $\ds\myCC{1}{2} $ 
& $ \rightarrow $ & $\ds\myCC{2}{2} $ 
& $ \rightarrow $ & $\ds\myCC{3}{2} $ 
    \\\hline 
$\dim$    & 6 && 30 && 45 && 21  \\\hline
$\ker\dim$& 0 && 10 && 30 && 21 \\\hline
Betti &     0 && 4 && 10 && 6 \\ \hline 
\end{tabular} 
\hfil 
\begin{tabular}{|c|*{7}{c}|} \hline 
wt=3 &              $\ds\myCC{0}{3} $
& $ \rightarrow $ & $\ds\myCC{1}{3} $ 
& $ \rightarrow $ & $\ds\myCC{2}{3} $ 
& $ \rightarrow $ & $\ds\myCC{3}{3} $ 
    \\\hline 
$\dim$    & 10 && 45 && 63 && 28  \\\hline
$\ker\dim$& 0 && 15 && 42 && 28 \\\hline
Betti &     0 && 5 && 12 && 7 \\ \hline 
\end{tabular} 
\hfil 
\begin{tabular}{|c|*{7}{c}|} \hline 
wt=4 &              $\ds\myCC{0}{4} $
& $ \rightarrow $ & $\ds\myCC{1}{4} $ 
& $ \rightarrow $ & $\ds\myCC{2}{4} $ 
& $ \rightarrow $ & $\ds\myCC{3}{4} $ 
    \\\hline 
$\dim$    & 15 && 63 && 84 && 36  \\\hline
$\ker\dim$& 0 && 21 && 56 && 36 \\\hline
Betti &     0 && 6 && 14 && 8 \\ \hline 
\end{tabular} 
\hfil 
\begin{tabular}{|c|*{7}{c}|} \hline 
wt=5 &              $\ds\myCC{0}{5} $
& $ \rightarrow $ & $\ds\myCC{1}{5} $ 
& $ \rightarrow $ & $\ds\myCC{2}{5} $ 
& $ \rightarrow $ & $\ds\myCC{3}{5} $ 
    \\\hline 
$\dim$    & 21 && 84 && 108 && 45  \\\hline
$\ker\dim$& 0 && 28 && 72 && 45 \\\hline
Betti &     0 && 7 && 16 && 9 \\ \hline 
\end{tabular} 
\end{center}

\bigskip 

\begin{kmRemark} 
In the concrete examples above, the Euler characteristic of Poisson
polynomial cohomology groups is zero except the case of weight is
minimum and the cochain complex is single.  And we expect that the Euler
characteristic of the Poisson polynomial cohomology groups of
$2$-homogeneous Poisson structure may be zero in general.   When the
case of  $3$-homogeneous Poisson (only depends on homogeneity 3 but not
depends on the structure itself) on $\ds\mR^3$, we have the distribution
of the Euler characteristic below.     
\begin{center} 
        \setlength{\tabcolsep}{6truept}
        \setlength{\extrarowheight}{0pt}
\begin{tabular}{c|*{14}{r}}
$h=3$, wt & $-6$ & $-5$ & $-4$ & $-3$ & $-2$ & $-1$ & 0 & 1 & 2 & 3& 4 & 5 &
        6\\\hline
        Euler & $-1$ & $-3$ & $-3$ & $ -1$ & 0 & 0 & 0  & 0 & 0 & 0  & 0 & 0 & 0 \end{tabular}
\end{center} 
The results only depends on homogeneity 3 but not depends on the Poisson
structure itself on $\ds\mR^3$.     Still we expect the Euler
characteristic may be zero for higher weights.  
\end{kmRemark}
\medskip

In fact, we have the following result including Proposition
\ref{prop:euler:hOne}.  
\begin{thm} \label{thm:euler:hAll} 
On $\mR^n$, for each $h$-homogeneous Poisson structure, the Euler
characteristic of Poisson polynomial cohomology groups is always zero
for each weight $w \geq 1-n $ when $h>0$ and 
for each weight $w \geq 1 $ when $h=0$.  
\end{thm}

In order to prove the theorem above, we follow binomial expansion
theorem twice.  

Let \(h\) be non-negative integer and \(x\) be an indeterminate
variable.  Let us start from the binomial expansion 

\begin{equation}\label{eqn:moto}
( (x+1)^{h-1} -1)^n = \sum_{m} \tbinom{n}{m} (-1)^{n-m} (x+1)^{(h-1)m}
\;. 
\end{equation}
Multiplying the above (\ref{eqn:moto})
by \( (x+1)^{n-1+w}\), and expand as follows: 
\begin{align}
        \label{eqn:last}
        &  (x+1)^{n-1+w}
( (x+1)^{h-1} -1)^n
=  \sum_{m} \tbinom{n}{m} (-1)^{n-m} 
(x+1)^{n-1+w}
(x+1)^{(h-1)m} \\ \notag
=& \sum_{m} \tbinom{n}{m} (-1)^{n-m} 
(x+1)^{n-1+w + (h-1) m } 
= \sum_{m} \tbinom{n}{m} (-1)^{n-m} \sum_{k}
\tbinom{ n-1+w + (h-1) m }{k} 
x^k  \\ \notag
=& (-1)^n \sum_{k} \sum_{m} (-1)^m \tbinom{n}{m} 
\tbinom{ n-1+w + (h-1) m }{k} 
x^k  \;. 
\end{align} 
When $h>0$,  
comparing the coefficients of \( x^{n-1}\) of the both sides, we
conclude that if \(n-1+w \geq 0\), then 
\begin{equation}
        \label{eqn:final}
%\text{If}\quad n-1+w \geq 0,\quad  \text{then}\quad 
 \sum_{m}  (-1)^{m} 
\tbinom{ n-1+w + (h-1) m }{n-1} \tbinom{n}{m}= 0  
\end{equation}
holds. When $h=0$, we rewrite 
the left-hand-side of 
(\ref{eqn:last}) and have $\ds (x+1)^{-1+w} (-x)^n$, so assuming
$-1+w\geq 0$, we have 
(\ref{eqn:final}) for $h=0$.  
\kmqed

% \young(123,455,6\null5)
% \young(\null\null\null,\null,\null)

%% file: Suppl-2.tex
\subsection{Discussion}

We remark that 
$w$ and $m$ have to satisfy $w+(h-1)m \geq 0$ and $0\leq m \leq n$. 

\begin{center}\setlength{\unitlength}{5mm}
        \fbox{
        \begin{picture}(10,6)(-2,-1) 
                \path(-2,0)(6,0) 
                \put(6,0){\vector(1,0){1}}
        \put(7.1,0){\makebox(0,0)[l]{$w$}}
                \path(0,-1)(0,3) 
        \put(0,3){\vector(0,1){1}}
        \put(0,4.1){\makebox(0,0)[b]{$m$}}
        \put(2,4.1){\makebox(0,0)[b]{$h=0$}}
        \put(-0.2,-0.2){\makebox(0,0)[tr]{$O$}}
        \path(-2,3)(6,3)
        \path(0,0)(4,4)
        \put(-0.2,2.8){\makebox(0,0)[tr]{$n$}}
        \shade{\path(0,0)(3,3)(5,3)(5,0)(0,0)}
        \end{picture}
}\hfil 
        \kmcomment{
        \begin{picture}(10,6)(-2,-1) 
                \path(-2,0)(6,0) 
                \put(6,0){\vector(1,0){1}}
        \put(7.1,0){\makebox(0,0)[l]{$w$}}
                \path(0,-1)(0,3) 
        \put(0,3){\vector(0,1){1}}
        \put(0,4.1){\makebox(0,0)[b]{$m$}}
        \put(2,4.1){\makebox(0,0)[b]{$h=1$}}
        \put(-0.2,-0.2){\makebox(0,0)[tr]{$O$}}
        \path(-2,3)(6,3)
 %       \path(0,0)(4,4)
        \put(-0.2,2.8){\makebox(0,0)[tr]{$n$}}
        \shade{\path(0,0)(0,3)(5,3)(5,0)(0,0)}
        \end{picture}
}%endOFkmcomment 
        \fbox{
        \begin{picture}(10,6)(-2,-1) 
                \path(-2,0)(6,0) 
                \put(6,0){\vector(1,0){1}}
        \put(7.2,0){\makebox(0,0)[l]{$w$}}
                \path(0,-1)(0,3) 
        \put(0,3){\vector(0,1){1}}
        \put(0,4.2){\makebox(0,0)[b]{$m$}}
        \put(2,4.1){\makebox(0,0)[b]{$h>1$}}
        \put(-0.2,-0.2){\makebox(0,0)[tr]{$O$}}
        \path(-2,3)(6,3)
        \path(0.5,-1)(-2,4) 
        \put(1.2,-0.6){\makebox(0,0)[l]{$w+(h-1)m=0$}}
        \put(-0.2,3.2){\makebox(0,0)[br]{$n$}}
        \shade{\path(0,0)(-1.5,3)(5,3)(5,0)(0,0)}
        \end{picture}
} 
\end{center}

If $h=1$ in (\ref{defn:cochain:poly}) then $w$ and $m$ are independent. 

\subsubsection{h=0}
Since 
$w$ and $m$ have to satisfy $w+(h-1)m\geq 0$ and $ 0 \leq m \leq n$, 
when $h=0$ then $m$ runs $ 0 \leq m \leq \min(n,w) $. 
\paragraph{w = 0, not satisfy $w >0$:}
\[
        \text{
        LHS of (\ref{eqn:final})} = \sum_{m=0}^0 (-1)^m 
\binom{n-1 + 0-m}{n-1} \binom{n}{m} = 
 \binom{n-1}{n-1} \binom{n}{0}   
= 1 
\]
\paragraph{w = 1:}
\begin{align*}
        \text{
        LHS of (\ref{eqn:final})} =&  \sum_{m=0}^1 (-1)^m 
                \binom{n-1 + 1-m}{n-1} \binom{n}{m}\\ = &  
 \binom{n}{n-1} \binom{n}{0}   
 - \binom{n-1}{n-1} \binom{n}{1}   
=  \binom{n}{n-1} - \binom{n}{1}   = 0 
\end{align*}

\paragraph{w = 2:}
\begin{align*}
        \text{
        LHS of (\ref{eqn:final})} =& \sum_{m=0}^2 (-1)^m 
                \binom{n-1 + 2-m}{n-1} \binom{n}{m}  =  
 \binom{n+1}{n-1} \binom{n}{0}   
 - \binom{n}{n-1} \binom{n}{1}   
 + \binom{n-1}{n-1} \binom{n}{2}   
 \\
 = &  \binom{n+1}{2} - \binom{n}{1}^2 + \binom{n}{2} = 0 
\end{align*}

\subsubsection{h=2}
Since 
$w$ and $m$ have to satisfy $w+(h-1)m\geq 0$ and $ 0 \leq m \leq n$, 
when $h=2$ then $m$ runs $\max(0,-w) \leq m \leq n$. 

\paragraph{n = 4 and w = -2:}
\[
        \text{
        LHS of (\ref{eqn:final})} = \sum_{m=2}^4 (-1)^m 
\binom{1+m}{3} \binom{4}{m} = 
\binom{3}{3} \binom{4}{2} -  
\binom{4}{3} \binom{4}{3} +  
\binom{5}{3} \binom{4}{4}   
= 
 \binom{4}{2} - 4^2  
 +  
\binom{5}{2} = 0 
\]

\paragraph{n = 4 and w = -3:}
\[
        \text{
        LHS of (\ref{eqn:final})} = \sum_{m=3}^4 (-1)^m 
\binom{m}{3} \binom{4}{m} = 
- 
\binom{3}{3} \binom{4}{3} +  
\binom{4}{3} \binom{4}{4}   
= 
- \binom{4}{3} + 
 +  
\binom{4}{3} = 0 
\]

\paragraph{n = 4 and w = -4,  not satisfy $w \geq 1-n$:}
\[
        \text{
        LHS of (\ref{eqn:final})} = \sum_{m=4}^4 (-1)^m 
\binom{m-1}{3} \binom{4}{m} = 
- 
\binom{3}{3} \binom{4}{4}   
= 
1 
\]

\subsubsection{h=3}
Since 
$w$ and $m$ have to satisfy $w+(h-1)m\geq 0$ and $ 0 \leq m \leq n$, 
when $h=3$ then $m$ runs $\max(0,-w/2) \leq m \leq n$. 
\paragraph{n = 4 and w = -2:}
\[
        \text{
        LHS of (\ref{eqn:final})} = \sum_{m=1}^4 (-1)^m 
\binom{1+2m}{3} \binom{4}{m} = 
- \binom{3}{3} \binom{4}{1}   
+ \binom{5}{3} \binom{4}{2}   
- \binom{7}{3} \binom{4}{3}   
+ \binom{9}{3} \binom{4}{4}   
= 0 
\]
\paragraph{n = 4 and w = -3:}
\[
        \text{
        LHS of (\ref{eqn:final})} = \sum_{m=2}^4 (-1)^m 
\binom{2m}{3} \binom{4}{m} = 
  \binom{4}{3} \binom{4}{2}   
- \binom{6}{3} \binom{4}{3}   
+ \binom{8}{3} \binom{4}{4}   
= 0 
\]
\paragraph{n = 4 and w = -4 not satisfy $w \geq 1-n$:}
\[
        \text{
        LHS of (\ref{eqn:final})} = \sum_{m=2}^4 (-1)^m 
\binom{2m-1}{3} \binom{4}{m} = 
  \binom{3}{3} \binom{4}{2}   
- \binom{5}{3} \binom{4}{3}   
+ \binom{7}{3} \binom{4}{4}   
= 1 
\]

\subsection{Stability of Poisson polynomial cohomology groups of special Lie Poisson structures on 3-space}
Two Lie Poisson structures due to Heisenberg Lie algebra and 
$\fraksp(\ds\mR^2)$ have special features. 

\subsubsection{Heisenberg case:} 
\begin{thm} The dimension of  
       the kernel subspaces and the Betti numbers of Heisenberg Lie
        Poisson structure are given by 
\begin{center} %\setlength{\tabcolsep}{2truept}
\begin{tabular}{|c|*{7}{c}|} 
        \hline 
$w>0$ &              $\ds\myCC{0}{w} $
& $ \rightarrow $ & $\ds\myCC{1}{w} $ 
& $ \rightarrow $ & $\ds\myCC{2}{w} $ 
& $ \rightarrow $ & $\ds\myCC{3}{w} $ 
    \\\hline 
    $\dim$    & $\tbinom{2+w}{2}$  && 
    $3\tbinom{2+w}{2}$  && 
    $3\tbinom{2+w}{2}$  && 
    $\tbinom{2+w}{2}$   \\\hline
    $\ker\dim$ & 1 && $\tbinom{3+w}{2}$ && $(w+3)(w+1)$ &&
    $\tbinom{w+2}{2}$ \\\hline
Betti &     1 && $w+3$ && $2w+3$ && $w+1$ \\ \hline\hline 
    $\ker\dim(w=0)$ & 1 && 2  && 3 && 1  \\\hline
    Betti ($w=0$)&     1 && 2 && 2 && 1 \\ \hline
\end{tabular} 
\end{center}
Thus, we see that 
\begin{align*} \dim 
        \myHH{0}{w+1} & -  \dim  \myHH{0}{w} = 0\; ( w> 0), \quad   
        \dim \myHH{1}{w+1} -  \dim  \myHH{1}{w} = 1\; ( w> 0), \\    
        \dim \myHH{2}{w+1} & -   \dim \myHH{2}{w} = 2\; ( w> 0), \quad   
        \dim \myHH{3}{w+1} -  \dim  \myHH{3}{w} = 1\; ( w> 0). \end{align*} 

\end{thm}

\textbf{Proof: } 
2-vector field is given by 
$\ds \Pkt{z}{x} = 0,  \Pkt{z}{y} =  0,  \Pkt{x}{y} = z$ i.e.,  
\[\pi =\frac{1}{2} \sum_{i,j}\Pkt{x_i}{x_j}\pdel_{x_i}\wedge\pdel_{x_j} 
= z \pdel_x \wedge \pdel_y  
\] 
and satisfies $\SbtS{\pi}{\pi}=0$ and so a Poisson structure on $\ds\mR^3$. 
\begin{alignat*}{3}
        \SbtS{\pi}{x} &= - z \pdel_y \;, &
        \SbtS{\pi}{y} &=  z \pdel_x \;, & 
\SbtS{\pi}{z} &=  0 \;, \\
        \SbtS{\pi}{\pdel_x} &= 0 \;, &
        \SbtS{\pi}{\pdel_y} &= 0 \;, &
\SbtS{\pi}{\pdel_z} &= -  \pdel_x \wedge  \pdel_y \;, 
\end{alignat*}
are the necessary information. We see that 
\[
        \SbtS{\pi}{\pdel_x \wedge \pdel y} = 
        \SbtS{\pi}{\pdel_x \wedge \pdel z} = 
\SbtS{\pi}{\pdel_y \wedge \pdel z} = 0 \] 
and 
\[ \SbtS{ \pi }{u^A} = -a_{1} x^{a_1 -1} y^{a_2} z^{a_3+1} \pdel_{y}
+a_{2} x^{a_1 } y^{a_2 -1} z^{a_3+1} \pdel_{x} \]
where $\ds u^A = x^{a_1} y^{a_2} z^{a_3}$ for $A=[a_1,a_2,a_3]\in\mN^3$.

Now we assume $w>0$. 
$\ds \myCC{0}{w}:$
\begin{align*} 
        0  &= \SbtS{\pi} { \sum_{|A|=w} c_A u^A } 
        = \sum_{|A|=w}  c_A ( 
        -a_{1} x^{a_1 -1} y^{a_2} z^{a_3+1} \pdel_{y} 
+a_{2} x^{a_1 } y^{a_2 -1} z^{a_3+1} \pdel_{x}) 
\end{align*} implies 
\begin{align}
        \label{al:one}
        & \sum_{|A|=w}  c_A  
        a_{1} x^{a_1 -1} y^{a_2} z^{a_3+1} = 0 \\
        \label{al:two}
        &\sum_{|A|=w}  c_A  
a_{2} x^{a_1 } y^{a_2 -1} z^{a_3+1} = 0 
\end{align} 
(\ref{al:one}) says that if $\ds a_1 \ne 0$ then $\ds c_A=0$ or 
(\ref{al:two}) says that if $\ds a_2 \ne 0$ then $\ds c_A=0$. Thus,  the
kernel subspace of 
$\ds \myCC{0}{w}$ is spanned by $\ds c_{[0,0,w]} z^w$ and 1-dimensional. 

\medskip

About $\ds \myCC{n}{w}$ ($n=3$), the whole space is the kernel subspace and 
\( \tbinom{n-1+w}{n-1}\)-dimensional. The difference of the dimension of
the kernel subspaces with the weight $w$ and weight $w+1$ is
$1+w$. 

\medskip

$\ds \myCC{1}{w}$:  
\begin{align*} 0 =& 
\SbtS{\pi}{ \sum_{A}( \alpha_A u^A \pdel_x + 
\beta_A u^A \pdel_y + \gamma_A u^A \pdel_z ) } \\
=& 
\sum_{A} \alpha_A \SbtS{\pi}{u^A } \wedge \pdel_x + 0 
+ \sum_{A} \beta_A \SbtS{\pi}{u^A } \wedge \pdel_y + 0 
%\\& 
+ \sum_{A} \gamma_A \SbtS{\pi}{u^A } \wedge \pdel_z  
- \sum_{A} \gamma_A u^A  \pdel_x \wedge \pdel_y 
\\
= &
\sum_{A} \alpha_A a_1 x^{a_1 -1} y^{a_2}  z^{a_3+1} \pdel_x \wedge \pdel_y 
%\\&
+ \sum_{A} \beta_A a_2 x^{a_1 } y^{a_2-1}  z^{a_3+1} \pdel_x \wedge \pdel_y 
%\\& 
+ \sum_{A} \gamma_A \SbtS{\pi}{u^A } \wedge \pdel_z  
\\& 
- \sum_{A} \gamma_A u^A  \pdel_x \wedge \pdel_y 
\end{align*}
Thus, we have 
\begin{align} 
        \label{ali:one}
        & 
\sum_{A} \alpha_A a_1 x^{a_1 -1} y^{a_2}  z^{a_3+1}  
 + \sum_{A} \beta_A a_2 x^{a_1 } y^{a_2-1}  z^{a_3+1} 
%\\ \notag &
- \sum_{A} \gamma_A u^A  = 0\;,  
\\ 
        \label{ali:two}
        &
\sum_{A} \gamma_A a_1 x^{a_1 -1} y^{a_2}  z^{a_3+1} =0 \;,  \\ 
        \label{ali:three}
        & 
\sum_{A} \gamma_A a_2 x^{a_1 } y^{a_2-1}  z^{a_3+1} =0 \;, 
        \\
        \noalign{
%\end{align}
                and 
(\ref{ali:two}) and (\ref{ali:three}) say that} 
%\begin{equation}
        \label{eqn:one} 
         & \gamma_A =0 \quad \text{ unless }
A = [0,0,w]\; . %\end{equation} 
\end{align}

(\ref{ali:one}) and 
(\ref{eqn:one}) imply that 
\[
\sum_{A} \alpha_A a_1 x^{a_1 -1} y^{a_2}  z^{a_3+1}  
 + \sum_{A} \beta_A a_2 x^{a_1 } y^{a_2-1}  z^{a_3+1} 
- \gamma_{[0,0,w]}  z^w   = 0 \;. 
\] 
Comparing $\ds z^w$, we have 
\begin{equation}
        \alpha_{[1,0,w-1]} + \beta_{[0,1,w-1]} = \gamma_{[0,0,w]}\; .  
\end{equation} 

Comparing $\ds z^j$ ($j < w$), we have ($w-j+1$) linear equations 
\[ \alpha_{[1+p,w-j-p, j-1]} (1+p) + \beta_{[p, w-j+1-p,j-1]} (w-j+1-p)
        = 0\;. \] 
Thus, the dim of kernel subspace of $\ds \myCC{1}{w}$ is 
\[ 3 \tbinom{n-1+w}{n-1} - \left( \tbinom{n-1+w}{n-1} -1\right) -
                \sum_{j=1}^w (w-j+1) = \frac{(w+2)(w+3)}{2}\;. \]
%\medskip

$\ds \myCC{2}{w}$:  
\begin{align*}
        0 &= \SbtS{\pi}{
        \sum_{A} ( \alpha_{A} u^A \pdel_x \wedge \pdel_y  + 
        \beta_{A} u^A \pdel_x \wedge \pdel_z  + 
 \gamma_{A} u^A \pdel_y \wedge \pdel_z ) }\\
& = \sum_{A} \SbtS{\pi}{u^A} \wedge (
        \alpha_{A}  \pdel_x \wedge \pdel_y  + 
        \beta_{A} \pdel_x \wedge \pdel_z  + 
 \gamma_{A} \pdel_y \wedge \pdel_z ) \\
 &= 
 \sum_{A}( \beta_{A} a_1 x^{a_1 -1} y^{a_2} z^{a_3+1} 
 + \gamma_{A} a_2 x^{a_1 } y^{a_2-1} z^{a_3+1} ) 
  \pdel_x \wedge \pdel_y \wedge \pdel_z  
 \end{align*}

By the almost same argument, we have linearly independent $(w-j)$ linear
equations for ($0 \leq j < w$). Thus, the dim of the kernel subspace of
$\ds \myCC{2}{w}$ is \[ 3 \tbinom{n-1+w}{n-1} - \sum_{j=0}^{w-1} (w-j) =
(w+1)(w+3)\;.\] Thus, we know the whole stuff and Betti \# as the table
in the theorem.  \kmqed

\subsubsection{$\ds\fraksp(\mR^2)$ Poisson case:}

\begin{thm} The dimension of 
the kernel subspaces and the Betti numbers of $\ds\fraksp(\mR^2)$ Lie
        Poisson structure are given as follows:   
        We have two kinds of 
lists of Betti numbers and they are modulo 2 periodic with respect to
the weight $w$. 

\begin{center} %\setlength{\tabcolsep}{2truept}
\begin{tabular}{|c|*{7}{c}|} 
        \hline 
 &              $\ds\myCC{0}{w} $
& $ \rightarrow $ & $\ds\myCC{1}{w} $ 
& $ \rightarrow $ & $\ds\myCC{2}{w} $ 
& $ \rightarrow $ & $\ds\myCC{3}{w} $ 
    \\\hline 
    $\dim$    & $\tbinom{2+w}{2}$  && 
    $3\tbinom{2+w}{2}$  && $3\tbinom{2+w}{2}$  && $\tbinom{2+w}{2}$\\\hline
    $\ker\dim\ (w \text{ odd }) $ & 0 && $\tbinom{2+w}{2}$ && 
    $2 \tbinom{2+w}{2}$ && $\tbinom{w+2}{2}$ \\\hline
    Betti ($w$ odd) &     0 &&  0 && 0 && 0 \\ \hline\hline 
    $\ker\dim (w \text{ even }) $ & 1 && $\tbinom{2+w}{2}-1$ && 
    $2 \tbinom{2+w}{2}+1$ && $\tbinom{w+2}{2}$ \\\hline
    Betti ($w$ even) &     1 &&  0 && 0 && 1 \\ \hline 
\end{tabular} 
\end{center} 

%\medskip

\end{thm}
\textbf{Proof: } 
2-vector field given by 
$\ds \Pkt{z}{x} = 2 x,  \Pkt{z}{y} = -2 y,  \Pkt{x}{y} = z$ is 
\[\pi =\frac{1}{2} \sum_{i,j}\Pkt{x_i}{x_j}\pdel_{x_i}\wedge\pdel_{x_j} 
= 2 x \pdel_z \wedge \pdel_x - 
 2 y \pdel_z \wedge \pdel_y + 
  z \pdel_x \wedge \pdel_y  
\] 
and satisfies $\SbtS{\pi}{\pi}=0$ and so a Poisson structure on $\ds\mR^3$. 
\begin{alignat*}{3}
\SbtS{\pi}{x} &= 2 x \pdel_z - z \pdel_y , &
\SbtS{\pi}{y} &= -2 y \pdel_z + z \pdel_x ,&
\SbtS{\pi}{z} &= -2 x \pdel_x + 2 y \pdel_y\;, \\
\SbtS{\pi}{\pdel_x} &= -2  \pdel_z \wedge  \pdel_x ,&
\SbtS{\pi}{\pdel_y} &= 2  \pdel_z \wedge  \pdel_y ,&
\SbtS{\pi}{\pdel_z} &= -  \pdel_x \wedge  \pdel_y , 
\end{alignat*}
are the necessary information to get Poisson polynomial cohomologies of
our Poisson structure $\pi$. 

Since $\ds \myCC{m}{w} = w\text{-homog} \otimes
\Lambda^m \frakX{0}$,  $\ds \myCC{m}{0} = \Lambda^m \frakX{0}$ and   
\[ 0= \SbtS{\pi}{\sum_{j}c_j \pdel_{x_j}} = c_1 
(-2  \pdel_z \wedge  \pdel_x )+ c_2 ( 
 2  \pdel_z \wedge  \pdel_y) + c_3 (
 -  \pdel_x \wedge  \pdel_y ) \] 
implies $\ds c_1=c_2=c_3=0$, i.e., 
$\ds\dim\ker \text{ subspace of } \myCC{1}{0} = 0$.  

Since \(\ds 
\SbtS{\pi}{\sum_{j}c_j \pdel_{x_{j+1}} \wedge \pdel_{x_{j+2}} 
} = 0\), we see that  
 $\dim$ of the kernel subspace of $\myCC{2}{0} = 3$, and   
 $\dim$ of the kernel subspace of $\myCC{3}{0} = 1$. We can continue this
 argument for weight 1,2,\ldots.  

%\medskip

\paragraph{Kernel subspace of 0-th cochain space:} 
About the kernel subspace of $\ds \myCC{0}{w}$, we see that if the weight
is odd then 0 and if the weight is even then 1 dimensional as follows.   
Take a cochain $\ds \sum_{A} c_A u^A\in \myCC{0}{w}$ and suppose 
$\ds 0= \SbtS{\pi} { \sum_{A} c_A u^A}$. 
Then we have 
\begin{align}
        \label{ii:one}        & 2 \sum_{A}(a_1-a_2)c_A u^A =0\;, \\ 
  \label{ii:two}       & \sum_{A}c_A ( z \frac{\pdel }{\pdel y}   -2
        x \frac{\pdel }{\pdel z} ) u^A =0 ,\; 
        \\& 
        \sum_{A}c_A ( z \frac{\pdel }{\pdel x}  -2 y \frac{\pdel }{\pdel z} ) u^A  =0\;.   
\end{align} 
(\ref{ii:one}) implies $\ds (a_1-a_2)c_A  =0$, so we see that 
if $\ds a_1\ne a_2$ then $\ds c_A=0$.  Thus, we have to care only about 
$\ds c_{[i,i,w-2i]}$ which we denote by $\ds \myP{i}$   
for $i$ from 0 to $w$.  
Using (\ref{ii:two}), we have 
\begin{align*} 0= & 
\sum_{0\leq i \leq w/2} \myP{i} ( 
        z \frac{\pdel }{\pdel y}  -2 x \frac{\pdel }{\pdel z} )  x^{i} y^{i} z^{w-2i}
        \\
        =& 
\sum_{0 < i \leq w/2} i \myP{i}  
         x^{i} y^{i-1} z^{w-2i+1}
        - 
         \sum_{0\leq i < w/2} 2(w-2i) \myP{i}  
         x^{i+1} y^{i} z^{w-2i-1}
         \\
         \kmcomment{
         =&          \sum_{0 \leq i \leq w/2-1} (i+1) \myP{i+1}  
         x^{i+1} y^{i} z^{w-2i-1}
         \\& 
        - 
         \sum_{0\leq i < w/2}2 (w-2i) \myP{i}  
         x^{i+1} y^{i} z^{w-2i-1} 
         \\
 }%endOFkmcomment
         =& 
         \sum_{0 \leq i \leq w/2-1}
         \left( (i+1) \myP{i+1}  - 2(w-2i) \myP{i} \right)  
         x^{i+1} y^{i} z^{w-2i-1}
%         \\& 
        -  \sum_{w/2-1 < i < w/2} 2 (w-2i) \myP{i}  
         x^{i+1} y^{i} z^{w-2i-1} 
        \end{align*}
If $w$ is even, then the last term does not happen and we see that
$\ds (i+1) \myP{i+1}  - 2(w-2i) \myP{i} =0$ for $i=0.. w/2-1$. Thus, the
freedom of $c_A$ is 1. 
If $w$ is odd, say $w=2\Omega+1$, then the last term tells that
$\ds\myP{\Omega} = 0$ and also 
$\ds (i+1) \myP{i+1}  - 2(w-2i) \myP{i} =0$ for $i=0.. \Omega-1$. Thus, 
$\ds \myP{i} =0$ for $i=0.. \Omega$ and the kernel subspace is
0-dimensional. Namely, 0-th Betti number is 1 if $w$ is even and 0
otherwise.

\paragraph{Kernel subspace of 2-nd cochain space:} 
We try now 3-rd Betti number. Since the kernel subspace of 
$\ds \myCC{3}{w}$ is itself, we search 
the kernel subspace of $\ds \myCC{2}{w}$.  
Take a general 2-cochain $\ds \sum_{A} u^A (
\alpha_A \pdel_x \wedge \pdel_y 
+\beta_A \pdel_x \wedge \pdel_z 
+\gamma _A \pdel_y \wedge \pdel_z)$ and take the Schouten bracket by
$\pi$. Then the coefficient of $\pdel_x \wedge \pdel_y \wedge \pdel_z$
is 
\begin{align} &  \sum_{A} \left(  2 \alpha_A (a_1-a_2)  +
\beta_A 
( z \frac{\pdel}{\pdel x} - 2 y \frac{\pdel}{\pdel z} )
%\\& \notag 
+ 
\gamma_A 
( z \frac{\pdel}{\pdel y} - 2 x \frac{\pdel}{\pdel z} )\right) u^A = 0 \; . 
\end{align}
The coefficient of $\ds x^{i} y^{i} z^{w-2i}$ is 
\begin{align}
        \label{align:one}
        & \beta_{[1+i,i]} (1+i) -2 \beta_{[i,i-1]} (w-2i+1) 
%        \\ \notag & 
        + \gamma_{[i,1+i]}(1+i) -2 \gamma_{[i-1,i]}(w-2i+1) = 0   
\end{align} 
Using the notation $\ds Q_{[i,j]} = \beta_{[i,j]} + \gamma_{[j,i]}$, 
(\ref{align:one}) is equivalent to 
\begin{align}
        \label{jj:one} & Q_{[1,0]} = 0 \\
        \label{jj:two} & (1+i) Q_{[i+1,i]} =  2 ( w-2i+1) Q_{[i,i-1]} \ 
        (0<i\leq \frac{w-1}{2}) \\ 
        \label{jj:three} & 2 ( w-2i+1) Q_{[i,i-1]} = 0  \ (\text{if } \frac{w-1}{2} < i \leq
\frac{w+1}{2} )  
\end{align}
Using (\ref{jj:one}) and (\ref{jj:two}), we have 
\begin{equation} Q_{[i+1,i]} = 0\;, \text{ i.e., } 
        \label{eqn::zz}
        \beta_{[i+1,i]} = - \gamma_{[i,i+1]} 
        \quad 
        (0 \leq i\leq \frac{w-1}{2}) \end{equation} 
(\ref{jj:three}) says nothing when $w$ is odd, and 
if $w$ is even 
(\ref{jj:three}) is a part of (\ref{eqn::zz}).   
We summarize this observation as follows: 
if $w= 2\Omega$  then we have linear independent $\Omega$ relations and  
if $w= 2\Omega+1 $ then we have ($\Omega+1$) linear independent relations.   

%\medskip

Checking the coefficient of $\ds x^{i} y^{j} z^{w-i-j}$ ($i\ne j$), we have  
\begin{align*}
        - 2 \alpha_{[i,j]}(i-j) = & 
         \beta_{[1+i,j]} (1+i) -2 \beta_{[i,j-1]} (w-i-j+1) 
        % \\ & 
         + \gamma_{[i,1+j]}(1+j) -2 \gamma_{[i-1,j]}(w-i-j+1)\;,    
\end{align*} 
thus, we say the freedom of $\ds \alpha_A$ is $\ds \alpha_{[i,i]}$ with
 $2i\leq w$ and the number is $\Omega + 1$ even 
 if $w= 2\Omega$ or  $w= 2\Omega+1$.   
 We conclude that the kernel subspace of $\ds \myCC{2}{w}$ is $\ds
 (\Omega+1) + \{ \tbinom{2+w}{2} - \Omega^{*}\} + \tbinom{2+w}{2}$ where 
 $\ds\Omega^{*} = \begin{cases} \Omega & \text{ if } w=2\Omega\;,
         \\[-2mm]
 1+\Omega  & \text{ if } w = 2\Omega+1\;, \end{cases} $ and the 
 third Betti number is given by 
 \[ \tbinom{2+w}{2} + 
         \left(
         (\Omega+1) + \{ \tbinom{2+w}{2} - \Omega^{*}\} +
 \tbinom{2+w}{2}\right) - 3 \tbinom{2+w}{2} = \Omega +1 - \Omega^{*}\;. \] 

\paragraph{Kernel subspace of 1-st cochain space:}

Take a general 1-cochain $\ds \sum_{A} u^A ( \alpha_A \pdel_x +\beta_A
\pdel_y +\gamma _A \pdel_z )$ and solve the equation of the Schouten
bracket by $\pi$ is 0. 

Then the coefficient of $\pdel_z \wedge \pdel_x $ is 
\begin{align}  &  \sum_{A} 
\left( 2(a_1-a_2 - 1)  \alpha_A 
- 
\gamma_A ( -2 x \frac{\pdel}{\pdel_z} + z \frac{\pdel}{\pdel_y} )\right) u^A   
\end{align}
and the coefficient of $\pdel_z \wedge \pdel_y $ is 
\begin{align} &  \sum_{A} \left(  2(a_1-a_2+1) \beta_A   -
 \gamma_A 
( 2 y \frac{\pdel}{\pdel z} -  z \frac{\pdel}{\pdel x} )\right) u^A  
\end{align}
and the coefficient of $\pdel_x \wedge \pdel_y $ is 
\begin{align} &  \sum_{A}  
\left( \beta_A  
( - 2 x \frac{\pdel}{\pdel z} +  z \frac{\pdel}{\pdel y} )
- \alpha_A 
(  2 y \frac{\pdel}{\pdel z} -  z \frac{\pdel}{\pdel x} )
-\gamma_A\right) u^A \;. 
\end{align}

The coefficient of $\ds x^{i} y^{j} z^{w-i-j}$ are  
\begin{align}
\label{blign:one}
        & 2 (i-j-1) \alpha_{[i,j]} +2(w-i-j+1) \gamma_{[i-1,j]} 
- (j+1) \gamma_{[i,1+j]} = 0\;, \\
\label{blign:two}
        & 2 (i-j+1) \beta_{[i,j]} -2(w-i-j+1) \gamma_{[i,j-1]} 
+ (i+1) \gamma_{[1+i,j]} = 0\;, \\
\label{blign:three}
& \gamma_{[i,j]}  = -2 (w-i-j+1) \beta_{[i-1,j]}  + (j+1) \beta_{[i,j+1]} 
\\& \notag\hspace{10mm}
-2(w-i-j+1) \alpha_{[i,j-1]} 
+ (i+1) \alpha_{[1+i,j]}\; .  
\end{align} 
In (\ref{blign:one}), if $i-j-1 \ne 0$ then $\alpha_{[i,j]} $ are
determined by $\{ \gamma_{[k,\ell]} \}$. 
Also in 
 (\ref{blign:two}), if $i-j+1 \ne 0$ then $\beta_{[i,j]} $ are
determined by $\{ \gamma_{[k,\ell]} \}$. 

If 
 $i-j-1 = 0$ in (\ref{blign:one})  or 
 if $i-j+1 =0 $ in (\ref{blign:two}), then we have  the same relations
 \begin{equation}
         \label{eqn:kk:one}
         2(w-2j) \gamma_{[j,j]} = (j+1) \gamma_{[j+1,j+1]}\;.  
         %\quad (0\leq j \leq w/2-1) 
\end{equation}
We write the above recursive formula concretely as follows: 
\begin{align*}
        2w \gamma_{[0,0]} &=  \gamma_{[1,1]}\;, \\
        2(w-2) \gamma_{[1,1]} &= 2 \gamma_{[2,2]}\;, \\
%%        2(w-4) \gamma_{[2,2]} &= 3 \gamma_{[3,3]} \\
                               & \vdots \\
2(w-2 \Omega+2) \gamma_{[\Omega-1,\Omega-1]} &= \Omega
        \gamma_{[\Omega,\Omega]}\;, \\  
        2(w-2\Omega) \gamma_{[\Omega,\Omega]} &= (\Omega+1)
        \gamma_{[\Omega+1,\Omega+1]} = 0 \;.  
\end{align*}
When $w = 2\Omega +1$,  the last equation above tells 
$\ds \gamma_{[\Omega,\Omega]} = 0 $ and we have 
$\ds \gamma_{[j,j]} =0 $ for $j=0 .. \Omega$.  
In (\ref{blign:three}), we put $i=j$ and have 
\begin{align}
\label{blign:four}
\gamma_{[i,i]}  =&  -2 (w-2i+1) \beta_{[i-1,i]}  + (i+1) \beta_{[i,i+1]} 
\\& -2(w-2i+1) \alpha_{[i,i-1]} + (i+1) \alpha_{[1+i,i]} \notag\\
=& 
(i+1)( \beta_{[i,i+1]} + \alpha_{[1+i,i]})
\notag %\\ \notag & 
-2 (w-2i+1) (  \beta_{[i-1,i]} +  \alpha_{[i,i-1]}) \;. 
\end{align}
Denoting 
$\ds   \beta_{[i,i+1]}  + \alpha_{[1+i,i]} 
$ by $\myPP{i}$, the above
says that 
\[ (i+1) \myPP{i} =  ( 2w+2 -4i) \myPP{i-1}\;. \]
Since $\ds\myPP{0}=0$, we have $\ds \myPP{i}=0$ and so 
$\ds  \alpha_{[1+i,i]} + \beta_{[i,i+1]} =0 $ for each $i$.  Now we
count the freedom of $\ds \alpha_A, \beta_A$ and $\ds \gamma_A$.  
The freedom of  $\ds \gamma_A$ is $\ds \tbinom{2+w}{2} - (1+\Omega)$
because $\ds \gamma_{[i,i]} = 0$. We see the freedom of $\ds
\beta_A$ is $1+\Omega$ by    
$\ds  \beta_{[i,i+1]}  $, and $\ds\alpha_A$ have no freedom because 
$\ds  \alpha_{[1+i,i]} + \beta_{[i,i+1]} =0 $ for each $i$.  Thus, the 
dimension of the kernel subspace of $\ds\myCC{1}{w}$ is $\ds
\tbinom{2+w}{2}$. 
The first Betti number is 0 because $\ds \tbinom{2+w}{2} + 0 -
\tbinom{2+w}{2}$. 
Since we already know 
$\dim \ker $ 
subspace of $\ds\myCC{2}{w}$ is $\ds
2\tbinom{2+w}{2}$ when $w=2\Omega+1$, the second Betti number is also 0.

\bigskip

When $\ds w = 2\Omega$ (even), the situation is slightly different. 
We will count the freedom of $3\Omega+1$ variables 
$\ds \{ \gamma_{[i,i]}\}$ ($i=0.. \Omega$), 
$\ds \{ \beta_{[i,i+1]}\}$ ($i=0.. (\Omega-1)$) and  
$\ds \{ \alpha_{[i+1,i]}\}$ ($i=0..(\Omega-1)$).   
We have ($2\Omega+1$) linear equations   
(\ref{eqn:kk:one}) and 
(\ref{blign:four}) of those variables by using  
$\ds \myPP{i} :=  \beta_{[i,i+1]} +  \alpha_{[1+i,i]} $ (where 
$\ds\myPP{\Omega}=0$)  as follows:
\begin{align*}  
        &     \gamma_{[i,i]} + \phi(i)  \myPP{i-1} - (i+1) \myPP{i} 
        = 0 
        \;, \text{ where } 
        \phi(i) := 2 w+2-4i\\\noalign{and }
        &      \psi(i)  \gamma_{[i,i]} - (i+1) \gamma_{[i+1,i+1]} = 0 \;, 
        \text{ where }
\psi(i) := 2 w-4i \;. \end{align*}

The matrix which express the linear equations is the following: 
\[ M = \begin{bmatrix} 
         \begin{matrix} 
        1 & 0 & 0 & 0 & 0 \\
        0 & 1 & 0 & 0 & 0\\
        0 & 0 & \ddots & 0  & 0 \\
        0 & 0 & 0 & \ddots & 0\\
        0 & 0 & 0 & 0 & 1
                \\\hline
                \psi(0)  & -1 & 0 & 0 & 0 \\
                0 & \psi(1)  & -2 & 0 & 0\\
        0 & 0 & \ddots  & \ddots & 0 \\
                0 & 0 & 0 & \psi(\Omega-1) & -\Omega\\
        \end{matrix} 
       \quad \vline 
        &
        \begin{matrix} 
                -1 & 0 & 0 & 0 \\
                \phi(1) & -2 & 0 & 0\\
        0 & \ddots  & \ddots  & 0\\
        0 & 0 & \ddots  & -\Omega\\
                0 & 0 & 0 &  \phi(\Omega) 
                \\\hline
        0 & 0 & 0 & 0  \\
        0 & 0 & 0 & 0 \\
        0 & 0 & \ddots  & 0  \\
        0 & 0 & 0 & 0 \\
        \end{matrix} 
\end{bmatrix}
\]
where the (1,1)-matrix is the identity matrix of size $\Omega+1$, the
(2,2)-matrix is the zero matrix of size  $\Omega$, 
the (1,2)-matrix is 
the size $(\Omega+1, \Omega)$ and the diagonal entries are $-i$ and the
below line consists of $\phi(i)$, and 
the (2,1)-matrix is 
the size $(\Omega, \Omega+1 )$ and the diagonal entries are $\psi(i-1)$ and the
above line consists of $-i$.

We see that the rank of $M$ is full ($2\Omega +1$) by the following Lemma.  
So, the freedom of those variables is $3\Omega +1 -
 (2\Omega+1) = \Omega$. Thus the dimension of the kernel subspace of $\ds\myCC{1}{w}$ is
 $\ds \left (\tbinom{2+w}{2} - (\Omega+1) \right) + \Omega = 
 \tbinom{2+w}{2} - 1$ when $w= 2\Omega$.  This shows that 
 the first Betti number is 0 and the second Betti number is also 0.
 \kmqed

\begin{Lemma}
        \[ \det M = \Omega! \sum_{j=1}^{\Omega+1} \psi(0)
        \cdots \psi(j-2)\phi(j)\cdots \phi(\Omega)\]   holds,  
        where   $ \phi(i) = 4 \Omega +2-4i$ and 
$\psi(i) = 4 \Omega-4i$.  
\end{Lemma}
\textbf{Proof of Lemma:} We follow the definition of determinant directly. Since the
(2,2)-submatrix  of $M$ is zero matrix, among
permutations of $1,2,\ldots, 2\Omega +1$, the candidate $\sigma$ which
contribute the determinant of $M$ must satisfy $\sigma(i) \leq \Omega+1$
for each $i > \Omega+1$ and so for some  $\ds j\leq \Omega+1$, $\sigma$
satisfies  
\[
        \sigma( \{\Omega+2,\ldots,2\Omega+1\} ) =
        \{1,\ldots, \widehat{j},\ldots,\Omega+1\},\quad 
        \sigma( \{1, \ldots,\Omega+1\} ) =
\{j,\Omega+2,\ldots,2\Omega+1\}\;.  
\]
Furthermore, \( i \leq \sigma(1+\Omega+i) \leq i+1 \) and $j$ is
specified, then $\sigma(1+\Omega+i) = i$ for $i < j$ and 
$\sigma(1+\Omega+i) = i+1$ for $i \geq  j$.   For $\sigma(i)$ ($i\leq
1+\Omega$), effective possibilities are 
\begin{alignat*}{2} 
        \sigma(1)& =1 \quad  \text{ or }\quad &  & \sigma(1) =1+\Omega+1 \\
        \sigma(2)& =2 \quad  \text{ or }\quad &  1+\Omega+1 \leq &
        \sigma(2)  \leq 1+\Omega+2 \\
        \vdots &  &  & \\
        \sigma(\Omega)& =\Omega \quad  \text{ or }\quad & 
        2 \Omega \leq & \sigma(\Omega)  \leq 1+\Omega+\Omega \\
        \sigma(\Omega+1)& =\Omega+1 \quad  \text{ or }\quad  &
        2\Omega+1 = & \sigma(\Omega+1)    
\end{alignat*}
and if $j$ is specified, then $\sigma(j)=j$ holds necessarily  and 
\[ 
        \sigma( 1) = 1+\Omega +1 , \; \ldots\;  
        \sigma(j- 1) = 1+\Omega + j-1 , \; 
        \sigma(j) = j , \; 
        \sigma(j+ 1) = \Omega + j+1 , \; \ldots \; 
\sigma(\Omega + 1) = 2 \Omega + 1\;. \]
Thus, when $j$ is specified, the permutation $\sigma$ is uniquely determined
as 
\[\small
        \setcounter{MaxMatrixCols}{20} 
\begin{pmatrix} 
1 & \cdots & j-1 & j & j+1 &\cdots & \Omega+1
        & \Omega +2 & \cdots & \Omega+j & \Omega+j+1 & &
        2\Omega+1         \\
\Omega+2 & & \Omega+j & j & \Omega+1+j& & 2\Omega+1 & 1 &&j-1&
        j+1 && \Omega+1
        \end{pmatrix} 
        \] 
and is decomposed as $(1,\Omega+2)\cdots (j-1,\Omega+j)
(j+1,\Omega+j+1) \cdots (\Omega+1, 2\Omega+1)$ of $\Omega$
transpositions and so the signature of $\sigma$ is $\ds (-1)^{\Omega}$.  
\kmcomment{ 
\[
        \setcounter{MaxMatrixCols}{20} 
        \left( 
\begin{matrix} 
1 & \cdots & j-1 & j & j+1 &\cdots & \Omega+1
        & 1+\Omega +1 & \cdots & 1+\Omega+(j-1) & 1+\Omega+j & &
        2\Omega+1         \\
\Omega+2 & & \Omega+j & j & \Omega+1+j& & 2\Omega+1 & 1 &&j-1&
        j+1 && \Omega+1
        \end{matrix} 
        \right)\] 
}%endOFkmcomment
Therefore, we conclude 
\begin{align*}
        \det M =& \sum_{j=1}^{\Omega+1} (-1)^{\Omega}(-1)\cdots(-j+1)
         \phi(j)\cdots \phi(\Omega) 
        \psi(0)\cdots\psi(j-2) (-j) (-j-1) \cdots (-\Omega) 
        \\= & 
        \Omega! \sum_{j=1}^{\Omega+1}
        \psi(0)\cdots\psi(j-2) \phi(j)\cdots \phi(\Omega) 
        \;. \hspace{60mm}  \blacksquare 
\end{align*}

%\kmqed 

\begin{kmRemark}
        We already know 
        the behavior of Betti numbers about   
        $SL(2)$ ($A_1$), 
        $SO(3)$ ($B_1$) or $\ds Sp(\mR^2)$ ($C_1$). Here  
we calculate 
        Betti numbers for 
        lower weight of $SO(4)$ ($D_2$) Lie Poisson structure: 
        \begin{center}
        \setlength{\extrarowheight}{-3pt}
        \begin{tabular}{|c|*{7}{c}|} \hline
        $w \backslash m$ & 0 & 1 & 2 & 3 & 4 & 5 & 6 \\\hline
        0 & 1 & 0 & 0 & 2 & 0 & 0 & 1 \\
        1 & 0 & 0 & 0 & 0 & 0 & 0 & 0 \\
        2 & 2 & 0 & 0 & 4 & 0 & 0 & 2 \\
        3 & 0 & 0 & 0 & 0 & 0 & 0 & 0 \\
        4 & 3 & 0 & 0 & 6 & 0 & 0 & 3 \\
5 & 0 & 0 & 0 & 0 & 0 & 0 & 0 \\\hline
\end{tabular}
\end{center}

Concerning $w$, we expect some mechanism, in some sense, combinations of
patterns of $\ds Sp(\mR^2)$ and Heisenberg Lie algebra. 

\medskip
%\end{kmRemark} \begin{kmRemark}
        $SL(3)$ ($A_2$) is 8 dim and for the weight 0, 1, 2 we have 
        \kmcomment{the following table:}
        \begin{center}
        \setlength{\extrarowheight}{-3pt}
        \begin{tabular}{|c|*{9}{c}|} \hline
$w \backslash  m$ & 0 & 1 & 2 & 3 & 4 & 5 & 6 &7&8 \\\hline
                0 & 1 & 0 & 0 & 1 & 0 & 1 & 0 & 0 & 1  \\
                1 & 0 & 0 & 0 & 0 & 0 & 0 & 0 & 0 & 0  \\
                2 & 1 & 0 & 0 & 1 & 0 & 1 & 0 & 0 & 1   
\\\hline
\end{tabular}
        \end{center}

        About $SO(5)$ ($B_2$) or $\ds Sp(\mR^4)$ ($C_2$): 
        \kmcomment{is 10 dimensional and for the weight 0, 1, 2
        we have the following table:}
        \begin{center}
        \begin{tabular}{|c|*{11}{c}|} \hline
$w \backslash  m$ & 0 & 1 & 2 & 3 & 4 & 5 & 6 &7&8&9&10 \\\hline
                0 & 1 & 0 & 0 & 1 & 0 & 0 & 0 & 1 & 0 & 0 & 1 \\
                1 & 0 & 0 & 0 & 0 & 0 & 0 & 0 & 0 & 0 & 0 & 0\\
                2 & 1 & 0 & 0 & 1 & 0 & 0 & 0 & 1 & 0 & 0 & 1 
\\\hline
\end{tabular}
        \end{center} 
\end{kmRemark} 

\bigskip

Looking at those examples, it seems to be interesting 
to investigate the Betti numbers associated
with Lie Poisson structure of simple Lie algebras $A_n$, $B_n$,
$C_n$,\ldots and to understand periodicity or
symmetries of their Betti numbers.